\documentclass[11pt,twoside]{article}
\newcommand{\ifims}[2]{#1} 
\newcommand{\ifAMS}[2]{#1}   
\newcommand{\ifau}[4]{#1}  
\newcommand{\ifbook}[2]{#1}   
\newcommand{\ifunivariate}[2]{#1} 
\newcommand{\ifsupnorm}[2]{#1} 
\newcommand{\ifLaplace}[2]{#1} 
\newcommand{\ifNL}[2]{#1}  
\newcommand{\ifapp}[2]{#1}  

\usepackage{ifthen}
\usepackage[ampersand]{easylist}
\usepackage{amsmath,amssymb,amsthm}
\usepackage{mathtools}
\usepackage{natbib}
\usepackage{epsfig,graphicx}
\usepackage{comment}
\usepackage{color}
\usepackage{srcltx}
\usepackage[mathscr]{eucal}
\usepackage[math]{easyeqn}
\usepackage{etoolbox}
\usepackage{hyperref}
\hypersetup{
            colorlinks,
            linkcolor=hookersgreen,
            linktoc=blue,
            citecolor=ultramarine,
            urlcolor=black,
            filecolor=black
            }
\usepackage{nameref}
\usepackage{multirow}

\usepackage{accents}
\usepackage{mathrsfs}
\usepackage{bbm}
\usepackage{algorithm}
\usepackage{algpseudocode}

\ifims{
\textheight=22cm
\textwidth=14.8cm
\topmargin=0pt
\oddsidemargin=1.0cm
\evensidemargin=1.0cm
\linespread{1.3}

}{ 
}

\numberwithin{equation}{section}
\numberwithin{figure}{section}
\newcounter{example}[section]
\numberwithin{example}{section}
\newcounter{remark}[section]
\numberwithin{remark}{section}
\newtheorem{theorem}{Theorem}[section]
\newtheorem{proposition}[theorem]{Proposition}
\newtheorem{lemma}[theorem]{Lemma}
\newtheorem{corollary}[theorem]{Corollary}

\newtheorem{exmp}[example]{Example}
\newtheorem{rmrk}[remark]{Remark}
\newenvironment{example}{\begin{exmp}\rm}{\end{exmp}}
\newenvironment{remark}{\begin{rmrk}\rm}{\end{rmrk}}

\ifbook{
    \newcommand{\Chapter}[1]{\section{#1}}
    \newcommand{\Section}[1]{\subsection{#1}}
    \newcommand{\Subsection}[1]{\subsubsection{#1}}
    \def\Chname{Section }
    \def\chname{section }
    
  }
  {
    \newcommand{\Chapter}[1]{\chapter{#1}}
    \newcommand{\Section}[1]{\section{#1}}
    \newcommand{\Subsection}[1]{\subsection{#1}}
    \def\Chname{Chapter}
    
  }

\renewcommand{\(}{$\,}
\renewcommand{\)}{\,$}

\def\nquad{\hspace{-1cm}}
\def\eqdef{\stackrel{\operatorname{def}}{=}}
\def\tod{\stackrel{d}{\longrightarrow}}

\DeclareMathAlphabet{\mathbbmsl}{U}{bbm}{bx}{sl}

\DeclareMathSymbol{\Alpha}{\mathalpha}{operators}{"41}
\DeclareMathSymbol{\Beta}{\mathalpha}{operators}{"42}
\DeclareMathSymbol{\Epsilon}{\mathalpha}{operators}{"45}
\DeclareMathSymbol{\Zeta}{\mathalpha}{operators}{"5A}
\DeclareMathSymbol{\Eta}{\mathalpha}{operators}{"48}
\DeclareMathSymbol{\Iota}{\mathalpha}{operators}{"49}
\DeclareMathSymbol{\Kappa}{\mathalpha}{operators}{"4B}
\DeclareMathSymbol{\Mu}{\mathalpha}{operators}{"4D}
\DeclareMathSymbol{\Nu}{\mathalpha}{operators}{"4E}
\DeclareMathSymbol{\Omicron}{\mathalpha}{operators}{"4F}
\DeclareMathSymbol{\Rho}{\mathalpha}{operators}{"50}
\DeclareMathSymbol{\Tau}{\mathalpha}{operators}{"54}
\DeclareMathSymbol{\Chi}{\mathalpha}{operators}{"58}
\DeclareMathSymbol{\omicron}{\mathord}{letters}{"6F}

\newcommand{\cc}[1]{\mathscr{#1}}
\newcommand{\bb}[1]{\boldsymbol{#1}}

\renewcommand{\bar}[1]{%
  \hbox{%
    \vbox{%
      \hrule height 0.5pt 
      \kern0.5ex
      \hbox{%
        \kern-0.1em
        \ensuremath{#1}%
        \kern-0.1em
      }%
    }%
  }%
} 

\renewcommand{\hat}[1]{\widehat{#1}}
\renewcommand{\tilde}[1]{\widetilde{#1}}

\makeatletter
\def\mathcenterto#1#2{\mathclap{\phantom{#1}\mathclap{#2}}\phantom{#1}}
\let\old@widetilde\widetilde
\def\widetildeto#1#2{\mathcenterto{#2}{\old@widetilde{\mathcenterto{#1}{#2\,}}}}
\let\old@widehat\widehat
\def\widehatto#1#2{\mathcenterto{#2}{\old@widehat{\mathcenterto{#1}{#2\,}}}}
\makeatother

\newcommand{\thankstitle}[1]{\ifthenelse{\equal{#1}{}}{}{\thanks{#1}}}
\newcommand{\thanksau}[1]{\ifthenelse{\equal{#1}{}}{}{\thanks{#1}}}

\newcommand{\aua}[6]
{\def\authora{#1}
\def\runauthora{#2}
\def\addressa{#3}
\def\emaila{#4}
\def\affiliationa{#5}
\def\thanksa{#6}}

\def\theauthors{
\ifau{ 
  \author{
    \authora
    \thanksau{\thanksa}
    \\[5.pt]
    \addressa \\
    \texttt{ \emaila}
  }
}
{  
  \author{
    \authora
    \thanksau{\thanksa}
    \\[5.pt]
    \addressa \\
    \texttt{ \emaila}
    \and
    \authorb
    \thanksau{\thanksb}
    \\[5.pt]
    \addressb \\
    \texttt{ \emailb}
  }
}
{   
  \author{
    \authora
    \thanksau{\thanksa}
    \\[5.pt]
    \addressa \\
    \texttt{ \emaila}
    \and
    \authorb
    \thanksau{\thanksb}
    \\[5.pt]
    \addressb \\
    \texttt{ \emailb}
    \and
    \authorc
    \thanksau{\thanksc}
    \\[5.pt]
    \addressc \\
    \texttt{ \emailc}
  }
} {   
  \author{
    \authora
    \thanksau{\thanksa}
    \\[5.pt]
    \addressa \\
    \texttt{ \emaila}
    \and
    \authorb
    \thanksau{\thanksb}
    \\[5.pt]
    \addressb \\
    \texttt{ \emailb}
    \and
    \authorc
    \thanksau{\thanksc}
    \\[5.pt]
    \addressc \\
    \texttt{ \emailc}
    \and
    \authord
    \thanksau{\thanksd}
    \\[5.pt]
    \addressd \\
    \texttt{ \emaild}
  }
}
}

\renewcommand{\Gamma}{\varGamma}
\renewcommand{\Pi}{\varPi}
\renewcommand{\Sigma}{\varSigma}
\renewcommand{\Delta}{\varDelta}
\renewcommand{\Lambda}{\varLambda}
\renewcommand{\Psi}{\varPsi}
\renewcommand{\Phi}{\varPhi}
\renewcommand{\Theta}{\varTheta}
\renewcommand{\Omega}{\varOmega}
\renewcommand{\Xi}{\varXi}
\renewcommand{\Upsilon}{\varUpsilon}

\def\argmax{\operatornamewithlimits{argmax}}

\def\av{\bb{a}}

\def\uv{\bb{u}}

\def\wv{\bb{w}}
\def\xv{\bb{x}}

\def\Av{\bb{A}}

\def\Sv{\bb{S}}

\def\Uv{\bb{U}}

\def\Xv{\bb{X}}
\def\Yv{\bb{Y}}

\def\deltav{\bb{\delta}}

\def\epsv{\bb{\varepsilon}}

\def\gammav{\bb{\gamma}}

\def\xiv{\bb{\xi}}

\def\Psiv{\bb{\Psi}}

\def\sumi{\sum_{i=1}^{n}}

\usepackage{color}

\definecolor{blue(pigment)}{rgb}{0.2, 0.2, 0.6}
\definecolor{ultramarine}{rgb}{0.07, 0.04, 0.56}
\definecolor{darkspringgreen}{rgb}{0.09, 0.45, 0.27}
\definecolor{hookersgreen}{rgb}{0.0, 0.44, 0.0}
\definecolor{plum(traditional)}{rgb}{0.56, 0.27, 0.52}
\definecolor{purple(html/css)}{rgb}{0.5, 0.0, 0.5}
\definecolor{magenta(dye)}{rgb}{0.79, 0.08, 0.48}

\def\IFL{\IF}
\def\HUL{\mathsf{H}}
\def\XXL{\mathcal{X}}
\def\dimLL{\dimL(\xvs)}

\def\lamb{\bb{\lambda}}

\def\fn{g}
\def\upsvn{\upsvd}
\def\rrn{\rr}

\def\dltwaa{\CONSTi}
\def\uvdd{\uv}

\def\DPU{\mathbb{M}}
\def\HL{\mathsf{m}}
\def\HLG{\mathbf{M}}

\def\HPck{\HP}
\def\weak{\circ}
\def\GPw{\GP_{\weak}}
\def\DPw{\DP_{\weak}}

\def\CDG{\CONSTi_{\GP}}

\def\errt{\tilde{\err}}

\def\rsmall{\varrho}
\def\Psimean{\bar{\Psiv}}

\def\rhos{\rho^{*}}

\def\hmax{\mathsf{c}}

\def\hL{h}
\def\HVB{\mathbbmss{V}}

\def\dimVG{\dimA_{\GP}}
\def\fGu{h}
\def\Ex{\E_{\xx}}

\def\HPck{H}
\def\wPsi{\lambda}
\def\rrVG{\rr_{\GP}}
\def\CGP{w}
\def\GPa{\GP_{0}}
\def\CGPa{\CGP_{0}}
\def\CGPs{\mmps}

\def\smpa{\smp_{0}}

\def\lgd{f}
\def\elll{\ell}
\def\PfL{\P_{\lgd}}
\def\lgdL{\elll}

\def\rrL{\rr}

\def\rrLt{\tilde{\rrL}}

\def\dimL{\dimA}
\def\dimLt{\tilde{\dimL}}

\def\lamb{\bb{\lambda}}
\def\dagg{\prime}

\def\proj{P}

\def\HU{\HV}

\def\xvs{\xv^{*}}

\def\HPck{\HV}

\def\upsvck{\hat{\upsv}}

\def\amax{\nu}
\def\Cond{\,\, \Big| \, }

\def\Matr{\mathfrak{M}}

\def\xvb{\bar{\xv}}

\def\weights{\weight^{*}}

\def\rhou{b}

\def\fba{\bar{f}}

\def\nano{\circ}
\def\dual{*}
\def\zqv{\bb{\zq}}

\def\vH{\vA}

\def\Idd{\Gamma}

\def\HUDe{\Delta}

\def\Rem{\mathcal{R}}

\def\Remi{\mathcal{Q}}

\def\Eta{\mathcal{H}}

\def\nbin{N}

\def\WS{\mathfrak{A}}

\def\lsh{\alpha}

\def\HV{H}

\def\HU{\mathsf{D}}

\def\fG{f_{\GP}}

\def\smp{s}
\def\mux{\mu_{\xx}}

\def\dltw{\delta}

\def\dltwb{\omega}

\def\dltwbs{\dltwb^{*}}
\def\dltwm{\kappa}

\def\dltwu{\tau}
\def\dltwd{\dltw^{\dagg}}
\def\dltwbd{\dltwb^{\dagg}}
\def\dltwbss{\dltwb^{+}}

\def\II{\mathcal{I}}

\def\R{\mathbbmsl{R}}
\def\E{\mathbbmsl{E}}

\def\P{\mathbbmsl{P}}

\def\kappa{\varkappa}

\def\Bernoulli{\operatorname{Bernoulli}}

\def\diag{\operatorname{diag}}

\def\Exp{\operatorname{Exp}}
\def\Fr{\operatorname{Fr}}

\def\ND{\mathcal{N}}

\def\oper{\operatorname{op}}

\def\Var{\operatorname{Var}}

\def\T{\top}
\def\tr{\operatorname{tr}}

\def\TV{\operatorname{TV}}

\def\cond{\, \big| \,}

\def\nsize{{n}}

\def\sumi{\sum_{i=1}^{\nsize}}

\def\ex{\mathrm{e}}

\def\Id{I\!\!\!I}
\def\Ind{\operatorname{1}\hspace{-4.3pt}\operatorname{I}}

\def\alp{\alpha}

\def\biasv{\bb{b}}

\def\B{\cc{B}}

\def\BB{I\!\!B}     
\def\BB{B}

\def\BBB{\cc{B}}

\def\BBH{W}

\def\cdens{\phi}

\def\CA{\mathcal{A}}

\def\CONST{\mathtt{C} \hspace{0.1em}}
\def\CONSTi{\mathtt{C}}

\def\CONSTIF{\CONSTi_{\IF}}
\def\CONSTgp{\CONSTi_{\gp}}

\def\CONSTPsi{\CONSTi_{\Psi}}

\def\CONSTphi{\hmax_{\cdens}}

\def\CONSTn{\CONSTi_{0}}

\def\CONSTbal{\CONSTi_{0}}
\def\CONSTnabl{\CONSTi_{\circ}}
\def\CONSTV{\CONSTi_{\VP}}
\def\CONSTGP{\CONSTi_{\GP}}

\def\CS{\cc{E}}

\def\DP{D}

\def\DPt{\tilde{\DPc}}
\def\DPf{\mathbbmsl{\DP}}
\def\DPfGP{\DPf_{\GP}}

\def\DPGP{\DP_{\GP}}

\def\DPGPt{\DPt_{\GP}}

\def\DPt{\tilde{\DP}}

\def\DPck{\hat{\DP}}

\def\DV{\mathsf{D}}

\def\DVL{\DV}

\def\DU{\mathbb{D}}

\def\deltav{\bb{\delta}}

\def\dist{d}

\def\dimH{\dimA}

\def\dimp{p}

\def\dimA{\mathtt{p}}

\def\dimG{\dimA_{\GP}}

\def\dimq{q}

\def\dimt{\tilde{\dimA}}

\def\dens{f}

\def\Eta{\cc{H}}

\def\err{\diamondsuit}

\def\eps{\epsilon}			
\def\eps{\varepsilon}

\def\fs{f}

\def\gaussv{\bb{\gauss}}
\def\gauss{\gamma}

\def\gm{\mathtt{g}}
\def\gmc{\gm_{c}}
\def\gmb{\gm}

\def\gp{g}

\def\GP{G}

\def\HP{\mathsf{H}}

\def\IF{\mathbbmsl{F}}
\def\IFGP{\IF_{\GP}}

\def\IFt{\mathbbmsl{I}}

\def\IFt{\tilde{\IF}}
\def\IFt{\tilde{\IF}}

\def\IFba{\accentset{\circ}{\IF}}

\def\jc{j'}

\def\Kappa{\cc{K}}

\def\kullb{\cc{K}} 

\def\LT{L}
\def\LGP{\LT_{\GP}}

\def\ldens{\ell}

\def\mm{m}

\def\mmps{\mm_{0}}

\def\muc{\mu_{c}}

\def\pent{\operatorname{pen}}

\def\Pone{P}
\def\Pdom{\mu_{0}}
\def\PDOM{\bb{\mu}_{0}}

\def\PG{\P'}

\def\Psimean{\bar{\Psi}}
\def\Proj{\Pi}

\def\Prior{\Pi}

\def\priorh{\pi}

\def\priord{\pi}

\def\QP{Q}

\def\rhot{t}

\def\riskt{\cc{R}}

\def\rr{\mathtt{r}}
\def\rrs{\rr^{*}}

\def\rru{\rr_{\circ}}

\def\rrGP{\rr_{\GP}}
\def\rups{\rr_{\circ}}

\def\rups{\rr_{0}}

\def\Spsi{S}

\def\score{\nabla}

\def\supA{\lambda}

\def\thetas{\theta^{*}}

\def\Tau{T}

\def\uvc{\uv^{c}}

\def\ups{\upsilon}
\def\upsv{\bb{\ups}}

\def\upsvd{\upsv^{\circ}}
\def\upsvs{\upsv^{*}}

\def\vupsv{\upsv}

\def\ups{\upsilon}
\def\upsv{\bb{\ups}}

\def\upsd{\ups^{\circ}}

\def\upss{\ups^{*}}

\def\upsvb{\bar{\upsv}}

\def\UV{\mathcal{U}}

\def\UVL{\UV}

\def\UVz{\UV}

\def\UVt{\tilde{\UV}}

\def\Ups{\varUpsilon}
\def\Upsd{\Ups^{\circ}}

\def\vA{\mathtt{v}}

\def\VP{V}

\def\VV{\mathbbmsl{V}}

\def\vp{\mathbf{v}}	

\def\vp{\mathsf{v}}

\def\Weight{W}
\def\weight{w}

\def\wv{\bb{w}}

\def\xivGP{\xiv_{\GP}}

\def\xx{\mathtt{x}}
\def\xxc{\xx_{c}}

\def\XX{\cc{X}}

\def\zq{z}
\def\zqc{\zq_{c}}

\def\zz{\mathfrak{z}}

\def\thetitle{Finite samples inference and critical dimension for stochastically linear models}
\def\theruntitle {Finite samples inference and critical dimension}

\def\theabstract{
The aim of this note is to state a couple of general results about the properties of the penalized 
maximum likelihood estimators (pMLE) and of the posterior distribution for parametric models
in a non-asymptotic setup and for possibly large or even infinite parameter dimension.
We consider a special class of stochastically linear smooth (SLS) models satisfying two major conditions: the stochastic component
of the log-likelihood is linear in the model parameter, while the expected log-likelihood is a smooth function.
The main results simplify a lot if the expected log-likelihood is concave.
For the pMLE, we establish a number of finite sample bounds about its concentration and large deviations as well as 
the Fisher and Wilks expansion.
The later results extend the classical asymptotic Fisher and Wilks Theorems about the MLE to the  
non-asymptotic setup with large parameter dimension which can depend on the sample size.
For the posterior distribution, our main result states a Gaussian approximation of the posterior
which can be viewed as a finite sample analog of the prominent Bernstein--von Mises Theorem. 
In all bounds, the remainder is given explicitly and can be evaluated in terms of the effective sample size
and effective parameter dimension.
The results are dimension and coordinate free. 
In spite of generality, all the presented bounds are nearly sharp and the classical asymptotic results can be obtained
as simple corollaries. 
Interesting cases of logit regression and of estimation of a log-density with smooth or truncation priors are used to specify the results and to explain the main notions.
}

\def\kwdp{62F15}
\def\kwds{62F25}

\def\thekeywords{penalized MLE, Fisher and Wilks expansion, posterior, Laplace approximation}

\def\thankstitle{}

\aua
{Vladimir Spokoiny}
{Spokoiny, V. }
{
    Weierstrass Institute and HU Berlin,  
    \\
    Mohrenstr. 39, 10117 Berlin, Germany
    }
{spokoiny@wias-berlin.de}
{Weierstrass-Institute and Humboldt University Berlin}
{
  Financial support by the German Research Foundation (DFG) through the Collaborative Research Center 1294 ``Data assimilation'' is gratefully acknowledged.
}

\begin{document}
\thispagestyle{empty}
{
\title{\thetitle}
\theauthors

\maketitle
\begin{abstract}
{\footnotesize \theabstract}
\end{abstract}

\ifAMS
    {\par\noindent\emph{AMS 2010 Subject Classification:} Primary \kwdp. Secondary \kwds}
    {\par\noindent\emph{JEL codes}: \kwdp}

\par\noindent\emph{Keywords}: \thekeywords
} 

\tableofcontents

\Chapter{Introduction}
\label{Sgenintr}


This paper presents some general results describing 
the properties 
of the \emph{penalized Maximum Likelihood Estimator} (pMLE) and of the \emph{posterior distribution}.


Our starting point is a parametric assumption about the distribution \( \P \) of the data \( \Yv \):
\( \P \) belongs a given parametric family \( (\P_{\upsv} \, , \upsv \in \Ups) \) 
dominated by a sigma-finite measure \( \PDOM \).
This assumption is usually an idealization of reality
and the true distribution \( \P \) is not an element of \( (\P_{\upsv}) \).
However, a parametric assumption, even being wrong, may appear to be very useful, because it yields 
the method of estimation.
Namely, the MLE \( \tilde{\upsv} \) is defined by maximizing the log-likelihood function 
\( L(\Yv,\upsv) = L(\upsv) = \log \frac{d\P_{\upsv}}{d\PDOM}(\Yv) \) over the parameter set \( \Ups \):
\begin{EQA}
	\tilde{\upsv}
	&=&
	\argmax_{\upsv \in \Ups} L(\upsv) .
\label{tuauv}
\end{EQA}
For a penalty function \( \pent_{\GP}(\upsv) \) on \( \Ups \), the \emph{penalized MLE} 
\( \tilde{\upsv}_{\GP} \) is defined by maximizing the penalized MLE 
\( \LGP(\upsv) = L(\upsv) - \pent_{\GP}(\upsv) \):
\begin{EQA}
	\tilde{\upsv}_{\GP}
	&=&
	\argmax_{\upsv \in \Ups} \LGP(\upsv)
	=
	\argmax_{\upsv \in \Ups} \bigl\{ L(\upsv) - \pent_{\GP}(\upsv) \bigr\} .
\label{tuGauULGuauUp}
\end{EQA}
The sub-index \( \GP \) in the penalty relies to its quadratic structure:
\begin{EQA}
	\pent_{\GP}(\upsv)
	&=&
	\frac{1}{2} \| \GP \upsv \|^{2}
\label{vfgiubhygtfdehwjkie}
\end{EQA}
for a symmetric \( \dimp \times \dimp \) positive definite matrix \( \GP \in \Matr_{\dimp} \).

In the Bayesian setup with a prior measure \( \Prior_{\GP} \) on \( \Ups \) having a density \( \priord_{\GP} \),
the posterior distribution is a random measure on \( \Ups \) defined by normalizing the product
\( \priord_{\GP}(\upsv) \exp L(\upsv) \):
\begin{EQA}
	\upsv_{\GP} \cond \Yv
	& \propto &
	\priord_{\GP}(\upsv) \exp L(\upsv) .
\label{exLYupGptuGY}
\end{EQA}
In particular, the posterior mode defined as the point of maximum of the posterior density 
coincides with the penalized MLE \( \tilde{\upsv}_{\GP} \) for the penalty 
\( \pent_{\GP}(\upsv) = - \log \priord_{\GP}(\upsv) \).
The sub-index \( \GP \) relies to the Gaussian prior with parameters \( \ND(0,\GP^{2}) \).
The use of the penalty and the prior centered at zero is only for notational convenience. 
Everything below can be extended to the case of an arbitrary center \( \upsv_{0} \).

\Section{Classical asymptotic theory}

The classical Fisher parametric theory assumes that 
\( \Ups \) is a subset of a finite-dimensional Euclidean space \( \R^{\dimp} \),
the underlying data distribution \( \P \) indeed belongs to the considered parametric family 
\( (\P_{\upsv}) \), that is, \( \Yv \sim \P = \P_{\upsvs} \) for some \( \upsvs \in \Ups \).
In addition, some regularity of the family \( (\P_{\upsv}) \), or, equivalently, 
of the log-likelihood function \( L(\upsv) \) is assumed. 
This, in particular, enables us to 
apply the third order Taylor expansion of \( L(\upsv) \) around the point of maximum \( \tilde{\upsv} \)
and to obtain a Fisher expansion
\begin{EQA}
	\tilde{\upsv} - \upsvs
	& \approx &
	\IF^{-1} \nabla L(\upsvs) .
\label{tumusaxFm1nLus}
\end{EQA}
Here  \( \IF = \IF(\upsvs) \) is the total Fisher information at \( \upsvs \) 
defined as a negative Hessian of the expected log-likelihood function \( \E L(\upsv) \):
\begin{EQA}
	\IF(\upsv)
	&=&
	- \nabla^{2} \E L(\upsv) .
\label{IFuvmn2ELuv}
\end{EQA}
Under standard parametric assumptions, \( \IF(\upsv) \) is symmetric positive definite, \( \IF(\upsv) \in \Matr_{\dimp} \).
Moreover, if the data \( \Yv \) is generated as a sample of independent random variables 
\( Y_{1},\ldots,Y_{n} \), then the log-likelihood has an additive structure:
\( L(\upsv) = \sumi \ell(Y_{i},\upsv) \).
This allows to establish asymptotic standard normality of the standardized score 
\( \xiv \eqdef \IF^{-1/2} \nabla L(\upsvs) \) and hence, to state Fisher and Wilks Theorems:
as \( n \to \infty \)
\begin{EQ}[rclcl]
	\IF^{1/2} \bigl( \tilde{\upsv} - \upsvs \bigr)
	& \approx &
	\xiv
	& \tod &
	\gaussv,
	\\
	2 L(\tilde{\upsv}) - 2 L(\upsvs)
	& \approx &
	\| \xiv \|^{2}
	& \tod &
	\| \gaussv \|^{2} \sim \chi^{2}_{\dimp} \, ,
\label{2Lm2Lusg2sc2p}
\end{EQ}
where \( \gaussv \) is a standard Gaussian vector in \( \R^{\dimp} \) and 
\( \chi^{2}_{\dimp} \) is a chi-squared distribution with \( \dimp \) degrees of freedom.
These results are fundamental and build the basis for the most statistical applications 
like analysis of variance, canonical and correlation analysis, 
uncertainty quantification and hypothesis testing etc.
We refer to \cite{vdV00} for a comprehensive discussion and an historical overview of the related results
including the general LAN theory by L. Le Cam.

In the Bayesian framework, the Laplace approximation of the integral  
\( \int \priord(\upsv) \ex^{ L(\upsv)} \, d\upsv \) leads to the prominent
Bernstein -- von Mises Theorem claiming asymptotic normality of the posterior with 
mean \( \tilde{\upsv} \) and the covariance \( \IF^{-1} \):
\begin{EQA}
	\Prior\bigl\{ \IF^{1/2} ( \upsv - \tilde{\upsv} ) \cond \Yv \bigr\}
	& \tod & 
	\ND(0,\Id_{\dimp}) ,
	\qquad 
	n \to \infty.
\label{PrF12utucYni}
\end{EQA}
Again, this is one of most fundamental results in statistical inference allowing 
to effectively construct elliptic credible sets in typical applications; see e.g. \cite{GhVa2017}
and references therein. 

\Section{Challenges and goals}
Modern statistical problems require to extend the classical results in several directions.

\Subsection{Model misspecification}
First, the model can be misspecified and the underlying data generating measure \( \P \) 
is not an element of the family \( (\P_{\upsv} \, , \upsv \in \Ups) \).
This means that the used log-likelihood function is not necessarily a true log-likelihood.
In particular, the condition \( \E \exp L(\upsv) = 1 \) does not hold true.
Also the target of estimation \( \upsvs \) has to be redefined 
as the maximizer of the expected log-likelihood:
\begin{EQA}
	\upsvs
	& \eqdef &
	\argmax_{\upsv \in \Ups} \E L(\upsv) 
\label{pd47787tyhrujeirfjk}
\end{EQA}
leading to some modelling bias as the distance between \( \P \) and \( \P_{\upsvs} \).
This also concerns the use of a penalty or an impact of the prior leading to some estimation bias.
When operating with the penalized log-likelihood \( \LGP(\upsv) \), the target of estimation 
becomes
\begin{EQA}
	\upsvs_{\GP}
	& \eqdef &
	\argmax_{\upsv \in \Ups} \E \LGP(\upsv) ,
\label{upssaruELuaiU}
\end{EQA}
which might be significantly different from \( \upsvs \).
Everywhere we use the notation 
\begin{EQA}
	\IF(\upsv) 
	&=& 
	- \nabla^{2} \E L(\upsv) ,
	\\
	\IF_{\GP}(\upsv) 
	&=& 
	- \nabla^{2} \E \LGP(\upsv) 
	= 
	\IF(\upsv) + \GP^{2} .
\label{e456yhgfr567uhgf}
\end{EQA}
We also write
\begin{EQA}
	\IF
	=
	\IF(\upsvs_{\GP}),
	\quad
	\IF_{\GP} 
	&=& 
	\IF_{\GP}(\upsvs_{\GP})
	=
	\IF + \GP^{2} .
\label{12093nvcfhyw0eiu}
\end{EQA}

\Subsection{I.i.d. of independent samples, effective sample size}
Another important issue is a possibility to relax the assumption of an i.i.d. of independent sample. 
Below we operate with the general likelihood, its structure does not need to be specified. 
In fact we can even proceed with just one observation.
However, in the case of an independent sample \( \Yv = (Y_{1},\ldots,Y_{n})^{\T} \in \R^{n} \), 
the function \( L(\upsv) \) is of additive structure; see \eqref{LuYiPsiGLM} below for the GLM case. 
The same applies to the information matrix 
\( \IF(\upsv) = - \nabla^{2} L(\upsv) = - \nabla^{2} \E L(\upsv) \); see \eqref{IFGusiphpPiuvT}. 
Moreover, if the \( Y_{i} \)'s are i.i.d. then \( \IF(\upsv) \) is proportional to \( n \).
Therefore, we use the value 
\( n = \| \IF^{-1} \|^{-1} \) as a proxy for the ``sample size''.
It can be called the \emph{effective sample size}.

\Subsection{Effective parameter dimension}
One more important issue is the parameter dimension \( \dimp \).
The classical theory assumes \( \dimp \) fixed and \( n \) large.
We aim at relaxing both conditions by allowing a large/huge/infinite parameter dimension
and a small \( n \).
It appears that all the results below rely on the so called 
\emph{effective dimension} \( \dimG \) defined as
\begin{EQA}
	\dimVG
	& \eqdef &
	\tr \bigl\{ \IFGP^{-1} \Var(\nabla L(\upsvs_{\GP})) \bigr\} .
\label{nswe45678lkjmnbvc}
\end{EQA}
This quantity replaces the original dimension \( \dimp \) and it can be small or moderate even for \( \dimp \)
infinite. 
One of the main intension of our study is to understand the range of applicability of 
the mentioned results in terms of the effective parameter dimension \( \dimG \) and the effective sample size \( n \).
It appears that the most of the results below apply under the condition \( \dimG \ll n \) which 
replaces the classical signal-to-noise relation: the effective number of parameters to be estimated is smaller in order
than the effective sample size.
The only one final result about a Gaussian approximation of the posterior requires 
a much stronger condition \( \dimL^{3}(\upsv) \ll n \) with 
\begin{EQA}
	\dimL(\upsv)
	& \eqdef &
	\tr\bigl\{ \IFGP^{-1}(\upsv) \, \IF(\upsv) \bigr\} .
\label{5cjf763hfyenefiwdhg}
\end{EQA}

\Subsection{Bias-variance trade-off}
Introducing a quadratic penalty \( \| \GP \upsv \|^{2}/2 \) or, equivalently, a Gaussian prior 
\( \ND(0,\GP^{-2}) \) allow to consider in a unified way the classical parametric situation 
with \( \dimp \) fixed and \( n \) large and the nonparametric situation with \( \dimp \) large or infinite.
The link is given by the value \( \dimG \) which smoothly decreases as \( \GP^{2} \) increases.
From the other side, an increase of the penalty yields an increase of the bias 
measured by the difference \( \upsvs_{\GP} - \upsvs \) in some norm.
For most of tasks in nonparametric statistics such as 
minimax risk estimation or model selection, etc,
one has to balance the variance term \( \dimG \) and the squared bias
yielding the optimal rate of estimation.

%


\Section{Setup. Main steps of study}
Now we briefly describe our setup and the main focus of analysis.

\Subsection{Stochastically linear smooth models}

Below 
we limit ourselves to a special class of \emph{stochastically linear smooth} (SLS) statistical models.  
The major condition on this class is that the stochastic component 
\( \zeta(\upsv) = L(\upsv) - \E L(\upsv) \) of the log-likelihood \( L(\upsv) \) is linear in parameter \( \upsv \).
We also assume that the expected log-likelihood is a smooth function of the parameter \( \upsv \).
This class includes popular Generalized Linear Models but it is much larger. 
In particular, by an extending the parameter space, one can consider 
many nonlinear models including nonlinear regression or 
nonlinear inverse problems as a special case of SLS; see \cite{Sp2019NLI}.
We also focus on the case of a quadratic penalization \( \pent_{\GP}(\upsv) = \| \GP \upsv \|^{2}/2 \)
or a Gaussian prior \( \ND(0,\GP^{-2}) \).
This would not affect the SLS conditions.
The assumption of stochastic linearity helps to avoid heavy tools of empirical process theory 
which is typically used in the analysis of pMLE \( \tilde{\upsv}_{\GP} \) or the posterior \( \upsv \cond \Yv \);
see e.g. \cite{vdV00}, \cite{GhVa2017}, \cite{nickl_2015}.
We only need some accurate deviation bounds for quadratic forms of the errors; see Section~\ref{Sdevboundgen} in the appendix.

Our aim is to establish possibly sharp and accurate results under realistic assumptions 
on a SLS model and the amount of data.
The study includes several steps.

\Subsection{Concentration of the pMLE}
The first step of our analysis is in establishing a concentration result 
for the pMLE \( \tilde{\upsv}_{\GP} \) defined by maximization of \( \LGP(\upsv) \).
If the expected log-likelihood \( \E \LGP(\upsv) \) is strictly concave and smooth in \( \upsv \) then 
\( \tilde{\upsv}_{\GP} \) well concentrates 
in a small elliptic vicinity \( \CA_{\GP} \) of the ``target'' \( \upsvs_{\GP} \) from \eqref{upssaruELuaiU}:
\begin{EQA}
	\P\bigl( \| \IF_{\GP}^{1/2} (\tilde{\upsv}_{\GP} - \upsvs_{\GP}) \| > \rrGP + \sqrt{2\xx}  \bigr)
	& \leq &
	3 \ex^{-\xx} ,
\label{uyhhjkigvbghjugfg}
\end{EQA}
where \( \rrGP^{2} \asymp \dimG \).
The result becomes sensible provided that \( \dimG \ll n \) with \( \nsize^{-1} \asymp \| \IF_{\GP}^{-1} \| \).

\Subsection{Fisher and Wilks expansions}
Having established the concentration of \( \tilde{\upsv}_{\GP} \in \CA_{\GP} \), 
we can restrict the analysis to this vicinity and 
 use the Taylor expansion of the penalized log-likelihood function 
\( \LGP(\upsv) \).
This helps to derive rather precise approximation of the values 
\( \tilde{\upsv}_{\GP} - \upsvs_{\GP} \) and \( \LGP(\tilde{\upsv}_{\GP}) - \LGP(\upsvs_{\GP}) \):
\begin{EQ}[rcl]
    \bigl\| \IF_{\GP}^{1/2} \bigl( \tilde{\upsv}_{\GP} - \upsvs_{\GP} \bigr) - \xivGP \bigr\|^{2}
    & \lesssim &
    \dltwb_{\GP} \, \bigl\| \xivGP \bigr\|^{2} ,
    \\
    \biggl| 
    \LGP(\tilde{\upsv}_{\GP}) - \LGP(\upsvs_{\GP}) 
    - \frac{1}{2} \bigl\| \xivGP \bigr\|^{2}
    \biggr|
    & \lesssim &
    \dltwb_{\GP} \, \bigl\| \xivGP \bigr\|^{2} \, ,
\label{3d3Af12DGi}
\end{EQ}
where \( \xivGP = \IF_{\GP}^{-1/2} \nabla \LGP(\upsvs_{\GP}) \).
The accuracy of approximation is controlled by the value \( \dltwb_{\GP} \) and the presented results 
require \( \dltwb_{\GP} \ll 1 \).
In typical examples \( \dltwb_{\GP} \asymp \sqrt{\dimG/n} \) and again we come to the condition \( \dimG \ll n \).
The first result in \eqref{3d3Af12DGi} about the pMLE \( \tilde{\upsv}_{\GP} \) will be referred to as 
\emph{the Fisher expansion}, while the second one about \( \LGP(\tilde{\upsv}_{\GP}) \) is called 
\emph{the Wilks expansion}.
These two expansions provide a finite sample analog of the asymptotic statements \eqref{2Lm2Lusg2sc2p}
and are informative even in the classical parametric situation.
In fact, under standard assumptions, the normalized score vector \( \xivGP \) is asymptotically normal 
\( \ND(0,\Sigma_{\GP}) \) with \( \Sigma_{\GP} = \IF_{\GP}^{-1/2} \VP^{2} \IF_{\GP}^{-1/2} \in \Matr_{\dimp} \) 
and \( \VP^{2} = \Var\bigl( \nabla L(\upsv) \bigr) \in \Matr_{\dimp} \).
Stochastic linearity implies that the matrix \( \VP^{2} \) does not depend on the point \( \upsv \).
If the model is correctly specified, then \( \Sigma_{\GP} \) approaches the identity as \( n \to \infty \),
and we obtain the classical results \eqref{2Lm2Lusg2sc2p}.
Note that the use of stochastic linearity allows to obtain much more accurate bounds than 
in \cite{SP2011} or \cite{SP2013_rough}. 

\Subsection{Concentration of the posterior}

The study of the posterior \( \upsv \cond \Yv \) also starts with a concentration result 
that claims that the dominating part of the posterior distribution
is located in an elliptic vicinity \( \CA \) of the pMLE \( \tilde{\upsv}_{\GP} \).
The set \( \CA \) has to be slightly larger than  \( \CA_{\GP} \).
A great advantage of stochastic linearity is that it allows to replace the log-likelihood \( \LGP(\upsv) \)
in the Bayes formula with the expected log-likelihood \( \E \LGP(\upsv) \).
Further one can use quadratic approximation of \( \E \LGP(\upsv) \)
within the vicinity \( \CA \) and concavity of \( \E \LGP(\upsv) \) outside of \( \CA \).
The concentration result for the posterior is much more involved than the concentration result for pMLE.
However, the required conditions are very similar.
In particular, we need a relation \( \dimL(\upsv) \ll n \), where \( \dimL(\upsv) \)
is the \emph{Laplace effective dimension} at \( \upsv \):
\begin{EQA}
	\dimL(\upsv)
	&=&
	\tr\bigl\{ \IF(\upsv) \IF_{\GP}^{-1}(\upsv) \bigr\} .
\label{cf73gfuyvg5t63jgu7refi}
\end{EQA}
Note that in typical situations \( \dimL(\upsv) \asymp \dimG \) for \( \upsv \) from a local vicinity of \( \upsvs_{\GP} \).

\Subsection{Gaussian approximation of the posterior}

The desired Gaussian approximation of the posterior is a relatively simple corollary 
of smoothness properties of the expected log-likelihood using the ideas of Laplace approximation
for the exponent of a smooth concave function.
Our main result claims that the posterior distribution for a SLS model can be well approximated by 
a Gaussian distribution \( \ND\bigl( \tilde{\upsv}_{\GP}, \IFt_{\GP}^{-1} \bigr) \), 
where \( \IFt_{\GP} = \IF_{\GP}(\tilde{\upsv}_{\GP}) \).
This is again a direct extension of the classical BvM result \eqref{PrF12utucYni}.
It is curious that the Gaussian approximation of the posterior does not rely on any CLT result,
it is a purely analytic result.
Moreover, one can establish a very careful approximation accuracy for centrally symmetric sets 
of order \( \dimt^{3}/n \), where \( \dimt = \dimL(\tilde{\upsv}_{\GP}) \).
However, the condition on critical dimension becomes much more strict: 
now we need \( \dimL^{3}(\upsv) \ll n \).

\medskip

We can summarize as follows:
all the classical statistical results established for linear models can be extended to SLS modeling.
For most of classical results in parametric statistics, one can state a finite sample SLS version.
For the majority of our results,
it is sufficient to ensure the relation \( \dimG \ll n \) or \( \dimL(\upsv) \ll n \).
For Gaussian approximation of the posterior, we need a stronger condition \( \dimL(\upsv) \ll n^{1/3} \).

\bigskip
The paper is organized as follows.
Section~\ref{SgenBounds} presents some important properties of the penalized MLE \( \tilde{\upsv}_{\GP} \) including 
concentration, Fisher and Wilks expansion, bias and risk. 
Section~\ref{SgenBvMs} describes the properties of the posterior \( \vupsv_{\GP} \cond \Yv \).
The main result describes the accuracy of Gaussian approximation for this posterior. 
Section~\ref{Saniligot} specifies the result to the case of logistic regression.
The appendix collects some very useful facts about Laplace approximation for the integral  
\( \int \ex^{f(\upsv)} \, d\upsv \) for a smooth concave function \( f \) as well as some deviation bounds for Gaussian 
and non-Gaussian quadratic forms.


\Chapter{Properties of the pMLE \( \tilde{\upsv}_{\GP} \)}
\label{SgenBounds}
This \chname collects general results about concentration and expansion of the pMLE.
The use of self-concordance conditions from Section~\ref{Slocalsmooth} helps to substantially improve 
the bounds from \cite{SpPa2019}.
%
We assume to be given a pseudo log-likelihood random function \( L(\upsv) \), \( \upsv \in \Ups \subseteq \R^{\dimp} \),
\( \dimp < \infty \).
Given a quadratic penalty \( \| \GP \upsv \|^{2}/2 \), define 
\begin{EQA}
	\LGP(\upsv) 
	&=& 
	L(\upsv) - \| \GP \upsv \|^{2}/2 .
\label{LGLum2Gu22}
\end{EQA}
Typical examples of choosing \( \GP^{2} \) are given in Section~\ref{Ssmoothprior}.
Consider the penalized MLE \( \tilde{\upsv}_{\GP} \) and its population counterpart \( \upsvs_{\GP} \)
\begin{EQA}[rcccl]
	\tilde{\upsv}_{\GP} 
	&=& 
	\argmax_{\upsv} \LGP(\upsv),
	\qquad
	\upsvs_{\GP} 
	&=& 
	\argmax_{\upsv} \E \LGP(\upsv). 
\label{tuGauLGususGE}
\end{EQA}
The corresponding Fisher information matrix \( \IF_{\GP}(\upsv) \) is given by 
\begin{EQA}[rcccl]
	\IF(\upsv) 
	&=& 
	- \nabla^{2} \E L(\upsv) ,
	\qquad
	\IF_{\GP}(\upsv) 
	&=& 
	- \nabla^{2} \E \LGP(\upsv) 
	= 
	\IF(\upsv) + \GP^{2} .
\label{e456yhgfr567uhgf}
\end{EQA}
We assume \( \IF_{\GP}(\upsv) \) to be positive definite for all considered \( \upsv \). 
By \( \DPGP(\upsv) \) we denote a positive symmetric matrix with \( \DPGP^{2}(\upsv) = \IF_{\GP}(\upsv) \), and
\( \IF_{\GP} = \IF_{\GP}(\upsvs_{\GP}) \), \( \DPGP = \IF_{\GP}^{1/2} \). 

\Section{Conditions}
\label{Scondgeneric}
Now we present our conditions.
The most important one is about linearity of the stochastic component 
\( \zeta(\upsv) = L(\upsv) - \E L(\upsv) = \LGP(\upsv) - \E \LGP(\upsv) \).

\begin{description}
    \item[\label{Eref} \( \bb{(\zeta)} \)]
      \textit{The stochastic component \( \zeta(\upsv) = L(\upsv) - \E L(\upsv) \) of the process \( L(\upsv) \) is linear in 
      \( \upsv \). We denote by \( \nabla \zeta \equiv \nabla \zeta(\upsv) \in \R^{\dimp} \) its gradient
      }.
  \end{description}

Below we assume some concentration properties of the stochastic vector \( \nabla \zeta \).
More precisely, we require that \( \nabla \zeta \) obeys the following condition;
see \eqref{PxivbzzBBroB3} of Theorem~\ref{Tdevboundgm}.
\begin{description}
\item[\label{EU2ref}\( \bb{(\nabla \zeta)} \)]
	\textit{Let \( \VP^{2} = \Var(\nabla \zeta) \) and 
	\( \DPGP^{2} = \DPGP^{2}(\upsvs_{\GP}) \).
	Then for any considered \( \xx > 0 \)
	}
\begin{EQA}
	\P\bigl( \| \DPGP^{-1} \nabla \zeta \| \geq \rr_{\GP}(\xx) \bigr)
	& \leq &
	3 \ex^{-\xx} ,
\label{2emxGPm12nz122}
\end{EQA}
where for \( \dimG = \tr(\DPGP^{-2} \VP^{2}) \) and \( \lambda_{\GP} = \| \DPGP^{-1} \VP^{2} \DPGP^{-1} \| \)
\begin{EQA}
	\rr_{\GP}(\xx)
	& \eqdef &
	\sqrt{\dimG} + \sqrt{2 \xx \lambda_{\GP}} \,  .
\label{34rtyghuioiuyhgvftid}
\end{EQA}
\end{description}

This condition can be effectively checked if the errors in the data exhibit sub-Gaussian or sub-exponential behaviour;
see Section~\ref{SdevboundnonGauss}.
The important value \( \dimG = \tr(\DPGP^{-2} \VP^{2}) \) can be called the \emph{effective dimension}; see \cite{SP2013_rough}.

We also assume that the deterministic part \( \E \LGP(\upsv) \) 
of the penalized log-likelihood is a concave function.
It can be relaxed using localization; see \cite{Sp2019NLI}.
\begin{description}
    \item[\label{LLref} \( \bb{(\mathcal{C}_{\GP})} \)]
      \textit{\( \Ups \) is an open and convex set in \( \R^{\dimp} \).
        The function \( \E \LGP(\upsv) \) is concave on \( \Ups \).        
      }
\end{description}

In Section~\ref{Spostconcentr} we consider a stronger condition of semi-concavity of \( \E L(\upsv) \). 
Further we will also need some smoothness conditions on the function \( f(\upsv) = \E L(\upsv) \).
The class of models satisfying the conditions \nameref{Eref}, 
\nameref{EU2ref}
with a smooth function \( f(\upsv) = \E L(\upsv) \) will be referred to as \emph{stochastically linear smooth} (SLS). 
This class includes regression, generalized linear models (GLM) and log-density models; see \cite{SpPa2019} or Section~\ref{SGBvM} later. 
However, this class is much larger.
For instance, nonlinear regression and nonlinear inverse problems can be adapted to the SLS framework 
by an extension of the parameter space; see \cite{Sp2019NLI}.

\Section{Concentration of the pMLE \( \tilde{\upsv}_{\GP} \)}
\label{SgenpMLE}
This section discusses some concentration properties of the pMLE \( \tilde{\upsv}_{\GP} = \argmax_{\upsv} \LGP(\upsv) \).

Given \( \xx \) and \( \rr_{\GP} = \rr_{\GP}(\xx) \) from \eqref{34rtyghuioiuyhgvftid}, define for some \( \amax < 1 \) the set \( \UVz_{\GP} \) by
\begin{EQA}
	\UVz_{\GP}
	& \eqdef &
	\bigl\{ \uv \colon \| \DPGP \uv \| \leq \amax^{-1} \rr_{\GP} \bigr\} .
\label{rm1ucDGu2r0DGu}
\end{EQA}
The result of this section states the concentration properties of the pMLE \( \tilde{\upsv}_{\GP} \)
in the local vicinity \( \CA_{\GP} \) of \( \upsvs_{\GP} \) of the form 
\begin{EQA}
	\CA_{\GP}
	& \eqdef &
	\upsvs_{\GP} + \UVz_{\GP}
	=
	\bigl\{ \upsv = \upsvs_{\GP} + \uv \colon \uv \in \UVz_{\GP} \bigr\}
	\subseteq \Upsd .
\label{UsuvcruvusGiV0}
\end{EQA}
Local Gateaux-regularity of \( f(\upsv) = \E L(\upsv) \) within \( \CA_{\GP} \)
will be measured by the error of the second order Taylor approximation
\begin{EQ}[rcl]
	\dltw_{3}(\upsv,\uv) 
	&=& 
	f(\upsv + \uv) - f(\upsv) - \langle \nabla f(\upsv), \uv \rangle 
	- \frac{1}{2} \langle \nabla^{2} f(\upsv), \uv^{\otimes 2} \rangle , 
	\\
	\dltwd_{3}(\upsv,\uv) 
	&=&
	\langle \nabla f(\upsv + \uv), \uv \rangle - \langle \nabla f(\upsv), \uv \rangle 
	- \langle \nabla^{2} f(\upsv), \uv^{\otimes 2} \rangle \, .
\label{dltw3vufuv12f2g}
\end{EQ}
More precisely, 
define
\begin{EQA}
    \dltwb_{\GP}
    & \eqdef &
    \sup_{\uv \in \UVz_{\GP}}
    \frac{2 |\dltw_{3}(\upsvs_{\GP},\uv)|}{\| \DPGP \uv \|^{2}} \,\, ,
    \qquad
    \dltwbd_{\GP}
    \eqdef 
    \sup_{\uv \in \UVz_{\GP}} \frac{2 |\dltwd_{3}(\upsvs_{\GP},\uv)|}{\| \DPGP \uv \|^{2}} \,\, . 
\label{dtb3u1DG2d3GP}
\end{EQA}
The quantities \( \dltwb_{\GP} \) and \( \dltwbd_{\GP} \) can be effectively bounded under smoothness
conditions \nameref{LL3tref} or \nameref{LLtS3ref} given in Section~\ref{Slocalsmooth}.
Under \nameref{LL3tref} at \( \upsv = \upsvs_{\GP} \) with \( \HLG^{2}(\upsvs_{\GP}) = \DPGP^{2} \)
and \( \rr = \rrGP \), by Lemma~\ref{LdltwLa3t}, it holds for a small constant \( \dltwu_{3} \)
\begin{EQA}
	\dltwbd_{\GP}
	& \leq &
	\dltwu_{3} \, \amax^{-1} \rrGP \, ,
	\qquad
	\dltwb_{\GP}
	\leq 
	\dltwu_{3} \, \amax^{-1} \rrGP/ 3 .
\label{swdy6fwqd6qtxcvbdyfdtw}
\end{EQA}
Furthermore, under \nameref{LLtS3ref}, the same bounds apply with \( \dltwu_{3} = \hmax_{3} n^{-1/2} \); see Lemma~\ref{LdltwLaGP}.

\begin{proposition}
\label{PconcMLEgenc}
Suppose
\nameref{Eref},
\nameref{EU2ref},
and 
\nameref{LLref}.
Let also 
\begin{EQA}
	1 - \amax - \dltwbd_{\GP}
	& \geq &
	0 ;
\label{rrm23r0ut3u}
\end{EQA}
see \eqref{dtb3u1DG2d3GP} and \eqref{rm1ucDGu2r0DGu}.
Then \( \tilde{\upsv}_{\GP} \in \CA_{\GP} \) on a set \( \Omega(\xx) \) with \( \P\bigl( \Omega(\xx) \bigr) \geq 1 - 3 \ex^{-\xx} \), 
i.e.
\begin{EQA}
	\| \DPGP (\tilde{\upsv}_{\GP} - \upsvs_{\GP}) \|  
	& \leq &
	\amax^{-1} \rrGP \, . 
\label{rhDGtuGmusGU0}
\end{EQA}
\end{proposition}

\begin{proof}
By \nameref{EU2ref}, on a the random set \( \Omega(\xx) \) with 
\( \P(\Omega(\xx)) \geq 1 - 3 \ex^{-\xx} \), it holds \( \| \DPGP^{-1} \nabla \zeta \| \leq \rr_{\GP} \).
Now the result follows from Proposition~\ref{Pconcgeneric} with 
\( \fs(\upsv) = \E \LGP(\upsv) \), \( \fn(\upsv) = \LGP(\upsv) \), \( \rrn = \amax^{-1} \rrGP \), and \( \Av = \nabla \zeta \).
\end{proof}

\ifsupnorm{

\Section{A general condition of the stochastic gradient}
\label{Snanonorm}
The results about pMLE \( \tilde{\upsv}_{\GP} \) strongly rely on the condition \nameref{EU2ref} 
describing the concentration set of \( \nabla \zeta \) in the form
\( \| \DPGP^{-1} \nabla \zeta \| \leq \rups \).
Here we discuss some possible extensions using other norms than \( \ell_{2} \).

\begin{description}
\item[\label{EU0ref}\( \bb{(\nabla \zeta)_{\nano}} \)]
	\textit{There exists a norm \( \| \cdot \|_{\nano} \) and for any \( \xx > 0 \) a radius \( \rups = \rups(\xx) \) 
	such that}
\begin{EQA}
	\P\bigl( \| \nabla \zeta \|_{\nano} > \rups
	\bigr)
	& \leq &
	\CONSTnabl \ex^{-\xx} .
\label{PnzniUV0C0emx}
\end{EQA}
\end{description}

Previous results were stated for \( \| \nabla \zeta \|_{\nano} = \| \DPGP^{-1} \nabla \zeta \| \).
Another example of choosing the norm \( \| \cdot \|_{\nano} \) is given in Section~\ref{Snanonorm}.

Let \( \| \cdot \|_{\dual} \) be the corresponding dual norm for \( \| \cdot \|_{\nano} \):
\begin{EQA}
	\| \uv \|_{\dual}
	& \eqdef &
	\sup_{\| \wv \|_{\nano} \leq 1} \bigl\langle \wv, \uv \bigr\rangle .
\label{uvduswna1wu}
\end{EQA}
Given \( \xx \) and \( \rr_{\GP} = \rr_{\GP}(\xx) \), define for some \( \amax < 1 \) the set \( \UVz_{\GP} \) by
\begin{EQA}
	\UVz_{\GP}
	& \eqdef &
	\bigl\{ \uv \colon \| \DPGP \uv \|^{2} \leq \amax^{-1} \rr_{\GP} \| \uv \|_{\dual} \bigr\} .
\label{UzGuDGu2rm1na}
\end{EQA}
The next result states the concentration properties of the pMLE \( \tilde{\upsv}_{\GP} \)
in the local vicinity \( \CA_{\GP} \) of \( \upsvs_{\GP} \) of the form 
\begin{EQA}
	\CA_{\GP}
	& \eqdef &
	\upsvs_{\GP} + \UVz_{\GP}
	=
	\bigl\{ \upsv = \upsvs_{\GP} + \uv \colon \uv \in \UVz_{\GP} \bigr\}
	\subseteq \Upsd .
\label{UsuvcruvusGiV0na}
\end{EQA}
The quantities \( \dltwb_{\GP} \) and \( \dltwb_{\GP} \) are defined by \eqref{dtb3u1DG2d3GP}
with the set \( \UVz_{\GP} \) from \eqref{UzGuDGu2rm1na}.

\begin{proposition}
\label{PconcMLEgen}
Suppose
\nameref{Eref},
\nameref{EU0ref},
and 
\nameref{LLref}.
Let also \( \dltwbd_{\GP} \) from \eqref{dtb3u1DG2d3GP} and
\( \amax \) from \eqref{UzGuDGu2rm1na} fulfill
\begin{EQA}
	1 - \amax - \dltwbd_{\GP}
	& \geq &
	0 .
\label{rrm23r0ut3una}
\end{EQA}
Then on a set \( \Omega(\xx) \) with \( \P\bigl( \Omega(\xx) \bigr) \geq 1 - \CONSTnabl \ex^{-\xx} \), 
it holds \( \tilde{\upsv}_{\GP} \in \CA_{\GP} \), that is,
\begin{EQA}
	\tilde{\upsv}_{\GP} - \upsvs_{\GP}  
	& \in &
	\UVz_{\GP} \, . 
\label{rhDGtuGmusGU0na}
\end{EQA}
\end{proposition}


\begin{remark}
\label{RIF0DG2PPC}
If \( \DPU^{2} \leq \DPGP^{2} \)
then \( \tilde{\upsv}_{\GP} \in \CA_{\DPU} = \upsvs_{\GP} + \UVz_{\DPU} \) on \( \Omega(\xx) \) with the local set 
\begin{EQA}
	\UVz_{\DPU}
	& \eqdef &
	\bigl\{ \upsv \colon \| \DPU \uv \|^{2} \leq \amax^{-1} \rr_{\GP} \| \uv \|_{\dual} \bigr\} .
\label{dta3u1DG2srSAI0na}
\end{EQA}
\end{remark}

\begin{proof}[Proof of Proposition~\ref{PconcMLEgen}]
Below we restrict ourselves to the random set \( \Omega(\xx) \) and assume that \( \| \nabla \zeta \|_{\nano} \leq \rr_{\GP} \).
Then on \( \Omega(\xx) \) for any \( \uv \)
\begin{EQA}
	\bigl| \langle \nabla \zeta, \uv \rangle \bigr|
	& \leq &
	\sup_{\| \wv \|_{\nano} \leq \rr_{\GP}} \bigl| \langle \wv, \uv \rangle \bigr|
	=
	\rr_{\GP} \| \uv \|_{\dual} \, .
\label{LLoDGm1nzu}
\end{EQA}
If \( \uv \not\in \UVz_{\GP} \), then by definition
\eqref{UzGuDGu2rm1na} it holds \( \rr_{\GP} \| \uv \|_{\dual} \leq \amax \| \DPGP \uv \|^{2} \). 
Further we proceed as in the case of condition \nameref{EU2ref}.
\end{proof}

\Section{A bound for a sup-norm}
\label{Emaxjlulsq}
Now we consider the case with bounding the vector \( \nabla \zeta \) in some \( \ell_{\infty} \)-norm.
Suppose that the components of \( \nabla \zeta \) are sub-Gaussian r.v.'s satisfying 
\begin{EQA}
	\P\bigl( \bigl| (\nabla \zeta)_{j} \bigr| > \zq_{j}(\xx) \bigr)
	& \leq &
	\CONSTi_{0} \ex^{-\xx} ,
\label{C0emxPnzjzj}
\end{EQA} 
where \( \zq_{j}(\xx) \) is typically \( \vp_{j} \sqrt{\xx} \) for \( \vp_{j}^{2} = \Var( \nabla \zeta)_{j} \).
Let also \( \priorh \) be a positive measure on the discrete set \( \bigl\{ 1,\ldots,\dimp \bigr\} \).
Given \( \xx \), define \( \xx_{j} = \xx  - \log \priorh_{j} \).
Using the Bonferroni device, we can bound
\begin{EQA}
	\P\biggl( \max_{j} \frac{\bigl| (\nabla \zeta)_{j} \bigr|}{\zq_{j}(\xx_{j})} > 1 \biggr)
	& \leq &
	\sum_{j} \P\bigl( \bigl| (\nabla \zeta)_{j} \bigr| > \zq_{j}(\xx_{j}) \bigr)
	\leq 
	\CONSTi_{0} \sum_{j} \ex^{-\xx_{j}}
	=
	\CONSTi_{0} \ex^{-\xx} .
\label{C0C0mxxsjPnzj}
\end{EQA}
Introduce the diagonal matrix \( \zqv = \diag\bigl( \zq_{j}(\xx_{j}) \bigr) \).
Then \nameref{EU0ref} is fulfilled with
\( \| \uv \|_{\nano} = \| \zqv^{-1} \uv \|_{\infty} \) and \( \rr_{\GP}  = 1 \). 
Further,
\begin{EQA}
	\| \uv \|_{\dual}
	&=&
	\sup_{\| \wv \|_{\nano} \leq 1} \bigl\langle \uv, \wv \bigr\rangle
	=
	\sup_{\| \zqv^{-1} \wv \|_{\infty} \leq 1} \bigl\langle \zqv \uv,  \zqv^{-1} \wv \bigr\rangle
	=
	\| \zqv \uv \|_{1} \, .
\label{zu1szm1in1zuzm1}
\end{EQA}
Suppose now that there exists a diagonal matrix \( \DPU^{2} \) such that \( \DPGP^{2} \geq \DPU^{2} \).
Consider the set \( \UVz_{\DPU} \) from \eqref{dta3u1DG2srSAI0na}.
This set \( \UVz_{\DPU} \) can be described explicitly.
Indeed, with \( a_{j} = (2\amax)^{-1} \zq_{j}(\xx_{j})/\IF_{j} \)
\begin{EQA}
	\UVz_{\DPU}
	&=&
	\biggl\{ \uv \colon \sum_{j} \IF_{j} u_{j}^{2} \leq \amax^{-1} \sum_{j} \zq_{j} |u_{j}| \Bigr\}
	=
	\biggl\{ \uv \colon \sum_{j} \IF_{j} (|u_{j}| - a_{j})^{2} \leq \sum_{j} \IF_{j} a_{j}^{2} \Bigr\}
	\\
	&=&
	\Bigl\{ \uv \colon \bigl\| \DPU (|\uv| - \av) \bigr\|^{2} \leq (2 \amax)^{-1} \| \DPU^{-1} \zqv \|^{2} \Bigr\} .
\label{WGuIduFjj2aj}
\end{EQA}
Therefore, the set \( \UVz_{\DPU} \) has an elliptic shape in each quadrant of \( \R^{\dimp} \) with the center shifted 
from zero to the point \( \uv_{m} \) from this quadrant with \( |\uv_{m}| = \av \).
For a uniform measure \( \priorh \) with \( \priorh_{j} = 1/\dimp \) and \( \zq_{j}(\xx) = \vp_{j} \sqrt{\xx} \), we derive
\begin{EQA}
	\| \DPU^{-1} \zqv \|^{2}
	&=&
	\sum_{j} \frac{\zq_{j}^{2}(\xx+\log \dimp)}{\DPU_{j}^{2}}
	=
	(\xx+\log \dimp) \sum_{j} \frac{\vp_{j}^{2}}{\DPU_{j}^{2}} \, .
\label{HVm1zqv2sjxlp}
\end{EQA} 
The sum \( \sum_{j} \vp_{j}^{2}/\DPU_{j}^{2} \) has the same meaning as the effective dimension 
\( \dimG = \tr \bigl\{ \Var(\nabla \zeta) \IF_{\GP}^{-1} \bigr\} \).
It will be assumed that this value is small relatively to each \( \DPU_{j}^{2} \).
This assures the desired concentration of the pMLE \( \tilde{\upsv}_{\GP} \).

}{}

\Section{Fisher and Wilks expansions}
\label{SFiWiexs}
This section presents some finite sample results about the behavior of the penalized MLE 
\( \tilde{\upsv}_{\GP} \) and the excess \( \LGP(\tilde{\upsv}_{\GP}) - \LGP(\upsvs_{\GP}) \).
Proposition~\ref{PconcMLEgenc} 
states the concentration properties of \( \tilde{\upsv}_{\GP} \) around \( \upsvs_{\GP} \).
Now we show that this concentration can be used to establish a version of the Fisher expansion for 
the estimation error \( \tilde{\upsv}_{\GP} - \upsvs_{\GP} \) and the Wilks expansion for the excess 
\( L_{\GP}(\tilde{\upsv}_{\GP}) - L_{\GP}(\upsvs_{\GP}) \).

\begin{theorem}
\label{TFiWititG}
Assume the conditions of 
Proposition~\ref{PconcMLEgenc}
with \( \amax = 2/3 \).
Then on \( \Omega(\xx) \!\) 
\begin{EQ}[rcl]
    2 \LGP(\tilde{\upsv}_{\GP}) - 2 \LGP(\upsvs_{\GP}) 
    - \bigl\| \DPGP^{-1} \nabla \zeta \bigr\|^{2}
    & \leq &
    \frac{\dltwb_{\GP}}{1 - \dltwb_{\GP}} \bigl\| \DPGP^{-1} \nabla \zeta \bigr\|^{2} \, ,
    \\
    2 \LGP(\tilde{\upsv}_{\GP}) - 2 \LGP(\upsvs_{\GP}) 
    - \bigl\| \DPGP^{-1} \nabla \zeta \bigr\|^{2}
    & \geq &
    - \dltwb_{\GP} \bigl\| \DPGP^{-1} \nabla \zeta \bigr\|^{2} .
\label{3d3Af12DGttG}
\end{EQ}
Also
\begin{EQ}[rcl]
    \bigl\| \DPGP \bigl( \tilde{\upsv}_{\GP} - \upsvs_{\GP} \bigr) - \DPGP^{-1} \nabla \zeta \bigr\|^{2}
    & \leq &
    \frac{3 \dltwb_{\GP}}{(1 - \dltwb_{\GP})^{2}} \, \bigl\| \DPGP^{-1} \nabla \zeta \bigr\|^{2} \, ,
    \\
    \bigl\| \DPGP \bigl( \tilde{\upsv}_{\GP} - \upsvs_{\GP} \bigr) \bigr\|
    & \leq &
    \frac{1 + \sqrt{2 \dltwb_{\GP}}}{1 - \dltwb_{\GP}} \, \bigl\| \DPGP^{-1} \nabla \zeta \bigr\| \, .
\label{DGttGtsGDGm13rG}
\end{EQ}
\end{theorem}

\begin{proof}
The result follows from Proposition~\ref{PFiWigeneric} similarly to Proposition~\ref{PconcMLEgenc}.
\end{proof}


\Section{Effective sample size and critical dimension in pMLE}
\label{ScritdimMLE}
This section discusses the important question of the critical parameter dimension 
still ensuring the validity of the presented results.
A very important feature of our results is their dimension free and coordinate free form.
The true parametric dimension \( \dimp \) can be very large, it does not show up in the error terms.
Neither do we use any spectral decomposition or sequence space structure, in particular, we do not require that 
the Fisher information matrix \( \IF \) and the penalty matrix \( \GP^{2} \) are diagonal or can be jointly diagonalized. 
The results are stated for the general data \( \Yv \) and a quasi log-likelihood function.
In particular, we do not assume independent or progressively dependent observations and additive structure 
of the log-likelihood. 
The \emph{effective sample size} \( \nsize \) can be defined via the smallest eigenvalue of the matrix 
\( \IF_{\GP} = \DPGP^{2} = - \nabla^{2} \E \LGP(\upsvs_{\GP}) \):
\begin{EQA}
	n^{-1}
	& \eqdef &
	\| \IF_{\GP}^{-1} \| .
\label{ikdfu4rdu7yr47gfurkj}
\end{EQA}
Our results apply as long as this value is sufficiently small. 
%
Below in 
\ifapp{\Chname \ref{Saniligot} and \Chname \ref{SGBvM}}
{\Chname \ref{SGBvM}}
we will see that such defined value is closely related to the sample size \( \nsize \).

For the concentration result of Proposition~\ref{PconcMLEgenc}
we need the basic conditions \nameref{Eref} and \nameref{LLref}.
Further, \nameref{EU2ref} identifies the radius \( \rrGP \) of the local vicinity \( \CA_{\GP} \).
The final critical condition is given by \eqref{rrm23r0ut3u}.
Essentially it says that the values \( \dltwb_{\GP} \) and \( \dltwbd_{\GP} \) are significantly smaller than 1. 
Under \nameref{LLtS3ref}, \( \dltwbd_{\GP} 
\leq \hmax_{3} \, \amax^{-1} \rr_{\GP} \, \nsize^{-1/2} \); see Lemma~\ref{LdltwLaGP}.
So,  
\eqref{rrm23r0ut3u} means \( \rr_{\GP}^{2} \ll \nsize \).
Moreover, definition \eqref{2emxGPm12nz122} of \( \rr_{\GP} \) yields that 
\( \rr_{\GP}^{2} \asymp \tr(\DPGP^{-2} \VP^{2}) = \dimVG \), where \( \dimVG \) is the \emph{effective dimension} of the problem.
We conclude that the main properties of the pMLE \( \tilde{\upsv}_{\GP} \) 
are valid under the condition \( \dimVG \ll n \) meaning sufficiently many observations 
per effective number of parameters.

%


\Section{The use of \( \DPGPt^{2} \) instead of \( \DPGP^{2} \)}
\label{SDPGPDPGPt}
The penalized information matrix 
\( \DPGP^{2} = \DPGP^{2}(\upsvs_{\GP}) = - \nabla^{2} \E \LGP(\upsvs_{\GP}) \) 
plays an important role in our results.
In particular, \( \DPGP \) describes the shape of the concentration set \( \CA_{\GP} = \upsvs_{\GP} + \UVz_{\GP} \).
However, this matrix is not available as it involves the unknown point \( \upsvs_{\GP} \).
If the matrix function \( \IF(\upsv) \) is locally constant in \( \CA_{\GP} \), one can replace \( \upsvs_{\GP} \) 
with its estimate \( \tilde{\upsv}_{\GP} \).
Variability of \( \IF(\upsv) \), or, equivalently, \( \IF_{\GP}(\upsv) = \IF(\upsv) + \GP^{2} \) 
can be measured 
under the Fr\'echet smoothness of \( f(\upsv) = \E \LGP(\upsv) \) 
by the value \( \dltwbss_{\GP} \) from \eqref{jcxhydtferyufgy7tfdsy7ft}
\( \upsv = \upsvs_{\GP} \), \( \HLG(\upsv) = \DPGP \), and \( \rr = \amax^{-1} \rrGP \).

\begin{proposition}
\label{PDPGPtDPGP}
Assume the conditions of 
Proposition~\ref{PconcMLEgenc}
and let \( \dltwbss_{\GP} \leq 1/2 \); see \eqref{jcxhydtferyufgy7tfdsy7ft}.
The random matrix \( \DPGPt^{2} = \IF_{\GP}(\tilde{\upsv}_{\GP}) \) fulfills on \( \Omega(\xx) \) 
for any \( \uv \in \R^{\dimp} \)
\begin{EQ}[rcccl]
    \bigl\| \DPGP^{-1} \DPGPt^{2} \, \DPGP^{-1} - \Id_{\dimp} \bigr\|
    & \leq &
    \dltwbss_{\GP} \, ,
    \qquad
    \bigl\| \DPGP \, \DPGPt^{-2} \, \DPGP - \Id_{\dimp} \bigr\|
    & \leq &
    \frac{\dltwbss_{\GP}}{1 - \dltwbss_{\GP}} \, ,
    \qquad
\label{DPGPm1Cd3rG}
	\\
	(1 - \dltwbss_{\GP}) \, \| \DPGP \uv \|^{2}
	& \leq &
	\| \DPGPt \uv \|^{2}
	\leq 
	(1 + \dltwbss_{\GP}) \, \| \DPGP \uv \|^{2} .
\label{jamnjhgfrderr5t6784eds}
\end{EQ}
\end{proposition}

\begin{proof}
The value \( \tilde{\upsv}_{\GP} - \upsvs_{\GP} \) belongs to \( \UVz_{\GP} \) on \( \Omega(\xx) \) and 
\eqref{DPGPm1Cd3rG} follows from \eqref{ghdrd324ee4ew222w3ew}.
\end{proof}

\Section{Smoothness and bias}
\label{Ssmoothbias}
Due to 
{Proposition~\ref{PconcMLEgenc}}, 
the penalized MLE  \( \tilde{\upsv}_{\GP} \) is in fact an estimator of the vector \( \upsvs_{\GP} \).
However, \( \upsvs_{\GP} \) depends on penalization which introduces some bias.
This section discusses whether one can use \( \tilde{\upsv}_{\GP} \) for estimating the underlying truth \( \upsvs \)
defined as the maximizer of the expected log-likelihood:
\( \upsvs = \argmax_{\upsv} \E L(\upsv) \).
First we describe the bias \( \biasv_{\GP} = \upsvs_{\GP} - \upsvs \) induced by penalization.
It is important to mention that the previous results about the properties of the pMLE \( \tilde{\upsv}_{\GP} \)
require strong concavity of the expected log-likelihood function \( \E \LGP(\upsv) \) at least in a vicinity 
of the point \( \upsvs_{\GP} \). 
In some sense, this strong concavity is automatically forced 
by the penalizing term in the definition of \( \upsvs_{\GP} \).
However, the underlying truth \( \upsvs = \argmax_{\upsv} \E L(\upsv) \) is the maximizer of the non-penalized 
expected log-likelihood,
and the corresponding Hessian \( \IF(\upsvs) = - \nabla^{2} \E L(\upsvs) \) can degenerate. 
This makes evaluation of the bias more involved.
To bypass this situation, we assume later in this section 
that the Hessian \( \nabla^{2} \E \LGP(\upsv) \) cannot change much in a reasonably large vicinity of \( \upsvs \).
This allows to establish an accurate quadratic approximation of \( f(\upsv) \) and to evaluate the bias
\( \biasv_{\GP} = \upsvs_{\GP} - \upsvs \).

Define \( \DPfGP \) by \( \DPfGP^{2} = \IF_{\GP}(\upsvs) \);
cf. \( \DPGP^{2} = \IF_{\GP}(\upsvs_{\GP}) \). 
Let also \( \QP \) be a symmetric matrix satisfying  \( \QP^{2} \leq \DPfGP^{2} \).
Typical examples include \( \QP^{2} = \IF_{\GP} \), \( \QP^{2} = \IF(\upsvs) \), and \( \QP^{2} = \nsize \Id_{\dimp} \).
Later we bound the norm \( \| \QP \biasv_{\GP} \| \).
Introduce two vectors 
\begin{EQA}
	\Av_{\GP} 
	&=& 
	- \GP^{2} \upsvs,
	\qquad
	\av_{\GP} = \QP \, \DPfGP^{-2} \Av_{\GP} = - \QP \, \DPfGP^{-2} \, \GP^{2} \upsvs \, .
\label{gv8reyuo0hdfyt3hjrfi9}
\end{EQA}
Also denote \( \rru = \amax^{-1} \| \av_{\GP} \| = \amax^{-1} \| \QP \, \DPfGP^{-2} \, \GP^{2} \upsvs \| \) for \( \amax = 2/3 \) and 
\begin{EQA}
	\dltwbs_{\GP} 
	& \eqdef &
	\sup_{\uv \colon \| \QP \uv \| \leq \rru} \,\, 
	\| \DPfGP^{-1} \, \IF_{\GP}(\upsvs + \uv) \, \DPfGP^{-1} - \Id_{\dimp} \| ;
\label{ghdrd324ee4ew222}
\end{EQA}
cf. \eqref{jcxhydtferyufgy7tfdsy7ft} and \eqref{ghdrd324ee4ew222w3ew} for \( f(\upsv) = \E \LGP(\upsv) \).
Note that the definition of \( \dltwbss_{\GP} \) in Proposition~\ref{PDPGPtDPGP} uses another \( \rrn = \amax^{-1} \rrGP \),
therefore, different notation.
Proposition~\ref{Pbiasgeneric} yields the following result.

\begin{proposition}
\label{Lvarusetb} 
Let \( \DPfGP^{2} = \IF_{\GP}(\upsvs) \),
\( \Av_{\GP} \) and \( \av_{\GP} \) be given by \eqref{gv8reyuo0hdfyt3hjrfi9},
\( \amax \leq 2/3 \), and \( \rru = \amax^{-1} \| \av_{\GP} \| \).
Let also \( \dltwbs_{\GP} \leq 1/3 \); see \eqref{ghdrd324ee4ew222}.
Then the bias \( \biasv_{\GP} = \upsvs_{\GP} - \upsvs \) fulfills
\begin{EQA}[rcccl]
	\| \QP \, \biasv_{\GP} \|
	& \leq &
	\frac{1}{1 - \dltwbs_{\GP}} \, \| \av_{\GP} \| 
	&=&
	\frac{1}{1 - \dltwbs_{\GP}} \, \| \QP \, \DPfGP^{-2} \, \GP^{2} \upsvs \|\, ,
\label{11ma3eaelDebQ}
	\\
	\| \QP (\biasv_{\GP} - \DPfGP^{-2} \Av_{\GP}) \|
	& \leq &
	\frac{\dltwbs_{\GP}}{1 - \dltwbs_{\GP}} \, \| \av_{\GP} \| 
	&=&
	\frac{\dltwbs_{\GP}}{1 - \dltwbs_{\GP}} \, \| \QP \, \DPfGP^{-2} \, \GP^{2} \upsvs \| .
\label{11ma3eaelDeb}
\end{EQA}
\end{proposition}



As special cases of \eqref{11ma3eaelDebQ}, consider \( \QP = \DPfGP \) and \( \QP = \Id_{\dimp} \).

\begin{corollary}
\label{CbiasGPb}
Assume the conditions of Proposition~\ref{Lvarusetb}. 
Then 
\begin{EQA}[rcccl]
	\| \DPfGP \, \biasv_{\GP} \|
	& \leq &
	\frac{1}{1 - \dltwbs_{\GP}} \, \| \DPfGP^{-1} \, \GP^{2} \upsvs \| ,
	\qquad
	\| \biasv_{\GP} \|
	& \leq &
	\frac{1}{1 - \dltwbs_{\GP}} \, \| \DPfGP^{-2} \, \GP^{2} \upsvs \| .
\label{0mkvhgjnrw3dfwe3u8gty}
\end{EQA}
The same bounds apply with \( \DPGP^{2} = \IF_{\GP}(\upsvs_{\GP}) \) in place of \( \DPfGP^{2} = \IF_{\GP}(\upsvs) \).
\end{corollary}
The last statement follows from Remark~\ref{Rbiasgeneric}.

\Section{Loss and risk of the pMLE}
\label{Slossrisksb}
Now we combine the previous results about the stochastic term \( \tilde{\upsv}_{\GP} - \upsvs_{\GP} \)
and the bias term \( \biasv_{\GP} = \upsvs_{\GP} - \upsvs \) to obtain the sharp bounds 
on the loss and risk of the pMLE \( \tilde{\upsv}_{\GP} \).

%
%

\begin{theorem}
\label{PFiWibias}
Assume the conditions of Proposition~\ref{PconcMLEgenc} and \ref{Lvarusetb}.
Then on \( \Omega(\xx) \) with \( \P\bigl( \Omega(\xx) \bigr) \geq 1 - 3 \ex^{-\xx} \),
it holds with \( \xiv_{\GP} = \DPGP^{-1} \, \nabla \zeta \), \( \rr_{\GP} \) from \eqref{34rtyghuioiuyhgvftid},
and \( \nsize^{-1} = \| \DPGP^{-2} \| \)
\begin{EQA}[rcccl]
	\| \DPGP \, (\tilde{\upsv}_{\GP} - \upsvs) \|
	& \leq &
	\frac{1 + \sqrt{2 \dltwb_{\GP}}}{1 - \dltwb_{\GP}} \| \xiv_{\GP} \|
	+ \frac{\| \DPGP^{-1} \GP^{2} \upsvs \|}{1 - \dltwbs_{\GP}}  
	& \leq &
	\frac{1 + \sqrt{2 \dltwb_{\GP}}}{1 - \dltwb_{\GP}} \rr_{\GP}
	+ \frac{\| \DPGP^{-1} \GP^{2} \upsvs \|}{1 - \dltwbs_{\GP}} \, ,
	\qquad \qquad
\label{rG1r1d3GbvGb}
	\\
	\| \tilde{\upsv}_{\GP} - \upsvs \|
	& \leq &
	\frac{1 + \sqrt{2 \dltwb_{\GP}}}{\sqrt{\nsize} \, (1 - \dltwb_{\GP})} \| \xiv_{\GP} \| 
	+ \frac{\| \DPGP^{-2} \GP^{2} \upsvs \|}{1 - \dltwbs_{\GP}} 
	& \leq &
	\frac{3 \, \rrGP}{\sqrt{\nsize}} + 3 \, \| \DPGP^{-2} \GP^{2} \upsvs \| \, .
\label{g25re9fjfregdndg}
\end{EQA}
\end{theorem}

\begin{proof}
Let \( \Omega(\xx) \) be the random set from \nameref{EU2ref} on which
with \( \| \xiv_{\GP} \| \leq \rrGP \).
It follows from \eqref{DGttGtsGDGm13rG} of Theorem~\ref{TFiWititG} 
that on \( \Omega(\xx) \) 
with \( \biasv_{\GP} = \upsvs_{\GP} - \upsvs \)
\begin{EQA}
	\bigl\| \DPGP\bigl( \tilde{\upsv}_{\GP} - \upsvs \bigr) + \DPGP \biasv_{\GP} \bigr\|
	& \leq &
	\frac{1 + \sqrt{2 \dltwb_{\GP}}}{1 - \dltwb_{\GP}} \, \bigl\| \xiv_{\GP} \bigr\| \, .
\label{2biGe233GDttGQ}
\end{EQA}
This and \eqref{0mkvhgjnrw3dfwe3u8gty} imply \eqref{rG1r1d3GbvGb}.
\end{proof}


Now we state the results about the risk of the pMLE \( \tilde{\upsv}_{\GP} \).
To avoid technical burden, 
we fix a large \( \xx \), \( \rr_{\GP} = \rr_{\GP}(\xx) \), and exclude an event 
\( \bigl\{ \| \xiv_{\GP} \| > \rr_{\GP} \bigr\} \) 
having an exponentially small probability; see 
condition \eqref{2emxGPm12nz122} of Proposition~\ref{PconcMLEgenc}.

\begin{theorem}
\label{PFiWi2risk}
Assume the conditions of Proposition~\ref{PconcMLEgenc} and Proposition~\ref{Lvarusetb}.
Then for a set \( \Omega(\xx) \) with \( \P\bigl( \Omega(\xx) \bigr) \geq 1 - 3 \ex^{-\xx} \),
it holds with \( \dimVG = \tr(\DPGP^{-2} \, \VP^{2}) \)
\begin{EQA}
	&& \nquad
	\bigl\| \E \, \bigl\{ \DPGP (\tilde{\upsv}_{\GP} - \upsvs) \Ind\bigl( \Omega(\xx) \bigr) \bigr\} \bigr\|
	\leq 
	\frac{1}{1 - \dltwbs_{\GP}} \| \DPGP^{-1} \GP^{2} \upsvs \| 
	+ \frac{\sqrt{3 \dltwb_{\GP}}}{1 - \dltwb_{\GP}} \, \E \| \xiv_{\GP} \| 
	+ \CONSTi_{1} \ex^{-\xx}
	\\
	& \leq &
	\frac{1}{1 - \dltwbs_{\GP}} \| \DPGP^{-1} \GP^{2} \upsvs \| 
	+ \frac{\sqrt{3\dltwb_{\GP}}}{1 - \dltwb_{\GP}} \, \sqrt{\dimVG} 
	+ \CONSTi_{1} \ex^{-\xx} \,\, ,
\label{EtuGus11md3G}
\end{EQA}
and
\begin{EQA}
	&& \nquad
	\E \bigl\{ \| \DPGP \, (\tilde{\upsv}_{\GP} - \upsvs) \|^{2} \Ind\bigl( \Omega(\xx) \bigr) \bigr\}
	\\
	& \leq &
	\left( 1 + \frac{\sqrt{2\dltwb_{\GP}}}{1 - \dltwb_{\GP}} \right)^{2} \, \dimA_{\GP}
	+ \, \left( \frac{1}{1 - \dltwbs_{\GP}} \| \DPGP^{-1} \GP^{2} \upsvs \| 
		+ \frac{\sqrt{3\dltwb_{\GP}}}{1 - \dltwb_{\GP}} \, \sqrt{\dimVG}  + \CONSTi_{1} \ex^{-\xx} 
	\right)^{2} \, .
	\qquad
	\qquad
\label{EQtuGmstrVEQtG}
\end{EQA}
\end{theorem}

\begin{remark}
For \( \dltwbs_{\GP} \), \( \dltwb_{\GP} \)  small, 
\eqref{EQtuGmstrVEQtG} yields classical bias-variance decomposition:
\begin{EQA}
	\E \bigl\{ \| \DPGP \, (\tilde{\upsv}_{\GP} - \upsvs) \|^{2} \Ind\bigl( \Omega(\xx) \bigr) \bigr\}
	& \leq &
	\bigl( \dimVG + \| \DPGP^{-1} \GP^{2} \upsvs \|^{2} \bigr) \bigl\{ 1 + o(1) \bigr\} .
\label{EQtuGuvus2tr1o1}
\end{EQA}
With \( \nsize^{-1} = \| \DPGP^{-2} \| \), we also obtain
\begin{EQA}
	\E \bigl\{ \nsize \| \tilde{\upsv}_{\GP} - \upsvs \|^{2} \Ind\bigl( \Omega(\xx) \bigr) \bigr\}
	& \leq &
	\bigl( \dimVG + \nsize \| \DPGP^{-2} \GP^{2} \upsvs \|^{2}\bigr) 
	\bigl\{ 1 + o(1) \bigr\} .
\label{EQtuGuvus2tr1o1n}
\end{EQA}
Moreover, under the small bias condition \( \| \DPGP^{-2} \GP^{2} \upsvs \|^{2} \ll \dimA_{\GP}/\nsize \), 
the impact of the bias induced by penalization is negligible. 
The relation \( \| \DPGP^{-2} \GP^{2} \upsvs \|^{2} \asymp \dimA_{\GP}/\nsize \) is usually referred to as ``bias-variance trade-off''.
Our bound is sharp in the sense that even for the special case of a linear models, 
\eqref{EQtuGuvus2tr1o1} becomes equality.
\end{remark}

\begin{proof}
Below we denote \( \Ex \eta = \E \bigl\{ \eta \Ind\bigl( \Omega(\xx) \bigr) \bigr\} \) for any r.v. \( \eta \).
As \( \E \xiv_{\GP} = 0 \), we derive 
\begin{EQA}
	\Ex \bigl\{ \DPGP (\tilde{\upsv}_{\GP} - \upsvs) \bigr\}
	&=&
	\Ex \DPGP (\tilde{\upsv}_{\GP} - \upsvs - \xiv_{\GP}) 
	- \E \xiv_{\GP} \Ind\bigl( \Omega^{c}(\xx) \bigr) .
\label{EQtuGusIOQPc}
\end{EQA}
For the first term we apply \eqref{DGttGtsGDGm13rG} and \eqref{0mkvhgjnrw3dfwe3u8gty} yielding
\begin{EQA}
	\| \Ex \DPGP (\tilde{\upsv}_{\GP} - \upsvs - \xiv_{\GP}) \|
	& \leq &
	\frac{\| \DPGP^{-1} \GP^{2} \upsvs \| }{1 - \dltwbs_{\GP}} + \frac{\sqrt{3 \dltwb_{\GP} \, \dimVG}}{1 - \dltwb_{\GP}} \, .
\label{7554gtnhbgu8bg87yftwe}
\end{EQA}
To show \eqref{EtuGus11md3G}, we also have to bound the tail moments of 
\( \| \xiv_{\GP} \| \):
\begin{EQA}
	\bigl\| \E \xiv_{\GP} \Ind\bigl( \Omega^{c}(\xx) \bigr) \bigr\|
	& \leq &
	\E \| \xiv_{\GP} \| \Ind\bigl( \Omega^{c}(\xx) \bigr) 
	\leq 
	\ex^{-\xx/2}.
\label{exmxm2E12x2}
\end{EQA}
This can be easily done using deviation bounds for the quadratic form \( \| \xiv_{\GP} \|^{2} \);
see Theorem~\ref{Tdevboundinf}.
Similarly one can bound the variance of \( \DPGP \, \tilde{\upsv}_{\GP} \).
With \( \BB_{\GP} = \DPGP^{-1} \VP^{2} \DPGP^{-1} \)
\begin{EQA}
	\Var_{\xx} \bigl( \DPGP \, \tilde{\upsv}_{\GP} \bigr)
	& \leq &
	\Ex \,\bigl\{ \DPGP \, (\tilde{\upsv}_{\GP} - \upsvs_{\GP}) \bigr\} 
	\bigl\{ \DPGP \, (\tilde{\upsv}_{\GP} - \upsvs_{\GP}) \bigr\}^{\T} 
	\\
	& \leq &
	\left( 1 + \frac{\sqrt{2\dltwb_{\GP}}}{1 - \dltwb_{\GP}} \right)^{2} 
		\E \bigl( \xiv_{\GP} \, \xiv_{\GP}^{\T} \bigr)
	= 
	\left( 1 + \frac{\sqrt{2\dltwb_{\GP}}}{1 - \dltwb_{\GP}} \right)^{2} \BB_{\GP} \, .
\label{VaQtuGd3G1md3Gr2}
\end{EQA}
This yields for the quadratic risk \( \E \| \DPGP \, (\tilde{\upsv}_{\GP} - \upsvs) \|^{2} \)
\begin{EQA}
	\Ex \, \| \DPGP \, (\tilde{\upsv}_{\GP} - \upsvs) \|^{2} 
	& \leq &
	\tr \Var_{\xx}\bigl\{ \DPGP \, (\tilde{\upsv}_{\GP} - \upsvs) \bigr\}
	+ \bigl\| \Ex \, \DPGP (\tilde{\upsv}_{\GP} - \upsvs) \bigr\|^{2}
\label{EQtuGmstrVEQtGp}
\end{EQA}
and \eqref{EQtuGmstrVEQtG} follows. 
\end{proof}

\Section{Confidence sets}
\label{Sconfsets}
This section discusses the issue of constructing some frequentist confidence sets (CS) for the true parameter \( \upsvs \).
We consider two possible approaches. 
The first one is based on the Fisher expansion \eqref{DGttGtsGDGm13rG} for the pMLE \( \tilde{\upsv}_{\GP} \) and leads 
to elliptic CS. 
For the construction, one has to only specify the central point and the scale matrix. 
The second approach is motivated by the Wilks expansion \eqref{3d3Af12DGttG} and leads to likelihood-based CS.
However, due to \eqref{d3GrGELGtsG12}, both constructions are close to each other.

\Subsection{Elliptic confidence sets}
The Fisher expansion \( \DPGP (\upsv - \tilde{\upsv}_{\GP}) \approx \xiv_{\GP} \)
from \eqref{DGttGtsGDGm13rG} for \( \xiv_{\GP} = \DPGP^{-1} \nabla \zeta \)
lead to elliptic confidence sets \( \CS_{\GP}(\zq_{\alp}) \):
\begin{EQA}
	\CS_{\GP}(\zq_{\alp})
	&=&
	\bigl\{ \upsv \colon \| \DPGP (\upsv - \tilde{\upsv}_{\GP}) \| \leq \zq_{\alp} \bigr\},
\label{u2dfdaaqwdderdedfh}
\end{EQA}
where the radius \( \zq_{\alp} \) is fixed to ensure \( \P\bigl( \| \xiv_{\GP} \| > \zq_{\alp} \bigr) \approx \alp \).
Unfortunately, there are several issues related to this construction.
First of all, it involves the matrix \( \DPGP^{2} = - \nabla^{2} \E \LGP(\upsvs_{\GP}) \) which is unknown unless 
the model is linear and the noise variance is precisely known.
Furthermore, the radius \( \zq_{\alp} \) relates to the \( (1-\alp) \)-quantile of the distribution of the norm \( \| \xiv_{\GP} \| \) 
which is unknown as well.
Finally, the Fisher expansion provides an approximation for the difference \( \tilde{\upsv}_{\GP} - \upsvs_{\GP} \)
while we are interested in covering the true value \( \upsvs \).
This requires to account for the bias \( \biasv_{\GP} = \upsvs_{\GP} - \upsvs \).

The first problem can be bypassed by using the matrix \( \DPGPt = \DPGP(\tilde{\upsv}_{\GP}) \) instead of \( \DPGP \); see Section~\ref{SDPGPDPGPt}.
The second problem is more severe and cannot be addressed in the whole generality. 
Under additional structural assumptions on the data \( \Yv \) like a sample \( \Yv = (Y_{1},\ldots,Y_{n})^{\T} \) of independent 
\( Y_{i} \)'s, one can use one or another version of bootstrap resampling technique sketched below; see also \cite{SpZh2014} and reference therein. 
The bias issue can be addressed under the ``small bias'' condition \( \| \DPGP^{-1} \GP^{2} \upsvs \| \ll \rr_{\GP} \).

\Subsection{Likelihood-based confidence sets}
Theorem~\ref{TFiWititG} 
suggests to consider the likelihood-based confidence sets of the form
\begin{EQA}
	\CS_{\GP}(\zz_{\alp})
	&=&
	\bigl\{ \upsv \colon \LGP(\tilde{\upsv}_{\GP}) - \LGP(\upsv) \leq \zz_{\alp} \bigr\} .
\label{u2dfdaaqwdderdedfhL}
\end{EQA}
A great advantage of this construction is that it is entirely based on the penalized log-likelihood 
\( \LGP(\upsv) \).
The main issue for practical applications is the choice of the value \( \zz_{\alp} \).
Due to the Wilks expansion from Theorem~\ref{TFiWititG}, the quantities \( \zq_{\alp} \) from \eqref{u2dfdaaqwdderdedfh} 
and \( \zz_{\alp} \) are related as \( 2 \zz_{\alp} \approx \zq_{\alp}^{2} \).
Again, a resampling technique can be used for a data-driven choice of \( \zz_{\alp} \)
under the small bias condition.

\Subsection{Bootstrap technique and data-driven confidence sets}
The resampling procedure from \cite{SpZh2014} can be sketched as follows.
Assume the additive structure of the log-likelihood for the sample \( \Yv = (Y_{i}) \):
\begin{EQA}
	L(\upsv)
	&=&
	\sumi \ell_{i}(Y_{i},\upsv) , 
\label{lo7ht53680yfcvzkjtt}
\end{EQA}
where \( \ell_{i} \) is the log-density of the \( i \)th observation \( Y_{i} \).

\begin{enumerate}
\item
Draw i.i.d. bootstrap weights \( \Weight^{(m)} = (\weight_{i}^{(m)}) \) from \( \Exp(1) \) or 
\( \ND(1,1) \), \( i=1,\ldots,n \), \( m =1,\ldots,M \), and build the corresponding bootstrap weighted log-likelihood 
\begin{EQA}
	L^{(m)}(\upsv)
	&=&
	\sumi \ell_{i}(Y_{i},\upsv) \weight_{i}^{(m)} \, ;
\label{lo7ht53680yfcvzkjttw}
\end{EQA}

\item 
compute for each \( m \) the corresponding bootstrap estimate 
\( \tilde{\upsv}^{(m)}_{\GP} = \argmax_{\upsv} \LGP^{(m)}(\upsv) \) for 
\( \LGP^{(m)}(\upsv) = L^{(m)}(\upsv) - \| \GP \upsv \|^{2}/2 \).

\item
\( \zz_{\alp} \) from \eqref{u2dfdaaqwdderdedfhL} is fixed as the minimal value ensuring
\begin{EQA}
	\frac{1}{M} \sum_{m=1}^{M} \Ind\bigl( \LGP^{(m)}(\tilde{\upsv}^{(m)}_{\GP}) - \LGP^{(m)}(\tilde{\upsv}_{\GP}) > \zz_{\alp} \bigr)
	& \leq &
	\alp .
\label{odftdfs3ds3erwdfgert}
\end{EQA}
Similarly, the value \( \zq_{\alp } \) from \eqref{u2dfdaaqwdderdedfh} can be fixed by
\begin{EQA}
	\frac{1}{M} \sum_{m=1}^{M} \Ind\bigl( \| \DPGPt(\tilde{\upsv}^{(m)}_{\GP} - \tilde{\upsv}_{\GP}) \| > \zq_{\alp} \bigr)
	& \leq &
	\alp .
\label{odftdfs3ds3erwdfgertq}
\end{EQA}

\end{enumerate}
The proof of validity of this procedure is quite involved, the main technical issue apart the mentioned Fisher and Wilks expansions 
is a Gaussian approximation of the score vector \( \xiv_{\GP} \) for the original and the bootstrap data
as well as some bounds from the random matrix theory; see e.g. \cite{SpZh2014} for more details. 
A further discussion lies beyond the scope of this note.


\Chapter{Laplace approximation of the posterior}
\label{SgenBvMs}

This \chname studies the properties of the posterior \( \vupsv_{\GP} \cond \Yv \).
Our main result states Gaussian approximation of the posterior by \( \ND(\tilde{\upsv}_{\GP} , \DPGPt^{-2}) \).
%
More specifically, our aim is,
for any bounded measurable function \( g \), to compare the conditional moments of \( g(\vupsv_{\GP} - \tilde{\upsv}_{\GP}) \) 
and of \( g(\DPGPt^{-1} \gaussv) \),   
where \( \gaussv \) is standard normal conditionally on \( \Yv \).
The use of \( \nabla \LGP(\tilde{\upsv}_{\GP}) = 0 \) yields
\begin{EQA}
    && \nquad
    \E \bigl\{ g(\vupsv_{\GP} - \tilde{\upsv}_{\GP}) \cond \Yv \bigr\}
    =
    \frac{\int g(\uv - \tilde{\upsv}_{\GP}) \, \ex^{\LGP(\uv)} \, d \uv}
    	 {\int \ex^{\LGP(\uv)} \, d \uv}
    =
    \frac{\int g(\uv) \, \ex^{\LGP(\tilde{\upsv}_{\GP} + \uv) - \LGP(\tilde{\upsv}_{\GP})}
		 d \uv}
    	 {\int \ex^{\LGP(\tilde{\upsv}_{\GP} + \uv) - \LGP(\tilde{\upsv}_{\GP})}
			d \uv}
    \\
    &=&
    \frac{\int g(\uv) \, \exp \bigl\{ 
    		\LGP(\tilde{\upsv}_{\GP} + \uv) - \LGP(\tilde{\upsv}_{\GP}) 
			- \bigl\langle \nabla \LGP(\tilde{\upsv}_{\GP}), \uv \bigr\rangle
		  \bigr\} d \uv}
    {\int \exp \bigl\{ 
    		\LGP(\tilde{\upsv}_{\GP} + \uv) - \LGP(\tilde{\upsv}_{\GP}) 
    		- \bigl\langle \nabla \LGP(\tilde{\upsv}_{\GP}), \uv \bigr\rangle 
		   \bigr\} d \uv}
    \, .
\label{fergywerfdsjmljhdtf3}
\end{EQA}
Now consider the Bregman divergence of the expected log-likelihood \( \fG(\upsv) = \E \LGP(\upsv) \)  
\begin{EQA}
	\fG(\upsv;\uv) 
	&=& 
	\fG(\upsv + \uv) - \fG(\upsv) - \bigl\langle \nabla \fG(\upsv), \uv \bigr\rangle,
	\quad
	\uv \in \R^{\dimp} \, .
\label{koi9876578909876asd}
\end{EQA}
As the stochastic term of \( L(\upsv) \) and thus, of \( \LGP(\upsv) \) is linear in 
\( \upsv \), it holds for any \( \upsv, \uv \)
\begin{EQA}
	\LGP(\upsv + \uv) - \LGP(\upsv) - \bigl\langle \nabla \LGP(\upsv), \uv \bigr\rangle 
	&=&
	\fG(\upsv + \uv) - \fG(\uv) - \bigl\langle \nabla \fG(\upsv), \uv \bigr\rangle
	=
	\fG(\upsv;\uv).
\label{34erdbkjloiu8y7t6r5}
\end{EQA}
Given \( \tilde{\upsv}_{\GP} = \upsv \), we derive from \eqref{fergywerfdsjmljhdtf3}
\begin{EQA}
	\E \bigl\{ g(\vupsv_{\GP} - \tilde{\upsv}_{\GP}) \cond \Yv \bigr\}
    &=&
	\E \bigl\{ g(\vupsv_{\GP} - \tilde{\upsv}_{\GP}) \cond \tilde{\upsv}_{\GP} = \upsv \bigr\}
    =
    \frac{\int g(\uv) \, \ex^{\fG(\upsv;\uv)} \, d\uv}{\int \ex^{\fG(\upsv;\uv)} \, d\uv} \, .
\label{wepodlfkjgdhjdgtfyd7}
\end{EQA}
This basic identity will be systematically used below.
Laplace's approximation means nothing but the use of the second order Taylor approximation of the function \( \fG(\cdot) \) at \( \upsv \).
Namely,
\( \fG(\upsv;\uv) \approx - \| \DPGP(\upsv) \, \uv \|^{2}/2 \) and
\begin{EQA}
	\frac{\int g(\uv) \, \ex^{\fG(\upsv;\uv)} \, d\uv}{\int \ex^{\fG(\upsv;\uv)} \, d\uv}
	& \approx &
	\frac{\int g(\uv) \, \ex^{- \| \DPGP(\upsv) \, \uv \|^{2}/2} \, d\uv}{\int \ex^{- \| \DPGP(\upsv) \, \uv \|^{2}/2} \, d\uv} \, .
\label{uhg432er56wedfgytfefw}
\end{EQA}
The analysis includes two major steps: posterior concentration and a Gaussian approximation of the posterior distribution.

\Section{Posterior concentration}
\label{Spostconcentr}
We start with the important technical result describing the concentration sets of the posterior.
In all our result, the value \( \xx \) is fixed to ensure that \( \ex^{-\xx} \) is negligible.

Proposition~\ref{PconcMLEgenc} enables us to restrict the study to the case 
with \( \tilde{\upsv}_{\GP} \in \CA_{\GP} \).
To describe the concentration properties of the posterior 
we need a slightly stronger concavity condition on \( \E L(\upsv) \),
concavity of \( \E \LGP(\upsv) \) is not sufficient.

\begin{description}
\item[\label{LLCref}\( \bb{(\mathcal{C})} \)]
The function \( \E L(\upsv) \) is concave.

This condition can be relaxed to \emph{weak concavity}.

\item[\label{LLwref}\( \bb{(\mathcal{C}_{\weak})} \)]
There exists \( \GP_{\weak}^{2} \leq \GP^{2} \) such that
for any \( \upsv \in \CA_{\GP} \), 
the function \( 2 \E L(\upsv + \uv) - \| \GPw \uv \|^{2} \) is concave in \( \uv \).
\end{description}

\nameref{LLCref} is a special case of \nameref{LLwref} with \( \GPw = 0 \).
In what follow we assume \nameref{LLCref}.
However, all the results apply under \nameref{LLwref} after replacing \( \DP^{2} \) with \( \DPw^{2} = \DP^{2} + \GPw^{2} \).
Define 
\begin{EQ}[rcccl]
	\dimL(\upsv)
	& \eqdef &
	\tr \bigl\{ \DP^{2}(\upsv) \DPGP^{-2}(\upsv) \bigr\},
	\qquad
	\rrL(\upsv)
	& \eqdef & 
	2 \sqrt{\dimL(\upsv)} + \sqrt{2 \xx} \, ;
\label{tdgtftyfyt4234w34resrr}
\end{EQ}
cf. \eqref{34rtyghuioiuyhgvftid} for \( \dimG \) and \( \rrGP \).
This ensures with \( \gaussv \) standard normal
\begin{EQA}
	\P\bigl( \| \DP(\upsv) \, \DPGP^{-1}(\upsv) \, \gaussv \| > \rrL(\upsv) \bigr)
	& \leq &
	\ex^{-\xx} ,
	\quad
\label{nbvcxdertyuytrewst}
\end{EQA}
see \eqref{Pxiv2dimAxx12} of Corollary~\ref{RsochpHsA}.
%
With some fixed \( \amax \leq 1 \), e.g. \( \amax = 2/3 \), define for any \( \upsv \in \CA_{\GP} \)
\begin{EQA}
	\UV(\upsv) 
	&=& 
	\bigl\{ \uv \colon \| \DP(\upsv) \, \uv \| \leq \amax^{-1} \rrL(\upsv) \bigr\} \, .
\label{kvrder64yhefrnug7st}
\end{EQA} 
With \( f(\upsv) = \E L(\upsv) \) and \( \dltw_{3}(\upsv,\uv) \) from \eqref{dltw3vufuv12f2g},
local smoothness of \( f(\cdot) \) at \( \upsv \) will be measured by the value \( \dltwb(\upsv) \):
\begin{EQA}
	\dltwb(\upsv)
	& \eqdef &
	\sup_{\uv \in \UVL(\upsv)} \frac{1}{\| \DP \uv \|^{2}/2} \bigl| \dltw_{3}(\upsv,\uv) \bigr| ;
\label{om3esuU1H02d3st}
\end{EQA}
cf. \eqref{dtb3u1DG2d3GP}.
Under \nameref{LLtS3ref}, it holds \( \dltwb(\upsv) \leq \amax^{-1} \hmax_{3} \, \rrL(\upsv) \, n^{-1/2}/3 \); see Lemma~\ref{LdltwLaGP}.

\begin{proposition}
\label{PGaussconcge}
Suppose \nameref{Eref},
\nameref{EU2ref},
and \nameref{LLCref}.
Let also \( \dimL(\upsv) \) and \( \rrL(\upsv) \) be defined by \eqref{tdgtftyfyt4234w34resrr} 
and \( \UV(\upsv) \) by \eqref{kvrder64yhefrnug7st}.
If \( \dltwb(\upsv) \) from \eqref{om3esuU1H02d3st} satisfies 
\begin{EQA}
	\dltwb(\upsv)
	& \leq &
	1/3,
	\qquad
	\upsv \in \CA_{\GP} \, ,
\label{d3GusuiAG2d3D2}
\end{EQA} 
then on \( \Omega(\xx) \), it holds with \( \DPt = \DP(\tilde{\upsv}_{\GP}) \) and 
\( \rrLt = \rrL(\tilde{\upsv}_{\GP}) \)
\begin{EQA}
	\P\Bigl( \vupsv_{\GP} - \tilde{\upsv}_{\GP} \not\in \UVt \Cond \Yv \Bigr)
	=
	\P\Bigl( \| \DPt (\vupsv_{\GP} - \tilde{\upsv}_{\GP}) \| > \rrLt \Cond \Yv \Bigr)
	& \leq &
	\ex^{-\xx} .
\label{poybf3679jd532ff2st}
\end{EQA}
\end{proposition}

\begin{proof}
Let us fix \( \tilde{\upsv}_{\GP} = \upsv \) and apply \eqref{wepodlfkjgdhjdgtfyd7} with 
\( g(\uv) = \Ind\bigl( \| \DP(\upsv) \uv \| \not\in \UV(\upsv) \bigr) \).
Then it suffices to bound uniformly in \( \upsv \in \CA_{\GP} \) the ratio
\begin{EQA}
	\rho(\upsv)
	& \eqdef &
    \frac{\int \Ind\bigl( \DP(\upsv) \uv \not\in \UV(\upsv) \bigr) \, \ex^{ \fG(\upsv;\uv) } d \uv}
    	 {\int \ex^{ \fG(\upsv;\uv) } d \uv} \, .
	\qquad
\label{rhopifUVG}
\end{EQA}
\ifLaplace{Bound \eqref{poybf3679jd532ff2} of Theorem~\ref{TLaplaceTV}} 
{Proposition~A.11 of \cite{SpLaplace2022}}
yields the result.
\end{proof}

\Section{Posterior contraction}
\label{Scontractgen}

Now we bring together all the previous results  
to bound the posterior deviations \( \upsv_{\GP} - \upsvs \).
The difference \( \upsv_{\GP} - \upsvs \) can be decomposed as
\begin{EQA}
	\upsv_{\GP} - \upsvs 
	& = &
	\bigl( \upsv_{\GP} - \tilde{\upsv}_{\GP} \bigr) + \bigl( \tilde{\upsv}_{\GP} - \upsvs \bigr) .
\label{vtvGtsGGtssgen}
\end{EQA}
Result \eqref{g25re9fjfregdndg} of
Theorem~\ref{PFiWibias} provides a deviation bound for \( \| \tilde{\upsv}_{\GP} - \upsvs \| \) 
while Proposition~\ref{PGaussconcge} claims concentration of the posterior on the set
\( \bigl\{ \| \DPt (\vupsv_{\GP} - \tilde{\upsv}_{\GP}) \| \leq \rrLt \bigr\} \).
We conclude by the following result.

\begin{proposition}
\label{Ccontactionrategen}
Assume the conditions of Theorem~\ref{PFiWibias} and Proposition~\ref{PGaussconcge} and let 
\( \| \DP^{-2}(\upsv) \| \leq \nsize^{-1} \) for \( \upsv \in \CA_{\GP} \).
It holds on \( \Omega(\xx) \) 
\begin{EQA}
	\P\Bigl( 
		\| \upsv_{\GP} - \upsvs \| 
		\geq 
		\| \tilde{\upsv}_{\GP} - \upsvs \| + \amax^{-1} \rrLt \,/ \sqrt{\nsize}
		\, \cond \Yv 
	\Bigr)
	& \leq &
	2 \ex^{-\xx} ,
\label{PDvtGtsCrGYgen}
\end{EQA}
and \( \| \tilde{\upsv}_{\GP} - \upsvs \| \) satisfies \eqref{g25re9fjfregdndg},
while \( \rrLt \leq (1 - \dltwbss_{\GP})^{-1/2} \, \rrL(\upsvs_{\GP}) \).
\end{proposition}

\begin{proof}
Bound \eqref{PDvtGtsCrGYgen} follows from decomposition \eqref{vtvGtsGGtssgen} and Proposition~\ref{PGaussconcge}.
Further, the use of \eqref{DPGPm1Cd3rG} of Proposition~\ref{PDPGPtDPGP} yields
\( \rrLt \leq (1 - \dltwbss_{\GP})^{-1/2} \, \rrL(\upsvs_{\GP}) \).
\end{proof}

\noindent
The use of \eqref{g25re9fjfregdndg} of Theorem~\ref{PFiWibias} 
implies that the most of posterior mass is concentrated in the root-\( n \) vicinity of \( \upsvs \):
\begin{EQA}
	\P\Bigl( 
		\| \upsv_{\GP} - \upsvs \| 
		\geq 
		3 \| \DPGP^{-2} \GP^{2} \upsvs \| + \frac{3}{\sqrt{\nsize}} \bigl( \rrGP + \rrLt \bigr)
		\, \cond \Yv 
	\Bigr)
	& \leq &
	2 \ex^{-\xx} .
\label{kjdfvgrdf43dfrdsegedr}
\end{EQA}
A prior ensuring the bias-variance trade-off leads to the optimal contraction rate which 
corresponds to the optimal penalty choice in penalized maximum likelihood estimation.

\Section{Gaussian approximation of the posterior}
\label{SgenGaussa}
This section presents our main results about the accuracy of Gaussian approximation of the posterior
\( \vupsv_{\GP} \cond \Yv \) in the total variation distance.
%
The use of self-concordance type conditions from Section~\ref{Slocalsmooth} helps to obtain very accurate and precise finite sample guarantees,
which gradually improve the bounds from \cite{SpPa2019}.


Let \( \BBB(\R^{\dimp}) \) be the \( \sigma \)-field of all Borel sets in \( \R^{\dimp} \),
while \( \BBB_{s}(\R^{\dimp}) \) stands for all centrally symmetric sets from \( \BBB(\R^{\dimp}) \).

\begin{theorem}
\label{TBvMgen}
Assume \nameref{Eref},
\nameref{EU2ref},
and \nameref{LLCref}.
Furthermore, let 
\begin{EQA}[rcl]
	\dltwb(\upsv) \, \dimL(\upsv) 
	& \leq &
	2/3,
	\qquad
	\upsv \in \CA_{\GP} \, ;
\label{om3esuU1H02d3u}
\end{EQA}
cf. \eqref{d3GusuiAG2d3D2}.
Then with
\begin{EQA}
	\err_{2}(\upsv)
	&=& 
	\frac{0.75 \, \dltwb(\upsv) \, \dimL(\upsv)}{1 - \dltwb(\upsv)} \, 
\label{juytr90f2dzaryjhfyfde}
\end{EQA}
and \( \errt = \err_{2}(\tilde{\upsv}_{\GP}) \),
it holds on \( \Omega(\xx) \) with 
\begin{EQ}[rcl]
	\sup_{A \in \BBB(\R^{\dimp})} 
	\left| \P\bigl( \vupsv_{\GP} - \tilde{\upsv}_{\GP} \in A \cond \Yv \bigr) 
		- \PG\bigl( \DPGPt^{-1} \gaussv \in A \bigr)
	\right|
	& \leq & 
	\frac{2 (\errt + \ex^{-\xx})}{1 - \errt - \ex^{-\xx}} 
	\leq 
	4 (\errt + \ex^{-\xx}) .
	\qquad
	\qquad
\label{PPmdi1pdimDGu22G}
\end{EQ}
Here \( \PG \) means a standard Gaussian distribution of \( \gaussv \) given \( \Yv \).
\end{theorem}

Now we present more advanced bounds on the error of Gaussian approximation under conditions \nameref{LL3tref} and \nameref{LL4tref} 
(resp. \nameref{LLtS3ref} and \nameref{LLtS4ref})
from Section~\ref{Slocalsmooth}
for \( \Upsd = \CA_{\GP} \) and \( \HLG(\upsv) = \DP(\upsv) \).

\begin{theorem}
\label{PBvMgen2}
\label{TBvMgen34}
Assume \nameref{Eref},
\nameref{EU2ref},
\nameref{LLCref}, \nameref{LL3tref}, and 
let \( \dltwu_{3} \, \amax^{-1} \rrL(\upsv) \leq 3/4 \) for \( \rrL(\upsv) \) from \eqref{tdgtftyfyt4234w34resrr} 
and all \( \upsv \in \CA_{\GP} \).
Then the concentration bound \eqref{poybf3679jd532ff2st} holds.
Moreover, let 
\begin{EQA}
	\dltwu_{3} \, \amax^{-1} \rrL(\upsv) \, \dimL(\upsv)
	& \leq &
	2 ,
	\qquad
	\upsv \in \CA_{\GP} \, .
\label{0hcde4dftedrf94yr5twew}
\end{EQA}
With \( \dltwb(\upsv) \eqdef \dltwu_{3} \, \rrL(\upsv) / 3 \leq 1/4 \), define 
\begin{EQA}
	\err_{3}(\upsv)
	& \eqdef &
	\frac{\dltwu_{3}}{4 \{ 1 - \dltwb(\upsv) \}^{3/2}} \, \{ \dimL(\upsv) + 1 \}^{3/2} \, .
	\qquad
\label{err3def33}
\end{EQA}
Then the result \eqref{PPmdi1pdimDGu22G} applies on \( \Omega(\xx) \) with \( \errt = \err_{3}(\tilde{\upsv}_{\GP}) \).
Moreover, under \nameref{LL4tref}, 
\begin{EQA}
	\sup_{A \in \BBB_{s}(\R^{\dimp})}
	\left| \P\bigl( \vupsv_{\GP} - \tilde{\upsv}_{\GP} \in A \cond \Yv \bigr) 
		- \PG\bigl( \DPGPt^{-1} \gaussv \in A \bigr)
	\right|
	& \leq & 
	\frac{2 (\errt_{4} + \ex^{-\xx})}{1 - \errt_{4} - \ex^{-\xx}} 
	\leq 
	4 (\errt_{4} + \ex^{-\xx}) 
	\qquad
\label{wd97239rb3ryegdc6swq2ww22}
\end{EQA}
with \( \errt_{4} = \err_{4}(\tilde{\upsv}_{\GP}) \) and
\begin{EQA}
\label{errdef3322Hm1}
	\err_{4}(\upsv)
	& \eqdef &
	\frac{1}{16 \{ 1 - \dltwb(\upsv) \}^{2}} 
	\Bigl[ \dltwu_{3}^{2} \, \bigl\{ \dimL(\upsv) + 2 \bigr\}^{3} + 2 \dltwu_{4} \bigl\{ \dimL(\upsv) + 1 \bigr\}^{2} \Bigr] .
\end{EQA}
The results continue to apply with \nameref{LLtS3ref} (resp. \nameref{LLtS4ref}) in place of \nameref{LL3tref}
(resp. \nameref{LL4tref}) and \( \hmax_{3} \, n^{-1/2} \) (resp. \( \hmax_{4} \, n^{-1} \)) in place of \( \dltwu_{3} \)
(resp. \( \dltwu_{4} \)).
\end{theorem}

\begin{proof}[Proof of Theorem~\ref{TBvMgen} (resp. Theorem~\ref{TBvMgen34})]
Similarly to the proof of Proposition~\ref{PGaussconcge}, we restrict ourselves to the event 
\( \tilde{\upsv}_{\GP} \in \CA_{\GP} \).
Then we fix any possible value \( \upsv \in \CA_{\GP} \) of \( \tilde{\upsv}_{\GP} \) and 
use \eqref{wepodlfkjgdhjdgtfyd7}
to represent the posterior probability of a set \( A \) in the form
\begin{EQA}
	\P\bigl( \vupsv_{\GP} - \tilde{\upsv}_{\GP} \in A \cond \Yv \bigr)
	&=& 
	\frac{\int_{A} \ex^{\fG(\upsv;\uv)} \, d \uv}
		 {\int \ex^{\fG(\upsv;\uv)} \, d \uv} \, .
		 \qquad
\label{rorLGttGuuGG}
\end{EQA}
Now the result follows by 
\ifLaplace{Theorem~\ref{TLaplaceTV} (resp. Theorem~\ref{TLaplaceTV2})}{Theorem~2.1 (resp. Theorem~2.2) of \cite{SpLaplace2022}}.
\end{proof}

Under self-concordance conditions \nameref{LLtS3ref} and \nameref{LLtS4ref}, constraint \eqref{0hcde4dftedrf94yr5twew} reads as
\begin{EQA}
	\sup_{\upsv \in \CA_{\GP}} \, \frac{\hmax_{3} \, \amax^{-1} \rrL(\upsv) \, \dimL(\upsv)}{n^{1/2}}
	& \leq &
	2 .
\label{c6w3gtewr5e4e4wgewrer}
\end{EQA}
As \( \dltwb(\upsv) \leq 1/4 \), Theorem~\ref{TBvMgen34} yields on \( \Omega(\xx) \)
with \( \dimLt = \dimA(\tilde{\upsv}_{\GP}) \)
\begin{EQ}[rl]
	\sup_{A \in \BBB(\R^{\dimp})} 
	\left| \P\bigl( \vupsv_{\GP} - \tilde{\upsv}_{\GP} \in A \cond \Yv \bigr) 
		- \PG\bigl( \DPGPt^{-1} \gaussv \in A \bigr)
	\right|
	& \leq 
	2 \, \hmax_{3} \, \sqrt{\frac{(\dimLt + 1)^{3}}{n}} + 4 \ex^{-\xx} \, ,
	\\
	\sup_{A \in \BBB_{s}(\R^{\dimp})} 
	\left| \P\bigl( \vupsv_{\GP} - \tilde{\upsv}_{\GP} \in A \cond \Yv \bigr) 
		- \PG\bigl( \DPGPt^{-1} \gaussv \in A \bigr)
	\right|
	&
	\leq 
	\frac{\hmax_{3}^{2} \, (\dimLt + 2)^{3} + 2 \hmax_{4} (\dimLt + 1)^{2}}{2n} + 4 \ex^{-\xx} \, .
	\qquad
	\qquad
\label{scdugfdwyd2wywy26e6de432}
\end{EQ}

\Section{Critical dimension in Bayesian inference}
\label{Scritdim}
Posterior concentration in Proposition~\ref{PGaussconcge} only requires 
\( \dltwb(\upsv) \ll 1 \) for all \( \upsv \in \CA_{\GP} \).
Under \nameref{LLtS3ref}, one can bound \( \dltwb(\upsv) \asymp \sqrt{\dimL(\upsv)/n} \)
yielding the condition \( \dimL(\upsv) \ll n \) on the critical dimension 
which is essentially the same as the condition \( \dimG \ll n \) for the pMLE.
This is an important finding and an essential improvement of \cite{SpPa2019}. 
The main result of Theorem~\ref{TBvMgen} requires 
\( \dltwb(\upsv) \, \dimL(\upsv) \ll 1 \) which is much stronger because of the multiplicative factor \( \dimL(\upsv) \). 
Under \nameref{LLtS3ref}, the remainder \( \err_{3} \) 
is of order \( \sqrt{\dimA^{3}(\upsv)/n} \) while under \nameref{LLtS4ref}, \( \err_{4} \asymp \dimA^{3}(\upsv)/n \), still requiring 
\( \dimA^{3}(\upsv) \ll n \).
In some cases, e.g. for additive structure of the log-likelihood, it can be relaxed.
However, it seems that the \( \dimA^{3}(\upsv) \ll n \) condition is inherent in the problem and cannot be relaxed in general situation.
We guess that in the region \( n^{1/3} \ll \dimL(\upsv) \ll n \), 
another non-Gaussian type of limiting behavior of the posterior is well possible.

\Section{Laplace approximation with inexact parameters} 
\label{Spostparamm}
Our main result of Theorem~\ref{TBvMgen} states an approximation of the posterior distribution 
by the Gaussian measure with parameters \( \tilde{\upsv}_{\GP} \) and \( \DPGPt^{-2} \).
However, the vector \( \tilde{\upsv}_{\GP} = \argmax_{\upsv} \LGP(\upsv) \) is typically hard to compute, 
because it solves a high dimensional optimization problem.
If \( \tilde{\upsv}_{\GP} \) and thus \( \DPGPt = \DPGP(\tilde{\upsv}_{\GP}) \) are not available, 
one would be interested to use something more simple in place of \( \tilde{\upsv}_{\GP} \). 
Suppose to be given a vector \( \upsvck \) close to \( \tilde{\upsv}_{\GP} \) and a matrix 
\( \HPck^{2} \) close to \( \DPck_{\GP}^{2} = \DPGP^{2}(\upsvck) \).
A typical example to keep in mind corresponds to \( \upsvck \) being the numerically evaluated posterior mean and 
\( \HPck^{2} \) being the posterior covariance, also evaluated numerically.
Below we aim at presenting some sufficient conditions that ensure a reasonable approximation of the posterior 
by \( \ND(\upsvck,\HPck^{-2}) \) using general results on Gaussian comparison; see \cite{GNSUl2017} and references therein.
Of course, for this result we need all the conditions of Theorem~\ref{TBvMgen} corresponding to the special case 
with \( \upsvck = \tilde{\upsv}_{\GP} \) and \( \HPck^{2} = \DPGPt^{2} \).
We write \( \DPck = \DP(\upsvck) \).
Also we restrict ourselves to the class \( \BBB_{el}(\R^{\dimp}) \) 
of elliptic sets \( A \) in \( \R^{\dimp} \) of the form
\begin{EQA}
	A
	&=&
	\bigl\{ \upsv \in \R^{\dimp} \colon \| \QP (\upsv - \upsvck) \| \leq \rr \bigr\}
\label{vhg4dfe5w3tfdf54wteg}
\end{EQA}
for some linear mapping \( \QP \colon \R^{\dimp} \to \R^{\dimq} \) and \( \rr > 0 \).
Given two symmetric \( \dimq \)-matrices \( \Sigma_{1},\Sigma_{2} \) and a vector \( \av \in \R^{\dimq} \), define 
\begin{EQA}
	\dist(\Sigma_{1},\Sigma_{2},\av)
	& \eqdef &
	\biggl( \frac{1}{\| \Sigma_{1} \|_{\Fr}} + \frac{1}{\| \Sigma_{2} \|_{\Fr}} \biggr)
	\biggl( \| \lamb_{1} - \lamb_{2} \|_{1} + \| \av \|^{2} \biggr),
\label{btegdhertrewfdgvfddffdpo}
\end{EQA}
where \( \| \Sigma \|_{\Fr}^{2} = \tr \Sigma^{2} \),
\( \lamb_{1} \) is the vector of eigenvalues of \( \Sigma_{1} \) arranged in the non-increasing order
and similarly for \( \lamb_{2} \); see \cite{GNSUl2017} for this and related definitions.
By the Weilandt--Hoffman inequality, 
\( \| \lamb_{1} - \lamb_{2} \|_{1} \leq \| \Sigma_{1} - \Sigma_{2} \|_{1} \) , see e.g. 
\cite{MarkusEng}.
Here \( \| M \|_{1} = \tr |M| = \sum_{j} |\lambda_{j}(M)| \) for a symmetric matrix \( M \) with eigenvalues \( \lambda_{j}(M) \).
Obviously, if \( \Sigma_{1} \geq \Sigma_{2} \), then
\( \| \lamb_{1} - \lamb_{2} \|_{1} = \| \Sigma_{1} - \Sigma_{2} \|_{1} = \tr (\Sigma_{1} - \Sigma_{2}) \).

\begin{theorem}
\label{TBvMinexact}
Suppose the conditions of Theorem~\ref{TBvMgen} 
to be fulfilled.
Let also \( \upsvck \in \CA_{\GP} \). 
Given \( \QP \colon \R^{\dimp} \to \R^{\dimq} \), define \( \Sigma_{1} = \QP \DPGPt^{-2} \QP^{\T} \), 
\( \Sigma_{2} = \QP \HPck^{-2} \QP^{\T} \),
\( \av = \QP (\upsvck - \tilde{\upsv}_{\GP}) \)
and suppose \( \| \Sigma_{j} \|^{2} \leq 3 \| \Sigma_{j} \|_{\Fr}^{2} \) for \( j=1,2 \).
Then 
\begin{EQA}[rcl]
	\sup_{\rr > 0}
	\left| \P\bigl( \| \QP (\vupsv_{\GP} - \upsvck) \| \leq \rr \cond \Yv \bigr)
		- \PG\bigl( \| \QP \HPck^{-1} \gaussv \| \leq \rr \bigr) 
	\right|
	& \leq & 
	\frac{2 (\err + \ex^{-\xx})}{1 - \err - \ex^{-\xx}} 
	+ \CONST \dist(\Sigma_{1},\Sigma_{2},\av)\,\, ,
\label{PPmdi1pdimDGu22Ga}
\end{EQA}
where \( \err \) is from Theorem~\ref{TBvMgen}, \( \dist(\cdot) \) from \eqref{btegdhertrewfdgvfddffdpo},
and \( \CONST \) is an absolute constant.
\end{theorem}

\begin{proof}
Use 
\ifLaplace{Theorem~\ref{TLaplaceTVin}}{Theorem~2.9 of \cite{SpLaplace2022}}
with \( f(\upsv) = \E \LGP(\upsv) \), 
\( \xvs = \tilde{\upsv}_{\GP} \), \( \xv = \upsvck \), and \( \DV = \DPGPt \).
\end{proof}

The result is particularly transparent if \( \HPck = \DPGPt \) or, if these two matrices are sufficiently close.
\ifLaplace{Theorem~\ref{TpostmeanLan}}{Theorem~2.10 of \cite{SpLaplace2022}}
yields the following bound.

\begin{corollary}
\label{CTBvMinexact}
Under the conditions of Theorem~\ref{TBvMinexact}, it holds on \( \Omega(\xx) \)
\begin{EQA}[rcl]
	\sup_{\rr > 0}
	\left| \P\bigl( \| \QP (\vupsv_{\GP} - \upsvck) \| \leq \rr \cond \Yv \bigr)
		- \PG\bigl( \| \QP \DPGPt^{-1} \gaussv \| \leq \rr \bigr) 
	\right|
	& \leq & 
	\frac{2 (\err + \ex^{-\xx})}{1 - \err - \ex^{-\xx}} 
	+ \frac{\CONST \| \QP (\upsvck - \tilde{\upsv}_{\GP}) \|^{2}}{\| \QP \DPGPt^{-2} \QP^{\T} \|_{\Fr}} \,\, .
\label{PPmdi1pdimDGu22Gac}
\end{EQA}
\end{corollary}

\Section{Laplace approximation and Bernstein--von Mises Theorem}
\label{SBvMla}
The prominent Bernstein--von Mises (BvM) Theorem claims asymptotic normality of the posterior distribution
with the mean corresponding to the standard MLE \( \tilde{\upsv} = \argmax_{\upsv} L(\upsv) \)
and the variance \( \DPt^{-2} = \DP^{-2}(\tilde{\upsv}) \).
In particular, the prior does not show up in this result, its impact on the posterior distribution becomes negligible
as the sample size \( n \) grows.
In our setup, the situation is different.
The main results of Theorems~\ref{TBvMgen} through \ref{TBvMgen34} state another Gaussian approximation of the posterior
with the mean \( \tilde{\upsv}_{\GP} \) and the variance \( \DPGPt^{-2} \) both depending on the prior
covariance \( \GP^{-2} \).
This dependence is important because the accuracy of approximation is given in terms of 
\( \dimLt = \dimA(\tilde{\upsv}_{\GP}) \) also depending on \( \GP^{-2} \).
It is of interest to describe a kind of phase transition from the classical BvM approximation
by \( \ND(\tilde{\upsv},\DPt^{-2}) \)
to the prior-dependent Laplace approximation by \( \ND(\tilde{\upsv}_{\GP},\DPGPt^{-2}) \). 
%
Intuitively it is clear that the prior impact can be measured by the relation 
between the model-based Fisher information matrix \( \DP^{2} \) and the prior precision matrix \( \GP^{2} \).
The main result below confirms this intuitive guess, however, the result is not trivial and requires a careful treatment
based on Theorem~\ref{TBvMgen} and 
the Gaussian comparison technique mentioned in Section~\ref{Spostparamm}.
Implicitly we assume that the conditions ensuring concentration of the MLE \( \tilde{\upsv} \) corresponding to 
\( \GP^{2} = 0 \) to be fulfilled.
In particular, we need that \( \DP^{2}(\upsv) \) is sufficiently large for all \( \upsv \) in the vicinity of \( \upsvs \).
The radius of this vicinity is given by the value \( \rr_{\xx} = \sqrt{\tr (\DP^{-2} \VP^{2})} + \sqrt{2\xx} \);
see \eqref{34rtyghuioiuyhgvftid} of \nameref{EU2ref}.
Under correct model specification, it holds \( \VP^{2} \leq \CONST \DP^{2} \) and \( \rr_{\xx}^{2} \leq \CONST \dimp \).

\begin{theorem}
\label{TBvMLa}
Under the conditions of Theorem~\ref{TBvMgen}, it holds on \( \Omega(\xx) \)
\begin{EQA}[rcl]
	&& \nquad
	\sup_{\rr > 0}
	\left| \P\bigl( \| \QP (\vupsv_{\GP} - \tilde{\upsv}) \| \leq \rr \cond \Yv \bigr)
		- \PG\bigl( \| \QP \DPt^{-1} \gaussv \| \leq \rr \bigr) 
	\right|
	\\
	& \leq & 
	\frac{2 (\err + \ex^{-\xx})}{1 - \err - \ex^{-\xx}} 
	+ \frac{\CONST \| \QP (\tilde{\upsv} - \tilde{\upsv}_{\GP}) \|^{2}}{\| \QP \DPt^{-2} \QP^{\T} \|_{\Fr}}
	+ \CONST \| \DPGPt^{-1} \GP^{2} \DPGPt^{-1} \| \,\, \frac{\tr (\QP \DPt^{-2} \QP^{\T})}{\| \QP \DPt^{-2} \QP^{\T} \|_{\Fr}} \, .
\label{ygwef78gtewdygfw25dgvfceytfd }
\end{EQA}
Moreover, with \( \QP = \DPt \), \( \rr_{\xx} = \sqrt{\tr (\DP^{-2} \VP^{2})} + \sqrt{2\xx} \leq \CONST \sqrt{\dimp} \),
\begin{EQA}[rcl]
	&& \nquad
	\sup_{\rr > 0}
	\left| \P\bigl( \| \DPt (\vupsv_{\GP} - \tilde{\upsv}) \| \leq \rr \cond \Yv \bigr)
		- \PG\bigl( \| \gaussv \| \leq \rr \bigr) 
	\right|
	\\ 
	& \lesssim & 
	\err + \ex^{-\xx}
	+ {\| \GP \upsvs \|^{2}}/{\sqrt{\dimp}} 
	+ \| \DPGP^{-1} \GP^{2} \DPGP^{-1} \|^{2} \sqrt{\dimp} \,\, .
\label{PPmdi1pdimDGu22GacBvMtt}
\end{EQA}
\end{theorem}

\begin{proof}
Theorem~\ref{TBvMinexact} yields in view of \( \DPt^{-2} \geq \DPGPt^{-2} \)
\begin{EQA}[rcl]
	&& \nquad
	\sup_{\rr > 0}
	\left| \P\bigl( \| \QP (\vupsv_{\GP} - \tilde{\upsv}) \| \leq \rr \cond \Yv \bigr)
		- \PG\bigl( \| \QP \DPt^{-1} \gaussv \| \leq \rr \bigr) 
	\right|
	\\
	& \leq & 
	\frac{2 (\err + \ex^{-\xx})}{1 - \err - \ex^{-\xx}} 
	+ \frac{\CONST \| \QP (\tilde{\upsv} - \tilde{\upsv}_{\GP}) \|^{2}}{\| \QP \DPt^{-2} \QP^{\T} \|_{\Fr}}
	+ \frac{\CONST \tr \{ \QP (\DPt^{-2} - \DPGPt^{-2}) \QP^{\T} \}}{\| \QP \DPt^{-2} \QP^{\T} \|_{\Fr}} \,\, .
	\qquad
\label{PPmdi1pdimDGu22GacBvM}
\end{EQA}
The last term here can easily be bounded:
\begin{EQA}
	&& \nquad
	\tr \{ \QP (\DPt^{-2} - \DPGPt^{-2}) \QP^{\T} \}
	=
	\tr \{ \QP \DPt^{-1} (\Id_{\dimp} - \DPt \DPGPt^{-2} \DPt) \DPt^{-1} \QP^{\T} \}
	\\
	& \leq &
	\| \Id_{\dimp} - \DPt \DPGPt^{-2} \DPt \| \, \tr (\QP \DPt^{-2} \QP^{\T}) 
	=
	\| \DPGPt^{-1} \GP^{2} \DPGPt^{-1} \| \, \tr (\QP \DPt^{-2} \QP^{\T}) .
\label{mf34fdssatf7y64523214y68i}
\end{EQA}
Here we used that 
\begin{EQA}
	\| \DPGPt^{-1} \GP^{2} \DPGPt^{-1} \|
	&=&
	\| \DPGPt^{-1} (\DPGPt^{2} - \DPt^{2}) \DPGPt^{-1} \|
	=
	\| \Id - \DPGPt^{-1} \DPt^{2} \DPGPt^{-1} \|
	=
	\| \Id_{\dimp} - \DPt \DPGPt^{-2} \DPt \| .
\label{ytsafdyqw786ydfwqe8ygujahsdcv}
\end{EQA}
For \( \QP = \DPt \), we use that \( \| \QP \DPt^{-2} \QP^{\T} \|_{\Fr} = \sqrt{\dimp} \),
\( \tr (\QP \DPt^{-2} \QP^{\T}) = \dimp \), and
apply the Fisher expansion \eqref{DGttGtsGDGm13rG} of Theorem~\ref{TFiWititG} to \( \tilde{\upsv} \)
and \( \tilde{\upsv}_{\GP} \).
On \( \Omega(\xx) \), it holds \( \| \DP^{-1} \score \| \leq \rr_{\xx} \) and with \( \dltwb \) from \eqref{dtb3u1DG2d3GP}
\begin{EQA}
	\| \DP (\tilde{\upsv} - \tilde{\upsv}_{\GP}) \|
	& \leq &
	\| \DP (\upsvs - \upsvs_{\GP}) \|
	+ \| (\Id - \DP \DPGP^{-2} \DP) \DP^{-1} \score \|
	\\
	&&
	+ \, \| \DP (\tilde{\upsv} - \upsvs) - \DP^{-1} \score \| 
	+ \| \DP (\tilde{\upsv}_{\GP} - \upsvs_{\GP} - \DPGP^{-2} \score) \| 
	\\
	& \leq &
	\| \DP (\upsvs - \upsvs_{\GP}) \|
	+ \CONST \| \DPGP^{-1} \GP^{2} \DPGP^{-1} \| \, \rr_{\xx} 
	+ \CONST \dltwb \, \rr_{\xx} \, .
\label{usyfgasgsadiufgiuyw}
\end{EQA}
By \eqref{11ma3eaelDebQ} of Proposition~\ref{Lvarusetb} 
\begin{EQA}
	\| \DP (\upsvs - \upsvs_{\GP}) \|^{2}
	& \leq &
	\| \DP \DPGP^{-2} \DP^{\T} \| \, \| \GP \upsvs \|^{2} 
	\leq 
	\| \GP \upsvs \|^{2} .
\label{usafgy42er86fwe8yghscd}
\end{EQA} 
As \( \rr_{\xx}^{2} \leq \CONST \dimp \), \( \dltwb \sqrt{\dimp} \leq \CONST \err_{3} \),
\( \DPt^{2} \leq 2 \DP^{2} \), and \( \DPGPt^{-2} \leq 2 \DPGP^{-2} \), the assertion follows.
\end{proof}

\begin{remark}
The use of BvM requires rather strong bounds on the penalizing matrix \( \GP^{2} \) and the related bias 
\( \| \GP \upsvs \| \).
We need the condition of ``light penalization'' \( \| \DPGP^{-1} \GP^{2} \DPGP^{-1} \| \ll \dimp^{-1/2} \) 
which is much stronger than \( \GP^{2} \ll \DP^{2} \).
Similarly, the ``light bias'' condition \( \| \GP \upsvs \|^{2} \ll \dimp^{1/2} \)
is more restrictive than the ``small bias'' or ``undersmoothing'' condition \( \| \GP \upsvs \|^{2} \ll {\dimp} \).
\end{remark}

\Section{Posterior mean}
\label{Spostmeanvar}
This section addresses an important question of using the posterior mean in place of the MAP 
\( \tilde{\upsv}_{\GP} \) for Bayesian inference.
Our main result justifies the use of the posterior mean in place of the MAP under the same critical dimension condition
\( \dimL(\upsv) \ll n^{1/3} \) which is required for the Gaussian approximation result.
First we quantify the deviation of the posterior mean \( \bar{\upsv}_{\GP} \) from 
\( \tilde{\upsv}_{\GP} \).
Then we apply Corollary~\ref{CTBvMinexact} to measure the impact of using \( \bar{\upsv}_{\GP} \) in place of 
\( \tilde{\upsv}_{\GP} \).
By definition
\begin{EQA}
	\bar{\upsv}_{\GP} - \tilde{\upsv}_{\GP}
	& \eqdef &
	\E \bigl( \vupsv_{\GP} \cond \Yv \bigr) - \tilde{\upsv}_{\GP}
	=
	\frac{\int (\upsv - \tilde{\upsv}_{\GP}) \, \ex^{\LGP(\upsv)} \, d\upsv}{\int \ex^{\LGP(\upsv)} \, d\upsv} \, .
\label{uexLGuduidu}
\end{EQA}
More precisely, we consider a linear mapping \( \QP \colon \R^{\dimp} \to \R^{\dimq} \)
and evaluate the value \( \bigl\| \QP (\bar{\upsv}_{\GP} - \tilde{\upsv}_{\GP}) \bigr\| \).
The choice of \( \QP \) is important.
In particular, we cannot take \( \QP = \DPGP \) because this choice makes the bound 
dependent and linearly growing with \( \dimp \).

\begin{theorem}
\label{Tpostmean}
Assume the conditions of Theorem~\ref{TBvMgen34} and let \( \QP^{\T} \QP \leq \DP^{2}(\upsv) \) for all \( \upsv \in \CA_{\GP} \).
Then it holds with some absolute constant \( \CONST \)
\begin{EQA}
	\| \QP (\bar{\upsv}_{\GP} - \tilde{\upsv}_{\GP}) \|
	& \leq &
	2.4 \, \hmax_{3} \, \| \QP \DPt^{-2} \QP^{\T} \|^{1/2} \,
	{(\dimLt + 1)^{3/2}}{\, n^{-1/2}} + \CONST \ex^{-\xx} .
\label{hcdtrdtdehfdewdrfrhgyjufgerst}
\end{EQA}
\end{theorem}

\begin{proof}
One can apply the same trick as before: by \( \nabla\LGP(\tilde{\upsv}_{\GP}) = 0 \)
\begin{EQA}
	\QP (\bar{\upsv}_{\GP} - \tilde{\upsv}_{\GP})
	&=&
	\frac{\int \QP \uv \, \exp\bigl\{ \LGP(\tilde{\upsv}_{\GP} + \uv) - \LGP(\tilde{\upsv}_{\GP}) - 
			\langle \nabla\LGP(\tilde{\upsv}_{\GP}),\uv \rangle \bigr\} \, d\uv}
		 {\int \exp\bigl\{ \LGP(\tilde{\upsv}_{\GP} + \uv) - \LGP(\tilde{\upsv}_{\GP}) - 
			\langle \nabla\LGP(\tilde{\upsv}_{\GP}),\uv \rangle \bigr\} \, d\uv} \, .
\label{huGtuGfiuexdi}
\end{EQA}
For any particular value \( \tilde{\upsv}_{\GP} = \upsv \), 
stochastic linearity allows to replace the Bregman divergence of the log-likelihood \( \LGP(\upsv) \) by the similar one 
for the expected log-likelihood \( \fG(\upsv) = \E \LGP(\upsv) \).
This yields 
\begin{EQA}
	\| \QP(\bar{\upsv}_{\GP} - \tilde{\upsv}_{\GP}) \|
	& \leq &
	\biggl\| \frac{\int \QP \uv \, \ex^{ \fG(\upsv;\uv)} \, d\uv}{\int \ex^{ \fG(\upsv;\uv)} \, d\uv} 
	\biggr\| \, .
\label{QhuGtuGuiCGdu}
\end{EQA}
Now we may apply 
\ifLaplace{\eqref{hcdtrdtdehfdewdrfrhgyjufger} of Theorem~\ref{TpostmeanLa}}{Theorem~2.7 of \cite{SpLaplace2022}}
with \( \DV = \DPGP(\upsv) \).
\end{proof}

Now we specify the result for the special choice \( \QP = \DPt \).

\begin{corollary}
\label{CTpostmeanLast}
Assume the conditions of Theorem~\ref{TBvMgen34}.
Then 
\begin{EQA}
	\| \DPt (\bar{\upsv}_{\GP} - \tilde{\upsv}_{\GP}) \|
	& \leq &
	2.4 \, \hmax_{3} \, {(\dimLt + 1)^{3/2}}{n^{-1/2}} + \CONST \ex^{-\xx} \, .
\label{klu8gitfdgregfkhj7ytst}
\end{EQA}
\end{corollary}

Now we put together the result of Theorem~\ref{Tpostmean} and  
the accuracy bound from Theorem~\ref{TBvMinexact}.
To make the result more transparent, assume \( \QP = \DPt \) and \( \HPck = \DPGPt \).

\begin{theorem}
\label{Tpostmeanmode}
Assume the conditions of Theorem~\ref{TBvMgen34}.
Then on \( \Omega(\xx) \)
\begin{EQA}[rcl]
	&& \nquad
	\sup_{\rr > 0}
	\left| \P\bigl( \| \DPt (\vupsv_{\GP} - \upsvb_{\GP}) \| \leq \rr \cond \Yv \bigr)
		- \PG\bigl( \| \DPt \DPGPt^{-1} \gaussv \| \leq \rr \bigr) 
	\right|
	\leq 
	\CONST \Bigl( \frac{(\dimLt + 1)^{3/2}}{n^{1/2}} + \ex^{-\xx} \Bigr) .
\label{PPmdi1pdimDGu22Gab}
\end{EQA}
\end{theorem}

We conclude that the use of posterior mean in place of posterior mode
is possible under the same condition  \( \dimLt^{3} \ll \nsize \).
This is a non-trivial result based on recent progress in Gaussian probability from \cite{GNSUl2017} 
and it is only valid if we limit ourselves to elliptic credible sets.

\Section{Non-Gaussian priors}

\label{SnonGaussprior}
All the results of the paper are stated for Gaussian priors \( \ND(0,\GP^{-2}) \).
Of course, everything can be extended to the case of a prior \( \ND(\upsv_{0},\GP^{2}) \) by a parameter shift.
Relaxing the assumption of a Gaussian prior requires a bit more work.
Indeed, the effective dimension \( \dimG \) - one of most important notions in the paper - 
is defined in terms of the prior precision matrix \( \GP^{-2} \).
Here we discuss how the case of a general prior can be covered.

Let the (quasi) log-likelihood \( L(\upsv) \) satisfy \nameref{Eref}, that is, the gradient \( \nabla \zeta \) 
of the stochastic component \( \zeta(\upsv) = L(\upsv) - \E L(\upsv) \) does not depend on \( \upsv \).
Let also \( \priord_{\GP}(\upsv) \) stand for the prior density.
Define the penalized log-likelihood
\begin{EQA}
	\LGP(\upsv)
	& \eqdef &
	L(\upsv) + \log \priord_{\GP}(\upsv) .
\label{6jcvye3jhfkoigtuy}
\end{EQA} 
We assume that the function \( \LGP(\upsv) \) is \emph{strongly concave} in the sense that
\begin{EQA}
	\DPGP^{2}(\upsv)
	\eqdef
	- \nabla^{2} \LGP(\upsv)
	& \geq &
	\GP^{2} 
\label{jvfiur4ejdkcwst,kdhyfcj}
\end{EQA}
for some fixed positive matrix \( \GP^{2} \).
This enables us to define in a unique way the pMLE \( \tilde{\upsv}_{\GP} \) and its population counterpart \( \upsvs_{\GP} \):
\begin{EQA}
	\tilde{\upsv}_{\GP}
	&=&
	\argmax_{\upsv} \LGP(\upsv),
	\\
	\upsvs_{\GP}
	&=&
	\argmax_{\upsv} \E \LGP(\upsv),
\label{jhgf7ywe3jdcgyte76ekjyy}
\end{EQA}
as well as positive matrices 
\begin{EQA}
	\DP^{2}(\upsv)
	& \eqdef &
	\DPGP^{2}(\upsv) - \GP^{2} > 0 \, .
\label{hfuejdhcyswjswmn}
\end{EQA}
Note that \( \DP^{2}(\upsv) \) does not coincide with \( - \nabla^{2} L(\upsv) \) unless 
\( - \nabla^{2} \log \priord_{\GP}(\upsv) = \GP^{2} \), that is, \( \priord_{\GP} \) corresponds to \( \ND(\upsv_{0},\GP^{2}) \).
From this point, we can proceed as in the case of a Gaussian prior.
In particular, the \emph{Laplace effective dimension} \( \dimA(\upsv) \) is defined exactly as in the case of a Gaussian prior:
\begin{EQA}
	\dimA(\upsv)
	& \eqdef &
	\tr\bigl\{ \DPGP^{-2}(\upsv) \, \DP^{2}(\upsv) \bigr\} 
	=
	\tr\bigl\{ \Id_{\dimp} - \DPGP^{-2}(\upsv) \, \GP^{2} \bigr\}.
\label{hjcvdtywedhsxcjsuvcgtwh}
\end{EQA}
All the smoothness conditions including \nameref{LL3tref}, \nameref{LL4tref},
\nameref{LLtS3ref}, \nameref{LLtS4ref} concern the function \( f(\upsv) = \E \LGP(\upsv) = \E L(\upsv) + \log \priord_{\GP}(\upsv) \).
In the contrary, condition \nameref{Eref} only relies on the stochastic component of log-likelihood \( L(\upsv) \). 
Condition \nameref{EU2ref} effectively requires exponential moments of \( \DPGP^{-1} \nabla \zeta \) for \( \DPGP = \DPGP(\upsvs_{\GP}) \).

One may conclude that the approach extends in a straightforward way to the case when the function 
\( \LGP(\upsv) = L(\upsv) + \log \priord_{\GP}(\upsv) \) is \emph{strongly concave}, \emph{smooth}, 
and \emph{stochastically linear} with bounded exponential moments.


\Chapter{Smooth priors and rate over Sobolev classes}
\label{Spriorexample}
This section presents some examples of choosing \( \GP^{2} \) for 
achieving the ``bias-variance trade-off''
and obtaining rate optimal results.
Define the background true parameter
\begin{EQA}
	\upsvs 
	&=& 
	\argmax_{\upsv} \E L(\upsv) .
\label{v43hrf5fjh9ehft6f}
\end{EQA}
Remind the notation \( \LGP(\upsv) \eqdef L(\upsv) - \| \GP \upsv \|^{2}/2 \),
\begin{EQA}[rcccl]
	\tilde{\upsv}_{\GP}
	& \eqdef &
	\argmax_{\upsv \in \Ups} \LGP(\upsv)
\label{v8j3gug5ejfgifkri}
	\qquad	
	\upsvs_{\GP}
	& \eqdef &
	\argmax_{\upsv \in \Ups} \E \LGP(\upsv) .
\label{upssaruELuaiUa}
\end{EQA}
Also, for any \( \upsv \in \Ups \),
\begin{EQA}[rcccl]
	\IF(\upsv) 
	&=& 
	- \nabla^{2} \E L(\upsv) ,
	\qquad
	\IF_{\GP}(\upsv) 
	&=& 
	- \nabla^{2} \E \LGP(\upsv) 
	= 
	\IF(\upsv) + \GP^{2} .
\label{e456yhgfr567uhgfa}
\end{EQA}
Posterior concentration relies on the \emph{Laplace effective dimension}
\begin{EQA}
	\dimA(\upsv)
	& \eqdef &
	\tr\bigl\{ \IF_{\GP}^{-1}(\upsv) \, \IF(\upsv) \bigr\} 
\label{0n2ig82jfirjqgdkgyehj}
\end{EQA}
for \( \upsv \) in a small vicinity \( \CA_{\GP} \).
The pMLE \( \tilde{\upsv}_{\GP} \) concentrates in a vicinity of \( \upsvs_{\GP} \)
and the concentration radius is described via the \emph{effective dimension}
\begin{EQA}
	\dimG
	& \eqdef &
	\tr\bigl( \IF_{\GP}^{-1} \VP^{2} \bigr) ,
\label{ufh3u3whdjhduwkfig2}
\end{EQA}
where for \( \IF_{\GP} = \IF_{\GP}(\upsvs_{\GP}) \) and \( \VP^{2} = \Var (\nabla \zeta ) \).
It is worth noting that under conditions of Proposition~\ref{Lvarusetb} 
\( \IF_{\GP}(\upsvs) \) can be used in place of \( \IF_{\GP} \).
Moreover,
\( \VP^{2} = \IF(\upsvs) \) under the correct model specification.

Theorem~\ref{PFiWi2risk} yields the following bound for the risk of \( \tilde{\upsv}_{\GP} \): 
\begin{EQA}
	\E \| \DPGP (\tilde{\upsv}_{\GP} - \upsvs) \|^{2}
	& \lesssim &
	\dimG + \| \DP_{\GP}^{-1} \GP^{2} \upsvs \|^{2} \, .
\label{bnjmb3ed6w2jh21fg7fj}
\end{EQA}
This suggest to select the operator \( \GP^{2} \) which ensures the ``bias-variance trade-off'' 
\( \| \DP_{\GP}^{-1} \GP^{2} \upsvs \|^{2} \asymp \dimG \).

\Section{Bias-variance trade-off under \( \GPa \)-smoothness}
Suppose that \( \upsvs \) is \( \GPa \)-smooth, that is,
for a given matrix \( \GPa^{2} \)
\begin{EQA}
	\| \GPa \upsvs \|^{2}
	& \leq &
	1 .
\label{d4hgergv76ehrfrtfge}
\end{EQA}
Consider the univariate family of penalizing/precision matrices \( \GP^{2} \) of the form \( \GP^{2} = \CGP \GPa^{2} \).
Later we discuss a choice of the value \( \CGP \) ensuring the bias-variance relation 
\( \| \DP_{\GP}^{-1} \GP^{2} \upsvs \|^{2} \asymp \dimG \).
Note first that
\begin{EQA}
	\| \DP_{\GP}^{-1} \GP^{2} \upsvs \|^{2}
	& \leq &
	\| \GP \upsvs \|^{2} =	\CGP \, \| \GPa \upsvs \|^{2} 
	= 
	\CGP .
\label{6kv9gvbjdfcvw2cdged8e}
\end{EQA}
Next we consider the variance term measured by the effective dimension \( \dimG \).
Under the conditions of Proposition~\ref{Lvarusetb} and correct model specification, 
we may slightly change the definition to
\begin{EQA}
	\dimG
	&=&
	\dimA(\CGP)
	=
	\tr (\DPf^{2} \, \DPfGP^{-2})
	=
	\tr\bigl\{ \DPf^{2} (\DPf^{2} + \CGP \GPa^{2})^{-1} \bigr\} ,
\label{u98gfw3wd7uewj4ert7yg}
\end{EQA}
where \( \DPf^{2} = \IF(\upsvs) \).
The definition is particularly transparent for \( \DPf^{2} = \nsize \Id_{\dimp} \):
\begin{EQA}
	\dimA(\CGP)
	&=&
	\tr\bigl( \Id_{\dimp} + \nsize^{-1} \CGP \, \GPa^{2} \bigr)^{-1} .
\label{c78ju43njvf6vf4e3tg9ythgh}
\end{EQA}
Define also 
\begin{EQA}
	\rr(\CGP)
	&=&
	\sqrt{\dimA(\CGP)} + \sqrt{2 \xx} .
\label{3gyvf76gv7r4ur45r56}
\end{EQA}
Theorem~\ref{PFiWibias}, \eqref{rG1r1d3GbvGb}, yields the following bound.

\begin{theorem}
\label{LbiasCGPw}
Assume \( \| \GPa \upsvs \| \leq 1 \).
Consider \( \GP^{2} = \CGP \, \GPa^{2} \), \( \DPfGP^{2} = \DPf^{2} + \CGP \GPa^{2} \).
Suppose the conditions of Proposition~\ref{PconcMLEgenc} and \ref{Lvarusetb} hold for this choice of \( \GP^{2} \).
Then on \( \Omega(\xx) \), the estimate \( \tilde{\upsv}_{\GP} \) fulfills 
for some fixed constant \( \CONST \)
\begin{EQA}
	\| \DPfGP (\tilde{\upsv}_{\GP} - \upsvs) \|
	& \leq &
	2 \rr(\CGP) + 	\amax^{-1} \CGP^{1/2} \, .
\label{0mkvhgjnrw3dfwe3u8CGPw}
\end{EQA}
\end{theorem}

Clearly \( \dimA(0) = \dimp \) and \( \dimA(\CGP) \) monotonously decreases to zero as \( \CGP \) grows to infinity
provided that \( \GPa^{2} \) is positive definite.
Therefore, for any fixed \( \CONSTbal \), the equation
\begin{EQA}
	\rr(\CGP)
	&=&
	\CONSTbal \, \CGP^{1/2}
\label{h0hit5ti85tjrtfg78}
\end{EQA}
has a unique solution \( \CGP = \CGPs \) leading to the desired ``bias-variance trade-off''. 

\begin{theorem}
\label{LbiasCGP}
Assume \( \| \GPa \upsvs \| \leq 1 \).
With a fixed \( \CONSTbal \), define \( \CGPs \) by \( \dimA(\CGPs) = \CONSTbal \, \CGPs \); see \eqref{u98gfw3wd7uewj4ert7yg}.
Define also \( \GP_{\CGPs}^{2} = \CGPs \, \GPa^{2} \), \( \DPf_{\CGPs}^{2} = \DPf^{2} + \CGPs \GPa^{2} \).
Suppose the conditions of Proposition~\ref{PconcMLEgenc} and \ref{Lvarusetb} hold for this choice of \( \GP^{2} = \GP_{\CGPs}^{2} \).
Then on \( \Omega(\xx) \), the estimate \( \tilde{\upsv}_{\CGPs} = \tilde{\upsv}_{\GP_{\CGPs}} \) fulfills 
for some fixed constant \( \CONST \)
\begin{EQA}
	\| \DPf_{\CGPs} (\tilde{\upsv}_{\CGPs} - \upsvs) \|
	& \leq &
	\CONST \CGPs^{1/2} \, .
\label{0mkvhgjnrw3dfwe3u8CGP}
\end{EQA}
With \( \DPf^{2} = \nsize \Id_{\dimp} \), this implies
\begin{EQA}
	\| \tilde{\upsv}_{\CGPs} - \upsvs \|^{2}
	& \leq &
	{\CONST \CGPs}/{\nsize} \, .
\label{0mkvhgjnrw3dfwe3u8CGPn}
\end{EQA}
\end{theorem}

\begin{remark}
The value \( \dimA(\CGP) \) from \eqref{u98gfw3wd7uewj4ert7yg} can be written as
\begin{EQA}
	\dimA(\CGP)
	&=&
	\tr (\Id_{\dimp} + \CGP \DPf^{-1} \, \GPa^{2} \, \DPf^{-1})^{-1} .
\label{78cvhufu7yejhedukfiuer}
\end{EQA}
The operator \( \proj(\CGP) = (\Id_{\dimp} + \CGP \DPf^{-1} \, \GPa^{2} \, \DPf^{-1})^{-1} \) is obviously
a sub-projector in \( \R^{\dimp} \).
In many situation, its trace can be well approximated by the dimension of the subspace in \( \R^{\dimp} \)
on which \( \CGP \DPf^{-1} \, \GPa^{2} \, \DPf^{-1} \leq \Id_{\dimp} \); see the example later for the case 
of Sobolev smoothness.
In particular, with \( \DPf^{2} = \nsize \Id_{\dimp} \), the corresponding subspace is spanned by the eigenvectors
of \( \GPa^{2} \) with eigenvalues \( \gp_{j}^{2} \leq n/\CGP \).
\end{remark}

The obtained bound can be extended to the case 
when the smoothness of \( \upsvs \) is measured using a different operator \( \GPa \) than \( \GP \).
An important example corresponds to a \( (\smp,\CGP) \)-smooth prior 
and a signal \( \upsvs \) is from a Sobolev ball with the smoothness degree \( \smpa \neq \smp \);
see 
\ifapp{Section~\ref{Ssmoothprior}}{Section~\ref{SGBvM}}.
Suppose \( \upsvs \) is \( \GPa \)-smooth, and \( \CGPs \) is the corresponding solution of \eqref{h0hit5ti85tjrtfg78}
leading to \( \GP_{\CGPs}^{2} = \CGPs \, \GPa^{2} \) and \( \DPf_{\CGPs}^{2} = \DPf^{2} + \CGPs \GPa^{2} \).
We now consider some other penalizing matrix \( \GP^{2} \).

\begin{theorem}
\label{Pbiassmoone}
Let \( \upsvs \) follow \eqref{d4hgergv76ehrfrtfge} and \( \CGPs \) solve \eqref{h0hit5ti85tjrtfg78}.
Fix \( \GP_{\CGPs}^{2} = \CGPs \, \GPa^{2} \) and \( \DPf_{\CGPs}^{2} = \DPf^{2} + \CGPs \GPa^{2} \).
Suppose for some other \( \GP^{2} \), it holds
\begin{EQA}
	\DPf_{\GP}^{-2} \GP^{2}
	& \leq &
	\CONSTGP \, \DPf_{\CGPs}^{-2} \, \GP_{\CGPs}^{2}
\label{hgjkjhyoie3r56dftrtyf}
\end{EQA}
with a fixed constant \( \CONSTGP \),
and the conditions of Proposition~\ref{PconcMLEgenc} and \ref{Lvarusetb} hold for this choice of \( \GP^{2} \).
Then on \( \Omega(\xx) \) with 
\( \nsize^{-1} = \| \DPfGP^{-2} \| \vee \| \DPf_{\CGPs}^{2} \| \)
\begin{EQA}[rcl]
	\sqrt{\nsize} \| \tilde{\upsv}_{\GP} - \upsvs \|
	& \leq &
	\CONST \, \CONSTGP \, \CGPs^{1/2} \, .
\label{vfgd32hvfbdgi8ghj89ne}
\end{EQA}
\end{theorem}

\begin{proof}
Apply \eqref{g25re9fjfregdndg} and \eqref{hgjkjhyoie3r56dftrtyf}.
\end{proof}

Some sufficient conditions for \eqref{hgjkjhyoie3r56dftrtyf} are given in Lemma~\ref{Lsmoothbiase}.

\Section{Smooth penalty and rate of estimation over Sobolev balls}
\label{SrateSobball}
This section illustrates the general results for the case of a smooth penalty \( \GP^{2} \)
for a smooth signal \( \upsvs \) in a Sobolev sense.
Our main intension is to demonstrate that the obtained finite sample bounds imply the well known results on the rate of estimation
under standard smoothness assumptions.
We follow the setup of Section~\ref{Spriors} and assume the diagonal structure of the penalizing matrix 
\( \GP^{2} = \diag(\gp_{1}^{2},\ldots,\gp_{\dimp}^{2}) \) with eigenvalues
\( \gp_{j}^{2} \) satisfying the condition of polynomial growth \eqref{sumjJgjm2C}.
We also assume that the information matrix \( \IF = \IF(\upsvs) \) satisfies \eqref{uVmGu2C1C2gen}:
\begin{EQ}[rcccl]
	\CONSTIF^{-1} \, \nsize \| \uv \|^{2} 
	& \leq &
	\bigl\langle \IF \uv, \uv \bigr\rangle
	& \leq &
	\CONSTIF \, \nsize \| \uv \|^{2},
	\qquad
	\uv \in \R^{\dimp} \, .
\label{uVmGu2C1C2genS}
\end{EQ}
while \( \VP^{2} \) satisfies \eqref{uVmGu2C1C2genV}:
\begin{EQ}[rcccl]
	\CONSTV^{-1} \, \| \IF \uv \|^{2} 
	& \leq &
	\| \VP \uv \|^{2}
	& \leq &
	\CONSTV \, \nsize \| \IF \uv \|^{2},
	\qquad
	\uv \in \R^{\dimp} \, .
\label{uVmGu2C1C2genVS}
\end{EQ}
Here \( \nsize \) means a sample size; see, e.g., Section~\ref{SGBvM}.

%

\medskip
To illustrate the obtained results, assume that \( \upsvs \) belongs to a Sobolev ball \( \BBB(\smpa,\CGPa) \):
\begin{EQA}
	\BBB(\smpa,\CGPa)
	& \eqdef &
	\Bigl\{ \upsv = (\ups_{j}) \colon \sum_{j \geq 1} j^{2\smpa} \ups_{j}^{2} \leq \CGPa \Bigr\} 
\label{BBBsttsjj1sj2}
\end{EQA}
with \( \smpa > 0 \) and \( \CGPa \asymp 1 \); see Section~\ref{Ssmoothprior}.
First we discuss a smoothness-aware choice of \( \GP^{2} = \diag(\gp_{1}^{2},\ldots,\gp_{\dimp}^{2}) \) assuming \( \smpa > 1/2 \) known.
By \( \CONST \) we denote some fixed constant possibly depending on the other constants in our conditions like
\( \CONSTIF \), \( \CONSTV \), \( \smpa \), \( \hmax_{3} \), \( \CONSTGP \).

\begin{theorem}
\label{CGLMlossGP}
Let also \eqref{uVmGu2C1C2genS} and \eqref{uVmGu2C1C2genVS} hold
and \( \upsvs \in \BBB(\smpa,\CGPa) \), \( \smpa > 1/2 \).
Define \( \GP^{2} = \diag(\gp_{1}^{2},\ldots,\gp_{\dimp}^{2}) \) with
\( \gp_{j}^{2} = \CGP \, \CGPa^{-1} j^{2\smpa} \)
for \( \CGP \geq \CONST \CGPa \, \nsize^{1 - 2 \smpa} \).
Suppose the conditions of Proposition~\ref{PconcMLEgenc} and \ref{Lvarusetb} hold for this choice of \( \GP^{2} \).
Then on \( \Omega(\xx) \)
\begin{EQA}
	\| \DPGP (\tilde{\upsv}_{\GP} - \upsvs) \|^{2}
	& \lesssim &
	( \CGP^{-1} \CGPa \, \nsize)^{1/(2\smpa)} + \CGP .
\label{bnjhuy0w9e8r7t6bcgvftr}
\end{EQA}
Furthermore, the choice 
\begin{EQA}
	\CGP 
	&=& 
	(\CGPa \nsize)^{1/(2\smpa+1)} \eqdef \mmps
\label{kvytft64cdfrtdr213ere}
\end{EQA} 
leads to
\begin{EQA}
	\GP^{2} 
	&=& 
	\mmps \, \CGPa^{-1} \diag\{ j^{2\smpa} \}_{j\geq 1} 
	=
	\CGPa^{- 2\smpa/(2\smpa+1)} \, \nsize^{1/(2\smpa+1)} \diag\{ j^{2\smpa} \}_{j\geq 1} 
\label{njbdfhedue437dfhgv6te}
\end{EQA}
and yields for \( \nsize \geq \CONSTn \)
\begin{EQ}[rcl]
	\| \DPGP (\tilde{\upsv}_{\GP} - \upsvs) \|^{2}
	& \lesssim &
	\mmps
	=
	(\CGPa \nsize)^{1/(2\smpa+1)} \, ,
	\\
	\| \tilde{\upsv}_{\GP} - \upsvs \|^{2}
	& \lesssim &
	\mmps \, n^{-1}
	=
	\CGPa^{1/(2\smpa+1)} \, \nsize^{-2\smpa/(2 \smpa + 1)} .
\label{bnjhuy0w9e8r7t6bcgvftrop}
\end{EQ}
\end{theorem}

\begin{proof}
Define \( \mm \) by the relation \( \gp_{\mm}^{2} \approx n \).
This leads to \( \mm \approx (\CGP^{-1} \CGPa \, \nsize)^{1/(2\smpa)} \). 
By Lemma~\ref{LeffdimG} \( \dimG \asymp \mm \).
The condition \( \CGP \geq \CONST \CGPa \, \nsize^{1 - 2 \smpa} \) with a properly selected \( \CONST \)
ensures \( \dimG \ll \nsize \) and hence, \( \rrGP \, n^{-1/2} \) is sufficiently small.
This enables us to apply Proposition~\ref{PconcMLEgenc} and Theorem~\ref{TFiWititG}.  
Furthermore, smoothness of \( \upsvs \in \BBB(\smpa,\CGPa) \) and the definition of 
\( \GP^{2} = \CGP \, \CGPa^{-1} \diag\{ j^{2\smpa} \}_{j\geq 1} \) imply \( \| \GP \upsvs \|^{2} \leq \CGP \)
and hence, \( \| \DP_{\GP}^{-1} \GP^{2} \upsvs \|^{2} \leq \CGP \).
This and \eqref{EQtuGuvus2tr1o1} prove \eqref{bnjhuy0w9e8r7t6bcgvftr}.
Optimizing its right hand-side w.r.t. \( \CGP \) leads to the \( \CGP \asymp \mmps = (\CGPa \nsize)^{1/(2\smpa+1)} \).
The condition \( \mmps \ll \nsize \) reads \( \nsize^{2\smpa} \geq \CONST \CGPa \) and it is automatically fulfilled 
for \( \nsize \) sufficiently large.
This yields the final result \eqref{bnjhuy0w9e8r7t6bcgvftrop}.
\end{proof}

With \( \CGPa \asymp 1 \), this result yields the optimal rate of estimation 
\( \| \tilde{\upsv}_{\GP} - \upsvs \|^{2} \lesssim\nsize^{-2\smpa/(2 \smpa + 1)} \) 
over the smoothness class \( \BBB(\smpa,\CGPa) \) from \eqref{BBBsttsjj1sj2}
for the pMLE \( \tilde{\upsv}_{\GP} \) with \( \GP^{2} \) from \eqref{njbdfhedue437dfhgv6te}.

Further we check the situation when \( \upsvs \in \BBB(\smpa,\CGPa) \) but 
we use \( \GP^{2} = \CGP^{-1} \diag\{ j^{2\smp} \}_{j\geq 1} \) for \( \smp \neq \smpa \).
The main question is a proper choice of the scaling factor \( \CGP \) to ensure the same in order accuracy as 
for the \( \smpa \)-aware choice.
Lemma~\ref{LeffdimG} provides the key argument: the index \( \mmps \) from \eqref{kvytft64cdfrtdr213ere} should be preserved.
Due to Lemma~\ref{Lsmoothbiase}, to make the construction rigorous we also need to ensure \( \smp \geq \smpa \).
This can be achieved by selecting a reasonably large \( \smp \), e.g. \( \smp = 4 \) or \( \smp = 5 \).
Now the approach can be summarized as follows. 
Under the assumption \( \upsvs \in \BBB(\smpa,\CGPa) \) and \( n^{-1} \asymp \| \DP^{-2} \| \),
we fix \( \mmps \) by \eqref{kvytft64cdfrtdr213ere}.
With \( \GP^{2} = \CGP^{-1} \diag\{ j^{2\smp} \} \), the index \( \mm \) corresponding to \( \gp_{\mm}^{2} \approx \nsize \)
is given by \( \mm \approx (\CGP \, \nsize)^{1/(2\smp)} \).
This leads to the relation
\begin{EQA}
	(\CGP \, \nsize)^{1/(2\smp)}
	& \approx &
	\mmps
	=
	(\CGPa \nsize)^{1/(2\smpa+1)}
\label{6drvjhggv74e3ujf6hgewt}
\end{EQA}
and to the solution 
\begin{EQA}
	\CGP
	& \approx &
	\nsize^{-1} \mmps^{2\smp}
	=
	\nsize^{-1} (\CGPa \nsize)^{2\smp/(2\smpa+1)} .
\label{nmgvctf5rtrtdw2gf56} 
\end{EQA}
Alternatively one can apply the \( \mmps \)-truncation prior with \( \mmps \) from \eqref{6drvjhggv74e3ujf6hgewt}.
The next result assumes again \( \upsvs \in \BBB(\smpa,\CGPa) \) but for any \( \smpa > 0 \).

\begin{theorem}
\label{PGLMlossGPg}
Let also \( \upsvs \in \BBB(\smpa,\CGPa) \), \( \smpa > 0 \).
Define 
\( \GP^{2} = \CGP^{-1} \diag\{ j^{2\smp} \}_{j\geq 1} \) for some \( \smp > 1/2 \), \( \smp \geq \smpa \),
and \( \CGP \) from \eqref{6drvjhggv74e3ujf6hgewt} or \eqref{nmgvctf5rtrtdw2gf56}.
Suppose the conditions of Proposition~\ref{PconcMLEgenc} and \ref{Lvarusetb} hold for this choice of \( \GP^{2} \).
Then bounds \eqref{bnjhuy0w9e8r7t6bcgvftrop} continue to hold on \( \Omega(\xx) \) provided \( \nsize \geq \CONSTn \).
\end{theorem}

\begin{proof}
Lemma~\ref{LeffdimG} yields \( \dimG \asymp (\CGP \nsize)^{1/(2\smp)} \asymp \mmps \).
Further, Lemma~\ref{Lsmoothbiase} ensures condition \eqref{hgjkjhyoie3r56dftrtyf} of Proposition~\ref{Pbiassmoone}
\( \DPGP^{-2} \GP^{2} \leq \CONSTGP \DP_{\GPa}^{-2} \GPa^{2} \) for a constant \( \CONSTGP \) depending on \( \smpa \) and \( \smp \) only.
The result follows now similarly to \eqref{bnjhuy0w9e8r7t6bcgvftrop}.
\end{proof}

\Section{Bayesian inference under Sobolev smoothness}
Here we continue the discussion of Section~\ref{SrateSobball} assuming \( \upsvs \in \BBB(\smpa,\CGPa) \) with 
\( \CGPa \asymp 1 \).
A rate-optimal choice of the prior covariance \( \GP^{2} \) is an issue to addrress.

\begin{theorem}
\label{TLaplaceGP}
Assume \nameref{Eref},
\nameref{EU2ref},
\nameref{LLCref}, \nameref{LLtS3ref}.
Let also \eqref{uVmGu2C1C2genS} hold
and \( \upsvs \in \BBB(\smpa,\CGPa) \), \( \smpa > 0 \).
Define 
\( \GP^{2} = \CGP^{-1} \diag\{ j^{2\smp} \}_{j\geq 1} \) for \( \smp > 1/2 \), \( \smp \geq \smpa \),
and \( \CGP \) from \eqref{6drvjhggv74e3ujf6hgewt} or \eqref{nmgvctf5rtrtdw2gf56}.
Then 
for \( \nsize \geq \CONSTn \), it holds on \( \Omega(\xx) \)
\begin{EQA}
	\P\Bigl( \| \vupsv_{\GP} - \tilde{\upsv}_{\GP} \| > \CONST \nsize^{-\smpa/(2\smpa+1)} \Cond \Yv \Bigr)
	& \leq &
	\ex^{-\xx} .
\label{0j7jds2fc56gh98hbt3w}
\end{EQA}
Moreover, if \( \smpa > 1 \), then 
on \( \Omega(\xx) \) for \( \nsize \geq \CONSTn \)
\begin{EQ}[rcl]
	\sup_{A \in \cc{B}(\R^{\dimp})}
	\left| 
		\P\bigl( \vupsv_{\GP} - \tilde{\upsv}_{\GP} \in A \cond \Xv \bigr)
		- \PG\bigl( \DPGPt^{-1} \gammav \in A \bigr)  
	\right|
	& \leq &
	\CONST\, \nsize^{(1 - \smpa) / (2\smpa + 1)} \, ,
	\\
\label{1p1m1a1emxrrde}
	\sup_{A \in \cc{B}_{\smp}(\R^{\dimp})}
	\left| 
		\P\bigl( \vupsv_{\GP} - \tilde{\upsv}_{\GP} \in A \cond \Xv \bigr)
		- \PG\bigl( \DPGPt^{-1} \gammav \in A \bigr)  
	\right|
	& \leq &
	\CONST\, \nsize^{(2 - 2 \smpa) / (2\smpa + 1)} \, .
	\qquad
\end{EQ}
\end{theorem}

\begin{proof}
We proceed similarly to Theorem~\ref{CGLMlossGP} and Theorem~\ref{PGLMlossGPg}. 
Let \( \mmps \) be given by \eqref{6drvjhggv74e3ujf6hgewt}.
Then \( \dimL(\upsv) \asymp \mmps = (\CGPa \nsize)^{1/(2\smpa+1)} \) yielding \eqref{poybf3679jd532ff2st} for \( \nsize \geq \CONSTn \).
For validity of Laplace approximation \eqref{scdugfdwyd2wywy26e6de432} in Theorem~\ref{TBvMgen34},
we have to only check an addition condition \eqref{0hcde4dftedrf94yr5twew}.
It can be spelled out as \( \mmps^{3} \ll \nsize \).
For \( \smpa > 1 \), it is fulfilled automatically if \( \nsize \) is sufficiently large.
Bounds in \eqref{1p1m1a1emxrrde} are obtained from \eqref{scdugfdwyd2wywy26e6de432} by using \( \dimLt \asymp \mmps \).
\end{proof}


\Chapter{Anisotropic logistic regression}
\label{Saniligot}
This \chname specifies the general results for a popular logistic regression model. 
It is widely used e.g. in binary classification in machine learning for binary classification or in binary response model in econometrics.
The results presented here can be viewed as extension of \cite{SpPa2019} to the case 
of non-regular models.
Such an extension is very important for practical applications where the underlying probability of success 
can be very close to zero or one.
Most of existing theoretical studies fail in such cases and the analysis becomes very involved. 
Moreover, another type of limiting behavior of the posterior may happen; see e.g. \cite{bochkina2014} for some examples
of model with non-Gaussian posterior approximation in non-regular regression.
Our results indicate that the use of a sufficiently strong Gaussian prior on the canonical parameter vector allows 
to eliminate such cases and to establish a uniform Gaussian approximation of the posterior.

\Section{Setup and conditions}
Suppose we are given a vector of independent observations/labels \( \Yv = (Y_{1},\ldots,Y_{n})^{\T} \) and a set of the corresponding features vectors
\( \Psiv_{i} \in \R^{\dimp} \).
Each binary label \( Y_{i} \) is modelled as a Bernoulli random variable with the parameter 
\( \thetas_{i} = \P(Y_{i} = 1) \).
Logit regression operates with the canonical parameter \( \upss_{i} = \log \frac{\thetas_{i}}{1-\thetas_{i}} \) and 
assumes a linear dependence \( \upss_{i} = \langle \Psiv_{i}, \upsvs \rangle \) for the parameter vector
\( \upsv \in \R^{\dimp} \).
The corresponding log-likelihood reads 
\begin{EQA}
	L(\upsv)
	&=&
	\sumi \Bigl\{ Y_{i} \, \langle \Psiv_{i}, \upsv \rangle - \cdens\bigl( \langle \Psiv_{i}, \upsv \rangle \bigr) \Bigr\}
\label{LuYiPsiGLM}
\end{EQA}
with \( \cdens(\ups) = \log\bigl( 1 + \ex^{\ups} \bigr) \).
For simplicity we assume that the \( \Psiv_{i} \) are deterministic, otherwise we need to condition on them.

A penalized MLE \( \tilde{\upsv}_{\GP} \) is defined by maximization
of the penalized log-likelihood \( \LGP(\upsv) = L(\upsv) - \| \GP \upsv \|^{2}/2 \)
for the quadratic penalty \( \| \GP \upsv \|^{2}/2 \):
\begin{EQA}
	\tilde{\upsv}_{\GP} 
	&=& 
	\argmax_{\upsv \in \R^{\dimp}} \LGP(\upsv) .
\label{vb8rj4ru7yfueow3e}
\end{EQA}
In the Bayesian setup we consider a Gaussian prior \( \ND(0,\GP^{-2}) \) on \( \upsv \).
The posterior density is given by
\begin{EQA}
	\priord(\upsv \cond \Yv)
	&=&
	\frac{\exp {\LGP(\upsv)}}{\int \exp \LGP(\upsv) \, d\ups} \, .
\label{taucYLtuiL}
\end{EQA}
Inference for the pMLE \( \tilde{\upsv}_{\GP} \) and for the posterior density
is well studied in
regular situations when the true parameter \( \upsvs \) belongs to a bounded set \( \Ups \);
see e.g. \cite{SP2013_rough, SpPa2019}.
Here we explain how the general results of \Chname\ref{SgenBounds} and \Chname\ref{SgenBvMs} apply in the GLM case.

\Section{Check of general conditions}
\label{ScondGLM}

The truth and the penalized truth are defined via the expected log-likelihood
\begin{EQA}
	\upsvs 
	&=& 
	\argmax_{\upsv \in \R^{\dimp}} \E L(\upsv) ,
	\\
	\upsvs_{\GP} 
	&=& 
	\argmax_{\upsv \in \R^{\dimp}} \E \LGP(\upsv) .
\label{vb8rj4ru7yfueow3eEE}
\end{EQA}
The Fisher information matrix \( \IF(\upsv) \) at \( \upsv \) is given by
\begin{EQA}
	\IF(\upsv)
	&=&
	\sumi \cdens^{(2)}\bigl( \langle \Psiv_{i}, \upsv \rangle \bigr) \, \Psiv_{i} \, \Psiv_{i}^{\T} .
\label{IFGusiphpPiuvT}
\end{EQA}
Note that the matrix \( \IF(\upsv) \) is not diagonal even if \( \Psiv \Psiv^{\T} \) is diagonal.
Indeed, \( j,\jc \)-element of \( \IF(\upsv) = - \nabla^{2} \E L(\upsv) \)
is \( \IF(\upsv)_{j,\jc} = \sumi \psi_{i,j} \, \psi_{i,\jc} \, \cdens^{(2)}\bigl( \langle \Psiv_{i}, \upsv \rangle \bigr) \).
The penalized Fisher information matrix \( \IF_{\GP}(\upsv) \) at \( \upsv \) is given by
\begin{EQA}
	\IF_{\GP}(\upsv)
	&=&
	\IF(\upsv) + \GP^{2} .
\label{IFGusiphpPiuvGP}
\end{EQA}
As usual, we write \( \DPGP^{2} = \DPGP^{2}(\upsvs_{\GP}) = \IF(\upsvs_{\GP}) \).

Let \( \zeta(\upsv) = L(\upsv) - \E L(\upsv) = \sumi \bigl( Y_{i} - \E Y_{i} \bigr) \langle \Psiv_{i},\upsv \rangle \)
be the stochastic component of \( L(\upsv) \).
It is obviously linear in \( \upsv \) with
\begin{EQA}
	\nabla \zeta
	&=&
	\sumi \bigl( Y_{i} - \E Y_{i} \bigr) \Psiv_{i} \, .
\label{nabzesiYiEYiPi}
\end{EQA}
Therefore, \nameref{Eref} is granted by construction.
Convexity of \( \cdens(\cdot) \) yields concavity of \( L(\upsv) \) and thus, \nameref{LLCref}.
Now we check \nameref{EU2ref}.
It holds with \( \eps_{i} = Y_{i} - \E Y_{i} \)
\begin{EQA}
	\VP^{2}
	&=&
	\Var(\nabla \zeta)
	=
	\sumi \Var(\eps_{i}) \, \Psiv_{i} \Psiv_{i}^{\T} 
	=
	\sumi \thetas_{i} (1 - \thetas_{i}) \, \Psiv_{i} \Psiv_{i}^{\T} .
\label{V2siVeiPiT2}
\end{EQA}
Note that under correct GLM specification
\( Y_{i} \sim \Bernoulli(\thetas_{i}) \)
with \( \thetas_{i} = \ex^{\langle \Psiv_{i},\upsvs \rangle}/(1 + \ex^{\langle \Psiv_{i},\upsvs \rangle}) \) for
the true parameter \( \upsvs \), it follows \( \VP^{2} = \IF(\upsvs) \).

Define the \emph{effective dimension}
\begin{EQA}
	\dimVG 
	&=& 
	\tr\bigl( \DPGP^{-2} \VP^{2} \bigr) .
\label{ujf8uejecy7whwejd}
\end{EQA}
Let \( \xx \) be such that 
\begin{EQA}
	1 / \wPsi_{\GP} 
	& \geq &
	\xx^{1/2}/2 + (\xx \, \dimVG/4)^{1/4}
\label{e3gtfetyhwsfsedserd}
\end{EQA} 
with 
\( \wPsi_{\GP} = \max_{i} \| \DPGP^{-1} \Psiv_{i} \| \).
By Lemma~\ref{LHDxibound}, on a set \( \Omega(\xx) \) with 
\( \P\bigl( \Omega(\xx) \bigr) \geq 1 - 3 \ex^{-\xx} \) 
\begin{EQA}
	\| \DPGP^{-1} \nabla \zeta \|
	& \leq &
	\rrVG 
	= 
	\sqrt{\dimVG} + \sqrt{2\xx} \, ,
\label{wWeightsde}
\end{EQA}
yielding \nameref{EU2ref}.

We need one more condition on regularity of the design \( \Psiv_{1},\ldots,\Psiv_{n} \).
For any \( \upsv \in \Ups \), define \( \weight_{i}(\upsv) = \cdens^{(2)}\bigl( \langle \Psiv_{i}, \upsv \rangle \bigr) \) so that
\begin{EQA}
	\DP^{2}(\upsv)
	&=&
	\sumi \Psiv_{i} \, \Psiv_{i}^{\T} \weight_{i}(\upsv)
\label{?}
\end{EQA}
and for any vector \( \uv \)
\begin{EQA}
	\sumi \langle \Psiv_{i} , \uv \rangle^{2} \weight_{i}(\upsv)
	&=&
	\| \DP(\upsv) \uv \|^{2} .
\label{nqpAGusitGDm2}
\end{EQA}
Design regularity usually treated in the sense that each summand 
\( \langle \Psiv_{i} , \uv \rangle^{2} \weight_{i}(\upsv) \)
is smaller in order that the whole sum.
We also need a bound on the fourth moments of the design distribution. 
Suppose a subset \( \Upsd \subset \Ups \) is fixed.

\begin{description}
\item[\label{Psiref}\( \bb{(\Psi)} \)]
\begin{itemize}
\item[(i)]
For some \( \CONSTi_{n} \) and all \( \upsv \in \Upsd \)
\begin{EQA}
	\max_{i \leq n} \bigl| \langle \Psiv_{i}, \uv \rangle \bigr|
	& \leq &
	\frac{\CONSTi_{n}}{n^{1/2}} \, \| \DP(\upsv) \uv \|,
	\qquad
	\uv \in \R^{\dimp} .
\label{maxiPsiu214DG21}
\end{EQA}
\item[(ii)]
For some \( \CONSTPsi \) and all \( \upsv \in \Upsd \),  
\begin{EQA}
	\frac{1}{n} \sumi \bigl\langle \Psiv_{i}, \uv \bigr\rangle^{4} \, \weight_{i}(\upsv)
	& \leq &
	\left( \frac{\CONSTPsi}{n} \, \bigl\| \DP(\upsv) \uv \bigr\|^{2} \right)^{2} ,
	\qquad
	\uv \in \R^{\dimp} .
\label{k23CPpHkuf2}
\end{EQA}
\end{itemize}
%
\end{description}

It is straightforward to see that \( \CONSTi_{n} \geq \CONSTPsi \).
However, this bound is usually too rough and it is useful to separate these constants.
 
Now we check \nameref{LLtS3ref} and \nameref{LLtS4ref}.

\begin{lemma}
\label{LcondcheckGLM}
Assume \nameref{Psiref}.
Consider 
\begin{EQ}[rcccl]
	\dimL(\upsv)
	& \eqdef &
	\tr \bigl\{ \DP^{2}(\upsv) \DPGP^{-2}(\upsv) \bigr\},
	\qquad
	\rrL(\upsv)
	& \eqdef & 
	2 \sqrt{\dimL(\upsv)} + \sqrt{2 \xx} \, ,
\label{tdgtftyfyt4234w34resGLM}
\end{EQ}
and assume for all \( \upsv \in \Upsd \)
\begin{EQ}[rcl]
\label{dw3Gmin2211eGLM}
	\amax^{-1} \, 
	\CONSTi_{n} \, \rr(\upsv) \, n^{-1/2}
	& \leq &
	{1}/{2} \, ,
	\\
	\amax^{-1} \, \sqrt{\ex} \, \CONSTPsi \, \rr(\upsv) \, n^{-1/2}
	& \leq & 
	1/3 .
\end{EQ}
Then \nameref{LLtS3ref} and \nameref{LLtS4ref} hold with \( \HL^{2}(\upsv) = \nsize^{-1} \DP^{2}(\upsv) \), 
\( \hmax_{3} = \sqrt{\ex} \, \CONSTPsi \), and
\( \hmax_{4} = \sqrt{\ex} \, \CONSTPsi^{2} \).
Moreover, \( \amax^{-1} \, \hmax_{3} \, \rr(\upsv) \, n^{-1/2} \leq 1/3 \).
\end{lemma}

\begin{proof}
We start with some technical statements.
First observe that 
the function \( \cdens(\ups) = \log(1 + \ex^{\ups}) \) satisfies for all \( \ups \in \R \)
\begin{EQA}
	|\cdens ^{(k)}(\ups)|
	& \leq &
	\cdens^{(2)}(\ups) ,
	\qquad
	k=3,4.
\label{phukulppu43}
\end{EQA}
Indeed, it holds
\begin{EQA}
	\cdens'(\ups)
	&=&
	\frac{\ex^{\ups}}{1 + \ex^{\ups}} \, ,
	\\
	\cdens^{(2)}(\ups)
	&=&
	\frac{\ex^{\ups}}{(1 + \ex^{\ups})^{2}} \, ,
	\\
	\cdens^{(3)}(\ups)
	&=&
	\frac{\ex^{\ups}}{(1 + \ex^{\ups})^{2}} - \frac{2 \ex^{2\ups}}{(1 + \ex^{\ups})^{3}} \, ,
	\\
	\cdens^{(4)}(\ups)
	&=&
	\frac{\ex^{\ups}}{(1 + \ex^{\ups})^{2}} - \frac{6 \ex^{2\ups}}{(1 + \ex^{\ups})^{3}}
	+ \frac{6 \ex^{3\ups}}{(1 + \ex^{\ups})^{4}} \, .
\label{6314412f2e21e}
\end{EQA}
It is straightforward to see that 
\( | \cdens^{(k)}(\ups)| \leq \cdens^{(2)}(\ups) \) for \( k = 3,4 \) and any \( \ups \).

Next we check local variability of \( \cdens^{(2)}(\ups) \).
Namely, for any \( \rhou > 0 \), any \( \upsd \in \R \), 
\begin{EQA}
	\sup_{|\ups - \upsd| \leq \rhou} \bigl| \cdens^{(2)}(\ups) \bigr|
	& \leq &
	\ex^{\rhou} \cdens^{(2)}(\upsd) .
\label{takrrXetmk}
\end{EQA}
To see this, suppose that \( \upsd < 0 \).
As the function \( \cdens^{(2)}(\ups) \) is monotonously increasing in \( \ups < 0 \), it holds 
\begin{EQA}
	\sup_{|\ups - \upsd| \leq \rhou} \frac{\cdens^{(2)}(\ups)}{\cdens^{(2)}(\upsd)}
	& = &
	\frac{\cdens^{(2)}(\upsd + \rhou)}{\cdens^{(2)}(\upsd)}
	\leq 
	\ex^{\rhou}
\label{exrhfrfppusg}
\end{EQA}
and the result follows.

Putting together \eqref{phukulppu43} and \eqref{exrhfrfppusg} leads to a bound on variability of
\( \DP(\upsv) \).
Let us fix \( \upsv \in \Upsd \) and any \( \uv \) with \( \| \DP(\upsv) \uv \| \leq \rr(\upsv) \).
By definition
\begin{EQA}
	\DP^{2}(\upsv + \uv)
	&=&
	\sumi \Psiv_{i} \, \Psiv_{i}^{\T} \, \cdens^{(2)}\bigl( \langle \Psiv_{i}, \upsv + \uv \rangle \bigr) .
\label{0bvlgfiu84ehdf6ew3nh}
\end{EQA}
By \eqref{dw3Gmin2211eGLM}, for each \( i \leq n \), it holds 
\( \bigl| \langle \Psiv_{i}, \uv \rangle \bigr| \leq 1/2 \) and by \eqref{exrhfrfppusg}
\begin{EQA}
	\cdens^{(2)}\bigl( \langle \Psiv_{i}, \upsv + \uv \rangle \bigr)
	& \leq &
	\sqrt{\ex} \,\, \cdens^{(2)}\bigl( \langle \Psiv_{i}, \upsv \rangle \bigr) .
\label{g3d45de5r3ted5tvuye}
\end{EQA}
This yields 
\begin{EQA}
	\DP^{2}(\upsv + \uv)
	& \leq &
	\sqrt{\ex} \, \DP^{2}(\upsv) .
\label{e12f2Pivf2Piiv12}
\end{EQA}
As the next step we evaluate the derivative \( \nabla^{k} f(\upsv) \) for \( f(\upsv) = \E L(\upsv) \) and \( k \geq 2 \).
Let \( f(\upsv) = \E L(\upsv) \).
For any \( \upsv \in \CA_{\GP} \) and any \( \uv_{1},\ldots,\uv_{k} \in \UV \), \( k \geq 2 \),
\begin{EQA}
	\bigl\langle \nabla^{k} f(\upsv), \uv_{1} \otimes \ldots \otimes \uv_{k} \bigr\rangle
	& = &
	- \sumi \prod_{j=1}^{k} \langle \Psiv_{i}, \uv_{j} \rangle \,\, \cdens^{(k)}\bigl( \langle \Psiv_{i}, \upsv \rangle \bigr),
\label{Cparfxa2IaGusd}
\end{EQA}
and, in particular,
\begin{EQA}
	\| \DP(\upsv) \uv \|^{2}
	&=&
	- \bigl\langle \nabla^{2} f(\upsv), \uv^{\otimes 2} \bigr\rangle
	=
	\sumi \langle \Psiv_{i}, \uv \rangle^{2} \, \cdens^{(2)}\bigl( \langle \Psiv_{i}, \upsv \rangle \bigr).
\label{D2fxa2IaGusd}
\end{EQA}
With \( \weight_{i}(\upsv) = \cdens^{(2)}\bigl( \langle \Psiv_{i}, \upsv \rangle \bigr) \),
for \( \upsv \in \Upsd \) and any \( \uv,\wv \in \UV \), we derive by (ii) of \nameref{Psiref}, \eqref{phukulppu43},
and \eqref{e12f2Pivf2Piiv12} for \( t \in [0,1] \)
\begin{EQA}
	&& \nquad
	\bigl| \bigl\langle \nabla^{3} f(\upsv + t \wv), \wv \otimes \uv^{\otimes 2} \bigr\rangle \bigr|
	\leq
	\sumi \bigl| \langle \Psiv_{i}, \wv \rangle \bigr| \, \langle \Psiv_{i}, \uv \rangle^{2} \,\,
	\cdens^{(2)}\bigl( \langle \Psiv_{i}, \upsv + t \wv \rangle \bigr)
	\\
	& \leq &
	\sqrt{\ex} \sumi \bigl| \langle \Psiv_{i}, \wv \rangle \bigr| \, \langle \Psiv_{i}, \uv \rangle^{2} \,\,
	\weight_{i} 
	\leq 
	\sqrt{\ex} \biggl( 
		\sumi \langle \Psiv_{i}, \wv \rangle^{2} \, \weight_{i}(\upsv)
	\biggr)^{1/2}
	\biggl( 
		\sumi \langle \Psiv_{i}, \uv \rangle^{4} \, \weight_{i}(\upsv)
	\biggr)^{1/2}
	\\
	& \leq &
	\sqrt{\ex} \, \frac{\CONSTPsi}{n^{1/2}} \, \| \DP(\upsv) \wv \| \, \| \DP(\upsv) \uv \|^{2} .
\label{ph2Piu2Pief2}
\end{EQA}
This yields \nameref{LLtS3ref} with \( \HL^{2}(\upsv) = \nsize^{-1} \DP^{2}(\upsv) \) and \( \hmax_{3} = \sqrt{\ex} \, \CONSTPsi \).
Similarly \nameref{LLtS4ref} holds with \( \hmax_{4} = \sqrt{\ex} \, \CONSTPsi^{2} \).
\end{proof}

\Section{General properties}
All the conditions of Proposition~\ref{PconcMLEgenc} and Proposition~\ref{Lvarusetb} have been already checked
for GLM.
The results of Section~\ref{SgenBounds} for the pMLE \( \tilde{\upsv}_{\GP} \) can be summarized as follows.

\begin{theorem}
\label{TBernGLM}
Let \( Y_{i} \) be independent Bernoulli, \( Y_{i} \sim \Bernoulli(\thetas_{i}) \).
Consider GLM \eqref{LuYiPsiGLM}.
Given \( \GP^{2} \), assume \nameref{Psiref}.
Let also with \( \amax = 2/3 \) and \( \rrL(\upsv) \) from \eqref{tdgtftyfyt4234w34resGLM},
conditions \eqref{dw3Gmin2211eGLM} hold true.
Then on a set \( \Omega(\xx) \) with \( \P(\Omega(\xx)) \geq 1 - 3 \ex^{-\xx} \),
the concentration bound \eqref{rhDGtuGmusGU0} of Proposition~\ref{PconcMLEgenc}
and the Fisher-Wilks expansions \eqref{3d3Af12DGttG} and \eqref{DGttGtsGDGm13rG} of Theorem~\ref{TFiWititG}
apply. 
Bound \eqref{rG1r1d3GbvGb} of Theorem~\ref{PFiWibias} for the loss \( \tilde{\upsv}_{\GP} - \upsvs \) 
and bound \eqref{EQtuGmstrVEQtG} of Theorem~\ref{PFiWi2risk} for the corresponding risk  apply similarly.
\end{theorem}

Now we continue with the properties of the posterior \( \vupsv_{\GP} \cond \Xv \) for a Gaussian prior \( \ND(0,\GP^{-2}) \).
Theorem~\ref{Tconcdens} ensures a concentration of \( \tilde{\upsv}_{\GP} \) on a local set \( \CA_{\GP} \) which is a subset 
of the local ball \( \BBB_{\rsmall}(\upsvs_{\GP}) \).
The next result is just a specification of the general statements of Theorem~\ref{TBvMgen34}.

\begin{theorem}
\label{TBernGLMLapl}
Assume the conditions of Theorem~\ref{TBernGLM}.
The the concentration result \eqref{poybf3679jd532ff2st} of Proposition~\ref{PGaussconcge} 
as well as the contraction bound \eqref{PDvtGtsCrGYgen} of Proposition~\ref{Ccontactionrategen}
hold on \( \Omega(\xx) \).
Moreover, under \eqref{c6w3gtewr5e4e4wgewrer}, Laplace approximation 
\eqref{scdugfdwyd2wywy26e6de432} of Theorem~\ref{TBvMgen34} continues to apply.
\end{theorem}

\Section{Truncation and smooth priors}
\label{SregbalanceGLM}
Here we specify the results to the case of a regular design when all the eigenvalues of 
the matrix information \( \IF(\upsv) = \DP^{2}(\upsv) \) for all \( \upsv \in \CA_{\GP} \) are of order \( n \). 
With smooth or truncation priors of Section~\ref{Spriors},
we now state a finite sample version of 
the standard nonparametric \emph{rate optimal} results about concentration of the pMLE
and posterior contraction; cf \cite{CaNi2014}, \cite{CaRo2015}.
We use that \( \DPGPt^{2} \geq \nsize \, \CONSTPsi^{2} \, \Id_{\dimp} \)
and \( \nsize^{-1} \zq_{\GP}^{2} \asymp \nsize^{-1} \mmps \asymp \CGPa^{1/(2\smpa+1)} \nsize^{-2\smpa/(2\smpa+1)} 
\to 0 \) as \( \nsize \to \infty \) for any \( \smpa > 0 \).
The error term is of order \( \nsize^{(2 - 2 \smpa) / (2\smpa + 1)} \) 
and it tends to zero only as \( \nsize \to \infty \) only if \( \smpa > 1 \).

\begin{theorem}
\label{Tconcdensro}
Assume conditions of Theorem~\ref{TBernGLM}.
Define \( \mmps = (\CGPa \nsize)^{1/(2\smpa+1)} \).
For the \( \mm \)-truncation prior with \( \mm = \mmps \)
or for the \( (\smp,\CGP) \)-smooth prior with \( (\CGP \, \nsize)^{1/(2\smp)} = (\CGPa \, \nsize)^{1/(2\smpa+1)} = \mmps \), 
it holds on the same set \( \Omega(\xx) \) 
for any \( \upsv \in \BBB(\smpa,\CGPa) \)
\begin{EQA}
    \P\Bigl( \| \vupsv_{\GP} - \upsvs \|^{2}  
    		&>& \CONST \CGPa^{1/(2\smpa+1)} \nsize^{-2\smpa/(2\smpa+1)}
    	\cond \Yv 
    \Bigr)
    \leq 
    3 \ex^{-\xx} \, .
\label{DGttGtsCzGns}
\end{EQA}
\end{theorem}


\def\nsizen{\nsize_{0}}

\Chapter{Log-density estimation}
\label{SGBvM}
Suppose we are given a random sample \( X_{1},\ldots,X_{\nsize} \) in \( \R^{d} \).
The density model assumes that all these random variables are independent 
identically distributed from some measure \( \Pone \) with a density 
\( \dens(\xv) \) with respect to a \( \sigma \)-finite measure \( \Pdom \) in \( \R^{d} \).
This density function is the target of estimation.
By definition, the function \( \dens \) is non-negative, measurable, and integrates to one:
\( \int \dens(\xv) \, d\Pdom(\xv) = 1 \).
%
Here and below, the integral \( \int \) without limits means the integral over the whole space \( \R^{d} \).
If \( \dens(\cdot) \) has a smaller support \( \XX \), one can restrict integration to this set.
%
Below we parametrize the model by a linear decomposition of the log-density function.
Let \( \bigl\{ \psi_{j}(\xv), \, j=1,\ldots,\dimp \bigr\} \) with \( \dimp \leq \infty \)
be a collection of functions in \( \R^{d} \) (a dictionary).
%
For each \( \upsv = (\ups_{j}) \in \R^{\dimp} \), define   
\begin{EQA}
    \ldens(\xv,\upsv)
    &\eqdef &
    \ups_{1} \psi_{1}(\xv) + \ldots + \ups_{\dimp} \psi_{\dimp}(\xv) - \cdens(\upsv)
    =
    \bigl\langle \Psiv(\xv), \upsv \bigr\rangle - \cdens(\upsv),
\label{logdenssumj}
\end{EQA}
where \( \Psiv(\xv) \) is a vector with components \( \psi_{j}(\xv) \)
and \( \cdens(\upsv) \) is given by 
\begin{EQA}
    \cdens(\upsv)
    & \eqdef &
    \log \int \ex^{\langle \Psiv(\xv), \upsv \rangle} \, d\Pdom(\xv) .
\label{gtdelointTPxm0}
\end{EQA}
It is worth stressing that the data point \( \xv \) only enters in the linear term \( \bigl\langle \Psiv(\xv), \upsv \bigr\rangle \)
of the log-likelihood \( \ldens(\xv,\upsv) \).
The function \( \cdens(\upsv) \) is entirely model-driven.
Below we restrict \( \upsv \) to a subset \( \Ups \) in \( \R^{\dimp} \) such that 
\( \cdens(\upsv) \) is well defined and the integral of 
\( \ex^{\langle \Psiv(\xv), \upsv \rangle} \) is finite.
Linear log-density modeling assumes  
\begin{EQA}
  	\log \dens(\xv)
  	&=&
  	\ldens(\xv,\upsvs)
	=
	\bigl\langle \Psiv(\xv), \upsvs \bigr\rangle - \cdens(\upsvs)
\label{dnesdxtsl}
\end{EQA}  
for some \( \upsvs \in \Ups \subseteq \R^{\dimp} \).
A nice feature of such representation is that the function \( \log \dens(\xv) \) in the contrary to the density itself does not need to be non-negative.
One more important benefit of using the log-density is that the stochastic part of the corresponding log-likelihood is \emph{linear} w.r.t. the parameter \( \upsv \).
With \( \Spsi = \sumi \Psiv(X_{i}) \), for a given penalty operator \( \GP^{2} \),
the penalized log-likelihood \( \LGP(\upsv) \) reads as
\begin{EQA}
    \LGP(\upsv) 
    &=& 
    \sumi \bigl\langle \Psiv(X_{i}), \upsv \bigr\rangle - \nsize \cdens(\upsv) - \frac{1}{2} \| \GP \upsv \|^{2}
    = 
    \langle \Spsi, \upsv \rangle - \nsize \cdens(\upsv) - \frac{1}{2} \| \GP \upsv \|^{2} .
\label{density_likelyhood}
\end{EQA}
The penalized MLE 
\( \tilde{\upsv}_{\GP} \) and its population counterpart \( \upsvs_{\GP} \) are defined as
\begin{EQ}[rcccl]
    \tilde{\upsv}_{\GP}
    &=&
    \argmax_{\upsv \in \Ups} L_{\GP}(\upsv) ,
    \qquad
    \upsvs_{\GP}
    &=&
    \argmax_{\upsv \in \Ups} \E L_{\GP}(\upsv) .
\label{tsatTqLGttT122}
\end{EQ}

\Section{Conditions}
For applying the general results of Section~\ref{SgenBounds} and Section~\ref{SgenBvMs}, it suffices to check
the general conditions of Section~\ref{SgenBounds} for the log-density model.
First note that the generalized linear structure of the model automatically yields 
conditions \nameref{LLCref} and \nameref{Eref}.
Indeed, convexity of \( \cdens(\cdot) \) implies that 
\( \E L(\upsv) = \langle \E \Spsi, \upsv \rangle - \nsize \cdens(\upsv) \) is concave. 
Further, for the stochastic component \( \zeta(\upsv) = L(\upsv) - \E L(\upsv) \), it holds 
\begin{EQA}
\label{stochDensDef}
    \nabla \zeta(\upsv)
    &=&
    \nabla \zeta
    =
    \Spsi - \E \Spsi
    = 
    \sumi \bigl[\Psiv(X_{i}) - \E \, \Psiv(X_{i})\bigr],
\label{ztLtEltnztSPXi}
\end{EQA}
and \nameref{Eref} follows.
Further,
the representation \( \E L(\upsv) = \langle \E \Spsi, \upsv \rangle - \nsize \cdens(\upsv) \) implies
\begin{EQA}
	\IF(\upsv)
	&=&
	- \nabla^{2} \E L(\upsv)
	=
	- \nabla^{2} L(\upsv)
	=
	n \nabla^{2} \cdens(\upsv) .
\label{IFttn2ELtnn2cdt}
\end{EQA}
To simplify our presentation, we assume that \( X_{1},\ldots,X_{\nsize} \) are indeed i.i.d.
and the density \( \dens(\xv) \) can be represented in the form \eqref{dnesdxtsl} 
for some parameter vector \( \upsvs \).
This can be easily extended to non i.i.d. case at cost of more complicated notations.
Then
\begin{EQA}
    \upsvs
    & = &
    \argmax_{\upsv \in \Ups} \E L(\upsv)
    =
    \argmax_{\upsv \in \Ups} 
      	\bigl\{ \langle \E \Spsi, \upsv \rangle - \nsize \cdens(\upsv) \bigr\} 
	=
    \argmax_{\upsv \in \Ups} 
      	\bigl\{ \langle \Psimean, \upsv \rangle - \cdens(\upsv) \bigr\} ,
	\qquad
\label{ttsT1dedef}
\end{EQA}
where \( \Psimean = \E \, \Psiv(X_{1}) \) and  \( \E \Spsi = \nsize \, \Psimean \).
This yields the identity
\begin{EQA}
	\nabla \cdens(\upsvs)
	&=&
	\Psimean .
\label{cgsuwywdwdy2whq21ef2w}
\end{EQA}
Moreover, by \eqref{gtdelointTPxm0},
\( \nabla^{2} \cdens(\upsvs) = \Var \bigl\{ \Psiv(X_{1}) \bigr\} \) and
\begin{EQA}
	\VP^{2}
	&=&
	\Var(\nabla \zeta)
	=
	\nsize \, \nabla^{2} \cdens(\upsvs) 
	=
	\IF(\upsvs) .
\label{3094857ythgfjdkleo}
\end{EQA}
Here we present our conditions.
For any \( \upsv \in \Ups \) and \( \rsmall > 0 \), define \( \HL(\upsv) \) by \( \HL^{2}(\upsv) = \nabla^{2} \cdens(\upsv) \) and 
consider the corresponding balls in \( \R^{\dimp} \)
\begin{EQA}
	\B_{\rsmall}(\upsv)
	& \eqdef &
	\bigl\{ \uv \in \R^{\dimp} \colon \| \HL(\upsv) \uv \| \leq \rsmall \bigr\} =
	\bigl\{ \uv \in \R^{\dimp} \colon \langle \nabla^{2} \cdens(\upsv), \uv^{\otimes 2} \rangle \leq \rsmall^{2} \bigr\} .
\label{f765re276tf7we6t8erw63}
\end{EQA}

\begin{description}
\item[\label{upsvsref} \( \bb{(\fs)} \)]
\( X_{1},\ldots,X_{\nsize} \) are i.i.d. from a density \( \fs \) satisfying
\( \log \fs(\xv) = \Psiv(\xv)^{\T} \upsvs - \cdens(\upsvs) \).

\item[\label{Upsref} \( \bb{(\Ups)} \)]
The set \( \Ups \) is open and convex, 
the value \( \cdens(\upsv) \) from \eqref{gtdelointTPxm0} is finite for all \( \upsv \in \Ups \),
\( \upsvs \) from \eqref{ttsT1dedef} 
is an internal point in \( \Ups \) such that 
\( \B_{2\rsmall}(\upsvs) \subset \Ups \) for a fixed \( \rsmall > 0 \).

\item[\label{cdensref} \( \bb{(\cdens)} \)]
For the Bregman divergence \( \cdens(\upsv;\uv) \eqdef \cdens(\upsv + \uv) - \cdens(\upsv) - \langle \nabla \cdens(\upsv),\uv \rangle \),
it holds
\begin{EQA}
	\sup_{\upsv \in \B_{\rsmall}(\upsvs)} \,\, \sup_{\uv \in \B_{2\rsmall}(\upsv)} \exp \cdens(\upsv;\uv)
	& \leq &
	\CONSTi_{\rsmall} \, .
\label{vgtcgdcbwcte32vfurdf}
\end{EQA}
\end{description}

Introduce a measure \( P_{\upsv} \) by the relation:
\begin{EQA}
    \frac{d P_{\upsv}}{d\Pdom}(\xv)
    &=&
    \exp\bigl\{ \bigl\langle \Psiv(\xv), \upsv \bigr\rangle - \cdens(\upsv) 
    \bigr\} .
\label{dPdlxeLtSn}
\end{EQA}
Identity \eqref{gtdelointTPxm0} ensures that \( P_{\upsv} \) is a probabilistic measure.
Moreover, under \eqref{dnesdxtsl}, the data generating measure \( \P \) coincides with \( P_{\upsvs}^{\otimes n} \).

\begin{description}
\item[\label{PsiX4ref} \( \bb{(\Psiv_{4})} \)]
	\emph{There are \( \CONSTi_{\Psi,3} \geq 0 \) and \( \CONSTi_{\Psi,4} \geq 3 \) such that 
	for all \( \upsv \in \B_{\rsmall}(\upsvs) \) and \( \gammav \in \R^{\dimp} \)
\begin{EQA}
	\bigl| E_{\upsv} \bigl\langle \Psiv(X_{1}) - E_{\upsv} \Psiv(X_{1}),\gammav \bigr\rangle^{3} \bigr|
	& \leq &
	\CONSTi_{\Psi,3} \, E_{\upsv}^{3/2} \bigl\langle \Psiv(X_{1}) - E_{\upsv} \Psiv(X_{1}),\gammav \bigr\rangle^{2} ,
	\\
	E_{\upsv} \bigl\langle \Psiv(X_{1}) - E_{\upsv} \Psiv(X_{1}),\gammav \bigr\rangle^{4}
	& \leq &
	\CONSTi_{\Psi,4} \,  E_{\upsv}^{2} \bigl\langle \Psiv(X_{1}) - E_{\upsv} \Psiv(X_{1}),\gammav \bigr\rangle^{2} .
\label{dcief8yfe7eegywddwfeyg}
\end{EQA}
	}
\end{description}

In fact, conditions \nameref{cdensref} and \nameref{PsiX4ref} follow from \nameref{Upsref} and can be considered as a kind of definition 
of important quantities \( \CONSTi_{\rsmall} \), \( \CONSTi_{\Psi,3} \), and \( \CONSTi_{\Psi,4} \) 
which will be used for describing the smoothness properties of \( \cdens(\upsv) \).
The matrix \( \nabla^{2} \cdens(\upsv) \) is supposed well conditioned for \( \upsv \in \B_{\rsmall}(\upsvs) \).
\begin{description}
\item[\label{cdens2ref} \( \bb{(\nabla^{2}\cdens)} \)]
For the information matrix \( \nabla^{2}\cdens(\upsv) \),
it holds with some \( \CONSTIF \geq 1 \)
\begin{EQA}
	\CONSTIF^{-1} \Id_{\dimp}
	& \leq & 
	\nabla^{2} \cdens(\upsv)
	\leq 
	\CONSTIF \Id_{\dimp} \, ,
	\qquad
	\upsv \in \B_{\rsmall}(\upsvs) \, .
\label{3ertfghdfu8fcu78edgsedkid}
\end{EQA}
\end{description}
Later we show that it suffices to check \eqref{3ertfghdfu8fcu78edgsedkid} at \( \upsvs \), then
it will be fulfilled in \( \B_{\rsmall}(\upsvs) \) with a slightly larger constant \( \CONSTIF \).

\Subsection{Check of conditions \nameref{LLtS3ref} and \nameref{LLtS4ref}}
Let \( P_{\upsv} \) be defined by \eqref{dPdlxeLtSn}.
It is straightforward to check that \( E_{\upsv} \Psiv(X_{1}) = \nabla \cdens(\upsv) \) and 
\( \Var_{\upsv}(\Psiv(X_{1})) = \nabla^{2} \cdens(\upsv) \).
Further, if \( \uv \in \B_{\rsmall}(\upsv) \) and  
\( \upsv + \uv \in \Ups \), then 
\begin{EQA}
	\cdens(\upsv + \uv) 
	&=& 
	\log E_{0} \exp\{ \langle \Psiv(X_{1}),\upsv + \uv \rangle \} 
	=
	\log E_{\upsv} \exp\bigl\{ \bigl\langle \Psiv(X_{1}), \uv \bigr\rangle + \cdens(\upsv) \bigr\} .
\label{dfor4u7tgjhkgtirt}
\end{EQA}
This yields in view of \( E_{\upsv} \Psiv(X_{1}) = \nabla \cdens(\upsv) \) that
\( \epsv = \Psiv(X_{1}) - E_{\upsv} \Psiv(X_{1}) \) fulfills
\begin{EQA}
	\log E_{\upsv} \exp ( \langle \epsv,\uv \rangle )
	&=&
	\cdens(\upsv + \uv) - \cdens(\upsv) - \langle E_{\upsv} \Psiv(X_{1}),\uv \rangle
	\\
	&=&	
	\cdens(\upsv + \uv) - \cdens(\upsv) - \langle \nabla \cdens(\upsv),\uv \rangle .
\label{kjhr8dfrikweyyweyerdd}
\end{EQA}

%

\begin{lemma}
\label{L34derivat}
The function \( \cdens(\upsv) \) satisfies for any \( \upsv \in \B_{\rsmall}(\upsvs) \) and \( \gammav \in \R^{\dimp} \)
\begin{EQA}
\label{0det53543wyf7rryem}
	| \langle \nabla^{3}\cdens(\upsv), \gammav^{\otimes 3} \rangle |
	& \leq &
	\dltwaa_{\Psi,3} \, \langle \nabla^{2} \cdens(\upsv), \gammav^{\otimes 2} \rangle^{3/2} \, ,
	\\
	| \langle \nabla^{4}\cdens(\upsv), \gammav^{\otimes 4} \rangle |
	& \leq &
	(\dltwaa_{\Psi,4} - 3) \, \langle \nabla^{2} \cdens(\upsv), \gammav^{\otimes 2} \rangle^{2} .
\label{0det53543wyf7rryedm}
\end{EQA}
Moreover, for any \( \gammav_{1},\gammav \in \R^{\dimp} \)
\begin{EQA}
	| \langle \nabla^{3}\cdens(\upsv), \gammav_{1} \otimes \gammav^{\otimes 2} \rangle |
	& \leq &
	\sqrt{\dltwaa_{\Psi,4} \, \langle \nabla^{2} \cdens(\upsv), \gammav_{1}^{\otimes 2} \rangle} \, \,
	\langle \nabla^{2} \cdens(\upsv), \gammav^{\otimes 2} \rangle .
\label{hjfy6w3e35brdfrw3geg}
\end{EQA}
\end{lemma}

\begin{proof}
Denote \( \epsv = X_{1} - E_{\upsv} X_{1} \).
By \eqref{kjhr8dfrikweyyweyerdd} with \( \uv = t \gammav \) for \( t \) sufficiently small
\begin{EQA}
	\chi(t)
	& \eqdef &
	\log E_{\upsv} \exp ( t \langle \epsv,\gammav \rangle )
	=
	\cdens(\upsv + t \gammav) - \cdens(\upsv) - \langle \nabla \cdens(\upsv), t \gammav \rangle ,
\label{fg8yw8wf78wq87weqf87we8f}
\end{EQA}
and by \nameref{PsiX4ref} with \( \CONSTi_{\Psi,4} \geq 3 \)
\begin{EQA}
	\bigl| \chi^{(3)}(0) \bigr|
	&=&
	\bigl| E_{\upsv} \langle \epsv,\gammav \rangle^{3} \bigr|
	\leq 
	\CONSTi_{\Psi,3} \, E_{\upsv}^{3/2} \langle \epsv,\gammav \rangle^{2} \, ,
	\\
	\bigl| \chi^{(4)}(0) \bigr|
	&=&
	\bigl| E_{\upsv} \langle \epsv,\gammav \rangle^{4} 
	- 3 E_{\upsv}^{2} \langle \epsv,\gammav \rangle^{2} \bigr| 
	\leq 
	(\CONSTi_{\Psi,4} - 3) E_{\upsv}^{2} \langle \epsv,\gammav \rangle^{2} .
\label{hdudyue2g2w5dgeh3etyd}
\end{EQA}
If \( \gammav_{1} \neq \gammav \) then we may proceed in a similar way with the bivariate function 
\( \chi(t_{1},t) = \log E_{\upsv} \exp \bigl\{ t_{1} \langle \epsv,\gammav_{1} \rangle + t \langle \epsv,\gammav \rangle \bigr\} \).
Its mixed derivative at zero satisfies
\begin{EQA}
	\left| \frac{\partial^{3}}{\partial t_{1} \partial t^{2}} \chi(0,0) \right|
	&=&
	\bigl| E_{\upsv} \, \langle \epsv,\gammav_{1} \rangle \, \langle \epsv,\gammav \rangle^{2} \bigr| 
	\leq 
	\bigl\{ E_{\upsv} \langle \epsv,\gammav_{1} \rangle^{2} \, E_{\upsv} \langle \epsv,\gammav \rangle^{4} \bigr\}^{1/2}
\label{w43we4e4ew54ek0i0}
\end{EQA}
and the result follows as well.
\end{proof}

\begin{lemma}
\label{L2derivat}
If \( \upsv \in \B_{\rsmall}(\upsvs) \) then with 
\( \CONSTphi = \sqrt{\CONSTi_{\Psi,4} \, \CONSTi_{\rsmall}} \)
\begin{EQA}
	\sup_{\uvdd \in \B_{\rsmall}(\upsv)} \,\,
	\sup_{\gammav \in \R^{\dimp}} \,\, 
		\frac{\langle \nabla^{2} \cdens(\upsv + \uvdd), \gammav^{\otimes 2} \rangle}
			 {\langle \nabla^{2} \cdens(\upsv), \gammav^{\otimes 2} \rangle} 
	& \leq &
	\CONSTphi \, .
\label{vwehuf8ew87fuwyfe2}
\end{EQA}
\end{lemma}
\begin{proof}
Let \( \langle \nabla^{2} \cdens(\upsv), \uvdd^{\otimes 2} \rangle \leq \rsmall^{2} \).
By \eqref{kjhr8dfrikweyyweyerdd} with \( \epsv = X_{1} - E_{\upsv} X_{1} \)
\begin{EQA}
	\nabla^{2} \cdens(\upsv + \uvdd)
	&=&
	\nabla^{2} \log E_{\upsv} \ex^{\langle \epsv,\uvdd \rangle} 
	=
	\frac{E_{\upsv} \{ \epsv \epsv^{\T} \ex^{\langle \epsv,\uvdd \rangle} \}}
		 {(E_{\upsv} \, \ex^{\langle \epsv,\uvdd \rangle})^{2}} 
	- \frac{E_{\upsv} \{ \epsv \, \ex^{\langle \epsv,\uvdd \rangle} \} \, 
			E_{\upsv} \{ \epsv \, \ex^{\langle \epsv,\uvdd \rangle} \}^{\T} }
		 {( E_{\upsv} \, \ex^{\langle \epsv,\uvdd \rangle} )^{2}}
\label{cskjoihnbve4r5ujiko}
\end{EQA}
and by \eqref{0det53543wyf7rryedm} and \eqref{vgtcgdcbwcte32vfurdf} in view of 
\( E_{\upsv} \, \ex^{\langle \epsv,\uvdd \rangle} \geq 1 \)
\begin{EQA}
	\bigl\langle \nabla^{2} \cdens(\upsv + \uvdd), \gammav^{\otimes 2} \bigr\rangle
	& \leq &
	E_{\upsv} \bigl\{ \langle \epsv,\gammav \rangle^{2} \ex^{\langle \epsv,\uvdd \rangle} \bigr\}
	\\
	& \leq &
	E_{\upsv}^{1/2} \langle \epsv,\gammav \rangle^{4} \, \, E_{\upsv}^{1/2} \ex^{2 \langle \epsv,\uvdd \rangle} 
	\leq 
	\sqrt{\CONSTi_{\Psi,4} \, \CONSTi_{\rsmall}} \,\, \bigl\langle \nabla^{2} \cdens(\upsv), \gammav^{\otimes 2} \bigr\rangle
\label{dpe034985tuykhlgpd}
\end{EQA}
and the assertion follows.
\end{proof}


\begin{lemma}
\label{LS3S4logd}
Let \( \upsv \in \B_{\rsmall}(\upsvs) \) and \( \rr \leq \rsmall \, \sqrt{n} \).
Then \( f(\upsv) = \E_{\upsvs} L(\upsv) \) satisfies \nameref{LLtS3ref} and \nameref{LLtS4ref}
with \( \hL(\upsv) = \langle \nabla \cdens(\upsvs), \upsv \rangle - \cdens(\upsv) \), 
\( \HL^{2}(\upsv) = \nabla^{2} \cdens(\upsv) \),
and constants \( \hmax_{3} \) and \( \hmax_{4} \) depending on 
\( \CONSTi_{\rsmall} \), \( \CONSTi_{\Psi,3} \), and \( \CONSTi_{\Psi,4} \) only.
\end{lemma}

\begin{proof}
Let \( \upsv \in \B_{\rsmall}(\upsvs) \).
For any \( \uv \) with \( \| \HL(\upsv) \uv \| \leq \rr/\sqrt{n} \leq \rsmall \), by \eqref{0det53543wyf7rryem} and \eqref{vwehuf8ew87fuwyfe2}
\begin{EQA}
    \frac{|\langle \nabla^{3} \cdens(\upsv + t \uv),\uv^{ \otimes 3} \rangle|}{\| \HL(\upsv) \uv \|^{3}}
    & \leq &	
    \frac{\dltwaa_{\Psi,3} \, \| \HL(\upsv + t \uv) \uv \|^{3}}{\| \HL(\upsv) \uv \|^{3}}
    \leq 
    \dltwaa_{\Psi,3} \, \CONSTphi^{3/2} \, ,
\label{cuhywe8w3endsd72teu2hwd}
\end{EQA}
and \nameref{LLtS3ref} follows with \( \hmax_{3} = \dltwaa_{\Psi,3} \, \CONSTphi^{3/2} \). 
The proof of \nameref{LLtS4ref} is similar. 
\end{proof}

\Subsection{Check of \nameref{EU2ref}}

Now we check the deviation bound for \( \nabla \zeta = \Sv - \E \Sv \)
under \nameref{upsvsref} and \nameref{Upsref}.
I.i.d. structure of \( \Sv = \sum_{i} X_{i} \) and \eqref{3094857ythgfjdkleo} yield 
\( \Var(\Sv) = \VP^{2} = \nsize \nabla^{2} \cdens(\upsvs) \).
Further, for any \( \uv \in \B_{\rsmall}(\upsvs) \), again by the i.i.d. assumption
and by \eqref{kjhr8dfrikweyyweyerdd}
\begin{EQA}
	\nsize^{-1} \log \E_{\upsvs} \exp\bigl\{ \langle \nabla \zeta ,\uv \rangle \bigr\}
	&=&
	\log \E_{\upsvs} \ex^{\langle \epsv, \uv \rangle} 
	=
	\cdens(\upsvs + \uv) - \cdens(\upsvs) - \langle \nabla \cdens(\upsvs),\uv \rangle .
\label{gw7fere9rh8uyegrh978e9}
\end{EQA}
Fix \( \rr \leq \rsmall \, \nsize^{1/2} \) and consider all \( \uv \) with
\( \nsize \langle \nabla^{2} \cdens(\upsvs), \uv^{\otimes 2} \rangle \leq \rr^{2} \). 
Then by \nameref{LLtS3ref} and \eqref{gtcdsftdfvtwdsefhfdvfrvsewseGP} of Lemma~\ref{LdltwLaGP}
\begin{EQA}
	\cdens(\upsvs + \uv) - \cdens(\upsvs) - \langle \nabla \cdens(\upsvs),\uv \rangle
	& \leq &
	\frac{1 + \hmax_{3} \, \rr \, \nsize^{-1/2}/3}{2} \, \langle \nabla^{2} \cdens(\upsvs), \uv^{\otimes 2} \rangle
	\leq 
	\langle \nabla^{2} \cdens(\upsvs), \uv^{\otimes 2} \rangle 
\label{d0if897erfe87y9u}
\end{EQA}
provided that \( \hmax_{3} \, \rr \leq 3 \nsize^{1/2} \).
This implies \eqref{expgamgm} with \( \gm = \rsmall \sqrt{\nsize} \) and thus, the deviation bound
\eqref{PxivbzzBBroB3} of Theorem~\ref{Tdevboundgm} implies \nameref{EU2ref} for \( \rsmall \sqrt{\nsize} \)
sufficiently large.

\Section{Smoothness, bias-variance trade-off}
To handle the bias term, we impose some smoothness conditions on the underlying density parameter \( \upsvs \);
see Section~\ref{Spriorexample}.
We also limit ourselves to the penalty matrices \( \GP^{2} \) which ensure a kind of bias-variance trade-off.


\begin{description}
	\item[\label{ThetaG0ref} \( \bb{(\GPa \upsvs)} \)]
	\( \| \GPa \upsvs \|^{2} \leq 1 \) for some fixed \( \GPa^{2} \).

\end{description}

Later we follow the suggestion of Section~\ref{Spriorexample} and apply the penalizing matrix \( \GP^{2} = \CGP \GP_{0}^{2} \).
Set \( \DPf^{2} = \nsize \nabla^{2} \cdens(\upsvs) \), \( \DPfGP^{2} = \DPf^{2} + \GP^{2} \).
The particular value \( \CGP = \CGPs \) can be selected by the bias-variance relation \eqref{h0hit5ti85tjrtfg78}.
First we evaluate the bias term. 
This is important to ensure that the point \( \upsvs_{\GP} \) is still in the local vicinity \( \B_{\rsmall}(\upsvs) \). 

\begin{proposition}
\label{Pbiasdens}
Assume \nameref{upsvsref}, \nameref{cdensref}, \nameref{PsiX4ref}, \nameref{Upsref},
\nameref{ThetaG0ref}.
Let \( \GP^{2} = \CGP \GP_{0}^{2} \), \( \DPfGP^{2} = \nsize \nabla^{2} \cdens(\upsvs) + \GP^{2} \), 
and \( \amax^{-1} \CGP^{1/2} \, \nsize^{-1/2} \leq \rsmall \).
Then \( \upsvs_{\GP} \in \B_{\rsmall}(\upsvs) \) and
\begin{EQA}
	\| \DPfGP (\upsvs_{\GP} - \upsvs) \|
	& \leq &
	\amax^{-1} \CGP^{1/2} .
\label{7vcobvgy6e3yghergv6r}
\end{EQA} 
\end{proposition}

\begin{proof}
We intend to apply Proposition~\ref{Lvarusetb} with \( \QP = \DPf \) and
\( \rru = \amax^{-1} \| \DPfGP^{-1} \GP^{2} \upsvs \| \) for \( \amax = 2/3 \).
It holds by \nameref{ThetaG0ref} in view of \( \GP^{2} \leq \DPfGP^{2} \) and \( \GP^{2} = \CGP \GP_{0}^{2} \)
\begin{EQA}
	\rru
	&=&
	\amax^{-1} \| \DPfGP^{-1} \GP^{2} \upsvs \|
	\leq 
	\amax^{-1} \CGP^{1/2} \| \GPa \upsvs \|
	= 
	\amax^{-1} \CGP^{1/2} \, .
\label{udfufujef6he3d6rehjde}
\end{EQA}
Further, \( \amax^{-1} \CGP^{1/2} \, \nsize^{-1/2} \leq \rsmall \)
ensures that the set 
\( \{ \upsv \colon \| \DPf (\upsv - \upsvs) \| \leq \rru \} \) belongs to the ball \( \B_{\rsmall}(\upsvs) \).
Now Lemma~\ref{LS3S4logd} yields \nameref{LLtS3ref} for all \( \upsv \in \B_{\rsmall}(\upsvs) \).
By Lemma~\ref{LfreTay}, it holds 
\( \dltwbs_{\GP} \leq \hmax_{3} \rru \leq \hmax_{3} \, \amax^{-1} \, \CGP^{1/2} \, \nsize^{-1/2} \leq 1/3 \) 
for \( \nsize \geq \nsizen \).
Now Proposition~\ref{Lvarusetb} yields \eqref{7vcobvgy6e3yghergv6r}
and \( \upsvs_{\GP} \in \B_{\rsmall}(\upsvs) \).
Also by Lemma~\ref{LfreTay}, the matrix \( \DPGP^{2} = \DPGP^{2}(\upsvs_{\GP}) \) satisfies 
\( (1 - \dltwbs_{\GP})\DPfGP^{2} \leq \DPGP^{2} \leq (1 + \dltwbs_{\GP}) \DPfGP^{2} \).
\end{proof}

For \( \GP^{2} = \CGP \GP_{0}^{2} \), the effective dimension \( \dimG \) is given by \eqref{u98gfw3wd7uewj4ert7yg}:
\begin{EQA}
	\dimG
	&=&
	\dimA(\CGP)
	=
	\tr (\DPf^{2} \, \DPfGP^{-2})
	=
	\tr\bigl\{ \DPf^{2} (\DPf^{2} + \CGP \GPa^{2})^{-1} \bigr\} .
\label{u98gfw3wd7uewj4ert7ygl}
\end{EQA}
A particular value \( \CGPs \) is defined by the bias-variance relation reads \( \dimA(\CGP) \asymp \CGP \); see \eqref{h0hit5ti85tjrtfg78}.
Fix some \( \xx \) and consider
\begin{EQA}
	\rrGP 
	&=&
	\rr(\CGP) 
	=
	\sqrt{\dimA(\CGP)} + \sqrt{2\xx} .
\label{5r4hffgv743ehft5egt3gty}
\end{EQA}
Now expansion \eqref{rG1r1d3GbvGb} of Theorem~\ref{PFiWibias} applies provided that the concentration set 
\( \| \DPGP (\upsv - \upsvs_{\GP}) \| \leq \amax^{-1} \rrGP \) is contained in \( \B_{\rsmall}(\upsvs_{\GP}) \).
Proposition~\ref{Pbiasdens} allows to use \( \DPfGP \) in place of \( \DPGP \) in this condition. 
In our results, \( \CONST \) stands for a fixed constant depending on the other constants in our conditions 
like \( \smpa \), \( \CONSTIF \), \( \CONSTphi \), \( \rsmall \), and \( \CONSTi_{\Psi,4} \).

\begin{theorem}
\label{Tconcdens}
Assume \nameref{upsvsref}, \nameref{cdensref}, \nameref{PsiX4ref}, \nameref{Upsref},
\nameref{ThetaG0ref}, and let \( \GP^{2} = \CGP \GP_{0}^{2} \).
If \( \amax^{-1} \CGP^{1/2} \, \nsize^{-1/2} \leq \rsmall \) with \( \amax = 2/3 \) and \( \rr(\CGP) \nsize^{-1/2} \leq \rsmall \),
then on \( \Omega(\xx) \)
\begin{EQA}
	\| \DPfGP (\tilde{\upsv}_{\GP} - \upsvs) \|
	& \leq &
	2 \amax^{-1} \rr(\CGP) + \amax^{-1} \CGP^{1/2} .
\label{pogpogtr57645ryrtfjede4e}
\end{EQA}
If \( \rr(\CGPs) = \CONSTbal \CGPs^{1/2} \) then it holds on \( \Omega(\xx) \)
for \( \GP^{2} = \CGPs \GPa^{2} \) and \( \DPf_{\CGPs}^{2} = \DPf^{2} + \CGPs \GPa^{2} \) 
\begin{EQA}
	\| \DPf_{\CGPs} (\tilde{\upsv}_{\CGPs} - \upsvs) \|
	& \leq &
	\CONST \, \CONSTbal \, \CGPs^{1/2} \, . 
\label{ichdfhyeherdhytye4dfnhfcy}
\end{EQA}
\end{theorem}

\noindent
Under the same conditions one can specify the results of Section~\ref{SgenBounds}
including 
the Fisher--Wilks expansions of Theorem~\ref{TFiWititG}.
Moreover, for a Gaussian prior \( \ND(0,\GP^{-2}) \) with \( \GP^{2} = \CGP \GPa^{2} \),
one can derive the results about concentration and contraction
of the posterior \( \vupsv_{\GP} \cond \Xv \).
Only Laplace's approximation of the posterior in Theorem~\ref{TBvMgen34} requires a stronger condition
on critical dimension: \( \hmax_{3} \, \amax^{-1} \rr(\CGP) \, \dimA(\CGP) \leq 2 \nsize^{1/2} \).


\Subsection{Rate optimality under Sobolev smoothness}
\label{Spostlogdens}

To state standard rate-optimal results and to compare our conclusions with the existing results
in the literature, we cinsider the univariate case \( d=1 \) and introduce the condition on Sobolev smoothness of the 
log-density \( \log \dens(\xv) \).

\begin{description}
\item[\label{Thetaref} \( \bb{(\smpa,\CGPa)} \)]
\( \upsvs \in \BBB(\smpa,\CGPa) = \bigl\{ \upsv = (\ups_{j}) \colon \sum_{j} \ups_{j}^{2} \, j^{2\smpa} \leq \CGPa \bigr\} \) 
for \( \smpa > 0 \).
\end{description}
The assumption \( \upsvs \in \BBB(\smpa,\CGPa) \) is standard in log-density estimation; cf. \cite{CaNi2014}
or \cite{SpPa2019}.
However, one usually requires \( \smpa > 1 \), while our results below are valid for \( \smpa > 0 \).
The only exception is the final Gaussian approximation result requiring \( \smpa > 1 \).
We also assume that \( \CGPa \asymp 1 \).

Now we state the results under the smoothness condition \nameref{Thetaref}.
The precision/penalization matrix \( \GP^{2} \) is taken of the \( (\smp,\CGP) \)-form: 
\( \GP^{2} = \diag(\gp_{1}^{2},\ldots,\gp_{\dimp}^{2}) \)
with \( \gp_{j}^{2} = j^{2\smp}/\CGP \).
One can take any degree \( \smp > 1/2 \), \( \smp \geq \smpa \), we recommend a large value like \( \smp = 4 \) or \( \smp = 5 \).
Only the factor \( \CGP \) should be fixed carefully to get the optimal accuracy of estimation
from the relation \( (\CGP \nsize)^{1/(2\smp)} \approx (\CGPa \nsize)^{1/(2\smpa+1)} \);
see \eqref{6drvjhggv74e3ujf6hgewt} or \eqref{nmgvctf5rtrtdw2gf56}.
Under \nameref{cdens2ref}, 
the corresponding effective dimension \( \dimG \) and the Laplace effective dimension 
\( \dimL(\upsv) \) are determined by the index \( \mm \) for which 
\( \gp_{\mm}^{2} \approx \nsize \).
Alternatively one can use a \( \mmps \)-truncation prior with \( \mmps \approx (\CGPa \nsize)^{1/(2\smpa+1)} \).
Such a choice of the prior parameters
is frequently used for nonparametric \emph{rate optimal} results about concentration of the pMLE
and posterior contraction; cf. \cite{CaNi2014}, \cite{CaRo2015}.
Note that the mentioned results require \( \smpa > 1 \), while our concentration 
and contraction results apply under \( \smpa > 0 \).
Theorem~\ref{CGLMlossGP} through Theorem~\ref{TLaplaceGP} yield the following results.

\begin{theorem}
\label{PlglossGPg}
Assume \nameref{upsvsref}, \nameref{Upsref}, \nameref{cdensref}, \nameref{PsiX4ref}, \nameref{cdens2ref} and \nameref{Thetaref}.
Fix \( \smp > 1/2 \), \( \smp \geq \smpa \) and define 
\( \gp_{j}^{2} = \CGP^{-1} j^{2\smp} \) with \( \CGP \) satisfying \( (\CGP \nsize)^{1/(2\smp)} \approx (\CGPa \nsize)^{1/(2\smpa+1)} \); see \eqref{6drvjhggv74e3ujf6hgewt} or \eqref{nmgvctf5rtrtdw2gf56}.
Then on \( \Omega(\xx) \) for \( \nsize \geq \nsizen \),
\begin{EQA}
	\| \tilde{\upsv}_{\GP} - \upsvs \|
	& \leq &
	\CONST \nsize^{- \smpa/(2 \smpa + 1)} ,
\label{bnjhuy0w9e8r7t6bcgvftropglg}
\end{EQA}
and the posterior measure \( \vupsv_{\GP} \cond \Xv \) satisfies
\begin{EQA}
	\P\Bigl( \| \vupsv_{\GP} - \tilde{\upsv}_{\GP} \| > \CONST \nsize^{-\smpa/(2\smpa+1)} \Cond \Xv \Bigr)
	& \leq &
	\ex^{-\xx} .
\label{0j7jds2fc56gh98hbt3wdens}
\end{EQA}
\end{theorem}

Our main result about Gaussian approximation requires more smoothness of the log-density \( \log \fs(\cdot) \).
Namely, we require \nameref{Thetaref} with \( \smpa > 1 \).

\begin{theorem}
\label{Tconcdensex}
Assume the conditions of Theorem~\ref{PlglossGPg} and \nameref{Thetaref} with \( \smpa > 1 \).
Define \( \mmps = (\CGPa \nsize)^{1/(2\smpa+1)} \).
For the \( \mm \)-truncation prior with \( \mm = \mmps \)
or for the \( (\smp,\CGP) \)-smooth prior with \( (\CGP \, \nsize)^{1/(2\smp)} = (\CGPa \, \nsize)^{1/(2\smpa+1)} = \mmps \), 
on \( \Omega(\xx) \) for \( \nsize \geq \nsizen \),
\begin{EQ}[rcl]
	\sup_{A \in \cc{B}(\R^{\dimp})}
	\left| 
		\P\bigl( \vupsv_{\GP} - \tilde{\upsv}_{\GP} \in A \cond \Xv \bigr)
		- \PG\bigl( \DPGPt^{-1} \gammav \in A \bigr)  
	\right|
	& \leq &
	\CONST\, \nsize^{(1 - \smpa) / (2\smpa + 1)} \, ,
	\\
\label{1p1m1a1emxrrdens}
	\sup_{A \in \cc{B}_{\smp}(\R^{\dimp})}
	\left| 
		\P\bigl( \vupsv_{\GP} - \tilde{\upsv}_{\GP} \in A \cond \Xv \bigr)
		- \PG\bigl( \DPGPt^{-1} \gammav \in A \bigr)  
	\right|
	& \leq &
	\CONST\, \nsize^{(2 - 2 \smpa) / (2\smpa + 1)} \, .
	\qquad
\end{EQ}
\end{theorem}


The error terms in \eqref{1p1m1a1emxrrdens} 
tend to zero as \( \nsize \to \infty \) because \( \smpa > 1 \).



\Chapter{Dimension free bounds for Laplace approximation}
\label{STaylor}
Here we present several issues related to Laplace approximation.
Section~\ref{SboundsLapl} states general results about the accuracy of Laplace approximation for finite samples.
Section~\ref{SLaplinexact} discusses inexact approximation and the use of posterior mean.
Technical assertions and proofs are collected in Section~\ref{SLapltools}.

\Section{Setup and conditions}
\label{SsetupLapl}
Let \( \lgd(\xv) \) be a function in a high-dimensional Euclidean space \( \R^{\dimp} \) such that
\( \int \ex^{\lgd(\xv)} \, d\xv = \CONST < \infty \),
where the integral sign \( \int \) without limits means the integral over the whole space \( \R^{\dimp} \).
Then \( \lgd \) determines a distribution \( \PfL \) with the density
\( \CONST^{-1} \ex^{\lgd(\xv)} \).
Let \( \xvs \) be a point of maximum:
\begin{EQA}
	\lgd(\xvs)
	&=&
	\sup_{\uv \in \R^{\dimp}} \lgd(\xvs + \uv) .
\label{scdygw7ytd7wqqsquuqydtdtd}
\end{EQA}
We also assume that \( \lgd(\cdot) \) is at least three time differentiable. 
Introduce the negative Hessian \( \IFL = - \nabla^{2} \lgd(\xvs) \) and assume \( \IFL \) strictly positive definite.
We aim at approximating the measure \( \PfL \) by a Gaussian measure \( \ND(\xvs,\IFL^{-1}) \).
Given a function \( g(\cdot) \), define its expectation w.r.t. \( \PfL \) after centering at \( \xvs \):
\begin{EQA}
	\II(g)
	& \eqdef &
	\frac{\int g(\uv) \, \ex^{\lgd(\xvs + \uv)} \, d\uv}{\int \ex^{\lgd(\xvs + \uv)} \, d\uv} \, .
\label{IIgfigefxududu}
\end{EQA}
A Gaussian approximation \( \II_{\IFL}(g) \) for \( \II(g) \) is defined as
\begin{EQA}
	\II_{\IFL}(g)
	& \eqdef &
	\frac{\int g(\uv) \, \ex^{- \| \IFL^{1/2} \uv \|^{2}/2} \, d\uv}{\int \ex^{- \| \IFL^{1/2} \uv \|^{2}/2} \, d\uv} 
	=
	\E g(\gaussv_{\IFL}) ,
	\qquad
	\gaussv_{\IFL} \sim \ND(0,\IFL^{-1}) \, .
\label{IgiLapgeHvdu}
\end{EQA}
The choice of the distance between \( \PfL \) and \( \ND(\xvs,\IFL^{-1}) \) specifies the considered class 
of functions \( g \).
The most strong total variation distance can be obtained as 
the supremum of \( |\II(g) - \II_{\IFL}(g)| \) over all measurable functions \( g(\cdot) \) with 
\( |g(\uv)| \leq 1 \):
\begin{EQA}
	\TV\bigl( \PfL,\ND(\xvs,\IFL^{-1}) \bigr)
	&=&
	\sup_{\|g\|_{\infty} \leq 1} \bigl| \II(g) - \II_{\IFL}(g) \bigr| \, .
\label{ghdvcgftdftgegdftefd3434545}
\end{EQA}
The results can be substantially improved if only centrally symmetric functions \( g(\cdot) \)
with \( g(\xv) = g(-\xv) \) are considered.
Obviously, for any \( g(\cdot) \)
\begin{EQA}
	\II(g)
	& = &
	\frac{\int g(\uv) \, \ex^{\lgd(\xvs + \uv) - \lgd(\xvs)} \, d\uv}{\int \ex^{\lgd(\xvs + \uv) - \lgd(\xvs)} \, d\uv} \, .
\label{IIgfigefxuudufxudu}
\end{EQA}
Moreover, as \( \xvs = \argmax_{\xv} \lgd(\xv) \), it holds \( \nabla \lgd(\xvs) = 0 \) and 
\begin{EQA}
	\II(g)
	& = &
	\frac{\int g(\uv) \, \ex^{\lgd(\xvs;\uv)} \, d\uv}{\int \ex^{\lgd(\xvs;\uv)} \, d\uv} \, ,
\label{IIgfifxutudufxutdu}
\end{EQA}
where \( \lgd(\xv;\uv) \) is the Bregman divergence 
\begin{EQA}
	\lgd(\xv;\uv)
	&=&
	\lgd(\xv + \uv) - \lgd(\xv) - \bigl\langle \nabla \lgd(\xv), \uv \bigr\rangle .
\label{fxufxpufxfpxu}
\end{EQA}
Implicitly we assume that the negative Hessian \( \IFL = - \nabla^{2} \lgd(\xvs) \) is sufficiently large
in the sense that the Gaussian measure \( \ND(0,\IFL^{-1}) \) concentrates on a small local set \( \UVL \).
This allows to use a local Taylor expansion for 
\( \lgd(\xvs;\uv) \approx - \| \IFL^{1/2} \uv \|^{2}/2 \) in \( \uv \) on \( \UVL \).
If \( \lgd(\cdot) \) is also strongly concave, then the \( \PfL \)-mass  
of the complement of \( \UVL \) is exponentially small yielding the desirable Laplace approximation.

Motivated by applications to statistical inference, we consider \( \lgd \) in a special form
\begin{EQA}
	\lgd(\xv)
	&=&
	\lgdL(\xv) - \| \GP (\xv - \xv_{0}) \|^{2} / 2
\label{jhdctrdfred4322edt7y}
\end{EQA}
for some \( \xv_{0} \) and a symmetric \( \dimp \)-matrix \( \GP^{2} \geq 0 \).
Here \( \lgdL(\cdot) \) stands for a log-likelihood function while the quadratic penalty 
\( \| \GP (\xv - \xv_{0}) \|^{2} / 2 \) corresponds to a Gaussian prior \( \ND(\xv_{0},\GP^{-2}) \).
We also assume that \( \lgdL(\cdot) \) is concave with \( \DV^{2} \eqdef - \nabla^{2} \lgdL(\xv) \geq 0 \).
Then  
\begin{EQA}
	- \nabla^{2} \lgd(\xv)
	&=&
	- \nabla^{2} \lgdL(\xv) + \GP^{2} 
	=
	\DVL^{2} + \GP^{2} .
\label{hydsf42wdsdtrdstdg}
\end{EQA} 
In typical asymptotic setups, the log-likelihood function \( \lgdL(\xv) \) grows with the sample size 
or inverse noise variance, while the prior precision matrix \( \GP^{2} \) is kept fixed.
Decomposition \eqref{jhdctrdfred4322edt7y} is of great importance for obtaining the dimension free results.
The main reason is that the quadratic penalty does not affect smoothness properties of the function \( \lgdL(\cdot) \)
but greatly improves the quadratic approximation term.

\Subsection{Concavity}
Below we implicitly assume decomposition \eqref{jhdctrdfred4322edt7y} with 
a \emph{weak\-ly concave} function \( \lgdL(\cdot) \).
More specifically, we assume the following condition.

\begin{description}
    \item[\label{LLf0ref} \( \bb{(\mathcal{C}_{0})} \)]
      \textit{ There exists an operator \( \GP^{2} \leq - \nabla^{2} \lgd(\xvs) \) in \( \R^{\dimp} \) such that the function 
\begin{EQA}
	\lgdL(\xvs + \uv)
	& \eqdef &
	\lgd(\xvs + \uv) + \| \GP \uv \|^{2} / 2
\label{fTvufupumDTpDfx}
\end{EQA}
is concave. 
      }
\end{description}

If \( \lgdL(\cdot) \) in decomposition \eqref{jhdctrdfred4322edt7y} is concave then this condition is obviously fulfilled. 
More generally, if \( \lgdL(\cdot) \) in \eqref{jhdctrdfred4322edt7y} is weakly concave, so that 
\( \lgdL(\xvs + \uv) - \| \GP_{0} \uv \|^{2}/2 \) is concave in \( \uv \) with \( \GP_{0}^{2} \leq \GP^{2} \), then 
\nameref{LLf0ref} is fulfilled with \( \GP^{2} - \GP_{0}^{2} \) in place of \( \GP^{2} \).

The operator \( \DV^{2} \) plays an important role in our conditions and results:
\begin{EQA}
	\DV^{2} 
	&=&
	- \nabla^{2} \lgd(\xvs) - \GP^{2}
	\quad
	\bigl( = - \nabla^{2} \lgdL(\xvs) \text{ under \eqref{jhdctrdfred4322edt7y}} \bigr) .
\label{kc7c322dgdf43djhd}
\end{EQA}

\begin{remark}
The condition of strong concavity of \( \lgd \) on the whole space \( \R^{\dimp} \)
can be too restrictive.%
This condition can be replaced by its local version: \( \lgd \) is concave on a set \( \XXL_{0} \) 
such that the Gaussian prior \( \ND(\xv_{0},\GP^{-2}) \) concentrates on \( \XXL_{0} \)
and the maximizer \( \xvs_{\GP} \) belongs to \( \XXL_{0} \).%
In all the results, the integral over \( \R^{\dimp} \) has to be replaced by the integral over \( \XXL_{0} \).
\end{remark}

\Subsection{Laplace effective dimension}
\label{SeffdimLa}
With decomposition \eqref{kc7c322dgdf43djhd} in mind, we use another notation for \( \IFL = - \nabla^{2} \lgd(\xvs) \):
\begin{EQA}
	\DV_{\GP}^{2} 
	&=& 
	- \nabla^{2} \lgd(\xvs)
	=
	\DVL^{2} + \GP^{2} .
\label{dscytf5w2edte5rw4e24gyd}
\end{EQA}
Also we write \( \II_{\GP}(g) \) instead of \( \II_{\IFL}(g) \) in \eqref{IgiLapgeHvdu}.
The \emph{Laplace effective dimension} \( \dimLL \) is given by
\begin{EQA}
	\dimLL
	& \eqdef & 
	\tr\bigl( \DVL^{2} \, \DV_{\GP}^{-2} \bigr) .
\label{dAdetrH02Hm2}
\end{EQA}
Of course, \( \dimLL \leq \dimp \) but a proper choice of the penalty \( \GP^{2} \) in \eqref{jhdctrdfred4322edt7y}
allows to avoid the ``curse of dimensionality'' issue and ensure a small effective dimension \( \dimLL \) even for \( \dimp \) large or infinite; see 
\ifapp{Section~\ref{Ssmoothprior}}{\cite{SpPa2019}} 
for more rigorous discussion.

Later we write \( \dimL \) instead of \( \dimLL \) without risk of confusion because the parameter dimension \( \dimp \)
does not show up anymore.
The value \( \dimL \) helps to describe a local vicinity \( \UVL \) around \( \xvs \) such that the most of mass of \( \PfL \)
concentrates on \( \UVL \); see Section~\ref{StailLa}.
Namely, let us fix some \( \amax < 1 \), e.g. \( \amax = 2/3 \), and some \( \xx > 0 \) ensuring that \( \ex^{-\xx} \)
is our significance level.
Define
\begin{EQ}[rcccl]
	\rrL
	& = &
	2 \sqrt{\dimL} + \sqrt{2 \xx} ,
	\qquad
	\UVL 
	&=& 
	\bigl\{ \uv \colon \| \DVL \uv \| \leq \amax^{-1} \rrL \bigr\} .
\label{UvTDunm12spT}
\end{EQ}

\Subsection{Local smoothness conditions}
Let \( \dimp \leq \infty \) and
let \( \lgd(\cdot) \) be a three times continuously differentiable function on \( \R^{\dimp} \).
We fix a reference point \( \xv \) and local region around \( \xv \) given by the local set 
\( \UVL \subset \R^{\dimp} \) from \eqref{UvTDunm12spT}. 
Consider the remainder of the second and third order Taylor approximation 
\begin{EQ}[rcl]
	\dltw_{3}(\xv,\uv)
	&=&
	\lgd(\xv;\uv) - 
	\bigl\langle \nabla^{2} \lgd(\xv) , \uv^{\otimes 2} \bigr\rangle/2 ,
	\\
	\dltw_{4}(\xv,\uv)
	&=&
	\lgd(\xv;\uv) - 
	\bigl\langle \nabla^{2} \lgd(\xv) , \uv^{\otimes 2} \bigr\rangle/2 
	- 
	\bigl\langle \nabla^{3} \lgd(\xv), \uv^{\otimes 3} \bigr\rangle / 6
\label{d4fuv1216303}
\end{EQ}
with \( \lgd(\xv;\uv) \) from \eqref{fxufxpufxfpxu}.
The use of the Taylor formula allows to bound
\begin{EQ}[rcl]
	\bigl| \dltw_{k}(\xv,\uv) \bigr|
	& \leq &
	\sup_{t \in [0,1]}
	\frac{1}{k!} 
	\Bigl| \bigl\langle \nabla^{k} \lgd(\xv + t \uv), \uv^{\otimes k} \bigr\rangle \Bigr|,
	\quad
	k\geq 3. 
\label{k3t01f12n3fvtu}
\end{EQ}
Note that the quadratic penalty \( - \| \GP (\xv - \xv_{0}) \|^{2}/2 \) in \( \lgd \) does not affect 
the remainders \( \dltw_{3}(\xv,\uv) \) and \( \dltw_{4}(\xv,\uv) \).
Indeed, 
with \( \lgd(\xv) = \lgdL(\xv) - \| \GP (\xv - \xv_{0}) \|^{2}/2 \), it holds 
\begin{EQA}
	\lgd(\xv;\uv) 
	& \eqdef & 
	\lgd(\xv + \uv) - \lgd(\xv) - \bigl\langle \nabla \lgd(\xv), \uv \bigr\rangle 
	=
	\lgdL(\xv;\uv) - \| \GP \uv \|^{2}/2 
\label{D2fuvGu2fG}
\end{EQA}
and the quadratic term in definition of the values \( \dltw_{k}(\xv,\uv) \) cancels, \( k \geq 3 \). 
Local smoothness of \( \lgd(\cdot) \) or, equivalently, of \( \lgdL(\cdot) \), at \( \xv \) will be measured by
the value \( \dltwb(\xv) \):
\begin{EQA}
	\dltwb(\xv)
	& \eqdef &
	\sup_{\uv \in \UVL} \frac{1}{\| \DVL \uv \|^{2}/2} \bigl| \dltw_{3}(\xv,\uv) \bigr| ;
\label{om3esuU1H02d3}
\end{EQA}
cf. \eqref{dtb3u1DG2d3GPg}.
We also denote \( \dltwb \eqdef \dltwb(\xvs) \).
Our results apply if \( \dltwb \ll 1 \).
Local concentration of the measure \( \PfL \) requires \( \dltwb \leq 1/3 \);
see Proposition~\ref{PlocconLa}.
The results about Gaussian approximation are valid under a stronger condition 
\( \dltwb \, \dimL \leq 2/3 \) with 
\( \dimL \) from \eqref{dAdetrH02Hm2}.


\Section{Error bounds for Laplace approximation}
\label{SboundsLapl}
Our first result describes the quality of approximation of the measure \( \PfL \) by the Gaussian measure \( \ND(\xvs,\DV_{\GP}^{-2}) \)
with mean \( \xvs \) and the covariance \( \DV_{\GP}^{-2} \) in total variation distance.
In all our result, the value \( \xx \) is fixed to ensure that \( \ex^{-\xx} \) is negligible.
First we present the general results which will be specified later under the self-concordance condition.


\begin{theorem}
\label{TLaplaceTV}
Suppose \nameref{LLf0ref}.
Let also \( \dimL \) be defined by \eqref{dAdetrH02Hm2} and \( \rrL \) and \( \UVL \) by \eqref{UvTDunm12spT}.
If \( \dltwb \) from \eqref{om3esuU1H02d3} satisfies \( \dltwb \leq 1/3 \), then
\begin{EQA}
	\PfL(\Xv - \xvs \not\in \UVL)
	& \leq &
	\ex^{-\xx} .
\label{poybf3679jd532ff2}
\end{EQA}
If \( \dltwb \, \dimL \leq 2/3 \),
then for any \( g(\cdot) \) with \( |g(\uv)| \leq 1 \), it holds for \( \II(g) \) from \eqref{IIgfigefxududu}
\begin{EQA}
	\bigl| \II(g) - \II_{\GP}(g) \bigr|
	& \leq &
	\frac{2 (\err + \ex^{-\xx})}{1 - \err - \ex^{-\xx}}
	\leq 
	4(\err + \ex^{-\xx}) 
\label{ufgdt6df5dtgededsxd23gjg}
\end{EQA}
with 
\begin{EQA}
	\err
	&=&
	\err_{2} 
	= 
	\frac{0.75 \, \dltwb \, \dimL}{1 - \dltwb} \, .
\label{juytr90f2dzaryjhfyf}
\end{EQA}
\end{theorem}

This section presents more advanced bounds on the error of Laplace approximation
under conditions \nameref{LL3tref} and \nameref{LL4tref} with \( \upsv = \xvs \) and \( \HLG(\xvs) = \DV \)
or \nameref{LLtS3ref} and \nameref{LLtS4ref} with \( \upsv = \xvs \) and \( \HL(\xvs) = n^{-1/2} \, \DV \).
%

\begin{theorem}
\label{TLaplaceTV2}
\label{TLaplaceTV34}
Suppose \nameref{LLf0ref} and \nameref{LL3tref} 
and let \( \dltwu_{3} \, \amax^{-1} \rrL \leq 3/4 \) for \( \rrL \) from \eqref{UvTDunm12spT}.
Then the concentration bound \eqref{poybf3679jd532ff2} holds.
Moreover, if 
\begin{EQA}
	\dltwu_{3} \, \amax^{-1} \rrL \, \dimL
	& \leq &
	2 ,
\label{0hcde4dft3igthg94yr5twew}
\end{EQA}
then the accuracy bound \eqref{ufgdt6df5dtgededsxd23gjg} applies with
\begin{EQA}
	\err
	&=&
	\err_{3}
	\eqdef
	\frac{\dltwu_{3} (\dimL+1)^{3/2}}{4 (1 - \dltwb)^{3/2}} 
	\leq 
	\frac{\dltwu_{3} \, (\dimL+1)^{3/2}}{2} \, ,
\label{5qw7dyf4e4354coefw9dufih}
\end{EQA}
where \( \dltwb \eqdef \dltwu_{3} \, \amax^{-1} \rrL / 3 \leq 1/4 \).
Furthermore, under \nameref{LL4tref}, 
for any symmetric function \( g(\uv) = g(-\uv) \), \( |g(\uv)| \leq 1 \),
the accuracy bound \eqref{ufgdt6df5dtgededsxd23gjg} applies with 
\begin{EQA}
	\err
	&=&
	\err_{4}
	=
	\frac{\dltwu_{3}^{2} \, (\dimL + 2)^{3} + 2 \dltwu_{4} (\dimL + 1)^{2}}{16 (1 - \dltwb)^{2}} 	
	\leq 
	\frac{\dltwu_{3}^{2} \, (\dimL + 2)^{3} + 2 \dltwu_{4} (\dimL + 1)^{2}}{8} \, \, .
\label{hg25t6mxwhydseg3hhfdr}
\end{EQA}
Under \nameref{LLtS3ref} and \nameref{LLtS4ref} instead of \nameref{LL3tref} and \nameref{LL4tref}, the results apply with 
\( \dltwu_{3} = \hmax_{3} \, n^{-1/2} \) and \( \dltwu_{4} = \hmax_{4} \, n^{-1} \).
\end{theorem}

Let \( \BBB(\R^{\dimp}) \) be the \( \sigma \)-field of all Borel sets in \( \R^{\dimp} \),
while \( \BBB_{s}(\R^{\dimp}) \) stands for all centrally symmetric sets from \( \BBB(\R^{\dimp}) \).
By \( \Xv \) we denote a random element with the distribution \( \PfL \), while \( \gaussv_{\GP} \sim \ND(0,\DV_{\GP}^{-2}) \).

\begin{corollary}
\label{CTLaplaceTV}
Under the conditions of Theorem~\ref{TLaplaceTV2}, it holds for \( \Xv \sim \PfL \)
\begin{EQA}
	\sup_{A \in \BBB(\R^{\dimp})} \bigl| \PfL(\Xv - \xvs \in A) - \P(\gaussv_{\GP} \in A) \bigr|
	& \leq &
	4(\err_{3} + \ex^{-\xx}),
	\\
	\sup_{A \in \BBB_{s}(\R^{\dimp})} \bigl| \PfL(\Xv - \xvs \in A) - \P(\gaussv_{\GP} \in A) \bigr|
	& \leq &
	4(\err_{4} + \ex^{-\xx}) .
\label{cbc5dfedrdwewwgerg}
\end{EQA}
\end{corollary}

\Subsection{Critical dimension}
Here we briefly discuss the important issue of \emph{critical dimension}.
Theorem~\ref{TLaplaceTV34} states concentration of \( \PfL \) under the condition \( \dltwu_{3} \rrL \leq 1 \).
Under \nameref{LLtS3ref}, we can use \( \dltwu_{3} = \hmax_{3} \, n^{-1/2} \).
Together with \( \rrL \approx \sqrt{\dimL} \), this yields the condition
\( \dimL \ll n \).
Gaussian approximation applies under \( \hmax_{3} \, \amax^{-1} \rrL \, \dimL \, n^{-1/2} \leq 2 \); see \eqref{0hcde4dft3igthg94yr5twew},
yielding \( \dimL^{3} \ll n \).
We see that there is a gap between these conditions. 
We guess that in the region \( n^{1/3} \lesssim \dimL \lesssim n \), 
non-Gaussian approximation of the posterior is possible; cf. 
\cite{bochkina2014}.

\ifNL{
\Subsection{Kullback-Leibler divergence}
Theorem~\ref{TLaplaceTV} through \ref{TLaplaceTV34} quantify the approximation \( \PfL \approx \ND(\xvs,\DV_{\GP}^{-2}) \)
in the total variation distance.
Another useful characteristic could be the Kullback-Leibler (KL) divergence between \( \PfL \) and  
\( \ND(\xvs,\DV_{\GP}^{-2}) \).
The KL divergence \( \kullb(\P_{1},\P_{2}) = \E_{1} \log(d\P_{1}/d\P_{2}) \) is asymmetric, 
\( \kullb(\P_{1},\P_{2}) \neq \kullb(\P_{2},\P_{1}) \) with few exceptions 
like the case of Gaussian measures \( \P_{1} \) and \( \P_{2} \).
Moreover, \( \kullb(\P_{1},\P_{2}) \) can explode if \( \P_{1} \) is not absolutely continuous w.r.t. \( \P_{2} \). 
We present two bounds for each ordering.
For ease of presentation, we limit ourselves to the case when either \nameref{LL3tref} or \nameref{LLtS3ref} meets.

\begin{theorem}
\label{TLaplaceKL}
Suppose \nameref{LLf0ref} and \nameref{LL3tref} and let \( \dltwb \, \dimL \leq 2/3 \).
With \( \P_{\GP} = \ND(\xvs,\DV_{\GP}^{-2}) \),
\begin{EQA}
	\kullb(\PfL,\P_{\GP})
	& \leq &
	4 \err_{3} + 4 \ex^{-\xx} 
	\leq 
	\frac{\E \dltwu_{3}(\gaussv_{\GP})}{(1 - \dltwb)^{3/2}} + 4 \ex^{-\xx} .
\label{vlgvi8ugu7tr4ry43et31}
\end{EQA}
Moreover, under \nameref{LLtS3ref}
\begin{EQA}
	\kullb(\PfL,\P_{\GP})
	& \leq &
	2 \hmax_{3} \, \sqrt{\frac{(\dimL+1)^{3}}{n}} + 4 \ex^{-\xx} .
\label{cheiufheurhrtdgwb3wesrtd}
\end{EQA}
\end{theorem}

Now we briefly discuss the value \( \kullb(\P_{\GP},\PfL) \).
We already know that \( \PfL \) concentrates on \( \UVL \) and can be well approximated by \( \P_{\GP} \) on \( \UVL \).
However, this does not guarantee a small value of \( \kullb(\P_{\GP},\PfL) \).
It can even explode if e.g. \( \PfL \) has a compact support.
In fact, 
the log-density of \( \P_{\GP} \) w.r.t. \( \PfL \) reads
\begin{EQA}
	\log \frac{d\P_{\GP}}{d\PfL}(\xv)
	&=&
	- \lgd(\xv) - \frac{1}{2} \| \DV_{\GP} (\xv - \xv_{0}) \|^{2} - \CDG
\label{fduhefuifhew2tdw6ww3}
\end{EQA}
for some constant \( \CDG \),
and an upper bound on \( \kullb(\P_{\GP},\PfL) \) requires that the integral of \( \lgd(\xv) \) w.r.t. the measure \( \P_{\GP} \) is finite.

\begin{theorem}
\label{TLaplaceKLi}
Suppose \nameref{LLf0ref} and \nameref{LLtS3ref} and let \( \dltwb \, \dimL \leq 2/3 \).
Let also \( \rho_{\GP} = 2\xx /\rrL^{2} \); see \eqref{UvTDunm12spT}.
If \( \lgdL(\xvs;\uv) = \lgdL(\xvs + \uv) - \lgdL(\xvs) - \langle \nabla \lgdL(\xvs), \uv \rangle \) fulfills
\begin{EQA}
	\int \bigl| \lgdL(\xvs;\uv) \bigr| \, \exp \bigl\{ - \| \DV_{\GP} \uv \|^{2}/2 + \rho_{\GP} \| \DVL \uv \|^{2}/2 \bigr\} \, d\uv
	& \leq &
	\CONSTi_{\lgdL} \, 
\label{dchbhwdhwwdgscsn2efty2162}
\end{EQA} 
for some fixed constant \( \CONSTi_{\lgdL} \) then
\begin{EQA}
	\kullb(\P_{\GP},\PfL)
	& \leq &
	\hmax_{3} \, \sqrt{\frac{(\dimL+1)^{3}}{n}} + (2 + \CONSTi_{\lgdL}) \ex^{-\xx} .
\label{cheiufheurhrtdgwb3wesrtd}
\end{EQA}
\end{theorem}
}{}

\Subsection{Mean and MAP}
Here we present the bound on \( |\II(g) - \II_{\GP}(g)| \) for the case of a linear vector function \( g(\uv) = \QP \uv \) 
with \( \QP \colon \R^{\dimp} \to \R^{\dimq} \), \( \dimq \geq 1 \). 
A special case of \( \QP = \Id_{\dimp} \) corresponds to the mean value \( \xvb \) of \( \PfL \).
The next result presents an upper bound for the value \( \QP (\xvb - \xvs) \) under the conditions of Theorem~\ref{TLaplaceTV34}
including \nameref{LLtS3ref}.

\begin{theorem}
\label{TpostmeanLa}
Assume the conditions of Theorem~\ref{TLaplaceTV34} and let \( \QP^{\T} \QP \leq \DVL^{2} \).
Then 
\begin{EQA}
	\| \QP (\xvb - \xvs) \|
	& \leq &
	2.4 \, \hmax_{3} \, \| \QP \DV_{\GP}^{-2} \QP^{\T} \|^{1/2} \,
	\sqrt{(\dimL + 1)^{3}/n} + \CONST \ex^{-\xx} .
\label{hcdtrdtdehfdewdrfrhgyjufger}
\end{EQA}
\end{theorem}

Now we specify the result for the special choice \( \QP = \DVL \).

\begin{corollary}
\label{CTpostmeanLa}
Assume the conditions of Theorem~\ref{TLaplaceTV34}.
Then 
\begin{EQA}
	\| \DVL (\xvb - \xvs) \|
	& \leq &
	2.4 \, \hmax_{3} \, \sqrt{(\dimL + 1)^{3}/{n}} + \CONST \ex^{-\xx} \, .
\label{klu8gitfdgregfkhj7yt}
\end{EQA}
\end{corollary}

\ifapp{
\begin{remark}
\label{RTpostmeanc}
An interesting question is whether the result of Theorem~\ref{TpostmeanLa} or Corollary~\ref{CTpostmeanLa} applies 
with \( \QP = \DV_{\GP} \).
This issue is important in connection to inexact Laplace approximation; see the next section.
The answer is negative. 
The problem is related to the last term \( \CONST \ex^{-\xx} \) in the right hand-side of \eqref{hcdtrdtdehfdewdrfrhgyjufger}.
The constant \( \CONST \) involves the moments of \( \| \QP (\Xv - \xvs) \|^{2} \)
which explode for \( \QP = \DV_{\GP} \) and \( \dimp = \infty \).
\end{remark}
}{}

\Section{Inexact approximation and the use of posterior mean}
\label{SLaplinexact}

Now we change the setup.
Namely, we suppose that the true maximizer \( \xvs \) of the function \( \lgd \) is not available, but 
\( \xv \) is somehow close to the point of maximum \( \xvs \).
Similarly, the negative Hessian \( \DV_{\GP}^{2}(\xvs) = - \nabla^{2} \lgd(\xvs) \) is hard to obtain 
and we use a proxy \( \HUL^{2} \).
We already know that \( \PfL \) can be well approximated by \( \ND(\xvs,\DV^{-2}) \) with 
\( \DV^{2} = \DV^{2}(\xvs) \).
This section addresses the question whether \( \ND(\xv,\HUL^{-2}) \) can be used instead.
Here we may greatly benefit from the fact that Theorem~\ref{TLaplaceTV} provides a bound in the total variation distance 
between \( \PfL \) and \( \ND(\xvs,\DV_{\GP}^{-2}) \) yielding
\begin{EQA}
	\TV\bigl( \PfL, \ND(\xv,\HUL^{-2}) \bigr)
	& \leq &
	\TV\bigl( \PfL, \ND(\xvs,\DV_{\GP}^{-2}) \bigr)
	+ \TV\bigl( \ND(\xv,\HUL^{-2}), \ND(\xvs,\DV_{\GP}^{-2}) \bigr) .
	\qquad
\label{iuhvfuggeyfgeyjhefhgdy}
\end{EQA}
Therefore, it suffices to bound the TV-distance between the Gaussian distribution 
\( \ND(\xvs,\DV_{\GP}^{-2}) \) naturally arising in Laplace approximation, 
and the one used instead.
Pinsker's inequality provides an upper bound: for any two measures \( P,Q \)
\begin{EQA}
	\TV(P,Q)
	& \leq &
	\sqrt{\kullb(P,Q)/2} ,
\label{gffwdygweufeuyfeygfdwyd}
\end{EQA}
where \( \kullb(P,Q) \) is the Kullback-Leibler divergence between \( P \) and \( Q \).
The KL-divergence between two Gaussians has a closed form:
\begin{EQA}
	\kullb\bigl( \ND(\xvs,\DV_{\GP}^{-2}), \ND(\xv,\HUL^{-2}) \bigr)
	&=&
	\frac{1}{2} \bigl\{ \| \DV_{\GP} (\xv - \xvs) \|^{2} + \tr (\HUL^{-2} \DV_{\GP}^{2} - \Id_{\dimp}) + \log \det (\HUL^{-2} \DV_{\GP}^{2}) \bigr\} .
\label{h0bede4bweve7yjdyeghy}
\end{EQA}
Moreover, if the matrix \( \BB_{\GP} = \HUL^{-1} \DV_{\GP}^{2} \HUL^{-1} - \Id_{\dimp} \) satisfies \( \| \BB_{\GP} \| \leq 2/3 \) then
\begin{EQA}
	\TV\bigl( \ND(\xvs,\DV_{\GP}^{-2}), \ND(\xv,\HUL^{-2}) \bigr)
	& \leq &
	\frac{1}{2} \Bigl( \| \DV_{\GP} (\xv - \xvs) \| + \sqrt{\tr \BB_{\GP}^{2}} \Bigr) .
\label{dferrfwbvf6nhdnghfkeiry}
\end{EQA}
However, 
Pinsker's inequality is only a general upper bound which is applied to any two distributions \( P \) and \( Q \).
If \( P \) and \( Q \) are Gaussian, it might be too rough.
Particularly the use of \( \tr \BB_{\GP}^{2} \) is disappointing, this quantity is full dimensional even 
if each of \( \DV_{\GP}^{-2} \) and \( \HUL^{-2} \) has a bounded trace.
Also dependence on \( \| \DV_{\GP} (\xv - \xvs) \| \) is very discouraging.
\cite{Devroy2022} provides much sharper results, however, limited to the case of the same mean.
Even stronger results can be obtained if we restrict ourselves to the class \( \BBB_{el}(\R^{\dimp}) \) 
of elliptic sets \( A \) in \( \R^{\dimp} \) of the form
\begin{EQA}
	A
	&=&
	\bigl\{ \uv \in \R^{\dimp} \colon \| \QP (\uv - \xv) \| \leq \rr \bigr\}
\label{vhg4dfe5w3tfdf54wteg}
\end{EQA}
for some linear mapping \( \QP \colon \R^{\dimp} \to \R^{\dimq} \), \( \xv \in \R^{\dimp} \), and \( \rr > 0 \).
Given two symmetric \( \dimq \)-matrices \( \Sigma_{1},\Sigma_{2} \) and a vector \( \av \in \R^{\dimq} \), define 
\begin{EQA}
	\dist(\Sigma_{1},\Sigma_{2},\av)
	& \eqdef &
	\biggl( \frac{1}{\| \Sigma_{1} \|_{\Fr}} + \frac{1}{\| \Sigma_{2} \|_{\Fr}} \biggr)
	\biggl( \| \lamb_{1} - \lamb_{2} \|_{1} + \| \av \|^{2} \biggr),
\label{btegdhertrewfdgvfddffd}
\end{EQA}
where \( \lamb_{1} \) is the vector of eigenvalues of \( \Sigma_{1} \) arranged in the non-increasing order
and similarly for \( \lamb_{2} \), and \( \| \Sigma \|_{\Fr}^{2} = \tr (\Sigma^{\T} \Sigma) \).
By the Weilandt--Hoffman inequality, 
\( \| \lamb_{1} - \lamb_{2} \|_{1} \leq \| \Sigma_{1} - \Sigma_{2} \|_{1} \) , see e.g. 
\cite{MarkusEng}.
\cite{GNSUl2017} stated the following bound for \( \gaussv_{1} \sim \ND(0,\Sigma_{1}) \) and \( \gaussv_{2} \sim \ND(0,\Sigma_{2}) \):
with an absolute constant \( \CONST \)
\begin{EQA}
	\sup_{\rr > 0}
	\Bigl| \P\bigl( \| \gaussv_{1} + \av \| \leq \rr \bigr) - \P\bigl( \| \gaussv_{2} \| \leq \rr \bigr) \Bigr|
	& \leq &
	\CONST  \, \dist(\Sigma_{1},\Sigma_{2},\av) 
\label{gttegderqwerwvfdhewfs}
\end{EQA}
provided that \( \| \Sigma_{k} \|^{2} \leq \| \Sigma_{k} \|_{\Fr}^{2}/3 \), \( k=1,2 \).
Here \( \| M \|_{1} = \tr |M| = \sum_{j} |\lambda_{j}(M)| \) for a symmetric matrix \( M \) with eigenvalues \( \lambda_{j}(M) \).

\begin{theorem}
\label{TLaplaceTVin}
Assume the conditions of Theorem~\ref{TLaplaceTV} with \( \xvs \) being the maximizer of \( \lgd \) and 
\( \DV_{\GP}^{2} = - \nabla^{2} \lgd(\xvs) \).
For any \( \xv \) and \( \HUL \),
it holds with \( \gaussv_{\HUL} \sim \ND(0,\HUL^{-2}) \)
\begin{EQA}
	&& \nquad
	\sup_{A \in \BBB(\R^{\dimp})} \bigl| \PfL(\Xv - \xv \in A) - \P(\gaussv_{\HUL} \in A) \bigr|\; 
	\\
	& \leq &
	4(\err + \ex^{-\xx}) + \TV\bigl( \ND(\xv,\HUL^{-2}), \ND(\xvs,\DV_{\GP}^{-2}) \bigr) ,
\label{nhjrwsdgdehgdftregwbdctsh}
\end{EQA}
where \( \err = \err_{2} \), see \eqref{juytr90f2dzaryjhfyf}, or \( \err = \err_{3} \), see \eqref{5qw7dyf4e4354coefw9dufih}.

Furthermore, for \( \Xv \sim \PfL \) and \( \gaussv \sim \ND(0,\Id_{\dimp}) \), any linear mapping 
\( \QP \colon \R^{\dimp} \to \R^{\dimq} \),
it holds 
under \( 3 \| \QP \, \DV_{\GP}^{-2} \QP^{\T} \|^{2} \leq \| \QP \, \DV_{\GP}^{-2} \QP^{\T} \|_{\Fr}^{2} \)
\begin{EQA}
	&& \nquad
	\sup_{\rr > 0}
	\left| \PfL\bigl( \| \QP (\Xv - \xv) \| \leq \rr \bigr) 
	- \P\bigl( \| \QP \, \HUL^{-1} \gaussv \| \leq \rr \bigr) 
	\right|
	\\
	& \leq &
	4(\err_{3} + \ex^{-\xx}) + \frac{\CONST}{\| \QP \, \DV_{\GP}^{-2} \QP^{\T} \|_{\Fr}}
	\left( 
		\| \QP (\DV_{\GP}^{-2} - \HUL^{-2}) \QP^{\T} \|_{1} + \| \QP (\xv - \xvs) \|^{2} 
	\right) .
\label{jrwguvyr23jbviufdsfsdgf6w}
\end{EQA}
\end{theorem}

As a special case, consider the use of the posterior mean \( \xvb \) instead of \( \xvs \): 
\begin{EQA}
	\xvb
	& \eqdef &
	\frac{\int \xv \, \ex^{\lgd(\xv)} \, d\xv}{\int \ex^{\lgd(\xv)} \, d\xv} \, .
\label{054rfgofe3rdgbrfty6ry}
\end{EQA}

\begin{theorem}
\label{TpostmeanLan}
Assume the conditions of Theorem~\ref{TpostmeanLa} and Theorem~\ref{TLaplaceTVin}. 
Then it holds for any linear mapping \( \QP \colon \R^{\dimp} \to \R^{\dimq} \) with \( \QP^{\T} \QP \leq \DVL^{2} \)
\begin{EQA}
	\sup_{\rr > 0}
	\left| \PfL\bigl( \| \QP (\Xv - \xvb) \| \leq \rr \bigr) 
	- \P\bigl( \| \QP \gaussv_{\GP} \| \leq \rr \bigr) 
	\right|
	& \leq & 
	4(\err_{3} + \ex^{-\xx}) + \frac{\CONST \| \QP (\xvb - \xvs) \|^{2}}{\| \QP \, \DV_{\GP}^{-2} \QP^{\T} \|_{\Fr}} \, ,
\label{jrwguvyr23jbviufdsfsdgf6wm}
\end{EQA}
where \( \| \QP (\xvb - \xvs) \| \) follows \eqref{hcdtrdtdehfdewdrfrhgyjufger} and \eqref{klu8gitfdgregfkhj7yt}.
\end{theorem}

The case \( \QP = \DVL \) is particularly transparent.
In view of \eqref{klu8gitfdgregfkhj7yt} of Corollary~\ref{CTpostmeanLa} and 
\( \| \QP \, \DV_{\GP}^{-2} \QP^{\T} \|_{\Fr}^{2} = \tr\bigl\{ (\DVL \, \DV_{\GP}^{-2} \DVL^{\T})^{2} \bigr\} \asymp \dimL \),
the following result holds.

\begin{corollary}
\label{CTpostmeanLab}
Under the conditions of Corollary~\ref{CTpostmeanLa}, it holds for \( \Xv \sim \PfL \)
\begin{EQA}
	\sup_{\rr > 0}
	\Bigl| \PfL\bigl( \| \DVL (\Xv - \xvb) \| \leq \rr \bigr) 
	- \P\bigl( \| \DVL \, \gaussv_{\GP} \| \leq \rr \bigr) 
	\Bigr|
	& \leq &
	\CONST \biggl( \sqrt{{\dimL^{3}}/{n}} + \ex^{-\xx} \biggr) \, .
\label{klu8gitfdgregfkhj7yt3n}
\end{EQA}
\end{corollary}

The same bound applies with \( \QP = n^{1/2} \Id_{\dimp} \) in place of \( \DVL \) provided that 
\( \DVL^{2} \geq \CONSTi_{0} \, n \, \Id_{\dimp} \) for some fixed \( \CONSTi_{0} > 0 \).
We may conclude that the use of posterior mean \( \xvb \) in place of the posterior mode \( \xvs \) is justified 
under the same condition on critical dimension \( \dimL^{3} \ll n \) as required for the main result about Gaussian 
approximation.


{
\Section{Tools and proofs}
\label{SLapltools}
\renewcommand{\Section}[1]{\subsubsection{#1}}

Here we collect the proofs of the main results and some useful technical statements 
about the error of Laplace approximation.
Below we write \( \xv \) instead of \( \xvs \).
After passing to representation \eqref{IIgfifxutudufxutdu}, 
many results below apply to any \( \xv \), not necessarily for \( \xv = \xvs \).
We only use \( \DV_{\GP}^{2} = - \nabla^{2} \lgd(\xv) \) and 
\( \dltwb \) instead of \( \dltwb(\xv) \).
Everywhere we assume the local set \( \UVL \) to be fixed by \eqref{UvTDunm12spT}.
We separately study the integrals over \( \UVL \) and over its complement. 
The local error of approximation is measured by
\begin{EQA}
	\err
	=
	\err(\UVL)
	& \eqdef &
	\biggl| 
	\frac{\int_{\UVL} \ex^{\lgd(\xv;\uv)} \, g(\uv) \, d\uv - \int_{\UVL} \ex^{- \| \DV_{\GP} \uv \|^{2}/2} \, g(\uv) \, d\uv} 
		 {\int \ex^{- \| \DV_{\GP} \uv \|^{2}/2} d\uv} 
	\biggr| \, .
\label{errdefdiUaHu}
\end{EQA}
As a special case with \( g(\uv) \equiv 1 \) we obtain an approximation of the denominator in \eqref{IIgfifxutudufxutdu}.
In addition, we have to bound the tail integrals
\begin{EQ}[rcccl]
	\rho
	&=&
	\rho(\UVL)
	& \eqdef &
	\frac{\int \Ind(\uv \not\in \UVL) \, \ex^{\lgd(\xv;\uv)} \, d\uv}{\int \ex^{- \| \DV_{\GP} \uv \|^{2}/2} \, d\uv} \, ,
	\\
	\rho_{\GP}
	&=&
	\rho_{\GP}(\UVL)
	& \eqdef &
	\frac{\int \Ind(\uv \not\in \UVL) \, \ex^{- \| \DV_{\GP} \uv \|^{2}/2} \, d\uv}{\int \ex^{- \| \DV_{\GP} \uv \|^{2}/2} \, d\uv} \, .
\label{rhfiUaceHu22m}
\end{EQ}
Everywhere later \( \gaussv_{\GP} \sim \ND(0,\DV_{\GP}^{-2}) \) is a Gaussian element in \( \R^{\dimp} \).
The analysis will be split into several steps.

\Section{Overall error of Laplace approximation}
First we show how to seam together the error \( \err \) of local approximation and the bounds for the tail integrals
\( \rho \) and \( \rho_{\GP} \); see \eqref{rhfiUaceHu22m}.

\begin{proposition}
\label{PunbintLapl}
Suppose that for a function \( g(\uv) \in [0,1] \) and some \( \err , \err_{g} \)
\begin{EQA}[rcl]
	\biggl| 
		\frac{\int_{\UVL} \ex^{\lgd(\xv;\uv)} \, d\uv - \int_{\UVL} \ex^{- \| \DV_{\GP} \uv \|^{2}/2} \, d\uv}			 
			 {\int \ex^{- \| \DV_{\GP} \uv \|^{2}/2} \, d\uv} 
	\biggr|
	& \leq &
	\err \, ,
\label{erifxudieHu22}
	\\
\label{erifxudieHu22g}
	\biggl| 
		\frac{\int_{\UVL} g(\uv) \, \ex^{\lgd(\xv;\uv)} \, d\uv - \int_{\UVL} g(\uv) \, \ex^{- \| \DV_{\GP} \uv \|^{2}/2} \, d\uv}			 {\int \ex^{- \| \DV_{\GP} \uv \|^{2}/2} \, d\uv} 
	\biggr|
	& \leq &
	\err_{g} \, .
\end{EQA}
Then with \( \rho \) and \( \rho_{\GP} \) from \eqref{rhfiUaceHu22m}
\begin{EQ}[rcl]
	\frac{\int g(\uv) \, \ex^{\lgd(\xv;\uv)} \, d\uv}{\int \ex^{\lgd(\xv;\uv)} \, d\uv} 
	& \leq &
	\frac{1}{1 - \rho_{\GP} - \err} \,\, 
	\frac{\int g(\uv) \, \ex^{- \| \DV_{\GP} \uv \|^{2}/2} \, d\uv}{\int \ex^{- \| \DV_{\GP} \uv \|^{2}/2} \, d\uv}
	+ \frac{\rho + \err_{g}}{1 - \rho_{\GP} - \err} \, ,
	\qquad
	\\
	\frac{\int g(\uv) \, \ex^{\lgd(\xv;\uv)} \, d\uv}{\int \ex^{\lgd(\xv;\uv)} \, d\uv} 
	& \geq &
	\frac{1}{1 + \rho + \err} \,\, 
	\frac{\int g(\uv) \, \ex^{- \| \DV_{\GP} \uv \|^{2}/2} \, d\uv}{\int \ex^{- \| \DV_{\GP} \uv \|^{2}/2} \, d\uv}
	- \frac{\rho_{\GP} + \err_{g}}{1 + \rho + \err} \, .
	\qquad
\label{igefxumiguexHu22}
\end{EQ}
\end{proposition}

\begin{proof}
It follows from \eqref{erifxudieHu22} 
\begin{EQA}
	\int \ex^{\lgd(\xv;\uv)} \, d\uv
	& \geq &
	\int_{\UVL} \ex^{\lgd(\xv;\uv)} \, d\uv
	\geq 
	\int_{\UVL} \ex^{- \| \DV_{\GP} \uv \|^{2}/2} \, d\uv - \err \int \ex^{- \| \DV_{\GP} \uv \|^{2}/2} \, d\uv
	\\
	& \geq &
	(1 - \err - \rho_{\GP}) \int \ex^{- \| \DV_{\GP} \uv \|^{2}/2} \, d\uv ,
\label{1erriexmH22du22}
	\\
	\int \ex^{\lgd(\xv;\uv)} \, d\uv
	& \leq &
	\int_{\UVL} \ex^{\lgd(\xv;\uv)} \, d\uv + \rho \int \ex^{- \| \DV_{\GP} \uv \|^{2}/2} \, d\uv
	\\
	& \leq &
	(1 + \err + \rho) \int \ex^{- \| \DV_{\GP} \uv \|^{2}/2} \, d\uv .
	\qquad
\label{1erriexmH22du22u}
\end{EQA}
Similarly for \( g(\uv) \geq 0 \)
\begin{EQA}
	\int g(\uv) \, \ex^{\lgd(\xv;\uv)} \, d\uv
	& \geq &
	\int_{\UVL} g(\uv) \, \ex^{- \| \DV_{\GP} \uv \|^{2}/2} \, d\uv - \err_{g} \int \ex^{- \| \DV_{\GP} \uv \|^{2}/2} \, d\uv
	\\
	& \geq &
	\int g(\uv) \, \ex^{- \| \DV_{\GP} \uv \|^{2}/2} \, d\uv
	- (\rho_{\GP} + \err_{g}) \int \ex^{- \| \DV_{\GP} \uv \|^{2}/2} \, d\uv ,
\label{jhdjsmswnqwyts5423nd}
\end{EQA}	
\begin{EQA}
	\int g(\uv) \, \ex^{\lgd(\xv;\uv)} \, d\uv
	& \leq &
	\int_{\UVL} g(\uv) \, \ex^{\lgd(\xv;\uv)} \, d\uv + \rho \int \ex^{- \| \DV_{\GP} \uv \|^{2}/2} \, d\uv
	\\
	& \leq &
	\int g(\uv) \, \ex^{- \| \DV_{\GP} \uv \|^{2}/2} \, d\uv + (\rho + \err_{g}) \int \ex^{- \| \DV_{\GP} \uv \|^{2}/2} \, d\uv \, .
\label{iguexmHu22dreieH}
\end{EQA}
Putting together all these bounds yields \eqref{igefxumiguexHu22}.
\end{proof}

The next corollary is straightforward.

\begin{corollary}
\label{CPunbintLapl}
Let \( \rho_{\GP} \leq \rhos \), \( \rho \leq \rhos \); see \eqref{rhfiUaceHu22m}.
Let also for a function \( g(\uv) \) with \( |g(\uv)| \leq 1 \), \eqref{erifxudieHu22}, \eqref{erifxudieHu22g} hold with
\( \err_{g} \leq \err \).
If \( \err + \rhos \leq 1/2 \) then 
\begin{EQA}
	\left| 
		\frac{\int g(\uv) \, \ex^{\lgd(\xv;\uv)} \, d\uv}{\int \ex^{\lgd(\xv;\uv)} \, d\uv} 
		- \frac{\int g(\uv) \, \ex^{- \| \DV_{\GP} \uv \|^{2}/2} \, d\uv}{\int \ex^{- \| \DV_{\GP} \uv \|^{2}/2} \, d\uv}
	\right|
	& \leq &
	\frac{2 (\rhos + \err)}{1 - \rhos - \err} 
	\leq 
	4 (\rhos + \err)	\, . 
\label{2rherIHg2re}
\end{EQA}
\end{corollary}


Sometimes we need an extension to the case of an unbounded function \( g \).
This particularly arises when evaluating the moment of the posterior; see Theorem~\ref{TpostmeanLa}.
The next result corresponds to estimation of posterior mean with a linear function \( g \) 
and posterior variance with \( g \) quadratic.

\begin{proposition}
\label{PunbintLapg}
Given a function \( g(\uv) \), assume \eqref{erifxudieHu22}, \eqref{erifxudieHu22g}, and define
\begin{EQA}
	\rho_{g}
	& \eqdef &
	\frac{\int \Ind(\uv \not\in \UVL) \, | g(\uv) | \, \ex^{\lgd(\xv;\uv)} \, d\uv}{\int \ex^{- \| \DV_{\GP} \uv \|^{2}/2} \, d\uv} 
	 \, ,
	\\
	\rho_{\GP,g}
	& \eqdef &
	\frac{\int \Ind(\uv \not\in \UVL) \, | g(\uv) | \, \ex^{- \| \DV_{\GP} \uv \|^{2}/2} \, d\uv}{\int \ex^{- \| \DV_{\GP} \uv \|^{2}/2} \, d\uv} \, ,
\label{rhfiUaceHu22mg}
\end{EQA}
while \( \rho \) and \( \rho_{\GP} \) are given in \eqref{rhfiUaceHu22m}.
Then for \( \II_{\GP}(g) = \E g(\gaussv_{\GP}) \), \( \gaussv_{\GP} \sim \ND(0,\DV_{\GP}^{-2}) \), 
\begin{EQA}
	&& \nquad
	\left| 
		\frac{\int g(\uv) \, \ex^{\lgd(\xv;\uv)} \, d\uv}{\int \ex^{\lgd(\xv;\uv)} \, d\uv} 
		- \frac{\int g(\uv) \, \ex^{- \| \DV_{\GP} \uv \|^{2}/2} \, d\uv}{\int \ex^{- \| \DV_{\GP} \uv \|^{2}/2} \, d\uv}
	\right|
	\leq 
	\frac{\rho_{g} + \rho_{\GP,g} + \err_{g}}{1 - \rho_{\GP} - \err}
	+ \frac{| \II_{\GP}(g) | \, (\rho + \err)}{1 - \rho_{\GP} - \err} \, . 
\label{2rherIHg2red}
\end{EQA}
In particular, if \( \int g(\uv) \, \ex^{-\| \DV_{\GP} \uv \|^{2}/2} \, d\uv = 0 \) then
\begin{EQA}[rcl]
	\left| \frac{\int g(\uv) \, \ex^{\lgd(\xv;\uv)} \, d\uv}{\int \ex^{\lgd(\xv;\uv)} \, d\uv} \right| 
	& \leq &
	\frac{\rho_{g} + \rho_{\GP,g} + \err_{g}}{1 - \rho_{\GP} - \err} \, .
\label{igefxumiguexHu22g}
\end{EQA}
\end{proposition}

\begin{proof}
Suppose that \( \II_{\GP}(g) \geq 0 \).
Then 
\begin{EQA}
	&& \nquad
	\left| 
		\frac{\int g(\uv) \, \ex^{\lgd(\xv;\uv)} \, d\uv}{\int \ex^{\lgd(\xv;\uv)} \, d\uv} 
		- \frac{\int g(\uv) \, \ex^{- \| \DV_{\GP} \uv \|^{2}/2} \, d\uv}{\int \ex^{- \| \DV_{\GP} \uv \|^{2}/2} \, d\uv}
	\right|
	\\
	& \leq &
	\left| 
		\frac{\int g(\uv) \, \ex^{\lgd(\xv;\uv)} \, d\uv}{\int \ex^{\lgd(\xv;\uv)} \, d\uv} 
		- \frac{\int g(\uv) \, \ex^{- \| \DV_{\GP} \uv \|^{2}/2} \, d\uv}{\int \ex^{\lgd(\xv;\uv)} \, d\uv}
	\right|
	+ \II_{\GP}(g)
	\left| 
		\frac{\int \ex^{- \| \DV_{\GP} \uv \|^{2}/2} \, d\uv}{\int \ex^{\lgd(\xv;\uv)} \, d\uv} - 1
	\right| .
\label{igefxvv22igh}
\end{EQA}
By definitions
\begin{EQA}
	&& \nquad
	\left| 
		\int g(\uv) \, \ex^{\lgd(\xv;\uv)} \, d\uv - \int g(\uv) \, \ex^{-\| \DV_{\GP} \uv \|^{2}/2} \, d\uv
	\right|
	\\
	& \leq &
	\left| 
		\int_{\UVL} g(\uv) \, \ex^{\lgd(\xv;\uv)} \, d\uv - \int_{\UVL} g(\uv) \, \ex^{-\| \DV_{\GP} \uv \|^{2}/2} \, d\uv
	\right|
	\\
	&&
	+ \, \left| 
		\int \Ind(\uv \not\in \UVL) \, g(\uv) \, \ex^{\lgd(\xv;\uv)} \, d\uv 
	\right|
	+ \left| 
		\int \Ind(\uv \not\in \UVL) \, g(\uv) \, \ex^{-\| \DV_{\GP} \uv \|^{2}/2} \, d\uv
	\right|	
	\\
	& \leq &
	\bigl( \rho_{g} + \rho_{\GP,g} + \err_{g} \bigr) \int \ex^{-\| \DV_{\GP} \uv \|^{2}/2} \, d\uv
\label{lrguegfxumHu22}
\end{EQA}
and the assertion follows by \eqref{1erriexmH22du22} and \eqref{1erriexmH22du22u}.
\end{proof}

\ifunivariate{}{
\Section{Local approximation. Univariate case}
This section presents a bound of local approximation for the univariate case.
Later we provide an independent study of the multivariate case.

Define \( \lgd(x;u) = \lgd(x+u) - \lgd(x) - f'(x) u \), 
\begin{EQA}
	\dltw_{3}(x,u)
	& \eqdef &
	\lgd(x+u) - \lgd(x) - f'(x) u - f''(x) \frac{u^{2}}{2} \, .
\label{dl3xufxufxufppxu2}
\end{EQA}

\begin{proposition}
\label{Lintapproxdu}
Let a function \( \lgd(u) \) be two times continuously differentiable and satisfy 
\( f''(x) = - \DU^{2} < 0 \).
If for some \( \zq > 0 \) 
\begin{EQA}
	\sup_{|u| \leq \zq} \frac{1}{\DU^{2} u^{2}/2}\bigl| \dltw_{3}(x,u) \bigr|
	& \leq &
	\dltwb  
	\leq 
	1/3 ,
\label{123stzf3tH2}
\end{EQA}
then for any function \( g(u) \) with \( |g(u)| \leq 1 \) 
\begin{EQA}
	\frac{\bigl| \int_{-\zq}^{\zq}  \ex^{\lgd(x;u)} \, g(u) \, du 
			- \int_{-\zq}^{\zq} \ex^{- \DU^{2} u^{2}/2} \, g(u) \, du \bigr|}
		 {\int \ex^{- \DU^{2} u^{2}/2} \, du}
	& \leq &
	\err_{2}
	\eqdef
	\frac{\dltwb}{2 (1 - \dltwb)} \, .
\label{HimzzeftgttmH3}
\end{EQA}
Moreover, if, in addition to \eqref{123stzf3tH2}, for some \( \fba^{(3)} = \fba^{(3)}(x) \)
\begin{EQA}
	\bigl| \dltw_{3}(x,u) \bigr|
	& \leq &
	\frac{1}{6} \bigl| \fba^{(3)} u^{3} \bigr|,
	\qquad
	|u| \leq \zq,
\label{123stzf3tH3uzq}
\end{EQA}
then
\begin{EQA}
	\frac{\bigl| \int_{-\zq}^{\zq}  \ex^{\lgd(x;u)} \, g(u) \, du 
			- \int_{-\zq}^{\zq} \ex^{- \DU^{2} u^{2}/2} \, g(u) \, du \bigr|}
		 {\int \ex^{- \DU^{2} u^{2}/2} \, du}
	& \leq &
	\err_{3}
	\eqdef
	\frac{0.7 \bigl| \fba^{(3)} \bigr|}{\DU^{3}} \, .
\label{HimzzeftgttmH1}
\end{EQA}
\end{proposition}

\begin{proof}
By \eqref{123stzf3tH2} and \( \dltwb \leq 1/3 \), it holds 
\begin{EQA}
	&& \nquad
	\frac{\DU}{\sqrt{2\pi}}
	\biggl| 
		\int_{-\zq}^{\zq}  \ex^{\lgd(x;u)} \, g(u) \, du - \int_{-\zq}^{\zq} \ex^{- \DU^{2} u^{2}/2} \, g(u) \, du 
	\biggr|
	\\
	& \leq &
	\frac{\DU}{\sqrt{2\pi}}
	\biggl| 
		\int_{-\zq}^{\zq}  \ex^{- (1 - \dltwb) \DU^{2} u^{2}/2} \, du 
		- \int_{-\zq}^{\zq} \ex^{- \DU^{2} u^{2}/2} \, du 
	\biggr|
	\\
	& \leq &
	\frac{1}{\sqrt{2\pi}}
	\biggl| 
		\int  \ex^{- (1 - \dltwb) u^{2}/2} \, du 
		- \int \ex^{- u^{2}/2} \, du 
	\biggr|
	\leq 
	(1 - \dltwb)^{-1/2} - 1 .
\label{Himzzeftgt}
\end{EQA}
We now use that \( (1 - \dltwb)^{-1/2} - 1 \leq  0.5 \dltwb/(1 - \dltwb) \) for \( \dltwb < 1 \),
and the result \eqref{HimzzeftgttmH3} follows.

Now we show \eqref{HimzzeftgttmH1}. It holds
\begin{EQA}
	\int_{-\zq}^{\zq} \ex^{\lgd(x;u)} \, g(u) \, du
	& = &
	\int_{-\zq}^{\zq} \ex^{- \DU^{2} u^{2}/2 + \dltw_{3}(x,u)} \, g(u) \, du .
\label{iUefgu22guu}
\end{EQA}
Define for \( t \geq 0 \)
\begin{EQA}
	\Rem(t)
	&=&
	\int_{-\zq}^{\zq} \ex^{- \DU^{2} u^{2}/2 + t \dltw_{3}(x,u)} \, g(u) \, du .
\label{PhtiUexmHugu}
\end{EQA}
Then for \( t \in [0,1] \) by \eqref{123stzf3tH2} and \eqref{123stzf3tH3uzq}
\begin{EQA}
	\bigl| \Rem'(t) \bigr|
	&=&
	\left| \int_{-\zq}^{\zq} \dltw_{3}(x,u) \ex^{- \DU^{2} u^{2}/2 + t \dltw_{3}(x,u)} \, g(u) \, du \right|
	\leq 
	\int_{-\zq}^{\zq} \bigl| \dltw_{3}(x,u) \bigr| \ex^{- (1 - \dltwb) \DU^{2} u^{2}/2} \, du 
	\\
	& \leq &
	\frac{1}{6} \int_{-\zq}^{\zq} \bigl| \fba^{(3)} u^{3} \bigr| \ex^{- (1 - \dltwb) \DU^{2} u^{2}/2} \, du
	\leq 
	\frac{\bigl| \fba^{(3)} \bigr| \, \DU^{-4}}{6 (1 - \dltwb)^{2}} \int |u|^{3} \ex^{- u^{2}/2} \, du
	= 
	\frac{2\bigl| \fba^{(3)} \bigr| \, \DU^{-4}}{3 (1 - \dltwb)^{2}}
\label{Phit12u2H02u32}
\end{EQA}
and it holds
\begin{EQA}
	\bigl| \Rem(1) - \Rem(0) \bigr|
	& \leq &
	\sup_{t \in [0,1]} \bigl| \Rem'(t) \bigr|
	\leq 
	\frac{2\bigl| \fba^{(3)} \bigr| \, \DU^{-4}}{3 (1 - \dltwb)^{2}} \, 
\label{Pd2d3H0wem22}
\end{EQA}
yielding for \( \dltwb \leq 1/3 \)
\begin{EQA}
	\frac{\bigl| \Rem(1) - \Rem(0) \bigr|}{\int \ex^{- \DU^{2} u^{2}/2} du}
	& \leq &
	\frac{2\bigl| \fba^{(3)} \bigr| \, \DU^{-4}}{3 (1 - \dltwb)^{2} \sqrt{2\pi} \, \DU^{-1}}
	\leq 
	\frac{0.27 \bigl| \fba^{(3)} \bigr|}{(1 - \dltwb)^{2} \DU^{3}} \, 
\label{Phd3Ph0iUmHu22}
\end{EQA}
and \eqref{HimzzeftgttmH1} follows.
\end{proof}

Now we present some results based on a higher order Taylor approximation.
Define 
\begin{EQA}
	\dltw_{4}(x,u)
	& \eqdef &
	\lgd(x+u) - \lgd(x) - f'(x) u + \frac{\DU^{2} u^{2}}{2} - \frac{f^{(3)}(x) u^{3}}{6} 
	\\
	&=&
	\lgd(x,u) + \frac{\DU^{2} u^{2}}{2} - \frac{f^{(3)}(x) u^{3}}{6} \, .
\label{124stzf4tH24}
\end{EQA}
With \( \dltw_{3}(x,u) = \lgd(x,u) + \DU^{2} u^{2}/2 \), it obviously holds
\begin{EQA}
	\dltw_{3}(x,u)
	& = &
	\dltw_{4}(x,u) + {f^{(3)}(x) u^{3}}/{6} .
\label{d3d4uf30u36}
\end{EQA}

\begin{proposition}
\label{Lintapproxdu4}
Let \( \lgd(u) \) be three times continuously differentiable and satisfy \eqref{123stzf3tH2}.
Define \( \dltwm_{3} = f^{(3)}(x) / \DU^{2} \).
If \( \dltw_{4}(x,u) \) from \eqref{124stzf4tH24} satisfies
\begin{EQA}
	\sup_{|u| \leq \zq} \frac{24}{\DU^{2} u^{4}} \bigl| \dltw_{4}(x,u) \bigr|
	& \leq &
	\dltwm_{4}  
\label{124stzf4tH2}
\end{EQA}
then for any function \( g(u) \) with \( |g(u)| \leq 1 \) and \( g(u) = g(-u) \)
\begin{EQA}
	\frac{
	\bigl| \int_{-\zq}^{\zq}  \ex^{\lgd(x;u)} \, g(u) \, du - \int_{-\zq}^{\zq} \ex^{- \DU^{2} u^{2}/2} \, g(u) \, du \bigr|}
	{\int \ex^{- \DU^{2} u^{2}/2} \, du}
	& \leq &
	\err_{4} \, .
\label{HimzzeftgttmH2}
\end{EQA}
with 
\begin{EQA}
	\err_{4}
	& \eqdef &
	\left\{ \frac{5 \dltwm_{3}^{2}}{12 (1 - \dltwb)^{7/2} }
	+ \frac{\dltwm_{4}}{4 (1 - \dltwb)^{5/2}} \right\} \DU^{-2} .
\label{e4l1d372441m3}
\end{EQA}
\end{proposition}

\begin{proof}
We write \( \dltw_{3}(u) \) in place of \( \dltw_{3}(x,u) \) and \( f^{(3)} \) in place of \( f^{(3)}(x) \).
It holds 
\begin{EQA}
	\int_{-\zq}^{\zq} \ex^{\lgd(x;u)} \, g(u) \, du
	& = &
	\int_{-\zq}^{\zq} \exp\Bigl\{ - \frac{\DU^{2} u^{2}}{2} + \dltw_{3}(u) \Bigr\} \, g(u) \, du .
\label{P3tHmzzmt22t3}
\end{EQA}
Define for \( t \in [0,1] \)
\begin{EQA}
	\Rem(t)
	& \eqdef &
	\int_{-\zq}^{\zq} \exp\Bigl\{ - \frac{\DU^{2} u^{2}}{2} + t \dltw_{3}(u) \Bigr\} \, g(u) \, du .
\label{P3tHmzzmt22t}
\end{EQA}
Symmetricity \( g(u) = g(-u) \) implies that 
\begin{EQA}
	\Rem'(0)
	&=&
	\frac{1}{2} \int_{-\zq}^{\zq} 
	\exp\Bigl( - \frac{\DU^{2} u^{2}}{2} \Bigr) \bigl\{ \dltw_{3}(u) + \dltw_{3}(-u) \bigr\} \, g(u) \, du
	=
	\int_{-\zq}^{\zq} \exp\Bigl( - \frac{\DU^{2} u^{2}}{2} \Bigr) \dltw_{4}(u) \, g(u) \, du \, .
\label{Pp012zzmH2u212}
\end{EQA}
Moreover, as  \( |\dltw_{3}(u)| \leq \dltwb \, \DU^{2} u^{2}/2 \),
it holds for \( t \in [0,1] \) 
\begin{EQA}
	|\Rem''(t)|
	&=&
	\left| \int_{-\zq}^{\zq} \dltw_{3}^{2}(u) 
		\exp\Bigl\{ - \frac{\DU^{2} u^{2}}{2} + t \dltw_{3}(u) \Bigr\} \, g(u) \, du 
	\right|
	\leq 
	\int_{-\zq}^{\zq} \dltw_{3}^{2}(u) 
		\exp\Bigl( - \frac{(1 - \dltwb) \DU^{2} u^{2}}{2} \Bigr) \, du \, .
\label{F3pptHm2t6}
\end{EQA}
As \( \dltw_{3}(u) = f^{(3)} u^{3}/6 + \dltw_{4}(u) \) and \( |\dltw_{4}(u)| \leq 1 \), it holds 
\begin{EQA}
	|\Rem''(t)| 
	& \leq & 
	2 \int_{-\zq}^{\zq} \bigl\{ \dltw_{4}^{2}(u) + \bigl| f^{(3)} u^{3}/6 \bigr|^{2} \bigr\} 
		\exp\Bigl( - \frac{(1 - \dltwb)\DU^{2} u^{2}}{2} \Bigr) \, du
	\\
	& \leq & 
	2 \int_{-\zq}^{\zq} \bigl\{ |\dltw_{4}(u)| + \bigl| f^{(3)} \bigr|^{2} u^{6}/36 \bigr\} 
	\exp\Bigl( - \frac{(1 - \dltwb) \DU^{2} u^{2}}{2} \Bigr) \, du .
\label{u221md3u636}
\end{EQA}
The use of \( \bigl| f^{(3)} \bigr| = \dltwm_{3} \, \DU^{2} \) and \eqref{124stzf4tH2} yields
in view of \( \dltwb \leq 1/3 \) and \( |\dltw_{4}(u)| \leq \dltwm_{4} \, \DU^{2} u^{4} \)
\begin{EQA}
	&& \nquad
	\bigl| \Rem(1) - \Rem(0) \bigr|
	\leq 
	\bigl| \Rem'(0) \bigr| + \frac{1}{2} \sup_{t \leq 1} |\Rem''(t)|
	\\
	& \leq &
	2 \int_{-\zq}^{\zq} |\dltw_{4}(u)| \, \exp\Bigl( - \frac{(1 - \dltwb) \DU^{2} u^{2}}{2} \Bigr) \, du
	+
	\frac{\dltwm_{3}^{2} \, \DU^{4}}{36} \int_{-\zq}^{\zq} u^{6} \, 
		\exp\Bigl( - \frac{(1 - \dltwb)\DU^{2} u^{2}}{2} \Bigr) \, du 
	\\
	& \leq &
	\frac{\dltwm_{4} \, \DU^{2}}{12} 
	\int_{-\zq}^{\zq} u^{4} \, \exp\Bigl( - \frac{(1 - \dltwb) \DU^{2} u^{2}}{2} \Bigr) \, du
	+ \frac{\dltwm_{3}^{2} \, \DU^{4}}{36} 
	\int_{-\zq}^{\zq} u^{6} \, \exp\Bigl( - \frac{(1 - \dltwb)\DU^{2} u^{2}}{2} \Bigr) \, du 
	\\
	& \leq &
	\frac{\dltwm_{4} (1 - \dltwb)^{-5/2}}{12 \, \DU^{3}}	
	\int u^{4} \, \exp\Bigl( - \frac{\DU^{2} u^{2}}{2} \Bigr) \, du
	+ \frac{\dltwm_{3}^{2} (1 - \dltwb)^{-7/2}}{36 \, \DU^{3}} 
	\int u^{6} \, \exp\Bigl( - \frac{\DU^{2} u^{2}}{2} \Bigr) \, du 
	\\
	&=&
	\frac{3 \sqrt{2\pi} \, \dltwm_{4} (1 - \dltwb)^{-5/2} }{12 \, \DU^{3}}	
	+ \frac{15 \sqrt{2\pi} \, \dltwm_{3}^{2} (1 - \dltwb)^{-7/2}}{36 \, \DU^{3}}
\label{C4Hm2Fm4tH2t4}
\end{EQA}
yielding \eqref{HimzzeftgttmH2} in view of 
\( \int \ex^{- \DU^{2} u^{2}/2} \, du = \sqrt{2\pi} \, \DU^{-1} \)
and \( \dltwm_{3} \leq 1/3 \).
\end{proof}

\begin{remark}
Suppose that \( \lgd(\cdot) \) is 
three times continuously differentiable and \( |f^{(3)}(x+u)| \leq \dltwu_{3} \DU^{2} \)
for \( |u| \leq \zq \).
Then the expansion \eqref{123stzf3tH2} applies with \( \dltwb \leq \dltwu_{3} \zq/3 \).
The use of \( \zq = \xx/\DU \) yields
\begin{EQA}
	\dltwb \leq \dltwu_{3} \zq/ 3
	& \leq &
	\CONST \, \dltwu_{3} \, \xx \, \DU^{-1} .
\label{Cw5x5Hm3120}
\end{EQA}
Similarly, if \( \lgd(\cdot) \) is four times continuously differentiable and \( |f^{(4)}(u)| \leq \dltwm_{4} \DU^{2} \)
for \( |u| \leq \zq \), then \eqref{124stzf4tH2} holds as well.
One can see that the use of a higher order approximation of \( \lgd \) allows to improve 
the accuracy of approximation  
from \( \DU^{-1} \) as in \eqref{HimzzeftgttmH1} of Proposition~\ref{Lintapproxdu} to \( \DU^{-2} \) as in \eqref{HimzzeftgttmH2} of Proposition~\ref{Lintapproxdu4}.

\end{remark}
}

\Section{Lower and upper Gaussian measures}
This section introduces the lower and upper Gaussian measure which locally sandwich the measure \( \PfL \)
using the decomposition from condition \nameref{LLf0ref}. 
Denote \( - \nabla^{2} \lgd(\xv) = \DV_{\GP}^{2} \).
Definition \eqref{om3esuU1H02d3} enables us to bound with \( \dltwb = \dltwb(\xv) \)
\begin{EQA}
	\frac{1}{2} (\| \DV_{\GP} \uv \|^{2} - \dltwb \| \DVL \uv \|^{2})
	& \leq &
	\lgd(\xv;\uv)
	\leq 
	\frac{1}{2} (\| \DV_{\GP} \uv \|^{2} + \dltwb \| \DVL \uv \|^{2}) \, 
\label{ysdtydsfrtswedftedswfhd}
\end{EQA}
yielding two Gaussian measures which bounds \( \PfL \) locally from above and from below.
The next technical result provides sufficient conditions for their contiguity.

\begin{proposition}
\label{LlocintgrL}
Let \( \dltwb \) from \eqref{om3esuU1H02d3} satisfy \( \dltwb \leq 1/3 \).
Then with \( \dimL \) from \eqref{dAdetrH02Hm2}
\begin{EQA}
\label{12Ip2t3t4Hp3}
	\det\bigl( \Id + \dltwb \DV_{\GP}^{-1} \DVL^{2} \, \DV_{\GP}^{-1} \bigr)
	& \leq &
	\exp (\dltwb \, \dimL) \, ,
	\\
	\det\bigl( \Id - \dltwb \DV_{\GP}^{-1} \DVL^{2} \, \DV_{\GP}^{-1} \bigr)^{-1/2}
	& \leq &
	\exp \bigl\{ 3/2 \log (3/2) \, \dltwb \, \dimL \bigr\} .
	\qquad
\label{12Im2t3t4Hm30}
\end{EQA}
\end{proposition}

\begin{proof}
W.l.o.g. assume that \( \DV_{\GP}^{-1} \DVL^{2} \, \DV_{\GP}^{-1} \) is diagonal with 
eigenvalues \( \lambda_{j} \in [0,1] \).
As \( - x^{-1}\log(1 - x) \leq 3 \log(3/2) \) for \( x \in [0,1/3] \), it holds by \eqref{om3esuU1H02d3}
\begin{EQA}
	&& \nquad
	\log \det\bigl( \Id - \dltwb \, \DV_{\GP}^{-1} \DVL^{2} \, \DV_{\GP}^{-1} \bigr)^{-1}
	=
	- \sum_{j=1}^{\dimp} \log\bigl( 1 - \dltwb \lambda_{j} \bigr)
	\leq 
	3 \log (3/2) \sum_{j=1}^{\dimp} \dltwb \lambda_{j}
	\\
	&=&
	3 \log (3/2) \, \dltwb \, \tr\bigl( \DV_{\GP}^{-1} \DVL^{2} \, \DV_{\GP}^{-1} \bigr)
	=
	3 \log (3/2) \, \dltwb \, \dimL \, 
\label{125d3p0ldm1}
\end{EQA}
yielding \eqref{12Im2t3t4Hm30}.
The proof of \eqref{12Ip2t3t4Hp3} is similar using \( \log(1 + x) \leq x \) for \( x \geq 0 \).
\end{proof}

\Subsection{Gaussian moments}
The presented bounds involve the Gaussian moments \( \E \| \DV \gaussv_{\GP} \|^{k} \) for \( k=3,4,6 \)
and \( \gaussv_{\GP} \sim \ND(0,\DPGP^{-2}) \).
We make use of the following lemma.

\begin{lemma}
\label{LdltwLa}
It holds for \( \gaussv_{\GP} \sim \ND(0,\DPGP^{-2}) \)
\begin{EQA}
	\E \| \DVL \, \gaussv_{\GP} \|^{3}
	& \leq &
	(\dimL + 1)^{3/2} \, ,
	\\
	\E \| \DVL \, \gaussv_{\GP} \|^{4}
	& \leq &
	(\dimL + 1)^{2} \, ,
	\\
	\E \| \DVL \, \gaussv_{\GP} \|^{6}
	& \leq &
	(\dimL + 2)^{3} \, .
\label{eciu8hef8h8we3hy87u3y7y3fe7}
\end{EQA}
\end{lemma}

\begin{proof}
Represent \( \| \DVL \, \gaussv_{\GP} \|^{2} = \| \DVL \, \DVL_{\GP}^{-1} \gaussv \|^{2} 
= \langle \BB_{\GP} \gaussv, \gaussv \rangle \)
with \( \BB_{\GP} = \DVL \, \DVL_{\GP}^{-2} \DVL \leq \Id_{\dimp} \) and \( \gaussv \sim \ND(0,\Id_{\dimp}) \).
By Lemma~\ref{Gaussmoments} 
\begin{EQA}
	\E \| \DVL \, \gaussv_{\GP} \|^{4}
	&=&
	\E \bigl\langle \BB_{\GP} \gaussv, \gaussv \bigr\rangle^{2}
	\leq 
	\bigl\{ \tr (\BB_{\GP}) + 1 \bigr\}^{2}
	=
	(\dimL + 1)^{2} ,
	\\
	\E \| \DVL \, \gaussv_{\GP} \|^{6}
	&=&
	\E \bigl\langle \BB_{\GP} \gaussv, \gaussv \bigr\rangle^{3}
	\leq 
	\bigl\{ \tr (\BB_{\GP}) + 2 \bigr\}^{3}
	=
	(\dimL + 2)^{3} ,
\label{lkvu7rerycfrrwt3ghdtndf}
\end{EQA}
and \( \E \| \DVL \, \gaussv_{\GP} \|^{3} \leq \E^{3/4} \| \DVL \, \gaussv_{\GP} \|^{4}	\leq (\dimL + 1)^{3/2} \).
\end{proof}

\Section{Local approximation}
\label{SLaplapprmu}

This section presents the bounds on the error \( \err \) of local approximation \eqref{errdefdiUaHu}.
The first result only uses \( \dltwb \, \dimL \leq 2/3 \).
More advanced bounds also assume \nameref{LL3tref} and \nameref{LL4tref} with \( \upsv = \xv \) and \( \HLG(\upsv) = \DV \).
We also present some extensions for the moments of \( \PfL \).

\begin{proposition}
\label{Lintfxupp3}
Let \( \dltwb = \dltwb(\xv) \) from \eqref{om3esuU1H02d3} and \( \dimL \) from \eqref{dAdetrH02Hm2} satisfy 
\begin{EQA}
	\dltwb \, \dimL
	& \leq &
	2/3 \, .
\label{w3p0le13fx}
\end{EQA}
Then for any function \( g(\uv) \) with \( |g(\uv)| \leq 1 \) 
\begin{EQA}
	\biggl| 
	\frac{\int_{\UVL} \ex^{\lgd(\xv;\uv)} \, g(\uv) \, d\uv - \int_{\UVL} \ex^{- \| \DV_{\GP} \uv \|^{2}/2} \, g(\uv) \, d\uv} 
		 {\int \ex^{- \| \DV_{\GP} \uv \|^{2}/2} d\uv} 
	\biggr|
	& \leq &
	\err \, 
\label{ed3le2d3pGd}
\end{EQA}
with
\begin{EQA}
	\err
	&=&
	\err_{2} 
	= 
	\frac{0.75 \, \dltwb \, \dimL}{1 - \dltwb} \, .
\label{8fjre34et6dehgweyt6wsh}
\end{EQA}
\end{proposition}

\begin{proof}
The condition \( \dltwb \, \dimL \leq 2/3 \) from \eqref{w3p0le13fx} and bound \eqref{12Im2t3t4Hm30} imply 
\begin{EQA}
	\det\bigl( \Id - \dltwb \DV_{\GP}^{-1} \DVL^{2} \, \DV_{\GP}^{-1} \bigr)^{-1/2}
	& \leq &
	\exp \bigl\{ 3/2 \log (3/2) \, \dltwb \, \dimL \bigr\} 
	\leq 
	3/2 \, .
	\qquad
\label{12Im2t3t4Hm3}
\end{EQA}
Define for \( t \geq 0 \)
\begin{EQA}
	\Rem(t)
	&=&
	\int_{\UVL} \ex^{- \| \DV_{\GP} \uv \|^{2}/2 + t \dltw_{3}(\xv,\uv)} \, g(\uv) \, d\uv .
\label{PhtiUexmHu2mtH0u2g}
\end{EQA}
Then for \( t \in [0,1] \) by \eqref{om3esuU1H02d3} 
\begin{EQA}
	\bigl| \Rem'(t) \bigr|
	&=&
	\left| \int_{\UVL} \dltw_{3}(\xv,\uv) \ex^{- \| \DV_{\GP} \uv \|^{2}/2 + t \dltw_{3}(\xv,\uv)} \, g(\uv) \, d\uv \right|
	\\
	& \leq &
	\int_{\UVL} \bigl| \dltw_{3}(\xv,\uv) \bigr| \ex^{- (\| \DV_{\GP} \uv \|^{2} - \dltwb \| \DVL \uv \|^{2})/2} \, d\uv .
\label{Pht12u212H02u22}
\end{EQA}
Now we make change of variable \( \Idd \uv \) to \( \wv \) with
\( \Idd^{2} = \Id - \dltwb \, \DV_{\GP}^{-1} \DVL^{2} \, \DV_{\GP}^{-1} \).
By \eqref{12Im2t3t4Hm3} \( \det \Idd^{-1} \leq 3/2 \) and also 
\( \| \Idd^{-1} \| \leq (1 - \dltwb)^{-1/2} \).
By \eqref{om3esuU1H02d3} and \eqref{Pht12u212H02u22}
\begin{EQA}
	&& \nquad
	\bigl| \Rem(1) - \Rem(0) \bigr|
	\leq 
	\sup_{t \in [0,1]} \bigl| \Rem'(t) \bigr|
	\leq 
	\frac{\dltwb}{2}
	\int_{\UVL} \| \DVL \uv \|^{2} \ex^{- (\| \DV_{\GP} \uv \|^{2} - \dltwb \| \DVL \uv \|^{2})/2} \, d\uv 
	\\
	& \leq &
	\frac{3 \dltwb}{4}
	\int \| \DVL \Idd^{-1} \uv \|^{2} \ex^{- \| \DV_{\GP} \uv \|^{2} /2} \, d\uv 
	\leq 
	\frac{3 \dltwb}{4(1 - \dltwb)}
	\int \| \DVL \uv \|^{2} \ex^{- \| \DV_{\GP} \uv \|^{2} /2} \, d\uv 	.
\label{ppt2iUH0w2edw}
\end{EQA}
In view of \( \E \| \DVL \, \gaussv_{\GP} \|^{2} = \tr\bigl( \DVL^{2} \, \DV_{\GP}^{-2} \bigr) \) for a standard normal
\( \gaussv \), we derive 
\begin{EQA}
	\frac{\bigl| \Rem(1) - \Rem(0) \bigr|}{\int \ex^{- \| \DV_{\GP} \uv \|^{2}/2} d\uv}
	& \leq &
	\frac{3 \dltwb }{4 (1 - \dltwb)}
	\frac{\int \| \DVL \uv \|^{2} \ex^{- \| \DV_{\GP} \uv \|^{2} /2} \, d\uv}
		{\int \ex^{- \| \DV_{\GP} \uv \|^{2} /2} \, d\uv}
	\leq 
	\frac{3 \dltwb \, \dimL}{4 (1 - \dltwb)}
\label{Phd3Ph0iUmHu22}
\end{EQA}
and \eqref{ed3le2d3pGd} follows.
\end{proof}

\begin{proposition}
\label{Lintfxupp3T}
Assume \nameref{LL3tref} with \( \upsv = \xv \) and \( \HLG(\xv) = \DVL \), 
and let \( \dltwu_{3} \, \amax^{-1} \rrL \, \dimL \leq 2 \). 
Define \( \dltwb \eqdef \dltwu_{3} \, \amax^{-1} \rrL / 3 \).
Then bound \eqref{ed3le2d3pGd} applies with 
\begin{EQA}
	\err
	&=&
	\err_{3}
	\eqdef
	\frac{\dltwu_{3} \E \| \DV \gaussv_{\GP} \|^{3}}{4 (1 - \dltwb)^{3/2}} 
	\leq 
	\frac{\dltwu_{3} \, (\dimL+1)^{3/2}}{4 (1 - \dltwb)^{3/2}} \, .
\label{5qw7dyf4e4354co9dufih}
\end{EQA}
\end{proposition}

\begin{proof}
The proof follows the same line as for Proposition~\ref{Lintfxupp3}.
Under \nameref{LL3tref}, it holds \( |\dltw_{3}(\xv,\uv)| \leq \dltwu_{3} \| \DV \uv \|^{3} / 6 \) for \( \uv \in \UVL \) and
\begin{EQA}
	\bigl| \Rem(1) - \Rem(0) \bigr|
	& \leq &
	\frac{\dltwu_{3} \, \det (\Idd^{-1})}{6} 
	\int \| \DV \, \Idd^{-1} \uv \|^{3} \, \ex^{- \| \DV_{\GP} \uv \|^{2} /2} \, d\uv 
	\\
	& \leq &
	\frac{\dltwu_{3}}{4(1 - \dltwb)^{3/2}}
	\int \| \DV \uv \|^{3} \, \ex^{- \| \DV_{\GP} \uv \|^{2} /2} \, d\uv 
\label{Pd3p2d3H0wem22g1}
\end{EQA}
yielding the statement in view of \( \E \| \DVL \, \gaussv_{\GP} \|^{3} \leq (\dimL+1)^{3/2} \); see Lemma~\ref{LdltwLa}.
It remains to note that 
by Lemma~\ref{LdltwLa3t} \( \dltwb \leq \dltwu_{3} \, \amax^{-1} \rrL \, /3\).
Hence, 
\( \dltwu_{3} \, \amax^{-1} \rrL \, \dimL \leq 2 \) implies \( \dltwb \dimL \leq 2/3 \).
\end{proof}

\noindent
The result can be extended to the case of a \( m \)-homogeneous function \( g(\uv) \).

\begin{proposition}
\label{Pd3t3hot3d3u}
Suppose the conditions of Proposition~\ref{Lintfxupp3} and \nameref{LL3tref}.
Then for \( m \geq 1 \) and any \( m \)-homogeneous function \( g(\cdot) \) with \( g(t\uv) = t^{m} g(\uv) \)
\begin{EQA}
	\biggl| 
	\frac{\int_{\UVL} \ex^{\lgd(\xv;\uv)} \, g(\uv) \, d\uv 
	- \int_{\UVL} \ex^{- \| \DV_{\GP} \uv \|^{2}/2} \, g(\uv) \, d\uv} 
		 {\int \ex^{- \| \DV_{\GP} \uv \|^{2}/2} d\uv} 
	\biggr|
	& \leq &
	\frac{\dltwu_{3} \E \bigl\{ |g(\gaussv_{\GP})| \, \| \DVL \gaussv_{\GP} \|^{3} \bigr\}}{4 (1 - \dltwb)^{(m+3)/2}} \, .
	\qquad
\label{e3dEf3xGHo3}
\end{EQA}
\end{proposition}

\begin{proof}
As in the proof of Proposition~\ref{Lintfxupp3}, under \nameref{LL3tref}, it holds for \( \Rem(t) \) from
\eqref{PhtiUexmHu2mtH0u2g}
\begin{EQA}
	\bigl| \Rem(1) - \Rem(0) \bigr|
	& \leq &
	\frac{\dltwu_{3} \, \det (\Idd^{-1})}{6} 
	\int \| \DV \, \Idd^{-1} \uv \|^{3} \, \bigl| g(\Idd^{-1} \uv) \bigr| \, \ex^{- \| \DV_{\GP} \uv \|^{2} /2} \, d\uv 
	\\
	& \leq &
	\frac{\dltwu_{3}}{4(1 - \dltwb)^{(m+3)/2}}
	\int \| \DV \uv \|^{3} \, \bigl| g(\uv) \bigr| \, \ex^{- \| \DV_{\GP} \uv \|^{2} /2} \, d\uv 
\label{Pd3p2d3H0wem22}
\end{EQA}
yielding similarly to \eqref{Pd3p2d3H0wem22g1}
\begin{EQA}
	\frac{\bigl| \Rem(1) - \Rem(0) \bigr|}{\int \ex^{- \| \DV_{\GP} \uv \|^{2}/2} d\uv}
	& \leq &
	\frac{\dltwu_{3}}{4(1 - \dltwb)^{(m+3)/2}}
	\frac{\int \| \DV \uv \|^{3} \, \bigl| g(\uv) \bigr| \, \ex^{- \| \DV_{\GP} \uv \|^{2} /2} \, d\uv}
		{\int \ex^{- \| \DV_{\GP} \uv \|^{2} /2} \, d\uv} \, .
\label{Phd3Ph0iUmHu22}
\end{EQA}
This yields \eqref{e3dEf3xGHo3}.
\end{proof}

Important special cases correspond to \( m=1 \).

\begin{proposition}
\label{Perrbound3}
Suppose the conditions of Proposition~\ref{Pd3t3hot3d3u}.
Then it holds 
for any linear mapping \( \QP \colon \R^{\dimp} \to \R^{\dimq} \) and any vector \( \av \in \R^{\dimq} \)
\begin{EQA}
	\frac{\bigl| 
		\int_{\UVL} \ex^{\lgd(\xv;\uv)} \, \langle \QP \uv,\av \rangle \, d\uv 
		- \int_{\UVL} \ex^{- \| \DV_{\GP} \uv \|^{2}/2} \, \langle \QP \uv,\av \rangle \, d\uv 	 
		\bigr|} 
		 {\int \ex^{- \| \DV_{\GP} \uv \|^{2}/2} d\uv} 
	& \leq &
	\frac{\dltwu_{3} \, \E \bigl\{ \bigl| \langle \QP \gaussv_{\GP},\av \rangle \bigr| \, \| \DP \gaussv_{\GP} \|^{3} \bigr\}}
		 {4 (1 - \dltwb)^{2}} \, .
	\qquad 
	\qquad
\label{iUafxvdvH2vk3}
\end{EQA}
\end{proposition}

Now we state a sharper result based on \nameref{LL4tref}.

\begin{proposition}
\label{Lintfxupp2}
Suppose the conditions of Proposition~\ref{Lintfxupp3} and \nameref{LL4tref}.
Then for any function \( g(\uv) \) with \( |g(\uv)| \leq 1 \) and 
\( g(\uv) = g(-\uv) \)
\begin{EQA}
	\biggl| 
	\frac{\int_{\UVL} \ex^{\lgd(\xv;\uv)} \, g(\uv) \, d\uv - \int_{\UVL} \ex^{- \| \DV_{\GP} \uv \|^{2}/2} \, g(\uv) \, d\uv} 
		 {\int \ex^{- \| \DV_{\GP} \uv \|^{2}/2} d\uv} 
	\biggr|
	& \leq &
	\err_{4} \, 
\label{efuiguduUem}
\end{EQA}
with
\begin{EQA}
	\err_{4}
	& \eqdef &
	\frac{1}{16 (1 - \dltwb)^{2}} \Bigl\{ \E \bigl\langle \nabla^{3} \lgd(\xv) , \gaussv_{\GP}^{\otimes 3} \bigr\rangle^{2} 
	+ 2 \dltwu_{4} \E \| \DV \gaussv_{\GP} \|^{4} \Bigr\}
	\\
	& \leq &
	\frac{1}{16 (1 - \dltwb)^{2}} \Bigl\{ \dltwu_{3}^{2} \, (\dimL + 2)^{3} + 2 \dltwu_{4} (\dimL + 1)^{2} \Bigr\}
	\, .
\label{errdef3322Hm2}
\end{EQA}
If the function \( g(\cdot) \) is not bounded by one but it is symmetric and \( 2m \)-homogeneous,
i.e. \( g(t \uv) = t^{2m} g(\uv) \), then
\eqref{efuiguduUem} still applies with 
\begin{EQA}
	\err_{4}
	& \eqdef &
	\frac{1}{16 (1 - \dltwb)^{2+m}} \E \Bigl\{ 
		\bigl\langle \nabla^{3} \lgd(\xv) , \gaussv_{\GP}^{\otimes 3} \bigr\rangle^{2} \, g(\gaussv_{\GP})
	+ 2 \dltwu_{4} \, \| \DV \gaussv_{\GP} \|^{4} \, g(\gaussv_{\GP}) \Bigr\}
	\, .
\label{errdef3322Hm2m}
\end{EQA}
\end{proposition}

\begin{proof}
We write \( f^{(3)} \) and \( \dltw_{k}(\uv) \) in place of \( \nabla^{3} \lgd(\xv) \) and
\( \dltw_{k}(\xv,\uv) \), \( k=3,4 \).
It holds 
\begin{EQA}
	\int_{\UVL} \ex^{\lgd(\xv;\uv)} \, g(\uv) \, d\uv
	& = &
	\int_{\UVL} \exp\Bigl\{ - \frac{\| \DV_{\GP} \uv \|^{2}}{2} + \dltw_{3}(\uv) \Bigr\} \, g(\uv) \, d\uv .
\label{P3tHmzzmt22t3v}
\end{EQA}
Define for \( t \in [0,1] \)
\begin{EQA}
	\Rem(t)
	& \eqdef &
	\int_{\UVL} \exp\Bigl\{ - \frac{\| \DV_{\GP} \uv \|^{2}}{2} + t \dltw_{3}(\uv) \Bigr\} \, g(\uv) \, d\uv .
\label{P3tHmzzmt22t}
\end{EQA}
Symmetricity of \( \UVL \) and \( g(\uv) = g(-\uv) \) implies that 
\begin{EQA}
	\Rem'(0)
	&=&
	\frac{1}{2} \int_{\UVL} 
	\exp\Bigl( - \frac{\| \DV_{\GP} \uv \|^{2}}{2} \Bigr) \bigl\{ \dltw_{3}(\uv) + \dltw_{3}(-\uv) \bigr\} \, g(\uv) \, d\uv
	\\
	&=&
	\int_{\UVL} \exp\Bigl( - \frac{\| \DV_{\GP} \uv \|^{2}}{2} \Bigr) \, \bar{\dltw}_{4}(\uv) \, g(\uv) \, d\uv \, 
\label{Pp012zzmH2u212v}
\end{EQA}
with \( \bar{\dltw}_{4}(\uv) = \bigl\{ \dltw_{4}(\uv) + \dltw_{4}(-\uv) \bigr\}/2 \).
Moreover, as  \( |\dltw_{3}(\uv)| \leq \dltwb \| \DVL \uv \|^{2}/2 \),
it holds for \( t \in [0,1] \) 
\begin{EQA}
	|\Rem''(t)|
	& \leq &
	\int_{\UVL} \dltw_{3}^{2}(\uv) 
		\exp\Bigl\{ - \frac{\| \DV_{\GP} \uv \|^{2}}{2} + t \dltw_{3}(\uv) \Bigr\} \, |g(\uv)| \, d\uv
	\\
	& \leq &
	\int_{\UVL} \dltw_{3}^{2}(\uv) 
		\exp\Bigl( - \frac{\| \DV_{\GP} \uv \|^{2} - \dltwb \| \DVL \uv \|^{2}}{2} \Bigr) \, d\uv \, .
\label{F3pptHm2t6v}
\end{EQA}
As \( \dltw_{3}(\uv) = \bigl\langle f^{(3)}, \uv^{\otimes 3} \bigr\rangle/6 + \dltw_{4}(\uv) \) and \( |\dltw_{4}(\uv)| \leq 1 \), one can bound for \( t \in [0,1] \)
\begin{EQA}
	|\Rem''(t)| 
	& \leq & 
	2 \int_{\UVL} \bigl\{ \bar{\dltw}_{4}^{2}(\uv) + \bigl| \bigl\langle f^{(3)}, \uv^{\otimes 3} \bigr\rangle/6 \bigr|^{2} \bigr\} 
		\exp\Bigl( - \frac{\| \DV_{\GP} \uv \|^{2} - \dltwb \| \DVL \uv \|^{2}}{2} \Bigr) \, d\uv
	\\
	& \leq & 
	2 \int_{\UVL} \bigl\{ |\bar{\dltw}_{4}(\uv)| + \bigl\langle f^{(3)}, \uv^{\otimes 3} \bigr\rangle^{2} /36 \bigr\} 
	\exp\Bigl( - \frac{\| \DV_{\GP} \uv \|^{2} - \dltwb \| \DVL \uv \|^{2}}{2} \Bigr) \, d\uv .
\label{u221md3u636v}
\end{EQA}
This and \eqref{Pp012zzmH2u212v} yield
\begin{EQA}
	&& \nquad
	\bigl| \Rem(1) - \Rem(0) \bigr|
	\leq 
	\bigl| \Rem'(0) \bigr| + \frac{1}{2} \sup_{t \in [0,1]} |\Rem''(t)|
	\leq 
	2 \int_{\UVL} |\bar{\dltw}_{4}(\uv)| \, 
	\ex^{- (\| \DV_{\GP} \uv \|^{2} - \dltwb \| \DVL \uv \|^{2})/{2} } \, d\uv
	\\
	&&
	+ \, \frac{1}{36} \int_{\UVL} \bigl\langle f^{(3)}, \uv^{\otimes 3} \bigr\rangle^{2} \, 
		\ex^{- (\| \DV_{\GP} \uv \|^{2} - \dltwb \| \DVL \uv \|^{2})/{2} } \, d\uv .
\label{C4Hm2Fm4tH2tv}
\end{EQA}
Change of variable \( \bigl( \Id - \dltwb \, \DV_{\GP}^{-1} \DVL^{2} \, \DV_{\GP}^{-1} \bigr)^{1/2} \uv \) to \( \wv \) yields by 
\eqref{12Im2t3t4Hm3} in view of \( \dltwb \leq 1/3 \)
\begin{EQA}
	&& \nquad
	\frac{1}{36} \int_{\UVL} \bigl\langle f^{(3)}, \uv^{\otimes 3} \bigr\rangle^{2} \, 
		\exp\Bigl( - \frac{\| \DV_{\GP} \uv \|^{2} - \dltwb \| \DVL \uv \|^{2}}{2} \Bigr) \, d\uv 
	\\
	& \leq &
	\frac{3/2}{36 (1 - \dltwb)^{3}} 
	\int \bigl\langle f^{(3)}, \wv^{\otimes 3} \bigr\rangle^{2} \, \exp\Bigl( - \frac{\| \DV_{\GP} \wv \|^{2}}{2} \Bigr) \, d\wv .
\label{C4Hm2Fm4tH2t3v}
\end{EQA}
Similarly by \nameref{LL4tref}
\begin{EQA}
	\int_{\UVL} |\bar{\dltw}_{4}(\uv)| \,\exp\Bigl( - \frac{\| \DV_{\GP} \uv \|^{2} - \dltwb \| \DVL \uv \|^{2}}{2} \Bigr) \, d\uv
	& \leq &
	\frac{3/2}{24(1 - \dltwb)^{2}} 
	\int  
	\dltwu_{4} \| \DV \wv \|^{4} \, \exp\Bigl( - \frac{\| \DV_{\GP} \wv \|^{2}}{2} \Bigr) \, d\wv .
\label{C4Hm2Fm4tH2t4v}
\end{EQA}
The use of \( \dltwb \leq 1/3 \) implies that
\begin{EQA}
	\frac{\bigl| \Rem(1) - \Rem(0) \bigr|}{\int_{\UVL} \ex^{- \| \DV_{\GP} \uv \|^{2}/2} d\uv}
	& \leq &
	\frac{3/2}{24 (1 - \dltwb)^{2}}
	\Bigl\{  
	\E \bigl\langle f^{(3)} , \gaussv_{\GP}^{\otimes 3} \bigr\rangle^{2}
		+ 2 \dltwu_{4} \E \| \DV \gaussv_{\GP} \|^{4}
	\Bigr\} 
	\leq 
	\err_{4} \, 
\label{dwHw22n3n432}
\end{EQA}
and \eqref{efuiguduUem} follows.
Further, \nameref{LLtS3ref} yields
\( \bigl\langle \nabla^{3} \lgd(\xv) , \uv^{\otimes 3} \bigr\rangle^{2} \leq \dltwu_{3}^{2} \| \DV \uv \|^{6} \).
Now \eqref{errdef3322Hm2} follows from Lemma~\ref{LdltwLa}.
The proof of \eqref{errdef3322Hm2m} is similar.
\end{proof}

\Section{Tail integrals}
\label{StailLa}
In this section we also write \( \xv \) in place of \( \xvs \).
Below we evaluate \( \rho \) from \eqref{rhfiUaceHu22m} which bounds the integral of \( \ex^{\lgd(\xv;\uv)} \) 
over the complement of the local set \( \UVL \) of a special form 
\( \UVL = \bigl\{ \uv \colon \| \DVL \uv \| \leq \amax^{-1} \rrL \bigr\} \) for \( \DVL \) from \nameref{LLf0ref}.
Our results help to understand how the radius \( \rrL \) should be fixed to ensure \( \rho \) sufficiently small.
The main tools for the analysis are deviation probability bounds for Gaussian quadratic forms; see Section~\ref{Sdevboundgen}.

\begin{proposition}
\label{PUVstarLLfx}
Suppose \nameref{LLf0ref}.
Given \( \amax < 1 \) and \( \xx > 0 \), 
let \( \UVL \) and \( \rrL \) be defined by \eqref{UvTDunm12spT}.
Let also \( \dltwb \) from \eqref{om3esuU1H02d3} satisfy \( \dltwb \leq 1 - \amax \).
Then 
\begin{EQA}
\label{fiinIHUniUef0}
	\frac{\int \Ind\bigl( \uv \not\in \UVL \bigr) \, \ex^{\lgd(\xv;\uv)} \, d\uv}
		 {\int \ex^{ - \| \DV_{\GP} \uv \|^{2}/2} \, d \uv}
	& \leq &
	4 \ex^{-\xx - \dimL/2} \, ,
	\\
	\frac{\int \Ind\bigl( \uv \not\in \UVL \bigr) \, \ex^{ - \| \DV_{\GP} \uv \|^{2}/2} \, d\uv}
		 {\int \ex^{ - \| \DV_{\GP} \uv \|^{2}/2} \, d \uv}
	& \leq &
	\ex^{-\xx - \dimL/2} \, .
\label{fiinIHUniUef}
\end{EQA}
\end{proposition}

\begin{proof}
Let \( \uv \not \in \UVL \), i.e. \( \| \DVL \uv \| > \rr  \) with \( \rr = \amax^{-1} \rrL \).
Define \( \uvc = \rr \| \DVL \uv \|^{-1} \uv \) yielding \( \| \DVL \uvc \| = \rr \).
We also write \( \uv = (1 + \tau) \uvc \) for \( \tau > 0 \).
By \eqref{om3esuU1H02d3} 
and \( \nabla^{2} \lgdL(0) = - \DVL^{2} \)
\begin{EQ}[rcl]
	\lgdL(\uvc) - \lgdL(0) - \bigl\langle \nabla \lgdL(0), \uvc \bigr\rangle 
	& \leq &
	- (1 - \dltwb) \| \DVL \uvc \|^{2}/2 ,
	\\
	\bigl\langle \nabla \lgdL(\uvc) - \nabla \lgdL(0), \uv - \uvc \bigr\rangle
	& \leq &
	- (1 - \dltwb) \bigl\langle \DVL^{2} \uvc, \uv - \uvc \bigr\rangle .
\label{2ud1mud3Hv1dfx}
\end{EQ}
Concavity of \( \lgdL(\uv) \) implies for \( \uv = (1 + \tau) \uvc \),
\begin{EQA}
	\lgdL(\uv) 
	& \leq &
	\lgdL(\uvc) + \bigl\langle \nabla \lgdL(\uvc), \uv - \uvc \bigr\rangle 
\label{fudfpudumudwuld}
\end{EQA}
yielding by \eqref{2ud1mud3Hv1dfx} in view of 
\( \bigl\langle \DVL \uvc, \DVL \uv \bigr\rangle = \| \DVL \uvc \| \, \| \DVL \uv \| \)
\begin{EQA}
	&& \nquad
	\lgdL(\uv) - \lgdL(0) - \bigl\langle \nabla \lgdL(0), \uv \bigr\rangle
	=
	\lgdL(\uv) - \lgdL(\uvc) - \bigl\langle \nabla \lgdL(\uvc), \uv - \uvc \bigr\rangle
	\\
	&&
	\qquad
	\qquad
	\qquad
	\qquad
	+ \, \lgdL(\uvc) - \lgdL(0) - \bigl\langle \nabla \lgdL(0), \uvc \bigr\rangle
	+ \bigl\langle \nabla \lgdL(\uvc) - \nabla \lgdL(0), \uv - \uvc \bigr\rangle
	\\
	& \leq &
	(1 - \dltwb) \| \DVL \uvc \|^{2}/2
	- (1 - \dltwb) \bigl\langle \DVL \uvc, \DVL \uv \bigr\rangle 
	\leq 
	- (1 - \dltwb) \| \DVL \uvc \| \, \| \DVL \uv \| /2.
\label{tnfudnfHm1Hm1d}
\end{EQA}
We now use that 
\( \| \DVL \uvc \| = \rr \),  
\( \uvc = \uv/(1 + \tau) \), and thus,
\begin{EQA}
	&& \nquad
	\lgd(\xv + \uv) - \lgd(0) - \bigl\langle \nabla \lgd(\xv), \uv \bigr\rangle
	=
	\lgdL(\uv) - \lgdL(0) - \bigl\langle \nabla \lgdL(0), \uv \bigr\rangle
	- \| \DV_{\GP} \uv \|^{2}/2 + \| \DVL \uv \|^{2}/2 
	\\
	& \leq &
	- (1 - \dltwb) \rr \| \DVL \uv \|/2 
	- \| \DV_{\GP} \uv \|^{2}/2 + \| \DVL \uv \|^{2}/2.
\label{DGa22DTa221mw3rT22}
\end{EQA}
This yields by \( \rrL = \amax \, \rr \leq (1 - \dltwb) \rr \) with \( \Tau = \DVL \DV_{\GP}^{-1} \)
\begin{EQA}
	&& \nquad
	\frac{\int \Ind\bigl( \uv \not\in \UVL \bigr) 
		\exp\bigl\{  \lgd(\xv + \uv) - \lgd(\xv) - \langle \nabla \lgd(\xv), \uv \rangle \bigr\} \, d\uv}
		{\int \exp \bigl( - \| \DV_{\GP} \uv \|^{2}/{2} \bigr) \, d \uv}
	\\
	& \leq &
	\frac{\int \Ind\bigl( \| \DVL \uv \| > \rr \bigr)
		\exp \bigl\{ - {(1 - \dltwb) \rr} \| \DVL \uv \|/2 - \| \DV_{\GP} \uv \|^{2}/2 + \| \DVL \uv \|^{2}/2 \bigr\} \, d\uv}
		{\int \exp \bigl( - \| \DV_{\GP} \uv \|^{2}/{2} \bigr) \, d \uv}
	\\
	& \leq &
	\E \exp \bigl\{ 
		- \rrL \| \Tau \gaussv \| / 2 
		+ \| \Tau \gaussv \|^{2} / 2 
		\bigr\} \Ind\bigl( \| \Tau \gaussv \| > \rrL \bigr) 
\label{2C2emx1212IganiU}
\end{EQA}
with \( \gaussv \) standard normal in \( \R^{\dimp} \).
Next, define
\begin{EQA}
	\riskt_{0}(\rr)
	& \eqdef &
	\E \exp\bigl( - \rr \| \DVL \gaussv \| / 2 + \| \DVL \gaussv \|^{2}/2 \bigr) 
		\Ind(\| \DVL \gaussv \| > \rr) .
\label{mp2mxr2g22T}
\end{EQA}
Integration by parts allows to represent the last integral as
\begin{EQA}
	\riskt_{0}(\rr)
	&=&
	- \int_{\rr}^{\infty} \exp\Bigl( - { \rr \, \zq}/{2} + {\zq^{2}}/{2} \Bigr) \, 
		d\P\bigl( \| \Tau \gaussv \| > \zq \bigr)
	\\
	&=&
	\P\bigl( \| \Tau \gaussv \| > \rr \bigr) 
	+ \int_{\rr}^{\infty} (\zq - \rr/2) \exp\Bigl( - {\rr \zq}/{2} + {\zq^{2}}/{2} \Bigr) \, 
		\P\bigl( \| \Tau \gaussv \| > \zq \bigr) \, d\zq \, .
\label{dzzTPzr2rT}
\end{EQA}
By Theorem~\ref{TexpbLGA}, for any \( \zq \geq \sqrt{\dimL} \) for \( \dimL = \tr(\Tau \, \Tau^{\T}) = \tr (\DVL^{2} \, \DV_{\GP}^{-2}) \)
\begin{EQA}
	\P\bigl( \| \Tau \gaussv \| > \zq \bigr)
	& \leq & 
	\exp\bigl\{ - (\zq - \sqrt{\dimL})^{2}/2 \bigr\}
\label{2emxPTgasps2d}
\end{EQA}
yielding for \( \zq \geq \rrL = 2 \sqrt{\dimL} + \sqrt{2\xx} \)
\begin{EQA}
	\P\bigl( \| \Tau \gaussv \| > \zq \bigr)
	& \leq &
	\exp\bigl\{ - (\zq - \sqrt{\dimL})^{2}/2 \bigr\}
	\leq 
	\ex^{-\xx - \dimL/2}
\label{2exmp22ezdp2}
\end{EQA}
and for \( \rr \geq 2 \sqrt{\dimL} + \sqrt{2\xx} \) and \( \xx \geq 2 \)
\begin{EQA}
	\riskt_{0}(\rr)
	& \leq &
	\ex^{-\xx - \dimL/2}
	+ \int_{\rr}^{\infty} (\zq - \rr/2) 
	\exp\Bigl\{ - \frac{\rr \zq}{2} + \frac{\zq^{2}}{2} - \frac{(\zq - \sqrt{\dimL})^{2}}{2} \Bigr\} \, d\zq
	\\
	& \leq &
	\ex^{-\xx - \dimL/2}
	+ \exp\Bigl( - \frac{(\rr - \sqrt{\dimL})^{2}}{2} \Bigr) \int_{0}^{\infty} \Bigl(\zq + \frac{\rr}{2} \Bigr) 
	\exp\Bigl\{ - \frac{(\rr - 2 \sqrt{\dimL}) \zq}{2} \Bigr\} \, d\zq
	\\
	& \leq &
	2 \ex^{-\xx - \dimL/2} .
\label{Cexmxmp222Tr22fx}
\end{EQA}
This completes the proof of the result \eqref{fiinIHUniUef0}.
Statement \eqref{fiinIHUniUef} is about Gaussian probability 
\( \P\bigl( \| \Tau \gaussv \| \geq \rrL \bigr) \) for a standard normal element \( \gaussv \), and we derive
\begin{EQA}
	\P\bigl( \| \Tau \gaussv \| \geq 2 \sqrt{\dimL} + \sqrt{2\xx} \bigr)
	& \leq &
	\exp\bigl\{ - \bigl( \sqrt{\dimL} + \sqrt{2\xx} \bigr)^{2}/2 \bigr\}
	\leq 
	\exp\bigl( -\xx - \dimL/2 \bigr) 
\label{fiinIHUniUefi}
\end{EQA}
and \eqref{fiinIHUniUef} follows.
\end{proof}

\medskip

The next result extends \eqref{fiinIHUniUef0}.
\begin{proposition}
\label{PUVstarLLgenlifx}
Assume the conditions of Proposition~\ref{PUVstarLLfx}
with
\begin{EQA}
	\rrL
	& \geq &
	2 \sqrt{\dimL} + \sqrt{2 \xx} + m
\label{1w3rT2spTs2xlifx}
\end{EQA}
for some \( m \geq 0 \).
Then 
\eqref{fiinIHUniUef0} can be extended to
\begin{EQA}
\label{fiinIHUniUef0m}
	\frac{\int \Ind\bigl( \uv \not\in \UVL \bigr) \, \| \DVL \uv \|^{m} \, \ex^{\lgd(\xv;\uv)} \, d\uv}
		 {\int \ex^{ - \| \DV_{\GP} \uv \|^{2}/2} \, d \uv}
	& \leq &
	4 \ex^{-\xx - \dimL/2} \, ,
	\\
	\frac{\int \Ind\bigl( \uv \not\in \UVL \bigr) \, \| \DVL \uv \|^{m} \, \ex^{ - \| \DV_{\GP} \uv \|^{2}/2} \, d\uv}
		 {\int \ex^{ - \| \DV_{\GP} \uv \|^{2}/2} \, d \uv}
	& \leq &
	\ex^{-\xx - \dimL/2} \, .
\label{fiinIHUniUefm}
\end{EQA}
\end{proposition}

\begin{proof}
The case \( m > 0 \) can be proved similarly to \( m=0 \) using 
\( m \log z \leq m z \).
\end{proof}

\Section{Local concentration}
Here we show that the measure \( \PfL \) well concentrates on the local set \( \UVL \) from \eqref{UvTDunm12spT}.
Again we fix \( \xv = \xvs \).

\begin{proposition}
\label{PlocconLa}
Assume \( \dltwb \leq 1/3 \).
Then
\begin{EQA}
	\int_{\UVL} \ex^{\lgd(\xv;\uv)} \, d\uv  
	& \geq &
	\ex^{ - \dltwb \, \dimL/2} \, \int_{\UVL} \ex^{- \| \DV_{\GP} \uv \|^{2}/2} d\uv \, .
\label{emdw3dp0iUfxu}
\end{EQA}
Moreover,
\begin{EQA}
	\frac{\int_{\UVL^{c}} \ex^{\lgd(\xv;\uv)} \, d\uv}{\int \ex^{\lgd(\xv;\uv)} \, d\uv}  
	& \leq &
	4 \ex^{ - \xx - (1 - \dltwb) \, \dimL/2} 
	\leq 
	\ex^{-\xx} \, .
\label{emdw3dp0iUfxuL}
\end{EQA}
\end{proposition}

\begin{proof}
By \eqref{om3esuU1H02d3}
\begin{EQA}
	\int_{\UVL} \ex^{\lgd(\xv;\uv)} \, d\uv
	& = &
	\int_{\UVL} \ex^{- \| \DV_{\GP} \uv \|^{2}/2 + \dltw_{3}(\xv,\uv)} \, d\uv 
	\geq 
	\int_{\UVL} \ex^{- \| \DV_{\GP} \uv \|^{2}/2 - \dltwb \| \DVL \uv \|^{2}/2} \, d\uv \, .
\label{PhtiUexmHu222}
\end{EQA}
Change of variable \( \bigl( \Id + \dltwb \, \DV_{\GP}^{-1} \DVL^{2} \, \DV_{\GP}^{-1} \bigr)^{1/2} \uv \) to \( \wv \) yields by \eqref{12Ip2t3t4Hp3} 
\begin{EQA}
	\int_{\UVL} \ex^{\lgd(\xv;\uv)} \, d\uv
	& \geq &
	\det \bigl( \Id + \dltwb \, \DV_{\GP}^{-1} \DVL^{2} \, \DV_{\GP}^{-1} \bigr)^{-1/2} 
	\int_{\UVL} \ex^{- \| \DV_{\GP} \wv \|^{2}/2 } \, d\wv
	\\
	& \geq &
	\ex^{- \dltwb \, \dimL / 2} \int_{\UVL} \ex^{- \| \DV_{\GP} \wv \|^{2}/2 } \, d\wv ,
\label{iUefxudetId3H121}
\end{EQA}
and \eqref{emdw3dp0iUfxu} follows.
This and \eqref{fiinIHUniUef0}, \eqref{fiinIHUniUef} of Proposition~\ref{PUVstarLLfx} imply
\begin{EQA}
	&& \nquad
	\frac{\int_{\UVL^{c}} \ex^{\lgd(\xv;\uv)} \, d\uv}{\int \ex^{\lgd(\xv;\uv)} \, d\uv}
	= 
	\frac{\int_{\UVL^{c}} \ex^{\lgd(\xv;\uv)} \, d\uv}{\int_{\UVL} \ex^{\lgd(\xv;\uv)} \, d\uv + \int_{\UVL^{c}} \ex^{\lgd(\xv;\uv)} \, d\uv}
	\\
	& \leq &
	\frac{4 \ex^{-\xx - \dimL/2} \int \ex^{- \| \DV_{\GP} \uv \|^{2}/2} d\uv}
		 {\ex^{ - \dltwb \, \dimL/2} \, \int_{\UVL} \ex^{- \| \DV_{\GP} \uv \|^{2}/2} d\uv 
			+ 4 \ex^{-\xx - \dimL/2} \int \ex^{- \| \DV_{\GP} \uv \|^{2}/2} d\uv}
	\leq 
	4 \ex^{-\xx - (1 - \dltwb) \dimL/2} 
\label{jyfftc434rdssdsaaa}
\end{EQA}
as required in \eqref{emdw3dp0iUfxuL}.
\end{proof}

\Section{Finalizing the proof of Theorem~\ref{TLaplaceTV} and \ref{TLaplaceTV2}}
These results are proved by compiling the previous technical statements.
Proposition~\ref{PUVstarLLfx} provides some upper bounds for the quantities \( \rho \) and \( \rho_{\GP} \),
while Proposition~\ref{Lintfxupp3}, Proposition~\ref{Lintfxupp3T}, and Proposition~\ref{Lintfxupp2} bound the local errors 
\( \err \) and \( \err_{g} \).
The final bound \eqref{ufgdt6df5dtgededsxd23gjg} follows from Corollary~\ref{CPunbintLapl}.

\ifNL{
\Section{Proof of Theorem~\ref{TLaplaceKL} and Theorem~\ref{TLaplaceKLi}}
Define
\begin{EQA}
	\CDG 
	& \eqdef &
	\log \int \ex^{- \| \DV_{\GP} \uv \|^{2}/2} \, d\uv 
	- \log \int \ex^{\lgd(\xvs;\uv)} \, d\uv 
	\, . 
\label{sihedyfw3ytw3ydqwdbhwd}
\end{EQA}
For \( \uv = \xv - \xvs \), it holds as in \eqref{IIgfifxutudufxutdu}
\begin{EQA}
	\PfL
	& \sim &
	\frac{\ex^{\lgd(\xv)}}{\int \ex^{\lgd(\xv)} \, d\uv}
	=
	\frac{\ex^{\lgd(\xv) - \lgd(\xvs)}}{\int \ex^{\lgd(\xv) - \lgd(\xvs)} \, d\uv}
	=
	\frac{\ex^{\lgd(\xvs;\uv)}}{\int \ex^{\lgd(\xvs;\uv)} \, d\uv} 
	=
	\frac{\ex^{\lgd(\xvs;\uv) + \CDG}}{\int \ex^{- \| \DV_{\GP} \uv \|^{2}/2} \, d\uv} \, .
\label{dfhuwffwt5wtcsbydyqyqh}
\end{EQA}
Further, with \( \P_{\GP} = \ND(\xvs,\DV_{\GP}^{-2}) \)
\begin{EQA}
	\log \frac{d\PfL}{d\P_{\GP}}(\xv) 
	&=&
	\lgd(\xvs;\uv) - \| \DV_{\GP} \uv \|^{2}/2 + \CDG
	=
	\dltw_{3}(\uv) + \CDG 
\label{sufuwdeuquqeygwyg22e24}
\end{EQA}
and
\begin{EQA}
	\kullb(\PfL,\P_{\GP})
	&=&
	\ex^{\CDG} \, \frac{\int \dltw_{3}(\uv) \, \ex^{\lgd(\xvs;\uv)} \, d\uv}
		 {\int \ex^{- \| \DV_{\GP} \uv \|^{2}/2} \, d\uv}
	+ \CDG \, .
\label{dvjnue37e37e7weedhe}
\end{EQA}
Similarly to \eqref{Pht12u212H02u22} and \eqref{Pd3p2d3H0wem22g1}, we can bound
\begin{EQA}
	\left| \int_{\UVL} \dltw_{3}(\uv) \, \ex^{\lgd(\xvs;\uv)} \, d\uv \right|
	& \leq &
	\int_{\UVL} |\dltw_{3}(\uv)| \, \ex^{- \| \DV_{\GP} \uv \|^{2}/2 - \dltwb \| \DVL \uv \|^{2}/2} \, d\uv 
	\\
	& \leq &
	\frac{1}{4(1 - \dltwb)^{3/2}}
	\int \dltwu_{3}(\uv) \, \ex^{- \| \DV_{\GP} \uv \|^{2} /2} \, d\uv \, .
\label{ucuwudhuwhi2wefhw78w}
\end{EQA}
On the complement \( \uv \in \UVL^{c} \), we use that
\begin{EQA}
	\dltw_{3}(\uv)
	=
	\lgd(\xvs;\uv) + \frac{1}{2} \| \DV_{\GP} \uv \|^{2}
	&=&
	\lgdL(\xvs;\uv) + \frac{1}{2} \| \DVL \uv \|^{2}
	\leq 
	\frac{1}{2} \| \DVL \uv \|^{2} .
\label{yuqt2sstwwtwbstdywhsyw}
\end{EQA}
The last inequality is based on concavity of \( \lgdL(\cdot) \) and local approximation 
\( \lgdL(\xvs;\uv) \approx - \| \DVL \uv \|^{2}/2 \) within \( \UVL \) yielding \( \lgdL(\xvs;\uv) < 0 \) for \( \uv \in \UVL^{c} \).
By Proposition~\ref{PUVstarLLgenlifx} with \( m = 2 \),
\begin{EQA}
	\int_{\UVL^{c}} \dltw_{3}(\uv) \, \ex^{\lgd(\xvs;\uv)} \, d\uv
	& \leq &
	\frac{1}{2} \int_{\UVL^{c}} \| \DVL \uv \|^{2} \, \ex^{\lgd(\xvs;\uv)} \, d\uv
	\leq 
	\ex^{-\xx} \int \ex^{- \| \DV_{\GP} \uv \|^{2}/2} \, d\uv \, .
\label{sywdywydwbwtwgqb1qt213ba}
\end{EQA}
We conclude that
\begin{EQA}
	\frac{\int \dltw_{3}(\uv) \, \ex^{\lgd(\xvs;\uv)} \, d\uv}
		 {\int \ex^{- \| \DV_{\GP} \uv \|^{2}/2} \, d\uv}
	& \leq &
	\frac{\E \dltwu_{3}(\gaussv_{\GP})}{4(1 - \dltwb)^{3/2}} + \ex^{- \xx} 
	\leq 
	\err_{3} + \ex^{-\xx} .
\label{0fmdfeyjdwhqkshyxiwck}
\end{EQA}
Similarly
\begin{EQA}
	&& \nquad
	\bigl| \ex^{\CDG} - 1 \bigr|
	=
	\left| \frac{\int \ex^{\lgd(\xvs;\uv)} \, d\uv - \int \ex^{- \| \DV_{\GP} \uv \|^{2}/2} \, d\uv}
		 {\int \ex^{- \| \DV_{\GP} \uv \|^{2}/2} \, d\uv} 
	\right|
	\\
	& \leq &
	\left| \frac{\int_{\UVL} \ex^{\lgd(\xvs;\uv)} \, d\uv - \int_{\UVL} \ex^{- \| \DV_{\GP} \uv \|^{2}/2} \, d\uv}
		 {\int \ex^{- \| \DV_{\GP} \uv \|^{2}/2} \, d\uv} 
	\right|
	+ \frac{\int_{\UVL^{c}} \ex^{\lgd(\xvs;\uv)} \, d\uv}{\int \ex^{- \| \DV_{\GP} \uv \|^{2}/2} \, d\uv}
	+ \frac{\int_{\UVL^{c}} \ex^{- \| \DV_{\GP} \uv \|^{2}/2} \, d\uv}
		 {\int \ex^{- \| \DV_{\GP} \uv \|^{2}/2} \, d\uv}
	\\
	& \leq &
	\frac{\E \dltwu_{3}(\gaussv_{\GP})}{4(1 - \dltwb)^{3/2}} + 2 \ex^{- \xx} 
	\leq 
	\err_{3} + 2 \ex^{-\xx} .
\label{hsyc4w2tde6dyfdgyww23}
\end{EQA}
This yields \( \ex^{\CDG} \leq 1 + \err_{3} + 2 \ex^{-\xx} \) and \( \CDG \leq \err_{3} + 2 \ex^{-\xx} \).
Putting all together results in
\begin{EQA}
	\kullb(\PfL,\P_{\GP})
	& \leq &
	\bigl( 1 + \err_{3} + 2 \ex^{-\xx} \bigr) \bigl( \err_{3} + \ex^{-\xx} \bigr) + \err_{3} + 2 \ex^{-\xx}
	<  
	4 \err_{3} + 4 \ex^{-\xx} 
\label{gcusgu255eetfghbhjjkkook}
\end{EQA}
provided that \( 4 \err_{3} + 4 \ex^{-\xx} \leq 1 \),
and \eqref{vlgvi8ugu7tr4ry43et31} follows.

The proof of Theorem~\ref{TLaplaceKLi} is similar and even simpler except one special part, namely, 
the bound for the tail integral of \( - \dltw_{3}(\uv) \).
By definition \( \| \DVL \uv \| \geq \rr_{\GP} \) for \( \uv \in \UVL^{c} \).
This implies by \eqref{dchbhwdhwwdgscsn2efty2162} similarly to \eqref{yuqt2sstwwtwbstdywhsyw} 
\begin{EQA}
	&& \nquad
	- \int_{\UVL^{c}} \dltw_{3}(\uv) \, \ex^{- \| \DV_{\GP} \uv \|^{2}/2} \, d\uv
	\leq
	\int_{\UVL^{c}} \bigl| \lgdL(\xvs;\uv) \bigr| \, \ex^{- \| \DV_{\GP} \uv \|^{2}/2 + \rho_{\GP} \| \DVL \uv \|^{2}/2} \,
		\ex^{- \rho_{\GP} \| \DVL \uv \|^{2}/2} \, d\uv
	\\
	& \leq &
	\ex^{- \rho_{\GP} \rr_{\GP}^{2}/2} \, 
	\int \bigl| \lgdL(\xvs;\uv) \bigr| \, \ex^{- \| \DV_{\GP} \uv \|^{2}/2 + \rho_{\GP} \| \DVL \uv \|^{2}/2} \,
		d\uv
	\leq \CONSTi_{\lgdL} \, \ex^{-\xx} ,
\label{sguswdw2e2j2wuhq2wdw}
\end{EQA}
and the result follows.
}{}

\Section{Finalizing the proof of Theorem~\ref{TpostmeanLa} and Corollary~\ref{CTpostmeanLa}}
For Theorem~\ref{TpostmeanLa}, we follow the same line as for Theorem~\ref{TLaplaceTV34}.
Note first that
\begin{EQA}
	\QP (\xvb - \xvs)
	=
	\frac{\int \QP(\xvs + \uv) \, \ex^{\lgd(\xvs + \uv)} \, d\uv}{\int \ex^{\lgd(\xvs + \uv)} \, d\uv} - \QP \xvs
	&=&
	\frac{\int \QP \uv \, \ex^{\lgd(\xvs;\uv)} \, d\uv}{\int \ex^{\lgd(\xvs;\uv)} \, d\uv} \, 
\label{dvhjt6efedfchsijcfte4ws}
\end{EQA}
and
\begin{EQA}
	\| \QP (\xvb - \xvs) \|
	&=&
	\sup_{\av \in \R^{\dimq} \colon \| \av \| = 1} \bigl| \langle \QP (\xvb - \xvs), \av \rangle \bigr|
	= 
	\sup_{\av \in \R^{\dimq} \colon \| \av \| = 1} 
		\left| \frac{\int \langle \QP \uv, \av \rangle \, \ex^{\lgd(\xvs;\uv)} \, d\uv}{\int \ex^{\lgd(\xvs;\uv)} \, d\uv} \right| \, .
\label{p3hb893dfvfdsiuw5whlh}
\end{EQA}
Now fix \( \av \in \R^{\dimq} \) with \( \| \av \| = 1 \) and \( g(\uv) = \langle \QP \uv, \av \rangle \).
\eqref{igefxumiguexHu22g} implies
\begin{EQA}[rcl]
	\left| \frac{\int g(\uv) \, \ex^{\lgd(\xvs;\uv)} \, d\uv}{\int \ex^{\lgd(\xvs;\uv)} \, d\uv} \right| 
	& \leq &
	\frac{\rho_{g} + \rho_{\GP,g} + \err_{3,g}}{1 - \rho_{\GP} - \err_{3}} \, .
\label{igefxumiguexHu22gp}
\end{EQA}
The bound \( 1 - \err_{3} - \rho_{\GP} \geq 1/2 \) has been already checked.
Proposition~\ref{PUVstarLLgenlifx} for \( m=1 \) helps to bound the values \( \rho_{g} \) and \( \rho_{\GP,g} \) 
by \( \CONST \ex^{-\xx} \).
Next we bound \( \err_{3,g} \).
Under \nameref{LLtS3ref},
\eqref{iUafxvdvH2vk3} of Proposition~\ref{Perrbound3} combined with Lemma~\ref{LdltwLa} and Lemma~\ref{Gaussmoments} yield
\begin{EQA}
	&& \nquad
	4 \err_{3,g}
	=
	\frac{1}{(1 - \dltwb)^{2}} \,
	\E \bigl\{  | \langle \QP \gaussv_{\GP},\av \rangle | \,\, \dltwu_{3}(\gaussv_{\GP}) \bigr\}
	= 
	\frac{\hmax_{3} \, n^{-1/2}}{(1 - \dltwb)^{2}} \, 
	\E \bigl\{ | \langle \QP \, \gaussv_{\GP},\av \rangle | \,\, \| \DVL \, \gaussv_{\GP} \|^{3} \, \bigr\} 
	\\
	& \leq &
	\frac{\hmax_{3} \, n^{-1/2}}{(1 - \dltwb)^{2}} \, 
	\E^{3/4} \| \DVL \, \gaussv_{\GP} \|^{4} \,\, \E^{1/4} \langle \QP \gaussv_{\GP},\av \rangle^{4} 
	\leq 
	\frac{3^{1/4} \, \hmax_{3} \, (\dimL + 1)^{3/2} }{(1 - \dltwb)^{2} \sqrt{n}} \, 
	\sqrt{ \av^{\T} \QP \, \DV_{\GP}^{-2} \QP^{\T} \av } .
\label{sfdhysdfyf11ewds3wsded}
\end{EQA}
Here we used that 
\( \E \langle \QP \gaussv_{\GP},\av \rangle^{4} = \E \langle \gaussv, \DPGP^{-1} \QP^{\T}\av \rangle^{4}
= 3 (\av^{\T} \QP \, \DV_{\GP}^{-2} \QP^{\T} \av)^{2} \).
Now \eqref{hcdtrdtdehfdewdrfrhgyjufger} follows from \( 3^{1/4} (1 - \dltwb)^{-2} \leq 2.4 \) and
\begin{EQA}
	\sup_{\av \in \R^{\dimq} \colon \| \av \| = 1} \av^{\T} \QP \, \DV_{\GP}^{-2} \QP^{\T} \av
	&=&
	\| \QP \, \DV_{\GP}^{-2} \QP^{\T} \| \, .
\label{hdrtesw5sdghww3r4d5tdy6gh}
\end{EQA}
With \( \QP = \DVL \), this implies \eqref{klu8gitfdgregfkhj7yt}.

}

\newpage
\appendix


\Chapter{Local smoothness conditions}
\label{Slocalsmooth}
This section discusses different local smoothness characteristics of a 
multivariate function \( f(\upsv) = \E L(\upsv) \), \( \upsv \in \R^{\dimp} \).

\Section{Smoothness in Gateaux sense and self-concordance}
Below we assume that
the function \( f(\upsv) \) is three or sometimes even four times Gateaux differentiable in \( \upsv \in \Ups \).
For any particular direction \( \uv \in \R^{\dimp} \), we consider the univariate function 
\( f(\upsv + t \uv) \) and measure its smoothness in \( t \).
Local smoothness of \( f \) will be described by the relative error of the Taylor expansion 
of the third or four order.
Namely, define
\begin{EQ}[rcl]
	\dltw_{3}(\upsv,\uv) 
	&=& 
	f(\upsv + \uv) - f(\upsv) - \langle \nabla f(\upsv), \uv \rangle 
	- \frac{1}{2} \langle \nabla^{2} f(\upsv), \uv^{\otimes 2} \rangle , 
	\\
	\dltwd_{3}(\upsv,\uv) 
	&=&
	\langle \nabla f(\upsv + \uv), \uv \rangle - \langle \nabla f(\upsv), \uv \rangle 
	- \langle \nabla^{2} f(\upsv), \uv^{\otimes 2} \rangle \, ,
\label{dltw3vufuv12f2ga}
\end{EQ}
and
\begin{EQA}
	\dltw_{4}(\upsv,\uv)
	& \eqdef &
	f(\upsv + \uv) - f(\upsv) - \langle \nabla f(\upsv), \uv \rangle 
	- \frac{1}{2} \langle \nabla^{2} f(\upsv), \uv^{\otimes 2} \rangle
	- \frac{1}{6} \langle \nabla^{3} f(\upsv), \uv^{\otimes 3} \rangle .
	\qquad
	\qquad
\label{hvcduywgedfuyg2y1y35e3wweg}
\end{EQA}
Now, for each \( \upsv \), suppose to be given a positive symmetric operator \( \HLG(\upsv) \in \Matr_{\dimp} \) 
defining a local metric and a local vicinity around \( \upsv \):
\begin{EQA}
	\UVz(\upsv)
	&=&
	\bigl\{ \uv \in \R^{\dimp} \colon \| \HLG(\upsv) \uv \| \leq \rr \bigr\}
\label{ed7sycf7wedwgedq2ftwdfgtv}
\end{EQA}
for some radius \( \rr \).
If \( f \) is strongly concave, then the Hessian \( \IF(\upsv) \eqdef - \nabla^{2} f(\upsv) \in \Matr_{\dimp} \) 
is positive definite 
and one can use \( \HLG(\upsv) = \IF^{1/2}(\upsv) \).

Local smoothness properties of \( f \) are given via the quantities
\begin{EQA}[rcccl]
    \dltwb(\upsv)
    & \eqdef &
    \sup_{\uv \colon \| \HLG(\upsv) \uv \| \leq \rr} \,
    \frac{2 |\dltw_{3}(\upsv,\uv)|}{\| \HLG(\upsv) \uv \|^{2}} 
    \,\, ,
    \qquad
    \dltwbd(\upsv)
    & \eqdef &
    \sup_{\uv \colon \| \HLG(\upsv) \uv \| \leq \rr} \, \frac{2 |\dltwd_{3}(\upsv,\uv)|}{\| \HLG(\upsv) \uv \|^{2}} \,\, . 
    \qquad
\label{dtb3u1DG2d3GPg}
\end{EQA}
The Taylor expansion yields for any \( \uv \) with \( \| \HLG(\upsv) \uv \| \leq \rr \)
\begin{EQ}[rcccl]
	\bigl| \dltw_{3}(\upsv,\uv) \rangle \bigr|
	& \leq &
	\frac{\dltwb(\upsv)}{2} \| \HLG(\upsv) \uv \|^{2} 
	\, ,
	\qquad
	\bigl| \dltwd_{3}(\upsv,\uv) \bigr|
	& \leq &
	\frac{\dltwbd(\upsv)}{2} \| \HLG(\upsv) \uv \|^{2}
	\, .
	\qquad
\label{dta3u1DG2d3GPa1g}
\end{EQ}
%
The introduced quantities \( \dltwb(\upsv) \), \( \dltwbd(\upsv) \) 
strongly depend on the radius \( \rr \) of the local vicinity \( \UVz(\upsv) \).
The results about Laplace approximation can be improved 
provided a homogeneous upper bound on the error of Taylor expansion. 
Assume a subset \( \Upsd \) of \( \Ups \) to be fixed.

\begin{description}
    \item[\label{LL3tref} \( \bb{(\mathcal{T}_{3})} \)]
      \textit{There exists \( \dltwu_{3} \) such that for all \( \upsv \in \Upsd \)}
\begin{EQA}
	\bigl| \dltw_{3}(\upsv,\uv) \bigr|
	& \leq &
	\frac{\dltwu_{3}}{6} \| \HLG(\upsv) \, \uv \|^{3} \, ,
	\qquad
	\bigl| \dltwd_{3}(\upsv,\uv) \bigr|
	\leq 
	\frac{\dltwu_{3}}{2} \| \HLG(\upsv) \, \uv \|^{3} \, ,
	\qquad
	\uv \in \UVz(\upsv).
\label{bd3xu16f3uo3st}
\end{EQA}
\end{description}
 
\begin{description}
    \item[\label{LL4tref} \( \bb{(\mathcal{T}_{4})} \)]
      \textit{ There exists \( \dltwu_{4} \) such that for all \( \upsv \in \Upsd \)}
\begin{EQA}
	\bigl| \dltw_{4}(\upsv,\uv) \bigr|
	& \leq &
	\frac{\dltwu_{4}}{24} \| \HLG(\upsv) \, \uv \|^{4} \, ,
	\qquad
	\uv \in \UVz(\upsv).
\label{1mffmxum5st}
\end{EQA}
\end{description}
%

\begin{lemma}
\label{LdltwLa3t}
Under \nameref{LL3tref},
the values \( \dltwb(\upsv) \) and \( \dltwbd(\upsv) \) from \eqref{dtb3u1DG2d3GPg} satisfy
\begin{EQA}[rcccl]
\label{gtcdsftdffrvsewsea}
	\dltwb(\upsv)
	& \leq &
	\frac{\dltwu_{3} \, \rr}{3 } \, ,
	\qquad
	\dltwbd(\upsv)
	& \leq &
	{\dltwu_{3} \, \rr} \, ,
	\qquad
	\upsv \in \Upsd .
\label{gtcdsftdfvtwdsefhfdvfrvsewseG}
\end{EQA}
\end{lemma}

\begin{proof}
For any \( \uv \in \UVz(\upsv) \) with \( \| \HLG(\upsv) \uv \| \leq \rr \)
\begin{EQA}
	\bigl| \dltw_{3}(\upsv,\uv) \bigr|
	& \leq &
	\frac{\dltwu_{3}}{6} \, \| \HLG(\upsv) \uv \|^{3} 
	\leq 
	\frac{\dltwu_{3} \, \rr}{6} \, \| \HLG(\upsv) \uv \|^{2},
\label{jrgeteteer2234587654}
\end{EQA}
and the bound for \( \dltwb(\upsv) \) follows.
The proof for \( \dltwbd(\upsv) \) is similar.
\end{proof}

The values \( \dltwu_{3} \) and \( \dltwu_{4} \) are usually very small.
Some quantitative bounds are given later in this section
under the assumption that the function \( f(\upsv) = \E L(\upsv) \) can be written in the form \( - f(\upsv) = n \hL(\upsv) \) 
for a fixed smooth function \( h(\upsv) \) with the Hessian \( \nabla^{2} \hL(\upsv) \). 
The factor \( n \) has meaning of the sample size; see \Chname \ref{ScritdimMLE} or \Chname \ref{SGBvM}.

\begin{description}
    \item[\label{LLtS3ref} \( \bb{(\mathcal{S}_{3})} \)]
      \emph{ \( - f(\upsv) = n \hL(\upsv) \) with \( \hL(\upsv) \) satisfying 
      \( |\nabla^{2} \hL(\upsv)| \leq \HL^{2}(\upsv) \) 
      and
\begin{EQA}
	\sup_{\uv \colon \| \HL(\upsv) \uv \| \leq \rr/\sqrt{n}} 
	\frac{\bigl| \langle \nabla^{3} \hL(\upsv + \uv), \uv^{\otimes 3} \rangle \bigr|}{\| \HL(\upsv) \uv \|^{3}}
	& \leq &
	\hmax_{3} \, .
\end{EQA}
}
    \item[\label{LLtS4ref} \( \bb{(\mathcal{S}_{4})} \)]
      \emph{ the function \( \hL(\cdot) \) satisfies \nameref{LLtS3ref} and  
\begin{EQA}
	\sup_{\uv \colon \| \HL(\upsv) \uv \| \leq \rr/\sqrt{n}}
	\frac{\bigl| \langle \nabla^{4} \hL(\upsv + \uv), \uv^{\otimes 4} \rangle \bigr|}{\| \HL(\upsv) \uv \|^{4}}
	& \leq &
	\hmax_{4} \, .
\end{EQA}
}
\end{description}

\noindent
\nameref{LLtS3ref} and \nameref{LLtS4ref}
are local versions of the so called self-concordance condition; see \cite{Ne1988}.
In fact, they require that each univariate function \( \hL(\upsv + t \uv) \) of \( t \in \R \)
is self-concordant with some universal constants \( \hmax_{3} \) and \( \hmax_{4} \).
Under \nameref{LLtS3ref} and \nameref{LLtS4ref}, we can use \( \HLG^{2}(\upsv) = n \, \HL^{2}(\upsv) \) 
and easily bound the values 
\( \dltw_{3}(\upsv,\uv) \), \( \dltw_{4}(\upsv,\uv) \), and \( \dltwb(\upsv) \), \( \dltwbd(\upsv) \).

\begin{lemma}
\label{LdltwLaGP}
Suppose \nameref{LLtS3ref}.
Set \( \HLG^{2}(\upsv) = n \, \HL^{2}(\upsv) \).
Then 
\nameref{LL3tref} follows with \( \dltwu_{3} = \hmax_{3} n^{-1/2} \).
Moreover, for \( \dltwb(\upsv) \) and \( \dltwbd(\upsv) \) from \eqref{dtb3u1DG2d3GPg}, it holds
\begin{EQA}[rcccl]
	\dltwb(\upsv)
	& \leq &
	\frac{\hmax_{3} \, \rr}{3 n^{1/2}} \, ,
	\qquad
	\dltwbd(\upsv)
	& \leq &
	\frac{\hmax_{3} \, \rr}{n^{1/2}} \, .
\label{gtcdsftdfvtwdsefhfdvfrvsewseGP}
\end{EQA}
Also \nameref{LL4tref} follows from \nameref{LLtS4ref} with \( \dltwu_{4} = \hmax_{4} n^{-1} \).
\end{lemma}

\begin{proof}
For any \( \uv \in \UVz(\upsv) \) and \( t \in [0,1] \), by the Taylor expansion of the third order
\begin{EQA}
	|\dltw(\upsv,\uv)|
	& \leq &
	\frac{1}{6} \bigl| \langle \nabla^{3} f(\upsv + t \uv), \uv^{\otimes 3} \rangle \bigr|
	=
	\frac{n}{6} \, \bigl| \langle \nabla^{3} \hL(\upsv + t \uv), \uv^{\otimes 3} \rangle \bigr|
	\leq 
	\frac{n \, \hmax_{3}}{6} \, \| \HL(\upsv) \uv \|^{3} 
	\\
	&=&
	\frac{n^{-1/2} \, \hmax_{3}}{6} \, \| \HLG(\upsv) \uv \|^{3}
	\leq 
	\frac{n^{-1/2} \, \hmax_{3} \, \rr}{6} \, \| \HLG(\upsv) \uv \|^{2} \, .
\label{jrgeteteer2234587654}
\end{EQA}
This implies \nameref{LL3tref} as well as \eqref{gtcdsftdfvtwdsefhfdvfrvsewseGP}; see \eqref{dta3u1DG2d3GPa1g}.
The statement about \nameref{LL4tref} is similar.
\end{proof}

\Section{Fr\'echet derivatives and smoothness of the Hessian}
\label{SlocsmooFr}
For evaluation of the bias, 
we also need stronger smoothness conditions in the Fr\'echet sense.
Let \( f \) be a strongly concave function. 
Essentially we need some continuity of the negative Hessian \( \IF(\upsv) = - \nabla^{2} f(\upsv) \).
For \( \upsv \in \Ups \), 
define \( \HLG(\upsv) = \IF^{1/2}(\upsv) \) and
\begin{EQA}
	\dltwbss(\upsv)
    & \eqdef & 
    \sup_{\uv \colon \| \HLG(\upsv) \uv \| \leq \rr} \,\, \sup_{\gammav \in \R^{\dimp}} \,\, 
    \frac{|\langle \IF(\upsv + \uv) - \IF(\upsv), \gammav^{\otimes 2} \rangle|}{\| \HLG(\upsv) \gammav \|^{2} } \, .
\label{jcxhydtferyufgy7tfdsy7ft}
\end{EQA}
This definition of \( \dltwbss(\upsv) \) is, of course, stronger than
the one-directional definition of \( \dltwb(\upsv) \) in \eqref{dtb3u1DG2d3GPg}.
However, in typical examples these quantities \( \dltwb(\upsv) \) and \( \dltwbss(\upsv) \) are similar; 
see the examples from 
\ifapp{\Chname \ref{Saniligot} and \Chname \ref{SGBvM}.}
{\Chname \ref{SGBvM}.}

We also present a Fr\'echet version of \nameref{LLtS3ref}.
\begin{description}

    \item[\label{LLsS3ref} \( \bb{(\mathcal{S}_{3}^{+})} \)]
      \emph{\( - f(\upsv) = n \hL(\upsv) \) with \( \hL(\cdot) \) strongly concave, 
      \( \HL^{2}(\upsv) = \nabla^{2} \hL(\upsv) \), 
      and
\begin{EQA}
	\sup_{\| \HL(\upsv) \uv \| \leq \rr/\sqrt{n}}  \,\, \sup_{\gammav \in \R^{\dimp}} \,\,
	\frac{\bigl| \langle \nabla^{3} \hL(\upsv + \uv), \uv \otimes \gammav^{\otimes 2} \rangle \bigr|}
		 {\| \HL(\upsv) \uv \|\, \| \HL(\upsv) \gammav \|^{2}}
	& \leq &
	\hmax_{3} \, .
\end{EQA}
}
\end{description}

\begin{lemma}
\label{LfreTay}
With \( \dltwbss(\upsv) \) from \eqref{jcxhydtferyufgy7tfdsy7ft}, it holds for any \( \uv \) with 
\( \| \IF^{1/2}(\upsv) \uv \| \leq \rr \)
\begin{EQA}
	\| \IF^{-1/2}(\upsv) \, \IF(\upsv + \uv) \, \IF^{-1/2}(\upsv) - \Id_{\dimp} \|
	& \leq &
	\dltwbss(\upsv) .
\label{ghdrd324ee4ew222w3ew}
\end{EQA}
Moreover, \nameref{LLsS3ref} yields \( \dltwbss(\upsv) \leq \hmax_{3} \, \rr \, n^{-1/2} \) and
for any \( \uv \) with \( \| \IF^{1/2}(\upsv) \uv \| \leq \rr \)
\begin{EQA}
	\| \IF^{-1/2}(\upsv) \, \IF(\upsv + \uv) \, \IF^{-1/2}(\upsv) - \Id_{\dimp} \|
	& \leq &
	\frac{\hmax_{3}}{n^{1/2}} \, \| \IF^{1/2}(\upsv) \uv \|
	\leq 
	\frac{\hmax_{3} \, \rr}{n^{1/2}} \, .
\label{ghdrd324ee4ew222w3ewss}
\end{EQA}
\end{lemma}

\begin{proof}
Denote 
\( \Delta(\uv) = \IF(\upsv + \uv) - \IF(\upsv) \).
Then by \eqref{jcxhydtferyufgy7tfdsy7ft} for any \( \gammav \in \R^{\dimp} \) with \( \deltav = \IF^{-1/2}(\upsv) \gammav \)
\begin{EQA}
	&& \nquad
	\bigl| \bigl\langle \IF^{-1/2}(\upsv) \, \IF(\upsv + \uv) \, \IF^{-1/2}(\upsv) - \Id_{\dimp}, \gammav^{\otimes 2} \bigr\rangle \bigr|
	=
	\bigl| \langle \Delta(\uv), \deltav^{\otimes 2} \rangle \bigr|
	\\
	& \leq &
	\dltwbss(\upsv) \, \| \IF^{1/2}(\upsv) \, \deltav \|^{2} 
	=
	\dltwbss(\upsv) \, \| \gammav \|^{2} .
\label{9487654r5tghasdfg}
\end{EQA}
This yields \eqref{ghdrd324ee4ew222w3ew}.
Bound \eqref{ghdrd324ee4ew222w3ewss} can be proved by arguments from Lemma~\ref{LdltwLaGP}.
\end{proof}

Define now for any \( \uv \)
\begin{EQA}
	\IFba(\upsv;\uv)
	& \eqdef &
	\int_{0}^{1} \IF(\upsv + t \uv) \, dt \, .
\label{vhudfg7sdfyt7s8hgbfuhbugrh9}
\end{EQA}

\begin{lemma}
\label{LfreTayb}
Under the conditions of Lemma~\ref{LfreTay}, for any \( \uv \) with 
\( \| \IF^{1/2}(\upsv) \uv \| \leq \rr \)
\begin{EQ}[rcl]
	\bigl\| \IF^{-1/2}(\upsv) \, \IFba(\upsv ; \uv) \, \IF^{-1/2}(\upsv) - \Id_{\dimp} \bigr\|
	& \leq &
	\dltwbss(\upsv) 
	\, .
\label{dsvfhgevrfgtwfyikoub}
\end{EQ}
Moreover,
\begin{EQA}[rcccl]
	\{ 1 - \dltwbss(\upsv) \} \IF(\upsv)
	& \leq &
	\IFba(\upsv ; \uv)
	& \leq &
	\{ 1 + \dltwbss(\upsv) \} \IF(\upsv)
\label{5chfdc7e3yvc5ededww}
\end{EQA}
\end{lemma}

\begin{proof}
For any \( \gammav \in \R^{\dimp} \) and \( t \in [0,1] \), the definition \eqref{jcxhydtferyufgy7tfdsy7ft} implies
\begin{EQA}
	\bigl| \bigl\langle \IF(\upsv + t \uv) - \IF(\upsv), \gammav^{\otimes 2} \bigr\rangle \bigr|
	& \leq &
	\dltwbss(\upsv) \, \| \IF^{1/2}(\upsv) \, \gammav \|^{2} .
\label{09vbfkmdior329iufe9o}
\end{EQA}
This obviously yields under \nameref{LLsS3ref}
\begin{EQA}
	\bigl| \bigl\langle \IFba(\upsv; \uv) - \IF(\upsv), \gammav^{\otimes 2} \bigr\rangle \bigr|
	& \leq &
	\frac{\hmax_{3} \, \rr}{n^{1/2}} \, \| \IF^{1/2}(\upsv) \, \gammav \|^{2} \int_{0}^{1} t \, dt
	=
	\frac{\hmax_{3} \, \rr}{2 n^{1/2}} \, \| \IF^{1/2}(\upsv) \, \gammav \|^{2} 
\label{09vbfkmdior329iufe9ob}
\end{EQA}
and \eqref{dsvfhgevrfgtwfyikoub} follows as in Lemma~\ref{LfreTay}.
\end{proof}



\Section{Quadraticity and near quadraticity}
\label{Squadnquad}

Let \( \fs(\upsv) \) be a smooth concave function, 
\begin{EQA}
	\upsvs
	&=&
	\argmax_{\upsv} \fs(\upsv),
\label{fg5hg3gf98tkj3dciryt}
\end{EQA}
and \( \DP^{2} = - \nabla^{2} \fs(\upsvs) \).
Later we study the question how the point of maximum and the value of maximum of \( \fs \) change if we add a linear or quadratic 
component to \( \fs \).

\Subsection{A linear perturbation}
This section studies the case of a linear change of \( \fs \).
More precisely, let another function \( \fn(\upsv) \) satisfy for some vector \( \Av \)
\begin{EQA}
	\fn(\upsv) - \fn(\upsvs) 
	&=&
	\bigl\langle \upsv - \upsvs, \Av \bigr\rangle + \fs(\upsv) - \fs(\upsvs) .
\label{4hbh8njoelvt6jwgf09}
\end{EQA}
A typical example corresponds to \( \fs(\upsv) = \E L(\upsv) \) and \( \fn(\upsv) = L(\upsv) \) 
for a random function \( L(\upsv) \) with a linear stochastic component \( \zeta(\upsv) = L(\upsv) - \E L(\upsv) \);
see \nameref{Eref}.
Then \eqref{4hbh8njoelvt6jwgf09} is satisfied with 
\begin{EQA}
	\Av 
	&=&
	\nabla \zeta .
\label{mdfwdgfwes3e4gfdyc7f}
\end{EQA}
The aim of the analysis is evaluate the values 
\begin{EQA}
	\upsvn
	& \eqdef &
	\argmax_{\upsv} \fn(\upsv),
	\qquad
	\fn(\upsvn)
	=
	\max_{\upsv} \fn(\upsv) .
\label{6yc63yhudf7fdy6edgehy} 
\end{EQA}
The results will be stated mainly in terms of the quantity \( \| \DP^{-1} \Av \| \).
First we consider the case of a quadratic function \( \fs \).

\begin{lemma}
\label{Pquadquad}
Let \( \fs(\upsv) \) be quadratic with \( \nabla^{2} \fs(\upsv) \equiv - \DP^{2} \).
If \( \fn(\upsv) \) satisfy \eqref{4hbh8njoelvt6jwgf09}, then 
\begin{EQA}
	\upsvn - \upsvs
	&=&
	\DP^{-2} \Av,
	\qquad
	\fn(\upsvn) - \fn(\upsvs)
	=
	\frac{1}{2} \| \DP^{-1} \Av \|^{2} .
\label{kjcjhchdgehydgtdtte35}
\end{EQA}
\end{lemma}

\begin{proof}
If \( \fs(\upsv) \) is quadratic, then, of course, under \eqref{4hbh8njoelvt6jwgf09}, \( \fn(\upsv) \) is quadratic as well
with \( \nabla^{2} \fn(\upsv) \equiv - \DP^{2} \).
This implies
\begin{EQA}
	\nabla \fn(\upsvs) - \nabla \fn(\upsvn)
	&=&
	\DP^{2} (\upsvn - \upsvs) .
\label{dcudydye67e6dy3wujhds7}
\end{EQA}
Further, \eqref{4hbh8njoelvt6jwgf09} and \( \nabla \fs(\upsvs) = 0 \) yield \( \nabla \fn(\upsvs) = \Av \).
Together with \( \nabla \fn(\upsvn) = 0 \), this implies
\( \upsvn - \upsvs = \DP^{-2} \Av \).
The Taylor expansion of \( \fn \) at \( \upsvn \) yields by \( \nabla \fn(\upsvn) = 0 \)
\begin{EQA}
	\fn(\upsvs) - \fn(\upsvn)
	&=&
	- \frac{1}{2} \| \DP (\upsvn - \upsvs) \|^{2}
	=
	- \frac{1}{2} \| \DP^{-1} \Av \|^{2} 
\label{8chuctc44wckvcuedje}
\end{EQA}
and the assertion follows.
\end{proof}
  
The next result describes the concentration properties of \( \upsvn \) from \eqref{6yc63yhudf7fdy6edgehy} in a local elliptic set
of the form
\begin{EQA}
	\CA(\rr)
	& \eqdef &
	\{ \upsv \colon \| \DP (\upsv - \upsvs) \| \leq \rr \} ,
\label{0cudc7e3jfuyvct6eyhgwe}
\end{EQA}
where \( \rr \) is slightly larger than \( \| \DP^{-1} \Av \| \).

\begin{proposition}
\label{Pconcgeneric}
Let \( \fs(\upsv) \) be a concave function with \( \fs(\upsvs) = \max_{\upsv} \fs(\upsv) \)  
and \( \DP^{2} = - \nabla^{2} \fs(\upsvs) \).
Let further \( \fn(\upsv) \) and \( \fs(\upsv) \) be related by \eqref{4hbh8njoelvt6jwgf09} with some vector \( \Av \).
Fix \( \amax \leq 2/3 \) and \( \rrn \) such that \( \| \DP^{-1} \Av \| \leq \amax \, \rrn \).
Suppose now that \( \fs(\upsv) \) satisfy \eqref{dtb3u1DG2d3GPg} for \( \upsv = \upsvs \) and \( \HLG(\upsvs) = \DP \), 
and \( \dltwbd \) such that 
\begin{EQA}
	1 - \amax - \dltwbd
	& \geq &
	0 .
\label{rrm23r0ut3ua}
\end{EQA}
Then  \( \upsvn \) from \eqref{6yc63yhudf7fdy6edgehy} satisfies 
\begin{EQA}
	\| \DP (\upsvn - \upsvs) \|  
	& \leq &
	\rrn \, . 
\label{rhDGtuGmusGU0a}
\end{EQA}
\end{proposition}

\begin{proof}
The bound \( \| \DP^{-1} \Av \| \leq \amax \, \rrn \) implies for any \( \uv \)
\begin{EQA}
	\bigl| \langle \Av, \uv \rangle \bigr|
	& = &
	\bigl| \langle \DP^{-1} \Av, \DP \uv \rangle \bigr|
	\leq 
	\amax \, \rrn \| \DP \uv \| \, .
\label{LLoDGm1nzua}
\end{EQA}
If \( \| \DP \uv \| > \rrn \), then \( \rrn \| \DP \uv \| \leq \| \DP \uv \|^{2} \). 
Therefore,
\begin{EQA}
	\bigl| \langle \Av, \uv \rangle \bigr|
	& \leq &
	\amax \| \DP \uv \|^{2} \, ,
	\qquad
	\| \DP \uv \| > \rrn \, .
\label{LLoDGm1nzun}
\end{EQA}
Let \( \upsv \) be a point on the boundary of the set \( \CA(\rrn) \) from \eqref{0cudc7e3jfuyvct6eyhgwe}.
We also write \( \uv = \upsv - \upsvs \).
The idea is to show that the derivative  \( \frac{d}{dt} \fn(\upsvs + t \uv) < 0 \) 
is negative for \( t > 1 \).
Then all the extreme points of \( \fn(\upsv) \) are within \( \CA(\rrn) \).
We use the decomposition
\begin{EQA}
	\fn(\upsvs + \rhot \uv) - \fn(\upsvs)
	&=&
	\bigl\langle \Av, \uv \bigr\rangle \, \rhot 
	+ \fs(\upsvs + \rhot \uv) - \fs(\upsvs) .
\label{LGtsGtuLGtsa}
\end{EQA}
With \( \fGu(t) = \fs(\upsvs + \rhot \uv) \), it holds
\begin{EQA}
	\frac{d}{d \rhot} \fn(\upsvs + \rhot \uv)
	&=&
	\bigl\langle \Av, \uv \bigr\rangle + \fGu'(\rhot) .
\label{frddtLtGstua}
\end{EQA}
By definition of \( \upsvs \), it also holds \( \fGu'(0) = 0 \).
The identity \( \nabla^{2} \fs(\upsvs) = - \DP^{2} \) yields \( \fGu''(0) = - \| \DP \uv \|^{2} \).
Bound \eqref{dta3u1DG2d3GPa1g} implies for \( | \rhot | \leq 1 \)
\begin{EQA}
	\bigl| \fGu'(\rhot) - \rhot \fGu''(0) \bigr|
	&=&
	\bigl| \fGu'(\rhot) - \fGu'(0) - \rhot \fGu''(0) \bigr|
	\leq 
	\rhot^{2} \, \bigl| \fGu''(0) \bigr| \, \dltwbd \, .
\label{fptfp0fpttfpp13a}
\end{EQA}
For \( \rhot = 1 \), we obtain by \eqref{rrm23r0ut3ua}
\begin{EQA}
	\fGu'(1) 
	& \leq &
	\fGu''(0) - \fGu''(0) \, \dltwbd
	=
	- \bigl| \fGu''(0) \bigr| (1 -  \dltwbd)
	< 0 .
\label{fp1fpp13d3rGa}
\end{EQA}
Moreover, concavity of \( \fGu(\rhot) \) and \( \fGu'(0) = 0 \) imply that \( \fGu'(\rhot) \) decreases in 
\( \rhot \) for \( \rhot > 1 \).
Further, summing up the above derivation yields 
\begin{EQA}
	\frac{d}{dt} \fn(\upsvs + \rhot \uv) \Big|_{\rhot=1}
	& \leq &
	- \| \DP \uv \|^{2} (1 - \amax - \dltwbd)
	< 0 .
\label{ddtLGtstu33a}
\end{EQA}
As \( \frac{d}{d \rhot} \fn(\upsvs + \rhot \uv) \) decreases with \( \rhot \geq 1 \) together with 
\( \fGu'(\rhot) \) due to \eqref{frddtLtGstua}, the same applies to all such \( \rhot \).
This implies the assertion.
\end{proof}

The result of Proposition~\ref{Pconcgeneric} allows to localize the point \( \upsvn = \argmax_{\upsv} \fn(\upsv) \)
in the local vicinity \( \CA(\rrn) \) of \( \upsvs \).
The use of smoothness properties of \( \fn \) or, equivalently, of \( \fs \), in this vicinity helps to obtain
rather sharp expansions for \( \upsvn - \upsvs \) and for \( \fn(\upsvn) - \fn(\upsvs) \); cf. \eqref{kjcjhchdgehydgtdtte35}.

\begin{proposition}
\label{PFiWigeneric}
Under the conditions of Proposition~\ref{Pconcgeneric}
\begin{EQ}[rcccl]
    - \frac{\dltwb}{1 + \dltwb} \| \DP^{-1} \Av \|^{2}
    & \leq &
    2 \fn(\upsvn) - 2 \fn(\upsvs) 
    - \| \DP^{-1} \Av \|^{2}
    & \leq &
    \frac{\dltwb}{1 - \dltwb} \| \DP^{-1} \Av \|^{2} \, .
    \qquad
\label{3d3Af12DGttGa}
\end{EQ}
Also
\begin{EQ}[rcl]
    \| \DP (\upsvn - \upsvs) - \DP^{-1} \Av \|^{2}
    & \leq &
    \frac{3 \dltwb}{(1 - \dltwb)^{2}} \, \| \DP^{-1} \Av \|^{2} \, ,
    \\
    \| \DP (\upsvn - \upsvs) \|
    & \leq &
    \frac{1 + \sqrt{2 \dltwb}}{1 - \dltwb} \, \| \DP^{-1} \Av \| \, .
\label{DGttGtsGDGm13rGa}
\end{EQ}
\end{proposition}

\begin{proof}
By \eqref{dtb3u1DG2d3GPg}, for any \( \upsv \in \CA(\rrn) \)
\begin{EQA}
	\Bigl| 
		\fs(\upsvs) - \fs(\upsv) - \frac{1}{2} \| \DP (\upsv - \upsvs) \|^{2} 
	\Bigr|
	& \leq &
	\frac{\dltwb}{2} \| \DP (\upsv - \upsvs) \|^{2} .
\label{d3GrGELGtsG12}
\end{EQA}
Further, 
\begin{EQA}[rcl]
	&& \nquad
	\fn(\upsv) - \fn(\upsvs) - \frac{1}{2} \| \DP^{-1} \Av \|^{2}
	\\
	&=&
	\bigl\langle \upsv - \upsvs, \Av \bigr\rangle
	+ \fs(\upsv) - \fs(\upsvs) - \frac{1}{2} \| \DP^{-1} \Av \|^{2} 
	\\
	&=&
	- \frac{1}{2} \bigl\| \DP (\upsv - \upsvs) - \DP^{-1} \Av \bigr\|^{2}
	+ \fs(\upsv) - \fs(\upsvs) + \frac{1}{2} \| \DP (\upsv - \upsvs) \|^{2} .
	\qquad 
\label{12ELGuELusG}
\end{EQA}
As \( \upsvn \in \CA(\rrn) \) and it maximizes \( \fn(\upsv) \), we derive by \eqref{d3GrGELGtsG12}
\begin{EQA}
	&& \nquad
	\fn(\upsvn) - \fn(\upsvs) - \frac{1}{2} \| \DP^{-1} \Av \|^{2}
	=
	\max_{\upsv \in \CA(\rrn)} 
	\Bigl\{ 
		\fn(\upsv) - \fn(\upsvs) - \frac{1}{2} \| \DP^{-1} \Av \|^{2} 
	\Bigr\}
	\\
	& \leq &
	\max_{\upsv \in \CA(\rrn)} 
	\Bigl\{ 
		- \frac{1}{2} \bigl\| \DP (\upsv - \upsvs) - \DP^{-1} \Av \bigr\|^{2} 
		+ \frac{\dltwb}{2} \| \DP (\upsv - \upsvs) \|^{2}
	\Bigr\} .
\label{d3G1212222B} 
\end{EQA}
Further, \( \max_{\uv} \bigl\{ \dltwb \| \uv \|^{2} - \| \uv - \xiv \|^{2} \bigr\} = \frac{\dltwb}{1 - \dltwb} \| \xiv \|^{2} \)
for \( \dltwb \in [0,1) \) and \( \xiv \in \R^{\dimp} \), 
yielding
\begin{EQA}
	\fn(\upsvn) - \fn(\upsvs) - \frac{1}{2} \| \DP^{-1} \Av \|^{2}
	& \leq &
	\frac{\dltwb}{2(1 - \dltwb)} \| \DP^{-1} \Av \|^{2} . 
\label{fd3G122B2} 
\end{EQA}
Similarly 
\begin{EQA}
	&& \nquad \quad
	\fn(\upsvn) - \fn(\upsvs) - \frac{1}{2} \| \DP^{-1} \Av \|^{2}
	\geq 
	\max_{\upsv \in \CA(\rrn)} 
	\Bigl\{ 
		- \frac{1}{2} \bigl\| \DP (\upsv - \upsvs) - \DP^{-1} \Av \bigr\|^{2} 
		- \frac{\dltwb}{2} \| \DP (\upsv - \upsvs) \|^{2}
	\Bigr\}
	\\
	&=&
	- \frac{\dltwb }{2(1 + \dltwb)} \, \| \DP^{-1} \Av \|^{2} . 
	\qquad \quad
\label{fd3G122B2m} 
\end{EQA}
These bounds imply 
\eqref{3d3Af12DGttGa}.

Now we derive similarly to \eqref{12ELGuELusG} that for \( \upsv \in \CA(\rrn) \)
\begin{EQA}
	\fn(\upsv) - \fn(\upsvs) 
	& \leq &
	\bigl\langle \upsv - \upsvs, \Av \bigr\rangle
	- \frac{1 - \dltwb}{2} \| \DP (\upsv - \upsvs) \|^{2} .
\label{LGvLGvsGf1d3G2}
\end{EQA}
A particular choice \( \upsv = \upsvn \) yields
\begin{EQA}
	\fn(\upsvn) - \fn(\upsvs) 
	& \leq &
	\bigl\langle \upsvn - \upsvs, \Av \bigr\rangle
	- \frac{1 - \dltwb}{2} \| \DP (\upsvn - \upsvs) \|^{2} .
\label{21GsvtvGDG3G2}
\end{EQA}
Combining with \eqref{fd3G122B2m} allows to bound
\begin{EQA}
	&& \nquad
	\bigl\langle \upsvn - \upsvs, \Av \bigr\rangle
	- \frac{1 - \dltwb}{2} \| \DP (\upsvn - \upsvs) \|^{2} 
	- \frac{1}{2} \| \DP^{-1} \Av \|^{2}
	\geq 
	- \frac{\dltwb}{2(1 + \dltwb)} \| \DP^{-1} \Av \|^{2} .
\label{2m1DGd3G123G}
\end{EQA}
Further, for \( \xiv = \DP^{-1} \Av \), 
\( \uv = \DP (\upsvn - \upsvs) \), and 
\( \dltwb \in [0,1/3] \), the inequality 
\begin{EQA}
	2 \bigl\langle \uv, \xiv \bigr\rangle - (1 - \dltwb) \| \uv \|^{2} - \| \xiv \|^{2}
	& \geq &
	- \frac{\dltwb}{1 + \dltwb} \| \xiv \|^{2} 
\label{dtxi2fd1d22}
\end{EQA}
implies 
\begin{EQA}
	\bigl\| \uv - \frac{1}{1-\dltwb} \xiv \bigr\|^{2}
	& \leq &
	\frac{2 \dltwb}{(1 + \dltwb) (1 - \dltwb)^{2}} \| \xiv \|^{2}
\label{uv11wxi22w1w}
\end{EQA}
yielding for \( \dltwb \leq 1/3 \)
\begin{EQA}
	\| \uv - \xiv \|
	& \leq &
	\biggl( \dltwb + \sqrt{\frac{2 \dltwb}{1 + \dltwb}} \biggr)
	\frac{\| \xiv \|}{1 - \dltwb}
	\leq 
	\frac{\sqrt{3 \dltwb} \| \xiv \|}{1 - \dltwb} \, ,
	\\
	\| \uv \|
	& \leq &
	\biggl( 1 + \sqrt{\frac{2 \dltwb}{1 + \dltwb}} \biggr)
	\frac{\| \xiv \|}{1 - \dltwb}
	\leq 
	\frac{1 + \sqrt{2 \dltwb} \| \xiv \|}{1 - \dltwb} \, ,
\label{uxiBws2w1w31w}
\end{EQA}
and \eqref{DGttGtsGDGm13rGa} follows.
\end{proof}

\Subsection{Quadratic penalization}
Here we discuss the case when \( \fn(\upsv) - \fs(\upsv) \) is quadratic.
The general case can be reduced to the situation with \( \fn(\upsv) = \fs(\upsv) - \| \GP \upsv \|^{2}/2 \).
To make the dependence of \( \GP \) more explicit, denote 
\begin{EQA}
	\fG(\upsv) 
	&=& 
	\fs(\upsv) - \| \GP \upsv \|^{2}/2 .
\label{8cvkfc9fujf6jnmcer4cd}
\end{EQA}
With \( \upsvs = \argmax_{\upsv} \fs(\upsv) \) and \( \upsvs_{\GP} = \argmax_{\upsv} \fG(\upsv) \),
we study the bias \( \upsvs_{\GP} - \upsvs \) induced by this penalization
under Fr\'echet-type smoothness conditions. 
To get some intuition, consider first the case of a quadratic function \( \fs(\upsv) \).

\begin{lemma}
\label{Lbiasquadgen}
Let \( \fs(\upsv) \) be quadratic with \( \IF \equiv - \nabla^{2} \fs(\upsv) \) and \( \Av_{\GP} \equiv - \GP^{2} \upsvs \).
Then it holds with \( \IF_{\GP} = \IF + \GP^{2} \)
\begin{EQA}[rcccl]
	\upsvs_{\GP} - \upsvs
	&=&
	\IF_{\GP}^{-1} \Av_{\GP} 
	&=&
	- \IF_{\GP}^{-1} \GP^{2} \upsvs,
	\\
	\fG(\upsvs_{\GP}) - \fG(\upsvs)
	&=&
	\frac{1}{2} \| \IF_{\GP}^{-1/2} \Av_{\GP} \|^{2}
	&=&
	\frac{1}{2} \| \IF_{\GP}^{-1/2} \GP^{2} \upsvs \|^{2} \, .
\label{kjfydf554dfwertdgwdf}
\end{EQA}
\end{lemma} 

\begin{proof}
Quadraticity of \( \fs(\upsv) \) implies quadraticity of \( \fG(\upsv) \) with \( \nabla^{2} \fG(\upsv) \equiv - \IF_{\GP} \).
This implies
\begin{EQA}
	\nabla \fG(\upsvs) - \nabla \fG(\upsvs_{\GP})
	&=&
	\IF_{\GP} (\upsvs_{\GP} - \upsvs) .
\label{dcudydye67e6dy3wujhdqu}
\end{EQA}
Further, \( \nabla \fs(\upsvs) = 0 \) yielding \( \nabla \fG(\upsvs) = \Av_{\GP} = - \GP^{2} \upsvs \).
Together with \( \nabla \fG(\upsvs_{\GP}) = 0 \), this implies
\( \upsvs_{\GP} - \upsvs = \IF_{\GP}^{-1} \Av_{\GP} \).
The Taylor expansion of \( \fG \) at \( \upsvs_{\GP} \) yields with \( \DPGP = \IF_{\GP}^{1/2} \) 
\begin{EQA}
	\fG(\upsvs) - \fG(\upsvs_{\GP})
	&=&
	- \frac{1}{2} \| \DPGP (\upsvs_{\GP} - \upsvs) \|^{2}
	=
	- \frac{1}{2} \| \DPGP^{-1} \Av_{\GP} \|^{2} 
\label{8chuctc44wckvcuedjequ}
\end{EQA}
and the assertion follows.
\end{proof}

Now we turn to the general case with \( \fs \) smooth in Fr\'echet sense.
Let \( \IF(\upsv) = - \nabla^{2} \fs(\upsv) \).
For \( \upsv \in \Ups \) and \( \uv \in \R^{\dimp} \), define
\begin{EQA}
	\IFba(\upsv;\uv)
	& \eqdef &
	\int_{0}^{1} \IF(\upsv + t \uv) \, dt \, ;
\label{vhudfg7hgbfuhbugr}
\end{EQA}
cf. \eqref{vhudfg7sdfyt7s8hgbfuhbugrh9}.
Similarly define 
\( \IF_{\GP}(\upsv) = - \nabla^{2} \fG(\upsv) = \IF(\upsv) + \GP^{2} \) and
\begin{EQA}
	\IFba_{\GP}(\upsv;\uv)
	& \eqdef &
	\int_{0}^{1} \IF_{\GP}(\upsv + t \uv) \, dt
	=
	\IFba(\upsv;\uv) + \GP^{2} .
\label{virjtfgy6f53e5fuhuyyuu}
\end{EQA}

\begin{lemma}
\label{LfreTgrad}
For the vector \( \biasv_{\GP} \eqdef \upsvs_{\GP} - \upsvs \), define \( \IFba = \IFba(\upsvs;\biasv) \).
Then
\begin{EQA}
	\upsvs_{\GP} - \upsvs
	&=&
	\IFba^{-1} \, \Av_{\GP} \, .
\label{wivufhu338fvjed3jfufb}
\end{EQA}
\end{lemma} 

\begin{proof}
First we show for any \( \upsv \in \Ups \) and \( \uv \in \R^{\dimp} \) that
\begin{EQA}
	\nabla \fG(\upsv + \uv) - \nabla \fG(\upsv) 
	& = &
	- \IFba_{\GP}(\upsv;\uv) \, \uv .
\label{dfhhjvqwbhewjmaqjmqjew}
\end{EQA}
Indeed, for any \( \gammav \in \R^{\dimp} \), consider the univariate function 
\( h(t) = \langle \nabla \fG(\upsv + t \uv) - \nabla \fG(\upsv), \gammav \rangle \).
Statement \eqref{dfhhjvqwbhewjmaqjmqjew} follows from definitions \eqref{vhudfg7hgbfuhbugr}, \eqref{virjtfgy6f53e5fuhuyyuu}, 
and the identity
\( h(1) - h(0) = \int_{0}^{1} h'(t) \, dt \).

Further, the definition \( \upsvs = \argmax_{\upsv} \fs(\upsv) \) yields \( \nabla \fs(\upsvs) = 0 \) and
\begin{EQA}
	\nabla \fG(\upsvs)
	&=& 
	\nabla f(\upsvs) - \GP^{2} \upsvs
	=
	\Av_{\GP} \, .
\label{dfse423q2esftcsdfgsygyb}
\end{EQA}
In view of \( \nabla \fG(\upsvs_{\GP}) = 0 \), we derive
\begin{EQA}
	\nabla \fG(\upsvs_{\GP}) - \nabla \fG(\upsvs)
	&=&
	- \Av_{\GP} \, .
\label{4hv78gfu7jh3ewo0ggt2lobhb}
\end{EQA}
Representation \eqref{dfhhjvqwbhewjmaqjmqjew} with \( \upsv = \upsvs \) and \( \uv = \biasv_{\GP} \) yields \eqref{wivufhu338fvjed3jfufb}.
\end{proof}

Representation \eqref{wivufhu338fvjed3jfufb} is very useful to bound from above the bias \( \upsvs - \upsvs_{\GP} \).
Indeed, assuming that \( \upsvs_{\GP} \) is in a local vicinity of \( \upsvs \) we may use 
Fr\'echet smoothness of \( \fs \) in terms of the value \( \dltwbss(\upsvs) \) from \eqref{jcxhydtferyufgy7tfdsy7ft} to approximate 
\( \IFba_{\GP} \approx \IF_{\GP}(\upsvs) \) and
\( \upsvs_{\GP} - \upsvs \approx - \IF_{\GP}^{-1}(\upsvs) \, \GP^{2} \upsvs \).
 
\begin{proposition}
\label{Pbiasgeneric} 
Define \( \DPfGP \) by
\begin{EQA}
	\DPfGP^{2}
	& = &
	\IF_{\GP}(\upsvs) ,
\label{95jioo00853rdvfyu7}
\end{EQA} 
and let \( \QP \leq \DPfGP \).
Fix 
\begin{EQA}
	\rrs 
	&=& 
	\| \QP \, \DPfGP^{-2} \, \GP^{2} \upsvs \| 
	=
	\| \QP \, \DPfGP^{-2} \, \Av_{\GP} \|
\label{fd9dfhy4ye6fuydfrerf}
\end{EQA} 
and \( \rru > \amax^{-1} \rrs \) with \( \amax = 2/3 \).
Assume
\begin{EQA}
	\dltwbs_{\GP}
	\eqdef
	\sup_{\uv \colon \| \QP \uv \| \leq \rru} \,\, 
	\| \DPfGP^{-1} \, \IF_{\GP}(\upsvs + \uv) \, \DPfGP^{-1} - \Id_{\dimp} \|
	& \leq &
	\frac{1}{3} \, ;
\label{ghdrd324ee4ew222gen}
\end{EQA}
cf. \eqref{ghdrd324ee4ew222w3ew}. 
Then \( \| \QP (\upsv - \upsvs) \| \leq \rru \) or, equivalently,
\begin{EQA}
	\upsvs_{\GP}
	 & \in &
	 \CA_{\GP}
	 \eqdef
	 \{ \upsv \colon \| \QP (\upsv - \upsvs) \| \leq \rru \} .
\label{odf6fdyr6e4deuewjug}
\end{EQA}
Moreover, 
\begin{EQA}
	\| \QP (\biasv_{\GP} - \DPfGP^{-2} \Av_{\GP}) \|
	& \leq &
	\dltwbs_{\GP} \, \rru \, .
\label{11ma3eaelDebbgen}
\end{EQA}
\end{proposition}

\begin{proof}
First we check that \( \upsvs_{\GP} \) concentrates in the local vicinity 
\( \CA_{\GP} \) from \eqref{odf6fdyr6e4deuewjug}.
Strong concavity of \( \fG \) implies that the solution \( \upsvs_{\GP} \) exists and unique.
Let us fix any direction \( \gammav \in \R^{\dimp} \) with \( \| \QP \, \gammav \| = \rru \).
Due to \eqref{wivufhu338fvjed3jfufb}, 
we are looking at the solution \( \IFba_{\GP}(\upsvs; s \gammav) \, s \gammav = \Av_{\GP} \) in \( s \) and \( \gammav \).
It suffices to ensure that 
\begin{EQA}
	s \, \QP \, \gammav 
	& = &
	\QP \, \IFba_{\GP}(\upsvs; s \gammav)^{-1} \Av_{\GP} 
\label{ihy666677erdft4e3te4yt}
\end{EQA}
is impossible for \( s \geq 1 \).
For \( s = 1 \), we can use that \( \IFba_{\GP}(\upsvs; \gammav) \geq (1 - \dltwbs_{\GP}) \DPfGP^{2} \);
see \eqref{5chfdc7e3yvc5ededww} of Lemma~\ref{LfreTayb}.
Therefore,
\begin{EQA}
	\| \QP \, \IFba_{\GP}^{-1}(\upsvs; \gammav) \Av_{\GP} \|
	& \leq &
	\| \QP \, \IFba_{\GP}^{-1}(\upsvs; \gammav) \DPfGP^{2} \QP^{-1} \| \, \| \QP \, \DPfGP^{-2} \, \Av_{\GP} \|
	\\
	& \leq &
	(1 - \dltwbs_{\GP})^{-1} \rrs
	<
	\rru
\label{v565tf5tfd56fd65tgtsydx}
\end{EQA}
and \eqref{ihy666677erdft4e3te4yt} with \( s=1 \) is impossible because \( \| \QP \gammav \| = \rru \).
It remains to note that the matrix \( s \, \IFba_{\GP}(\upsvs; s \gammav) \) grows with \( s \) as
\begin{EQA}
	s \, \IFba_{\GP}(\upsvs; s \gammav)
	&=&
	s \int_{0}^{1} \IF_{\GP}(\upsvs + t \, s \, \gammav) \, dt
	=
	\int_{0}^{s} \IF_{\GP}(\upsvs + t \, \gammav) \, dt .
\label{r0ofigvjbh3wd9eorikcjv}
\end{EQA}
Now we bound \( \| \QP \, \biasv_{\GP} \| \) assuming that \( \| \QP (\upsvs_{\GP} - \upsvs) \| \leq \rru \)
and \eqref{ghdrd324ee4ew222gen} applies.  
Statement \eqref{5chfdc7e3yvc5ededww} of Lemma~\ref{LfreTayb} implies 
\( \IFba^{-1} \leq (1 - \dltwbs_{\GP})^{-1} \DPfGP^{2} \) and
\begin{EQA}
	\| \QP \, \biasv_{\GP} \|
	&=&
	\| \QP \, \IFba^{-1} \, \DPfGP^{2} \, \QP^{-1} \, \QP \, \DPfGP^{-2} \, \Av_{\GP} \|
	\leq 
	\| \QP \, \IFba^{-1} \, \DPfGP^{2} \, \QP^{-1} \| \, \| \QP \, \DPfGP^{-2} \Av_{\GP} \|
	\\
	& \leq &
	\frac{1}{1 - \dltwbs_{\GP}} \, \| \QP \, \DPfGP^{-2} \Av_{\GP} \| 
	=
	\frac{\rrs}{1 - \dltwbs_{\GP}} \, .
\label{jhfdt3dtyiu89jgrert6b}
\end{EQA}
In the same way we derive
\begin{EQA}
	\| \QP \, (\biasv - \DPfGP^{-2} \Av_{\GP}) \|
	&=&
	\bigl\| \QP \, (\IFba^{-1} - \DPfGP^{-2}) \Av_{\GP} \bigr\|
	\leq 
	\frac{\dltwbs_{\GP}}{1 - \dltwbs_{\GP}} \, \| \QP \, \DPfGP^{-2} \Av_{\GP} \| \, ,
\label{ae3DeuseDeutb}
\end{EQA}
and \eqref{11ma3eaelDebbgen} follows as well. 
\end{proof}

\begin{remark}
\label{Rbiasgeneric}
Inspection of the proofs of Proposition~\ref{Pbiasgeneric} indicates that 
the results \eqref{odf6fdyr6e4deuewjug} through \eqref{11ma3eaelDebbgen} can be restated 
with \( \DPGP^{2} = \IF_{\GP}(\upsvs_{\GP}) \) in place of \( \DPfGP^{2} = \IF_{\GP}(\upsvs) \).
\end{remark}


\Chapter{Examples of priors}
\label{Spriors}
This section presents two typical examples of priors and some properties 
including the bounds for \emph{effective dimension} and \emph{Laplace effective dimension}.

\Section{Truncation and smooth priors}
\label{Ssmoothprior}
Below we consider two non-trivial examples of Gaussian priors: 
truncation and smooth priors.
To make the presentation clear, we impose some assumptions on the considered setup.
Most of them are non-restrictive and can be extended to more general situations.
We assume to be given a growing sequence of nested linear approximation subspaces 
\( \VV_{1} \subset \VV_{2} \subset \ldots \subset \R^{\dimp} \) of dimension 
\( \dim(\VV_{\mm}) = \mm \).
Below \( \Proj_{\mm} \) is the projector on \( \VV_{\mm} \) and 
\( \VV_{\mm}^{c} \) is the orthogonal complement of \( \VV_{\mm} \).
A \emph{smooth prior} is described by a self-adjoint operator \( \GP \) such that
\( \| \GP \uv \| / \| \uv \| \) becomes large for \( \uv \in \VV_{\mm}^{c} \) and \( \mm \) large.
One can write this condition in the form
\begin{EQ}[rcl]
	\| \GP \uv \|^{2}
	& \leq &
	\gp_{\mm}^{2} \| \uv \|^{2} ,
	\qquad 
	\uv \in \VV_{\mm} \, ,
	\\
	\| \GP \uv \|^{2}
	& \geq &
	\gp_{\mm}^{2} \| \uv \|^{2} \, ,
	\qquad 
	\uv \in \VV_{\mm}^{c} \, .
\label{uniVmgm2u2uiVm}
\end{EQ}
Often one assumes that \( \VV_{\mm} \) is spanned by the eigenvectors
of \( \GP^{2} \) corresponding to its smallest eigenvalues 
\( \gp_{1}^{2} \leq \gp_{2}^{2} \leq \ldots \leq \gp_{\mm}^{2} \).
We only need \eqref{uniVmgm2u2uiVm}.
A typical example is given by \( \GP^{2} = \diag(\gp_{j}^{2}) \) with \( \gp_{j}^{2} = \CGP^{-1} j^{2\smp} \) for \( \smp > 1/2 \)
and some window parameter \( \CGP \).
Below we refer to this case as \( (\smp,\CGP) \)-\emph{smooth prior}.

A \( \mm \)-\emph{truncation prior} assumes that the prior distribution is restricted to \( \VV_{\mm} \). 
This formally corresponds to a covariance operator \( \GP_{\mm}^{-2} \) with 
\( \GP_{\mm}^{-2} \bigl( \Id - \Proj_{\mm} \bigr) = 0 \).
Equivalently, we set \( \gp_{\mm+1} = \gp_{\mm+2} = \ldots = \infty \) in \eqref{uniVmgm2u2uiVm}.

\Section{Effective dimension}
\label{Seffdima}
This section explains how the \emph{effective dimension} and \emph{Laplace effective dimension}
can be evaluated for some typical situations.  

Let \( \IF \) be a generic information matrix while \( \GP^{2} \) a penalizing matrix. 
With \( \IF_{\GP} \eqdef \IF + \GP^{2} \), define a sub-projector \( \proj_{\GP} \) in \( \R^{\dimp} \) by
\begin{EQA}[c]
	\proj_{\GP}
	\eqdef
	\IF_{\GP}^{-1} \IF
	=
	(\IF + \GP^{2})^{-1} \IF .
\label{9bvgu7fh54rfgy7byrfuer4}
\end{EQA}
Also define the \emph{Laplace effective dimension}
\begin{EQA}
	\dimA(\GP)
	& \eqdef &
	\tr \proj_{\GP}
	=
	\tr \bigl\{ (\IF + \GP^{2})^{-1} \IF \bigr\} .
\label{pjnv6yfkjwwyfvjvikekm}
\end{EQA}
The penalizing matrix \( \GP^{2} \) will be supposed diagonal, 
\( \GP^{2} = \diag\{ \gp_{1}^{2},\ldots,\gp_{\dimp}^{2} \} \).
Moreover, we implicitly assume that the values \( \gp_{j}^{2} \) \emph{grow} with \( j \) at some rate, 
\emph{polynomial} or exponential, yielding for all \( \mm \geq 1 \)
\begin{EQA}
	\sum_{j > \mm} \gp_{j}^{-2} 
	& \leq &
	\CONSTgp \, \mm \, \gp_{\mm}^{-2} \, .
\label{sumjJgjm2C}
\end{EQA}
Our leading example is given by \( \gp_{j}^{2} = \CGP^{-1} j^{2\smp} \) for \( \smp > 1/2 \).
Then \eqref{sumjJgjm2C} holds with \( \CONSTgp = (2 \smp -1)^{-1} \).

Concerning the matrix \( \IF \), we assume 
\begin{EQ}[rcccl]
	\CONSTIF^{-1} \, \nsize \| \uv \|^{2} 
	& \leq &
	\bigl\langle \IF \uv, \uv \bigr\rangle
	& \leq &
	\CONSTIF \, \nsize \| \uv \|^{2},
	\qquad
	\uv \in \R^{\dimp} \, ,
\label{uVmGu2C1C2gen}
\end{EQ}
for some \( \CONSTIF \geq 1 \).
It appears that the value \( \dimA(\GP) \) is closely related to the index \( \mm \) for which 
\( \gp_{\mm}^{2} \approx \nsize \).

\begin{lemma}
\label{LeffdimG}
Let \( \IF \) satisfy \eqref{uVmGu2C1C2gen}. 
Let also \( \GP^{2} = \diag\{ \gp_{1}^{2},\ldots,\gp_{\dimp}^{2} \} \) with \( \gp_{j}^{2} \) satisfying
\eqref{sumjJgjm2C}.
Define the index \( \mm \) as the smallest index \( j \) with \( \gp_{j}^{2} \geq \nsize \):
\begin{EQA}
	\mm
	& = &
	\mm(\GP)
	\eqdef
	\min \{ j \colon \gp_{j}^{2} \geq n \} .
\label{gklghi43mjd7skje8vhr}
\end{EQA}
Then 
\begin{EQA}
	\frac{1}{\CONSTIF + 1}
	& \leq &
	\frac{\dimA(\GP)}{\mm}
	\leq 
	1 + \CONSTIF \, \CONSTgp \, .
\label{4hhyc6vfjhfgyyrheudj}
\end{EQA}
\end{lemma}

\begin{proof}
By \eqref{sumjJgjm2C}
\begin{EQA}
	\tr (\proj_{\GP})
	& \leq &
	\sum_{j \geq 1} \frac{\CONSTIF \, \nsize}{\CONSTIF \, \nsize + \gp_{j}^{2}}
	\leq 
	\mm + \sum_{j > \mm} \frac{\CONSTIF \, \nsize}{\CONSTIF \, \nsize + \gp_{j}^{2}}
	\leq 
	\mm + \CONSTIF \, \nsize \sum_{j > \mm} \gp_{j}^{-2}
	\\
	& \leq &
	\mm + \CONSTIF \, \nsize \, \CONSTgp \, \mm \, \gp_{\mm}^{-2}
	\leq 
	\mm (1 + \CONSTIF \, \CONSTgp) .
\label{pendf6ycvfhewkevtye}
\end{EQA}
Similarly
\begin{EQA}
	\tr (\proj_{\GP})
	& \geq &
	\sum_{j=1}^{\mm} \frac{\CONSTIF^{-1} \, \nsize}{\CONSTIF^{-1} \, \nsize + \gp_{j}^{2}}
	\geq 
	\mm \, \frac{\CONSTIF^{-1} \, \nsize}{\CONSTIF^{-1} \, \nsize + \gp_{\mm}^{2}}
	\geq 
	\mm \, \frac{\CONSTIF^{-1}}{\CONSTIF^{-1} + 1} \, ,
\label{2fdvfg9rfmvvhywhsh}
\end{EQA}
and the assertion follows.
\end{proof}

This result yields an immediate corollary.

\begin{corollary}
\label{CLeffdimG}
Let \( \GP_{1}^{2} \) and \( \GP_{2}^{2} \) be two different penalizing matrices 
satisfying \eqref{sumjJgjm2C} and such that 
\( \mm(\GP_{1}) = \mm(\GP_{2}) \); see \eqref{gklghi43mjd7skje8vhr}.
Then 
\begin{EQA}
	\frac{\dimA(\GP_{1})}{\dimA(\GP_{2})}
	& \leq &
	( 1 + \CONSTIF \, \CONSTgp ) (1 + \CONSTIF) .
\label{c98dfjme478fcdje3otwe}
\end{EQA}
\end{corollary}

Now we evaluate the \emph{effective dimension} \( \dimG = \tr(\IF_{\GP}^{-1} \VP^{2}) \), 
where the \emph{variance matrix} \( \VP^{2} \) satisfies the condition 
\begin{EQ}[rcccl]
	\CONSTV^{-1} \, \| \IF \uv \|^{2} 
	& \leq &
	\| \VP \uv \|^{2}
	& \leq &
	\CONSTV \, \nsize \| \IF \uv \|^{2},
	\qquad
	\uv \in \R^{\dimp} \, 
\label{uVmGu2C1C2genV}
\end{EQ}
with some constant \( \CONSTV \geq 1 \); cf. \eqref{uVmGu2C1C2gen}.

\begin{lemma}
\label{LeffdimGV}
Assume \eqref{uVmGu2C1C2gen} for \( \IF \) and \eqref{uVmGu2C1C2genV} for \( \VP^{2} \). 
Let also \( \GP^{2} = \diag\{ \gp_{1}^{2},\ldots,\gp_{\dimp}^{2} \} \) with \( \gp_{j}^{2} \) satisfying
\eqref{sumjJgjm2C} and let \( \mm \) be given by \eqref{gklghi43mjd7skje8vhr}.
Then \( \dimG = \tr(\IF_{\GP}^{-1} \VP^{2}) \) satisfies
\begin{EQA}
	\frac{\CONSTV^{-1}}{\CONSTIF + 1}
	& \leq &
	\frac{\dimG}{\mm}
	\leq 
	\CONSTV (1 + \CONSTIF \, \CONSTgp) \, .
\label{4hhyc6vfjhfgyyrheudjV}
\end{EQA}
\end{lemma}

\begin{proof}
It follows from \eqref{uVmGu2C1C2gen} and \eqref{uVmGu2C1C2genV} that 
\begin{EQA}
	\tr (\CONSTIF \Id_{\dimp} + \GP^{2})^{-1} \CONSTIF \CONSTV^{-1}
	& \leq &
	\tr ( \IF_{\GP}^{-1} \VP^{2} )
	\leq 
	\tr (\CONSTIF^{-1} \Id_{\dimp} + \GP^{2})^{-1} \CONSTIF^{-1} \CONSTV \, .
\label{iuhevopedfijwu8efjfw3}
\end{EQA}
Further we may proceed as in the proof of Lemma~\ref{LeffdimG}.
\end{proof}

\Section{Sobolev classes and smooth priors}
This section illustrates the introduced notions and results for a typical situation of 
a \( (\smp,\CGP) \)-smooth prior with \( \gp_{j}^{2} = \CGP j^{2\smp} \).

\Section{Properties of the sub-projector \( \proj_\GP \)}
\label{SsubprojPG}
For the sub-projector \( \proj_{\GP} \) from \eqref{9bvgu7fh54rfgy7byrfuer4}, this section analyzes 
the operator \( \Id_{\dimp} - \proj_{\GP} \) which naturally appears 
in the evaluation of the bias \( \upsvs_{\GP} - \upsvs \); see Section~\ref{Ssmoothbias}.
It turns out that the main characteristic of \( \proj_{\GP} \) is the index \( \mm \) defined by \eqref{gklghi43mjd7skje8vhr}.
The sub-projector \( \proj_{\GP} \) approximates the projector on the space \( \VV_{\mm} \).
The quality of approximation is controlled by the growth rate of the eigenvalues \( \gp_{j}^{2} \):
the faster is this rate the better is the approximation \( \proj_{\GP} \approx \Proj_{\mm} \).
To illustrate this point, we consider the situation with two operators 
\( \Id_{\nsize} - \proj_{\GP} \) and \( \Id_{\nsize} - \proj_{\GP_{0}} \) for two different penalizing matrices 
\( \GP^{2} \) and \( \GP_{0}^{2} \) with the same characteristic \( \mm \).
We slightly change the notations and assume that 
\begin{EQA}
	\GP^{2}
	&=&
	\CGP^{-1} \diag\{ \gp_{1}^{2}, \ldots, \gp_{\dimp}^{2} \},
	\qquad
	\GP_{0}^{2}
	=
	\CGP_{0}^{-1} \diag\{ \gp_{1,0}^{2}, \ldots, \gp_{\dimp,0}^{2} \} ,
\label{946e7f4e346dyeehdy}
\end{EQA}
with some fixed constants \( \CGP \) and \( \CGP_{0} \)
and growing sequences \( (\gp_{j}^{2}) \) and \( (\gp_{j,0}^{2}) \)
satisfying 
\begin{EQA}
	\CGP^{-1} \gp_{\mm}^{2}
	& \approx &
	\CGP_{0}^{-1} \gp_{\mm,0}^{2}
	\approx 
	\nsize .
\label{ywhsywshxetwtsgdew2w}
\end{EQA}
To simplify the presentation we later assume that these relations in \eqref{ywhsywshxetwtsgdew2w} are precisely fulfilled. 
We also assume that 
\begin{EQA}
	\gp_{j,0}^{2}/\gp_{\mm,0}^{2}
	& \leq &
	\gp_{j}^{2}/\gp_{\mm}^{2} \, ,
	\qquad
	j \leq \mm ,
\label{ohcybd7cyr6duyhcyf}
\end{EQA}
meaning that \( \gp_{j,0}^{2} \) grows faster than \( \gp_{j}^{2} \).

\begin{lemma}
\label{Lsmoothbiase}
Let \( \IF \) satisfy \eqref{uVmGu2C1C2gen}, and let \( \GP^{2} \) and \( \GP_{0}^{2} \)
be diagonal penalizing matrices satisfying
\eqref{ywhsywshxetwtsgdew2w} and \eqref{ohcybd7cyr6duyhcyf} for some \( \mm \leq \dimp \).
Then
\begin{EQA}
	(\Id_{\dimp} - \proj_{\GP}) \Proj_{\mm}
	& \leq &
	\CONSTIF^{2} (\Id_{\dimp} - \proj_{\GP_{0}}) \Proj_{\mm} \, ,
	\\
	(\CONSTIF + 1)^{-1} \Proj_{\mm}^{c}
	& \leq &
	(\Id_{\dimp} - \proj_{\GP}) \Proj_{\mm}^{c}
	\leq 
	\Proj_{\mm}^{c} \, ,
\label{8h8h8g5g77hduejwyddyh}
\end{EQA}
yielding
\begin{EQA}
	\Id_{\dimp} - \proj_{\GP}
	& \leq &
	\CONST (\Id_{\dimp} - \proj_{\GP_{0}}) ,
	\qquad
	\CONST = \CONSTIF^{2} \vee (\CONSTIF + 1) .
\label{4t66he5heudshwyseuw}
\end{EQA}
\end{lemma}

\begin{proof}
The definition implies \( \Id_{\dimp} - \proj_{\GP} = (\IF + \GP^{2})^{-1} \GP^{2} \)
and by \eqref{uVmGu2C1C2gen}
\begin{EQA}
	(\CONSTIF \, n \, \Id_{\dimp} + \GP^{2})^{-1} \GP^{2}
	\leq 
	(\IF + \GP^{2})^{-1} \GP^{2}
	& \leq &
	(\CONSTIF^{-1} \, n \, \Id_{\dimp} + \GP^{2})^{-1} \GP^{2} .
\label{6dtetwdtdywywddtetw}
\end{EQA}
Further, for \( j \leq \mm \), \eqref{ywhsywshxetwtsgdew2w} and
\( \gp_{j}^{2}/\gp_{\mm}^{2} \leq \gp_{j,0}^{2}/\gp_{\mm,0}^{2} \) imply
\( \CGP^{-1} \gp_{j}^{2} \leq \CGP_{0}^{-1} \gp_{j,0}^{2} \) that is, 
\begin{EQA}
	\GP^{2} \, \Proj_{\mm}
	& \leq &
	\GP_{0}^{2} \, \Proj_{\mm} \, .
\label{8efdvjdieewwedhfhdvh}
\end{EQA}
Therefore,
\begin{EQA}
	(\Id_{\dimp} - \proj_{\GP}) \Proj_{\mm}
	&=&
	\IF_{\GP}^{-1} \GP^{2} \, \Proj_{\mm}
	\leq 
	(\CONSTIF^{-1} \, n \, \Id_{\dimp} + \GP^{2})^{-1} \GP^{2} \, \Proj_{\mm}
	\\
	& \leq &
	(\CONSTIF^{-1} \, n \, \Id_{\dimp} + \GP_{0}^{2})^{-1} \GP_{0}^{2} \, \Proj_{\mm}
	\leq 
	\CONSTIF^{2} (\IF + \GP_{0}^{2})^{-1} \GP_{0}^{2} \Proj_{\mm} \, .
\label{bdgdigeifgwiegfeydycgs}
\end{EQA}
After restricting to the orthogonal complement \( \VV_{\mm}^{c} \), both operators 
\( \Id_{\dimp} - \proj_{\GP} \) and \( \Id_{\dimp} - \proj_{\GP_{0}} \) behave nearly as 
projectors:
in view of \( \gp_{j}^{2} \geq \nsize \) for \( j > \mm \)
\begin{EQA}
	(\CONSTIF + 1)^{-1} \Proj_{\mm}^{c}
	& \leq &
	(\Id_{\dimp} - \proj_{\GP}) \Proj_{\mm}^{c}
	\leq 
	\Proj_{\mm}^{c}
\label{jhh6dtd5eghdgdetegke5}
\end{EQA}
and similarly for \( \Id_{\dimp} - \proj_{\GP_{0}} \).
\end{proof}

Finally we evaluate the quantity \( \| \QP (\Id - \proj_{\GP}) \upsv \| \)
assuming \( \| \GP_{0} \upsv \| \) bounded.

\begin{lemma}
\label{Lsmoothbias}
It holds for any \( \QP \colon \R^{\dimp} \to \R^{\dimq} \) and any \( \GP^{2} \)
\begin{EQA}
	\| \QP (\Id - \proj_{\GP}) \upsv \|^{2}
	& \leq &
	\| \QP \IF_{\GP}^{-1} \QP^{\T} \| \, \| \GP \upsv \|^{2} .
\label{yc65e5eftw3dtwfefdtw}
\end{EQA}
Moreover, let \( \IF \) satisfy \eqref{uVmGu2C1C2gen}, and let \( \GP^{2} \) and \( \GP_{0}^{2} \)
be diagonal penalizing matrices satisfying
\eqref{ywhsywshxetwtsgdew2w} and \eqref{ohcybd7cyr6duyhcyf} for some \( \mm \leq \dimp \).
Then
\begin{EQA}
	\| \QP (\Id - \proj_{\GP}) \upsv \|
	& \leq &
	\CONST \| \QP \IF_{\GP_{0}}^{-1} \QP^{\T} \|^{1/2} \, \| \GP_{0} \upsv \|
\label{6ew6td366cetyeyegey}
\end{EQA}
with \( \CONST \) from \eqref{4t66he5heudshwyseuw}.
\end{lemma}

\begin{proof}
It holds \( (\Id - \proj_{\GP}) \upsv = \IF_{\GP}^{-1} \GP^{2} \upsv \) and in view of \( \GP^{2} \leq \IF_{\GP} \)
\begin{EQA}
	\| \QP (\Id - \proj_{\GP}) \upsv \|
	& = &
	\| \QP \IF_{\GP}^{-1} \GP^{2} \upsv \|
	\leq 
	\| \QP \IF_{\GP}^{-1/2} \| \, \| \IF_{\GP}^{-1/2} \, \GP^{2} \upsv \|
	\leq 
	\| \QP \IF_{\GP}^{-1} \QP^{\T} \|^{1/2} \, \| \GP \upsv \| .
\label{uchfuehe6fr3ehdjd}
\end{EQA}
For the second statement, we apply \eqref{4t66he5heudshwyseuw} of Lemma~\ref{Lsmoothbiase}.
As \( \IF_{\GP}^{-1} \, \GP^{2} \leq \CONST \IF_{\GPa}^{-1} \, \GPa^{2} \) implies
\( \IF_{\GPa}^{1/2} \, \IF_{\GP}^{-1} \, \GP^{2} \leq \CONST \IF_{\GPa}^{-1/2} \, \GPa^{2} \), it follows in a similar way
\begin{EQA}
	\| \QP \IF_{\GP}^{-1} \GP^{2} \upsv \|
	& \leq &
	\| \QP \IF_{\GPa}^{-1/2} \| \, \| \IF_{\GPa}^{1/2} \, \IF_{\GP}^{-1} \, \GP^{2} \upsv \|
	\\
	& \leq &
	\CONST \| \QP \IF_{\GPa}^{-1/2} \| \, \| \IF_{\GPa}^{-1/2} \, \GPa^{2} \upsv \|
	\leq 
	\CONST \| \QP \IF_{\GPa}^{-1/2} \| \, \| \GPa \upsv \|
\label{uchfuehe6fr3ehdjda}
\end{EQA}
as required.
\end{proof}

For the case of \( \QP = \IF_{\GP_{0}}^{1/2} \), we obtain a corollary of \eqref{6ew6td366cetyeyegey}
\begin{EQA}
	\| \IF_{\GP_{0}}^{1/2} (\Id - \proj_{\GP}) \upsv \|
	& \leq &
	\CONST \| \GP_{0} \upsv \| .
\label{p0p0h9h6f3w33rfhjk}
\end{EQA}

\def\Eg{\E_{\gammav \sim \ND(0,\Id)}}
\def\Pg{\P_{\gammav \sim \ND(0,\Id)}}
\def\HPi{\QP}
\def\dimH{\dimA}
\def\vH{\vA}
\def\BBH{W}
\def\HVB{\mathcal{V}}
\def\mux{\mu_{\xx}}

\Chapter{Deviation bounds for quadratic forms}
\label{Sdevboundgen}
Here we collect some useful results from probability theory mainly concerning 
Gaussian and non-Gaussian quadratic forms.

\Section{Moments of a Gaussian quadratic form}
Let \( \gaussv \) be standard normal in \( \R^{\dimp} \) for \( \dimp \leq \infty \).
Given a self-adjoint trace operator \( \BB \), consider a quadratic form 
\( \bigl\langle \BB \gaussv, \gaussv \bigr\rangle \).

\begin{lemma}
\label{Gaussmoments}
It holds
\begin{EQA}
	\E \bigl\langle \BB \gaussv, \gaussv \bigr\rangle 
	&=& 
	\tr \BB ,
	\\ 
	\Var \bigl\langle \BB \gaussv, \gaussv \bigr\rangle 
	&=& 
	2 \tr \BB^{2} .
\label{EAarAtrA2trA2}
\end{EQA}
Moreover, 
\begin{EQA}
	\E \bigl( \bigl\langle \BB \gaussv, \gaussv \bigr\rangle - \tr \BB \bigr)^{2}
	&=&
	2 \tr \BB^{2}  ,
	\\
	\E \bigl( \bigl\langle \BB \gaussv, \gaussv \bigr\rangle - \tr \BB \bigr)^{3}
	&=&
	8 \tr \BB^{3} ,
	\\
	\E \bigl( \bigl\langle \BB \gaussv, \gaussv \bigr\rangle - \tr \BB \bigr)^{4}
	&=&
	48 \tr \BB^{4} + 12 (\tr \BB^{2})^{2} ,
\label{2pG2trD2DGm22m2}
\end{EQA}
and
\begin{EQA}
	\E \bigl\langle \BB \gaussv, \gaussv \bigr\rangle^{2}
	&=&
	(\tr \BB)^{2} + 2 \tr \BB^{2},
	\\
	\E \bigl\langle \BB \gaussv, \gaussv \bigr\rangle^{3}
	& = &
	(\tr \BB)^{3} + 6 \tr \BB \,\, \tr \BB^{2} + 8 \tr \BB^{3} ,
	\\
	\E \bigl\langle \BB \gaussv, \gaussv \bigr\rangle^{4}
	& = &
	(\tr \BB)^{4} + 12 (\tr \BB)^{2} \tr \BB^{2}
	+ 32 (\tr \BB) \tr \BB^{3}
	+ 48 \tr \BB^{4} + 12 (\tr \BB^{2})^{2} ,
\label{2pG2trD2DGm22m2}
	\\
	\Var \bigl\langle \BB \gaussv, \gaussv \bigr\rangle^{2}
	& = &
	8 (\tr \BB)^{2} \tr \BB^{2}
	+ 32 (\tr \BB) \tr \BB^{3}
	+ 48 \tr \BB^{4} + 8 (\tr \BB^{2})^{2} .
\label{2pG2trD2DGm22m4}
\end{EQA}
Moreover, if \( \BB \leq \Id_{\dimp} \) and \( \dimA = \tr \BB \), then \( \tr \BB^{m} \leq \dimA \) for 
\( m \geq 1 \) and
\begin{EQA}[rcccl]
	\E \bigl\langle \BB \gaussv, \gaussv \bigr\rangle^{2}
	& \leq &
	\dimA^{2} + 2 \dimA
	&\leq &
	(\dimA + 1)^{2},
	\\
	\E \bigl\langle \BB \gaussv, \gaussv \bigr\rangle^{3}
	& \leq &
	\dimA^{3} + 6 \dimA^{2} + 8 \dimA 
	&\leq &
	(\dimA + 2)^{3},
	\\
	\E \bigl\langle \BB \gaussv, \gaussv \bigr\rangle^{4}
	& \leq &
	\dimA^{4} + 12 \dimA^{3} 
	+ 44 \dimA^{2}
	+ 48 \dimA  
	&\leq &
	(\dimA + 3)^{4},
\label{2pG2trD2DGm22m2}
	\\
	\Var \bigl\langle \BB \gaussv, \gaussv \bigr\rangle^{2}
	& \leq &
	8 \dimA^{3} + 40 \dimA^{2} + 48 \dimA .
\label{2pG2trD2DGm22m4}
\end{EQA}
\end{lemma}

\begin{proof}
Let \( \chi = \gauss^{2} - 1 \) for \( \gauss \) standard normal.
Then \( \E \chi = 0 \), \( \E \chi^{2} = 2 \), \( \E \chi^{3} = 8 \), \( \E \chi^{4} = 60 \).
Without loss of generality assume \( \BB \) diagonal: \( \BB = \diag(\lambda_{1},\lambda_{2},\ldots,\lambda_{\dimp}) \).
Then 
\begin{EQA}
	\xi
	\eqdef
	\bigl\langle \BB \gaussv, \gaussv \bigr\rangle - \tr \BB
	&=&
	\sum_{j=1}^{\dimp} \lambda_{j} (\gauss_{j}^{2} - 1) ,
\label{j1ljgj2m1}
\end{EQA}
where \( \gauss_{j} \) are i.i.d. standard normal. 
This easily yields
\begin{EQA}
	\E \xi^{2}
	&=&
	\sum_{j=1}^{\dimp} \lambda_{j}^{2} \E (\gauss_{j}^{2} - 1)^{2}
	=
	\E \chi^{2} \, \tr \BB^{2} 
	=
	2 \tr \BB^{2}  ,
	\\
	\E \xi^{3}
	&=&
	\sum_{j=1}^{\dimp} \lambda_{j}^{3} \E (\gauss_{j}^{2} - 1)^{3}
	=
	\E \chi^{3} \, \tr \BB^{3} 
	=
	8 \tr \BB^{3} ,
	\\
	\E \xi^{4}
	&=&
	\sum_{j=1}^{\dimp} \lambda_{j}^{4} (\gauss_{j}^{2} - 1)^{4}
	+ \sum_{i\neq j} \lambda_{i}^{2} \lambda_{j}^{2} \E (\gauss_{i}^{2} - 1)^{2} \E (\gauss_{j}^{2} - 1)^{2}
	\\
	&=&
	\bigl( \E \chi^{4} - 3 (\E \chi^{2})^{2} \bigr) \tr \BB^{4} + 3 (\E \chi^{2} \, \tr \BB^{2})^{2}
	=
	48 \tr \BB^{4} + 12 (\tr \BB^{2})^{2} ,
\label{2pG2trD2DGm22m2}
\end{EQA}
ensuring
\begin{EQA}
	\E \bigl\langle \BB \gaussv, \gaussv \bigr\rangle^{2}
	&=&
	\bigl( \E \bigl\langle \BB \gaussv, \gaussv \bigr\rangle \bigr)^{2} 
	+ \E \xi^{2}
	= 
	(\tr \BB)^{2} + 2 \tr \BB^{2},
	\\
	\E \bigl\langle \BB \gaussv, \gaussv \bigr\rangle^{3}
	& = &
	\E \bigl( \xi + \tr \BB \bigr)^{3}
	=
	(\tr \BB)^{3} + \E \xi^{3}
	+ 3 \tr \BB \,\, \E \xi^{2}
	\\
	&=&
	(\tr \BB)^{3} + 6 \tr \BB \,\, \tr \BB^{2} + 8 \tr \BB^{3} ,
\label{2pG2trD2DGm22m2}
\end{EQA}
and 
\begin{EQA}
	\Var \bigl\langle \BB \gaussv, \gaussv \bigr\rangle^{2}
	& = &
	\E \bigl( \xi + \tr \BB \bigr)^{4}
	- \bigl( \E \bigl\langle \BB \gaussv, \gaussv \bigr\rangle \bigr)^{2}
	\\
	&=&
	\bigl( \tr \BB \bigr)^{4} + 6 (\tr \BB)^{2} \E \xi^{2} + 4 \tr \BB \, \E \xi^{3} + \E \xi^{4}
	- \bigl( (\tr \BB)^{2} + 2 \tr \BB^{2} \bigr)^{2}
	\\
	&=& 
	8 (\tr \BB)^{2} \tr \BB^{2}
	+ 32 (\tr \BB) \tr \BB^{3}
	+ 48 \tr \BB^{4} + 8 (\tr \BB^{2})^{2} .
\label{2pG2trD2DGm22m4}
\end{EQA}
This implies the results of the lemma.
\end{proof}

Now we compute the exponential moments of centered and non-centered quadratic forms.

\begin{lemma}
\label{Lqfexpmom}
Let \( \| \BB \|_{\oper} \leq 1 \) and \( \gaussv \sim \ND(0,\Id_{\dimp}) \).
Then for any \( \mu \in (0,1) \), 
\begin{EQA}
	\E \exp \Bigl\{ \frac{\mu}{2} \bigl( \langle \BB \gaussv, \gaussv \rangle - \dimA \bigr) \Bigr\}
	&=&
	\det(\Id - \mu \BB)^{-1/2} \, .
\label{m2v241m41m}
\end{EQA}
Moreover, with \( \dimA = \tr \BB \) and \( \vA^{2} = \tr \BB^{2} \)
\begin{EQA}
	\log \E \exp \Bigl\{ \frac{\mu}{2} \bigl( \langle \BB \gaussv, \gaussv \rangle - \dimA \bigr) \Bigr\}
	& \leq &
	\frac{\mu^{2} \vA^{2}}{4 (1 - \mu)} \, .
\label{m2v241m41mb}
\end{EQA}
If \( \BB \) is positive semidefinite, \( \lambda_{j} \geq 0 \), then 
\begin{EQA}
	\log \E \exp \Bigl\{ - \frac{\mu}{2} \bigl( \langle \BB \gaussv, \gaussv \rangle - \dimA \bigr) \Bigr\}
	& \leq &
	\frac{\mu^{2} \vA^{2}}{4} \, .
\label{m2v241m41mbn}
\end{EQA}
\end{lemma}

\begin{proof}
Let \( \lambda_{j} \) be the eigenvalues of \( \BB \), 
\( |\lambda_{j}| \leq 1 \).
By an orthogonal transform, one can reduce the statement to the case of a diagonal matrix 
\( \BB = \diag\bigl( \lambda_{j} \bigr) \). 
Then \( \langle \BB \gaussv, \gaussv \rangle = \sum_{j=1}^{\dimp} \lambda_{j} \gauss_{j}^{2} \) and 
by independence of the \( \gauss_{j} \)'s
\begin{EQA}
	&& \nquad
	\E \Bigl\{ \frac{\mu}{2} \langle \BB \gaussv, \gaussv \rangle  \Bigr\}
	=
	\prod_{j=1}^{\dimp} \E \exp \Bigl( \frac{\mu}{2} \lambda_{j} \eps_{j}^{2} \Bigr)
	=
	\prod_{j=1}^{\dimp} \frac{1}{\sqrt{1 - \mu \lambda_{j}}} 
	=
	\det \bigl( \Id - \mu \BB \bigr)^{-1/2} .
\label{dOImuBm12EB}
\end{EQA}
Below we use the simple bound: 
\begin{EQ}[rcl]
\label{lo1uusk2iukkp}
	- \log(1 - u) - u
	&=&
	\sum_{k=2}^{\infty} \frac{u^{k}}{k}
	\leq 
	\frac{u^{2}}{2} \sum_{k=0}^{\infty} u^{k} 
	=
	\frac{u^{2}}{2 (1 - u)} \, ,
	\qquad 
	u \in (0,1),
	\\
	- \log(1 - u) + u
	&=&
	\sum_{k=2}^{\infty} \frac{u^{k}}{k}
	\leq 
	\frac{u^{2}}{2} \, ,
	\qquad \qquad
	u \in (-1,0).
\label{lo1uusk2iukk}
\end{EQ}
Now it holds 
\begin{EQA}
	&& \nquad
	\log \E \Bigl\{ \frac{\mu}{2} \bigl( \langle \BB \gaussv, \gaussv \rangle - \dimA \bigr) \Bigr\}
	=
	\log \det(\Id - \mu \BB)^{-1/2} - \frac{\mu \, \dimA}{2}
	\\
	&=&
	- \frac{1}{2} \sum_{j=1}^{\dimp} \bigl\{ \log(1 - \mu \lambda_{j}) + \mu \lambda_{j} \bigr\}
	\leq 
	\sum_{j=1}^{\dimp} \frac{\mu^{2} \lambda_{j}^{2}}{4 (1 - \mu)} 
	=
	\frac{\mu^{2} \vA^{2}}{4 (1 - \mu)} \, .
\label{m2v241m4mj1pd}
\end{EQA}
The last statement can be proved similarly.
\end{proof}

Now we consider the case of a non-centered quadratic form
\( \langle \BB \gaussv,\gaussv \rangle/2 + \langle \Av,\gaussv \rangle \) for a fixed vector \( \Av \).

\begin{lemma}
\label{Lexpmomnoncen}
Let \( \lambda_{\max}(\BB) < 1 \). 
Then for any \( \Av \)
\begin{EQA}
	\E \exp\Bigl\{ \frac{1}{2}\langle \BB \gaussv,\gaussv \rangle + \langle \Av,\gaussv \rangle \Bigr\}
	&=&
	\exp\Bigl\{ \frac{\| (\Id - \BB)^{-1/2} \Av \|^{2}}{2} \Bigr\} \, \det(\Id - \BB)^{-1/2} .
\label{EeBf12BggA}
\end{EQA}
Moreover, for any \( \mu \in (0,1) \)
\begin{EQA}
	&& \nquad
	\log \E \exp\Bigl\{ 
		\frac{\mu}{2} \bigl( \langle \BB \gaussv,\gaussv \rangle - \dimA \bigr) + \langle \Av,\gaussv \rangle 
	\Bigr\}
	\\
	&=&
	\frac{\| (\Id - \mu \BB)^{-1/2} \Av \|^{2}}{2} + \log \det(\Id - \mu \BB)^{-1/2} - \mu \dimA 
	\\
	& \leq &
	\frac{\| (\Id - \mu \BB)^{-1/2} \Av \|^{2}}{2} + \frac{\mu^{2} \vA^{2}}{4 (1 - \mu)} \, .
\label{EeBf12BggAmu}
\end{EQA}
\end{lemma}

\begin{proof}
Denote \( \av = (\Id - \BB)^{-1/2} \Av \). 
It holds by change of variables \( (\Id - \BB)^{1/2} \xv = \uv \) for \( \CONSTi_{\dimp} = (2\pi)^{-\dimp/2} \)
\begin{EQA}
	&& \nquad
	\E \exp\Bigl\{ \frac{1}{2}\langle \BB \gaussv,\gaussv \rangle + \langle \Av,\gaussv \rangle \Bigr\}
	=
	\CONSTi_{\dimp}
	\int \exp\Bigl\{ - \frac{1}{2}\langle (\Id - \BB) \xv,\xv \rangle + \langle \Av,\xv \rangle \Bigr\} d\xv
	\\
	&=&
	\CONSTi_{\dimp}
	\det(\Id - \BB)^{-1/2}
	\int \exp\Bigl\{ - \frac{1}{2} \| \uv \|^{2} + \langle \av,\uv \rangle \Bigr\} d\uv
	=
	\det(\Id - \BB)^{-1/2} \, 	\ex^{\| \av \|^{2}/2}  	.
\label{EeBf12BggAp}
\end{EQA}
The last inequality \eqref{EeBf12BggAmu} follows by \eqref{m2v241m41mb}.
\end{proof}

\Section{Deviation bounds for Gaussian quadratic forms}
\label{SdevboundGauss}
The next result explains the concentration effect of 
\( \langle \BB \xiv, \xiv \rangle \)
for a centered Gaussian vector \( \xiv \sim \ND(0,\HVB^{2}) \) and a symmetric trace operator \( \BB \) in \( \R^{\dimp} \),
\( \dimp \leq \infty \).
We use a version from \cite{laurentmassart2000}.
For completeness, we present a simple proof of the upper bound.

\begin{theorem}
\label{TexpbLGA}
\label{Lxiv2LD}
\label{Cuvepsuv0}
Let \( \xiv \sim \ND(0,\HVB^{2}) \) be a Gaussian element in \( \R^{\dimp} \) and \( \BB \) be symmetric 
such that \( \BBH = \HVB \BB \HVB \) is a trace operator in \( \R^{\dimp} \).
Then with \( \dimH = \tr(\BBH) \), \( \vH^{2} = \tr(\BBH^{2}) \), and 
\( \supA = \| \BBH \| \), it holds for each \( \xx \geq 0 \)
\begin{EQA}
\label{Pxiv2dimAvp12}
	\P\Bigl( \langle \BB \xiv, \xiv \rangle - \dimH > 2 \vH \, \sqrt{\xx} + 2 \supA \xx \Bigr)
	& \leq &
	\ex^{-\xx} .
\end{EQA}
It also implies 
\begin{EQA}
	\P\bigl( \bigl| \langle \BB \xiv, \xiv \rangle - \dimH \bigr| > \zq_{2}(\BBH,\xx) \bigr)
	& \leq &
	2 \ex^{-\xx} ,
\label{PxivTBBdimA2vp}
\end{EQA}
with
\begin{EQA}
	\zq_{2}(\BBH,\xx)
	& \eqdef &
	2 \vH \, \sqrt{\xx} + 2 \supA \xx \,\, .
\label{zqdefGQF}
\end{EQA}
%
\end{theorem}

\begin{proof}
W.l.o.g. assume that \( \supA = \| \BBH \| = 1 \).
We use the identity \( \langle \BB \xiv, \xiv \rangle = \langle \BBH \gaussv, \gaussv \rangle \) with
 \( \gaussv \sim \ND(0,\Id_{\dimp}) \).
We apply the exponential Chebyshev inequality: with \( \mu > 0 \)
\begin{EQA}
	\P\Bigl( \langle \BBH \gaussv, \gaussv \rangle > \zq^{2} \Bigr)
	& \leq &
	\E \exp \Bigl( \mu \langle \BBH \gaussv, \gaussv \rangle / 2 - \mu \zq^{2} / 2 \Bigr) \, .
\label{PBggiz2E2mz2}
\end{EQA}
Given \( \xx > 0 \), fix \( \mu < 1 \) by the equation
\begin{EQA}
	\frac{\mu}{1 - \mu} 
	&=&
	\frac{2 \sqrt{\xx}}{\vH} \, 
	\quad \text{ or } \quad
	\mu^{-1} 
	=
	1 + \frac{\vH}{2 \sqrt{\xx}} \, .
\label{1v2sxm12m1m}
\end{EQA}
Let \( \lambda_{j} \) be the eigenvalues of \( \BBH \), 
\( |\lambda_{j}| \leq 1 \).
It holds with \( \dimH = \tr \BBH \) in view of \eqref{m2v241m41mb}
\begin{EQA}
	&& \nquad
	\log \E \Bigl\{ \frac{\mu}{2} \bigl( \langle \BBH \gaussv, \gaussv \rangle - \dimH \bigr) \Bigr\}
	\leq 
	\frac{\mu^{2} \vH^{2}}{4 (1 - \mu)} \, .
\label{m2v241m4mj1p}
\end{EQA}
It remains to check that the choice \( \mu \) by \eqref{1v2sxm12m1m} and 
\( \zq = \zq(\BBH,\xx) \) yields
\begin{EQA}
	\frac{\mu^{2} \vH^{2}}{4 (1 - \mu)} - \frac{\mu (\zq^{2} - \dimH)}{2}
	& = &
	\frac{\mu^{2} \vH^{2}}{4 (1 - \mu)} - \mu \bigl( \vH \sqrt{\xx} + \xx \bigr)
	=
	\mu \Bigl( \frac{\vH \sqrt{\xx}}{2} - \vH \sqrt{\xx} - \xx \Bigr)
	=
	- \xx
	\qquad
	\qquad
\label{m2vA241muz2}
\end{EQA}
as required.
The last statement \eqref{PxivTBBdimA2vp} is obtained by applying this inequality twice to \( \BBH \) and \( - \BBH \).
\end{proof}

\begin{corollary}
\label{CTexpbLGAd}
Assume the conditions of Theorem~\ref{TexpbLGA}.
Then for \( \zq > \vH \)
\begin{EQA}
	\P\bigl( \bigl| \langle \BB \xiv, \xiv \rangle - \dimH \bigr| \ge \zq \bigr)
	& \leq &
	2 \exp\biggl\{ - \frac{\zq^{2}}{\bigl( \vH + \sqrt{\vH^{2} + 2 \supA \zq} \bigr)^{2}} \biggr\}
	\leq 
	2 \exp\biggl( - \frac{\zq^{2}}{4\vH^{2} + 4 \supA \zq} \biggr) .
	\qquad
\label{3z2spsp2z3z2}
\end{EQA}
\end{corollary}

\begin{proof}
Given \( \zq \), define \( \xx \) by 
\( 2 \vH \sqrt{\xx} + 2 \supA \xx = \zq \) or 
\( 2 \supA \sqrt{\xx} = \sqrt{\vH^{2} + 2 \supA \zq} - \vH \).
Then
\begin{EQA}
	\P\bigl( \langle \BB \xiv, \xiv \rangle - \dimH \ge \zq \bigr)
	& \leq &
	\ex^{-\xx} 
	=
	\exp\biggl\{ - \frac{\bigl( \sqrt{\vH^{2} + 2 \supA \zq} - \vH \bigr)^{2}}{4 \supA^{2}} \biggr\}
	=
	\exp\biggl\{ - \frac{\zq^{2}}{\bigl( \vH + \sqrt{\vH^{2} + 2 \supA \zq} \bigr)^{2}} \biggr\}.
\label{3emzmsp22z2c}
\end{EQA}
This yields \eqref{3z2spsp2z3z2} by direct calculus.
\end{proof}

Of course, bound \eqref{3z2spsp2z3z2} is sensible only if \( \zq \gg \vH \).

\begin{corollary}
\label{RsochpHsA}
Assume the conditions of Theorem~\ref{TexpbLGA}.
If also \( \BB \geq 0 \), then 
\begin{EQA}
\label{Pxiv2dimAxx12}
	\P\Bigl( \langle \BB \xiv, \xiv \rangle \geq \zq^{2}(\BB,\xx) \Bigr)
	& = &
	\P\bigl( \| \BB^{1/2} \xiv \| \geq \zq(\BB,\xx) \bigr)
	\leq 
	\ex^{-\xx} 
\end{EQA}
with 
\begin{EQA}
	\zq(\BB,\xx)
	& \eqdef &
	\sqrt{\dimH + 2 \vH \, \sqrt{\xx} + 2 \supA \xx} 
	\leq 
	\sqrt{\dimH} + \sqrt{2 \supA \xx} \, .
\label{zzxxppdBlroBB}
\end{EQA}
Also
\begin{EQA}
	\P\Bigl( \langle \BB \xiv, \xiv \rangle - \dimH < - 2 \vH \, \sqrt{\xx} \Bigr)
	& \leq &
	\ex^{-\xx} .
\label{Pxiv2dimAvp12d}
\end{EQA}
\end{corollary}

\begin{proof}
The definition implies \( \vH^{2} \leq \dimH \supA \).
One can use a sub-optimal choice of the value 
\( \mu(\xx) = \bigl\{ 1 + 2 \sqrt{\supA \dimH/\xx} \bigr\}^{-1} \) yielding the statement of the corollary.
\end{proof}

As a special case, we present a bound for the chi-squared distribution 
corresponding to \( \BB = \HVB^{2} = \Id_{\dimp} \), \( \dimp < \infty \).
Then \( \tr (\BBH) = \dimp \), \( \tr(\BBH^{2}) = \dimp \) and \( \supA(\BBH) = 1 \).

\begin{corollary}
\label{Cchi2p}
Let \( \gaussv \) be a standard normal vector in \( \R^{\dimp} \).
Then for any \( \xx > 0 \)
\begin{EQA}[ccl]
\label{Pxi2pm2px}
	\P\bigl( \| \gaussv \|^{2} \geq \dimp + 2 \sqrt{\dimp \, \xx} + 2 \xx \bigr)
	& \leq &
	\ex^{-\xx},
	\\
	\P\bigl( \| \gaussv \| \,\,  \geq \sqrt{\dimp} + \sqrt{2 \xx} \bigr)
	& \leq &
	\ex^{-\xx} ,
\label{Pxi2pm2px12}
	\\
	\P\bigl( \| \gaussv \|^{2} \leq \dimp - 2 \sqrt{\dimp \, \xx} \bigr)
	& \leq &
	\ex^{-\xx}	.
\label{Pxi2pm2px22}
\end{EQA}
\end{corollary}

The bound of Theorem~\ref{TexpbLGA} 
can be represented as a usual deviation bound.

\begin{theorem}
\label{CTexpbLGA}
Assume the conditions of Theorem~\ref{TexpbLGA} with \( \BB \geq 0 \).
Then for \( \zq > \sqrt{\dimH} + 1 \)
\begin{EQA}
	\P\bigl( \langle \BB \xiv, \xiv \rangle \ge \zq^{2} \bigr)
	& \leq &
	\exp\Bigl\{ - \frac{(\zq - \sqrt{\dimH})^{2}}{2 \supA} \Bigr\} ,
\label{3emzmsp22z2}
	\\
	\E \bigl\{ \langle \BB \xiv, \xiv \rangle^{1/2} \Ind\bigl( \langle \BB \xiv, \xiv \rangle \ge \zq^{2} \bigr) \bigr\}
	& \leq &
	\exp\Bigl\{ - \frac{(\zq - \sqrt{\dimH})^{2}}{2 \supA} \Bigr\} ,
\label{3emzmsp22z21}
	\\
	\E \bigl\{ \langle \BB \xiv, \xiv \rangle \Ind\bigl( \langle \BB \xiv, \xiv \rangle \ge \zq^{2} \bigr) \bigr\}
	& \leq &
	\frac{2 \zq}{\zq - \sqrt{\dimH}} \exp\Bigl\{ - \frac{(\zq - \sqrt{\dimH})^{2}}{2 \supA} \Bigr\} .
\label{3emzmsp22z2e}
\end{EQA}
\end{theorem}

\begin{proof}
Bound \eqref{3emzmsp22z2} follows from 
\eqref{Pxiv2dimAxx12}.
It obviously suffices to check the bound for the excess risk for \( \supA = 1 \).
It follows with \( \eta = \| \BB^{1/2} \xiv \| \) for \( \zq \geq \sqrt{\dimH} + 1 \)
\begin{EQA}
	\E \bigl\{ \eta \Ind(\eta > \zq) \bigr\}
	&=&
	\int_{\zq}^{\infty} \P(\eta \geq \zq) \, d\zq
	\leq 
	\int_{\zq}^{\infty} \exp\bigl\{ - \frac{(x - \sqrt{\dimH})^{2}}{2} \bigr\} \, dx
	\leq 
	\exp\bigl\{ - \frac{(\zq - \sqrt{\dimH})^{2}}{2} \bigr\}.
\label{zEe2Iezz2c2H23}
\end{EQA} 
Similarly
\begin{EQA}
	\E \bigl\{ \eta^{2} \Ind(\eta^{2} > \zq^{2}) \bigr\}
	&=&
	\int_{\zq^{2}}^{\infty} \P(\eta^{2} \geq \zz) \, d\zz
	\leq 
	\int_{\zq^{2}}^{\infty} \exp\bigl\{ - \frac{(\sqrt{\zz} - \sqrt{\dimH})^{2}}{2} \bigr\} \, d\zz .
\label{zEe2Iezz2c2H23d}
\end{EQA} 
By change of variables \( \sqrt{\zz} - \sqrt{\dimH} = u \) for \( \zq > \sqrt{\dimH} + 1 \)
\begin{EQA}
	&& \nquad
	\int_{\zq^{2}}^{\infty} \exp\bigl\{ - \frac{(\sqrt{\zz} - \sqrt{\dimH})^{2}}{2} \bigr\} \, d\zz
	\leq 
	2 \int_{\zq - \sqrt{\dimH}}^{\infty} (u + \sqrt{\dimH}) \, \exp\{ - u^{2}/2 \} \, du
	\\
	& \leq &
	2 \left( 1 + \frac{\sqrt{\dimH}}{\zq - \sqrt{\dimH}} \right) 
	\exp\bigl\{ - (\zq - \sqrt{\dimH})^{2}/2 \bigr\} 
	=
	\frac{2 \zq}{\zq - \sqrt{\dimH}} \exp\bigl\{ - (\zq - \sqrt{\dimH})^{2}/2 \bigr\}\, .
\label{dz21fsHzsHe22}
\end{EQA}
\end{proof}

\Section{Deviation bounds for non-Gaussian quadratic forms}
\label{Sprobabquad}
\label{SdevboundnonGauss}
This section collects some probability bounds for non-Gaussian quadratic forms
starting from the subgaussian case.
Then we extend the result to the case of exponential tails. 
%
Let \( \xiv \) be a random vector in \( \R^{\dimp} \), \( \dimp \leq \infty \)
satisfying \( \E \xiv = 0 \).
We suppose that there exists an operator \( \HVB \) in \( \R^{\dimp} \) such that
\begin{EQA}
	\log \E \exp \bigl( \langle \uv, \HVB^{-1} \xiv \rangle \bigr)
	& \leq &
	\frac{\| \uv \|^{2}}{2} \, ,
	\qquad 
	\uv \in \R^{\dimp} .
\label{devboundinf}
\end{EQA}
In the Gaussian case, one obviously takes \( \HVB^{2} = \Var(\xiv) \).
In general, \( \HVB^{2} \geq \Var(\xiv) \).
We consider a quadratic form \( \langle \BB \xiv, \xiv \rangle \), where 
\( \xiv \) satisfies \eqref{devboundinf} and \( \BB \) is a given symmetric non-negative 
operator in \( \R^{\dimp} \) such that \( \BB \leq \HVB^{-2} \) and
\( \BBH = \HVB \BB \HVB \) is a trace operator:
\begin{EQA}
	\dimH
	&=&
	\tr\bigl( \BBH \bigr) 
	< \infty.
\label{dimAtHm2Bpp}
\end{EQA}
Denote also 
\begin{EQA}[c]
    \vH^{2}
    \eqdef
    \tr(\BBH^{2}) .
\label{BBrddB}
\end{EQA}   
We show that under these conditions, the quadratic form \( \langle \BB \xiv, \xiv \rangle \)
follows the same deviation bound 
\( \P\bigl( \langle \BB \xiv, \xiv \rangle \ge \zq^{2}(\BBH,\xx) \bigr) \leq \ex^{-\xx} \) 
with \( \zq^{2}(\BBH,\xx) \) from \eqref{zqdefGQF}
as in the Gaussian case.

\begin{theorem}
\label{Tdevboundinf}
Suppose \eqref{devboundinf}. 
Let \( \dimH = \tr \BBH < \infty \) for \( \BBH = \HVB \BB \HVB \).
Then
\begin{EQA}
	\P\bigl( \langle \BB \xiv, \xiv \rangle > \zq^{2}(\BBH,\xx) \bigr)
	& \leq &
	\ex^{-\xx} .
\label{PxivbzzBBroBinf}
\end{EQA}
The bounds \eqref{3emzmsp22z2} through \eqref{3emzmsp22z2e} of Theorem~\ref{CTexpbLGA} continue to apply as well.
\end{theorem}

\begin{proof}
For any \( \mu < 1 \), we use the identity
\begin{EQA}
	\E \exp\bigl( \mu \langle \BB \xiv, \xiv \rangle / 2 \bigr)
	&=&
	\E \Eg \exp\bigl( \mu^{1/2} \langle \HVB \BB^{1/2} \gammav, \HVB^{-1} \xiv \rangle \bigr)
\label{Egexmu12lHm1B12}
\end{EQA}
Application of Fubini's theorem and \eqref{devboundinf} yields
\begin{EQA}
	\E \exp\bigl( \mu \langle \BB \xiv, \xiv \rangle / 2 \bigr)
	& \leq &
	\exp \Bigl( \frac{\mu^{2} \tr \BBH^{2}}{4 (1 - \mu)} + \frac{\mu \tr \BBH}{2} \Bigr) .
\label{wBmu2vA241mmutrBi}
\end{EQA}
Further we proceed as in the Gaussian case.
\end{proof}

Now we turn to the main case of light exponential tails of \( \xiv \).
Namely, we suppose that \( \E \xiv = 0 \) and 
for some fixed \( \gmb > 0 \) 
\begin{EQA}[c]
    \log \E \exp\bigl( \langle \uv, \HVB^{-1} \xiv \rangle \bigr)
    \le
    \frac{\| \uv \|^{2}}{2} \, ,
    \qquad
    \uv \in \R^{\dimp}, \, \| \uv \| \le \gmb ,
\label{expgamgm}
\end{EQA}
for some self-adjoint operator \( \HVB \) in \( \R^{\dimp} \), \( \HVB \geq \Id_{\dimp} \).
In fact, it suffices to assume that 
\begin{EQA}
	\sup_{\| \uv \| \leq \gmb} \E \exp\bigl( \langle \uv, \HVB^{-1} \xiv \rangle \bigr)
	& \leq &
	\CONST .
\label{sgagEexlgHx}
\end{EQA}
Then one can use the fact that existence of the exponential moment \( \E \ex^{\lambda_{0} \xi} \) for a centered random variable \( \xi \) and some fixed \( \lambda_{0} \) implies that the moment generating function 
\( f_{\xi}(\lambda) \eqdef \log \E \ex^{\lambda \xi} \) is analytic in \( \lambda \in (0, \lambda_{0}) \) with 
\( f_{\xi}(0) = f'_{\xi}(0) = 0 \) and hence, it can be well majorated by a quadratic function in a smaller interval 
\( [0,\lambda_{1}] \) for \( \lambda_{1} < \lambda_{0} \); see~\cite{GolSpo2009}.

Remind the notation \( \BBH = \HVB \BB \HVB \).
By normalization, one can easily reduce the study to the case \( \| \BBH \| = 1 \).
Let \( \dimH = \tr(\BBH) \), \( \vH^{2} = \tr(\BBH^{2}) \), and
\( \mu(\xx) \) be defined by \( \mu(\xx) = \bigl( 1 + \frac{\vH}{2 \sqrt{\xx}} \bigr)^{-1} \);
see \eqref{1v2sxm12m1m}. 
Obviously \( \mu(\xx) \) grows with \( \xx \).
Define the value \( \xxc \) as the root of the equation
\begin{EQA}
	\frac{\gmb - \sqrt{\dimH \, \mu(\xx)}}{\mu(\xx)} 
	& = &
	\zq(\BBH,\xx) + 1 .
\label{gmsqpAHmxzBx1}
\end{EQA}
The left hand-side here decreases with \( \xx \), while the right hand-side is increasing in \( \xx \) to infinity.
Therefore, the solution exists and is unique.
Also denote \( \muc = \mu(\xxc) \) and
\begin{EQA}
	\gmc
	&=&
	\gmb - \sqrt{\dimH \muc} \, ,
\label{gcgbsqpHmc}
\end{EQA}
so that
\begin{EQA}
	\gmc/\muc
	& = &
	\zq(\BBH,\xxc) + 1 .
\label{sxcgbspHspA}
\end{EQA}

\begin{theorem}
\label{Tdevboundgm}
\label{LLbrevelocroB}   
Let \eqref{expgamgm} hold and let \( \BB \) be such that \( \BBH = \HVB \BB \HVB \)
satisfies \( \| \BBH \| = 1 \)
and \( \dimH = \tr(\BBH) < \infty \).
Define \( \xxc \) by \eqref{gmsqpAHmxzBx1} and \( \gmc \) by \eqref{gcgbsqpHmc}, and suppose
\( \gmc \geq 1 \).
Then for any \( \xx > 0 \)
\begin{EQA}
    \P\bigl( \langle \BB \xiv, \xiv \rangle \ge \zqc^{2}(\BBH,\xx) \bigr)
    & \le &
    2 \ex^{-\xx} + \ex^{-\xxc} \Ind(\xx < \xxc) 
    \leq 
    3\ex^{-\xx},
\label{PxivbzzBBroB}
\end{EQA}    
where \( \zqc(\BBH,\xx) \) is defined by
\begin{EQA}
\label{PzzxxpBroB}
    \zqc(\BBH,\xx)
    & \eqdef &
    \begin{cases}
      \sqrt{ \dimH + 2 \vH \, \xx^{1/2} + 2 \xx } \, , &  \xx \le \xxc \, , 
      \\
      \gmc/\muc + 2 (\xx - \xxc)/\gmc \, , & \xx > \xxc \, ,
    \end{cases}
    \\
    & \leq &
    \begin{cases}
      \sqrt{\dimH} + \sqrt{2 \xx} \, , &  \xx \le \xxc \, , 
      \\
      \gmc/\muc + 2 (\xx - \xxc)/\gmc \, , & \xx > \xxc \, .
    \end{cases}
\label{zzxxppdBlroB}
\end{EQA} 
Moreover, if, given \( \xx \), it holds
\begin{EQA}
	\gm
	& \geq &
	\xx^{1/2}/2 + (\dimH \xx/4)^{1/4} ,
\label{g2x1214px14}
\end{EQA}
then
\begin{EQA}
    \P\bigl( \| \BB^{1/2} \xiv \| \ge \sqrt{\dimH} + \sqrt{2\xx} \bigr)
    & \le &
    3 \ex^{-\xx} .
\label{PxivbzzBBroB3}
\end{EQA}    
\end{theorem}


\begin{remark}
Depending on the value \( \xx \), we have two types of tail behavior of the 
quadratic form \( \langle \BB \xiv, \xiv \rangle \). 
For \( \xx \le \xxc \), we have essentially the same deviation bounds as in the Gaussian case
with the extra-factor two in the deviation probability.
For \( \xx > \xxc \), we switch to the special regime driven by the exponential moment
condition \eqref{expgamgm}.
Usually \( \gmb^{2} \) is a large number (of order \( n \) in the i.i.d. setup)
yielding \( \xxc \) also large, and the second term in \eqref{PxivbzzBBroB} can be simply ignored. 
The function \( \zqc(\BBH,\xx) \) is discontinuous at the point \( \xxc \).
Indeed, \( \zqc(\BBH,\xx) = \zq(\BBH,\xx) \) for \( \xx < \xxc \), while by \eqref{gmsqpAHmxzBx1}, 
it holds \( \gmc/\muc = \zq(\BBH,\xxc) + 1 \).
However, the jump at \( \xxc \) is at most one. 
\end{remark}

As a corollary, we state the result for the norm of \( \xiv \in \R^{\dimp} \) corresponding to the case
\( \HVB^{-2} = \BB = \Id_{\dimp} \) and \( \dimp < \infty \). 
Then 
\begin{EQA}
	\dimA
	&=&
	\vH^{2}
	=
	\dimp .
\label{pApHpvA2l1}
\end{EQA}

\begin{corollary}
\label{LLbrevelocro}   
Let \eqref{expgamgm} hold with \( \HVB = \Id_{\dimp} \).
Then for each \( \xx > 0 \)
\begin{EQA}
    \P\bigl( \| \xiv \| \ge \zqc(\dimp,\xx) \bigr)
    & \le &
    2 \ex^{-\xx} + \ex^{-\xxc } \Ind(\xx < \xxc) ,
\label{PxivbzzBBro}
\end{EQA}    
where \( \zqc(\dimp,\xx) \) is defined by
\begin{EQA}
\label{PzzxxpBro}
    \zqc(\dimp,\xx)
    & \eqdef &
    \begin{cases}
      \bigl( \dimp + 2 \sqrt{\dimp \, \xx} + 2 \xx\bigr)^{1/2}, &  \xx \le \xxc  , \\
      \gmc/\muc + 2 \gmc^{-1} (\xx - \xxc)   , & \xx > \xxc .
    \end{cases}
\label{zzxxppdBlro}
\end{EQA}    
If \( \gm \geq \xx^{1/2}/2 + (\dimp \xx/4)^{1/4} \), then 
\begin{EQA}
    \P\bigl( \| \xiv \| \ge \zq(\dimp,\xx) \bigr)
    & \le &
    3 \ex^{-\xx} .
\label{PxivbzzBBro3}
\end{EQA}    
\end{corollary}

\begin{proof}[Proof of Theorem~\ref{LLbrevelocroB}]
%
First we consider the most interesting case \( \xx \leq \xxc \). 
We expect to get Gaussian type deviation bounds for such \( \xx \).
%
The main tool of the proof is the following lemma.
\begin{lemma}
\label{LGDBqfexpB}
Let \( \mu \in (0,1) \) and \( \zz(\mu) = \gmb / \mu - \sqrt{\dimH/\mu} > 0 \).
Then \eqref{expgamgm} implies
\begin{EQA}
	\E \exp\bigl( \mu \langle \BB \xiv, \xiv \rangle / 2 \bigr)
	\Ind\bigl( \| \HVB \BB \xiv \| \leq \zz(\mu) \bigr)
	& \leq &
	2 \exp \Bigl( \frac{\mu^{2} \vH^{2}}{4 (1 - \mu)} + \frac{\mu \, \dimH}{2} \Bigr) .
\label{wBmu2vA241mmutrB}
\end{EQA}
\end{lemma}

\begin{proof}
Let us fix for a moment some \( \xiv \in \R^{\dimp} \) and \( \mu < 1 \) and
define 
\begin{EQA}
	\av = \HVB^{-1} \xiv, 
	&\qquad & 
	\Sigma = \mu \BBH = \mu \HVB \BB \HVB .
\label{aSigmHm1xivmuc}
\end{EQA}
Consider the Gaussian measure \( \P_{\av,\Sigma} = \ND(\av,\Sigma^{-1}) \), and 
let \( \Uv \sim \ND(0,\Sigma^{-1}) \).
By the Girsanov formula
\begin{EQA}
	\log \frac{d\P_{\av,\Sigma}}{d\P_{0,\Sigma}}(\uv)
	&=&
	\langle \Sigma \av, \uv \rangle - \frac{1}{2} \bigl\langle \Sigma \av, \av \bigr\rangle
\label{GirslogddPaSSm1}
\end{EQA}
and for any set \( A \in \R^{\dimp} \)
\begin{EQA}
	\P_{\av,\Sigma}(A)
	&=&
	\P_{0,\Sigma}(A - \av)
	=
	\E_{0,\Sigma} \Bigl[ 
		\exp\Bigl\{ \langle \Sigma \Uv, \av \rangle 
			- \frac{1}{2} \bigl\langle \Sigma \av, \av \bigr\rangle \Bigr\} \Ind(A) 
	\Bigr] .
\label{PaSAP0SAav}
\end{EQA}
Now we select \( A = \bigl\{ \uv \colon \| \Sigma \uv \| \leq \gmb \bigr\} \).
Under \( \P_{0,\Sigma} \), one can represent 
\( \Sigma \Uv = \Sigma^{1/2} \gammav 
= \mu^{1/2} \HVB \BB^{1/2} \gammav \) with a standard Gaussian 
\( \gammav \).
Therefore,
\begin{EQA}
	\P_{0,\Sigma}(A - \av)
	&=&
	\Pg\bigl( \| \Sigma^{1/2} (\gammav - \Sigma^{1/2} \av) \| \leq \gmb \bigr)
	\\
	& \geq &
	\Pg\bigl( \| \Sigma^{1/2} \gammav \| \leq \gmb - \| \Sigma \av \| \bigr) .
\label{P0SAmaPgSm12}
\end{EQA}
We now use that \( \Pg\bigl( \| \Sigma^{1/2} \gammav \|^{2} \leq \tr (\Sigma) \bigr) \geq 1/2 \)
with \( \tr(\Sigma) = \mu \tr(\BBH ) = \mu \, \dimH \).
Therefore, the condition \( \| \Sigma \av \| + \sqrt{\mu \, \dimH} \leq \gmb \) implies
in view of \( \langle \Sigma \av, \av \rangle = \mu \langle \BB \xiv, \xiv \rangle \)
\begin{EQA}
	1/2
	\leq 
	\P_{\av,\Sigma}(A)
	&=&
	\E_{0,\Sigma} \Bigl[ 
		\exp\Bigl\{ \langle \Sigma \Uv, \HVB^{-1} \xiv \rangle
			- \mu \langle \BB \xiv, \xiv \rangle / 2 \Bigr\} 
			\Ind( \| \Sigma \Uv \| \leq \gmb) 
	\Bigr]
\label{12PavSAE0SISm1}
\end{EQA}
or
\begin{EQA}
	&& \nquad
	\exp\bigl( \mu \langle \BB \xiv, \xiv \rangle / 2 \bigr)
	\Ind\bigl( \| \Sigma \HVB^{-1} \xiv \| \leq \gmb - \sqrt{\mu \, \dimH} \bigr)
	\\
	& \leq &
	2 \E_{0,\Sigma} \Bigl[ 
		\exp\Bigl\{ \langle \Sigma \Uv, \HVB^{-1} \xiv \rangle
			\Ind( \| \Sigma \Uv \| \leq \gmb) 
	\Bigr] .
\label{emNxx22E0Sm1u}
\end{EQA}
We now take the expectation of the each side of this equation w.r.t. \( \xiv \),
change the integration order, and use \eqref{expgamgm} yielding
\begin{EQA}
	&& \nquad
	\E \exp\bigl( \mu \langle \BB \xiv, \xiv \rangle / 2 \bigr)
	\Ind\bigl( \| \Sigma \HVB^{-1} \xiv \| \leq \gmb - \sqrt{\mu \, \dimH} \bigr)
	\leq 
	2 \E_{0,\Sigma} \exp\bigl( \| \Sigma \Uv \|^{2} / 2 \bigr)
	\\
	&=&
	2 \Eg \exp\bigl( \mu \| \BBH^{1/2} \gammav \|^{2} / 2 \bigr)
	=
	2 \det\bigl( \Id - \mu \BBH \bigr)^{-1/2} .
\label{IHm1B12EgEgg}
\end{EQA}
We also use that for any \( \mu > 0 \)
\begin{EQA}
	\log \det\bigl( \Id - \mu \BBH \bigr)^{-1/2} - \frac{\mu \tr \BBH}{2}
	& \leq &
	\frac{\mu^{2} \tr \BBH^{2}}{4 (1 - \mu)} \, ;
\label{mu2v241mmiulogIm12}
\end{EQA}
see \eqref{m2v241m4mj1p}, 
and the first statement follows in view of \( \Sigma \HVB^{-1} \xiv = \mu \HVB \BB \xiv \).
\end{proof}

The use of \( \mu = \mux \) from \eqref{1v2sxm12m1m} in \eqref{wBmu2vA241mmutrB} yields 
similarly to the proof of Theorem~\ref{TexpbLGA}
\begin{EQA}
	\P\Bigl( \langle \BB \xiv, \xiv \rangle > \zq^{2}(\BBH,\xx), \,
		\| \HVB \BB \xiv \| \leq \zz(\mux)
	 \Bigr)
	 & \leq &
	 2 \ex^{-\xx} .
\label{2emxPblrHm1B}
\end{EQA}
It remains to consider the probability of large deviation
\( \P\bigl( \| \HVB \BB \xiv \| > \zz(\mux) \bigr) \).

\begin{lemma}
\label{Ldvbetagmb}
For any \( \xxc > 0 \) such that \( \zq(\BBH,\xxc) + 1 \leq \gmc/\muc \), 
it holds with \( \muc = \bigl\{ 1 + \vH/(2\sqrt{\xxc}) \bigr\}^{-1} \) 
and \( \zqc = \zz(\muc) = \gmb/\muc - \sqrt{\dimH/\muc} \)
\begin{EQA}
	\P\bigl( \| \HVB \BB \xiv \| > \zqc \bigr) 
	& \leq &
	\P\bigl( \langle \BB \xiv, \xiv \rangle > \zqc^{2} \bigr)
	\leq 
	\ex^{-\xxc} .
\label{exmxxlPHm1Bxigm}
\end{EQA}
\end{lemma}

\def\Pmuvp{\Phi}
\begin{proof}
Define 
\begin{EQA}
	\Pmuvp(\mu)
	& \eqdef &
	\frac{\mu^{2} \vH^{2}}{4 (1 - \mu)} + \frac{\mu \, \dimH}{2} \, .
\label{sdfw6rewqrwdwtffwet66}
\end{EQA}
It follows due to \eqref{1v2sxm12m1m} and \eqref{m2vA241muz2} for any \( \mu \leq \muc \)
\begin{EQA}
	\Pmuvp(\mu)
	\leq 
	\Pmuvp(\muc)
	& \leq &
	\frac{\muc \zq^{2}(\BBH,\xxc)}{2} - \xxc ,
\label{mumucmuz2Bxmx24}
\end{EQA}
where the right hand-side does not depend on \( \mu \).
Denote \( \eta^{2} = \langle \BB \xiv, \xiv \rangle \) and use that 
\( \| \HVB \BB \xiv \| \leq \| \BB^{1/2} \xiv \| = \eta \).
Then by \eqref{wBmu2vA241mmutrB}
\begin{EQA}
	\E \exp(\mu \eta^{2}/2) \Ind\bigl( \eta \leq \zz(\mu) \bigr)
	& \leq &
	2 \exp \Pmuvp(\mu)
	\leq 
	2 \exp \Pmuvp(\muc) .
\label{me22Iezmemcz2}
\end{EQA}
Define the inverse function \( \mu(\zz) = \zz^{-1}(\mu) \).
For any \( \zz \geq \zqc \), it follows from \eqref{me22Iezmemcz2} with 
\( \mu = \mu(\zz) \)
\begin{EQA}
	\E \exp\bigl\{ \mu(\zz) (\zz-1)^{2} /2 \bigr\} 
		\Ind\bigl( \zz - 1 \leq \eta \leq \zz \bigr)
	& \leq &
	2 \exp \Pmuvp(\muc)  
\label{emcz2NBxm2m1I}
\end{EQA}
yielding
\begin{EQA}
	\P\bigl( \zz - 1 \leq \eta \leq \zz \bigr)
	& \leq &
	2 \exp \Bigl( 
		- \mu(\zz) \, (\zz - 1)^{2}/2 + \Pmuvp(\muc) 
	\Bigr)
\label{mzm12mcz2Bxm2x}
\end{EQA}
and hence,
\begin{EQA}
	\P\bigl( \eta > \zz \bigr)
	& \leq &
	2 \int_{\zz}^{\infty} \exp\bigl\{ - \mu(\zq) (\zq-1)^{2}/2 + \Pmuvp(\muc) \bigr\} d\zq .
\label{zzzzcm1Pezmcl22}
\end{EQA}
Further, \( \mu \, \zz(\mu) = \gmb - \sqrt{\dimH \mu} \) and
\begin{EQA}
	\gmc
	=
	\muc \, \zqc
	& \leq &
	\mu \, \zz(\mu)
	\leq 
	\gmb,
	\quad
	\mu \leq \muc .
\label{muxzxmzmglex}
\end{EQA}
This implies the same bound for the inverse function:
\begin{EQA}
	\gmc
	& \leq &
	\zz \, \mu(\zz)
	\leq 
	\gmb,
	\quad
	\zz \geq \zqc \, ,
\label{muxzxmzmglexi}
\end{EQA}
and for \( \zz \geq 2 \)
\begin{EQA}
	\P\bigl( \eta > \zz \bigr)
	& \leq &
	2 \int_{\zz}^{\infty} 
		\exp\bigl\{ - \mu(\zq) \bigl( \zq^{2}/2 - \zq \bigr) + \Pmuvp(\muc) \bigr\} d\zq
	\\
	& \leq &
	2 \int_{\zz}^{\infty} 
		\exp\bigl\{ - \gmc \, (\zq/2  - 1) + \Pmuvp(\muc) \bigr\} d\zq
	\\
	& \leq &
	\frac{4}{\gmc}
		\exp\bigl\{ - \gmc \, (\zz / 2 - 1) + \Pmuvp(\muc) \bigr\} .
\label{Perzxinzinmmuz2gmz2}
\end{EQA}
Conditions \( \gmc \zqc = \muc^{-1} \gmc^{2} \geq \muc \bigl\{ \zq(\BBH,\xxc) + 1 \bigr\}^{2} \)
and \( \gmc \geq 1 \) 
ensure that
\( \P\bigl( \eta > \zqc \bigr) \leq \ex^{-\xxc} \).
\end{proof}

Remind that \( \xxc \) is the largest \( \xx \)-value ensuring the condition
\( \gmc \geq \zq(\BBH,\xxc) + 1 \).
We also use that for \( \xx \leq \xxc \), it holds \( \zz(\mux) \geq \zz(\muc) = \zqc \).
Therefore, by \eqref{2emxPblrHm1B} and Lemma~\ref{Ldvbetagmb} 
\begin{EQA}
	\P\bigl( \langle \BB \xiv, \xiv \rangle \geq \zq^{2}(\BBH,\xx) \bigr)
	& \leq &
	\P\bigl( \langle \BB \xiv, \xiv \rangle \geq \zq^{2}(\BBH,\xx), \| \HVB \BB \xiv \| \leq \zz(\mux) \bigr)
	+ \P\bigl( \langle \BB \xiv, \xiv \rangle \geq \zqc^{2} \bigr)
	\\
	& \leq &
	2 \ex^{-\xx} + \ex^{-\xxc} \, .
\label{2emxpemxcPBxx}
\end{EQA}
Finally we consider \( \xx > \xxc \).
Applying \eqref{Perzxinzinmmuz2gmz2} yields by \( \zz \geq \zqc \)
\begin{EQA}
	\P\bigl( \eta > \zz \bigr)
	& \leq &
	\frac{2}{\muc \, \zqc}
		\exp\bigl\{ - \muc \, \zqc^{2} / 2 + \gmb + \muc \, \zq^{2}(\BBH,\xxc)/2 - \xxc \bigr\}
		\exp\bigl\{ - \muc \, \zqc (\zz - \zqc) / 2  \bigr\}
	\\
	& \leq &
	\ex^{- \xxc} \exp\bigl\{ - \gmc (\zz - \zqc) / 2  \bigr\} .
\label{emxcmcmzczzc2}
\end{EQA}
The choice \( \zz \) by 
\begin{EQA}
	\gmc (\zz - \zqc) / 2 
	&=&
	\xx - \xxc
\label{mzzc2xxczz}
\end{EQA}
ensures the desired bound.

Now, for a prescribed \( \xx \), we evaluate the minimal value \( \gm \) ensuring the bound \eqref{PxivbzzBBroB} with \( \xxc \geq \xx \). 
For simplicity we apply the sub-optimal choice 
\( \mu(\xx) = \bigl( 1 + 2 \sqrt{\dimH/\xx} \bigr)^{-1} \); see Remark~\ref{RsochpHsA}. 
Then for any \( \xx \geq 1 \)
\begin{EQA}
	\mu(\xx) \, \bigl\{ \zq(\BBH,\xx) + 1 \bigr\}
	& \leq &
	\frac{\sqrt{\xx}}{\sqrt{\xx} + 2 \sqrt{\dimH}} \, \Bigl( \sqrt{\dimH + 2(\xx \dimH)^{1/2} + 2\xx} + 1 \Bigr) 
	\, ,
	\\
	\dimH \, \mu(\xx)
	& = &
	\frac{\sqrt{\xx} \, \dimH }{\sqrt{\xx} + 2 \sqrt{\dimH}} 
	\, .
\label{mindH2sHx2x2x12}
\end{EQA} 
It is now straightforward to check that 
\begin{EQA}
	\mu(\xx) \, \bigl\{ \zq(\BBH,\xx) + 1 \bigr\} + \sqrt{\dimH \, \mu(\xx)}
	& \leq &
	\sqrt{\xx}/2 + (\xx \, \dimH/4)^{1/4} .
\label{x4p14sx2mx1}
\end{EQA}
Therefore, if \eqref{g2x1214px14} holds for the given \( \xx \), 
then \eqref{gmsqpAHmxzBx1} is fulfilled with \( \xxc \geq \xx \)
yielding \eqref{PxivbzzBBroB3}. 
%
\end{proof}

\Section{Weighted sums of Bernoulli r.v.'s: univariate case}
Let \( Y_{i} \) be independent \( \Bernoulli(\thetas_{i}) \) and \( \weight_{i} \in [0,1] \). 
First we state a deviation bound for a centered sum on non i.i.d. Bernoulli random variables.
\begin{lemma}
\label{LBeBvM}
Let \( Y_{i} \) be independent \( \Bernoulli(\thetas_{i}) \)
and \( \weight_{i} \in \R \). 
Define 
\begin{EQA}
	S
	&=&
	\sum_{i=1}^{\nsize} Y_{i} \weight_{i} \, ,
	\\
	\VP^{2}
	&=&
	\Var(S)
	=
	\sum_{i=1}^{\nsize} \thetas_{i}(1-\thetas_{i}) \weight_{i}^{2} \, ,
	\\
	\weights 
	&=& 
	\max_{i} |\weight_{i}| .
\label{Het2Fetetp}
\end{EQA}
Then it holds  
\begin{EQA}
	\log \E \exp\Bigl\{ \frac{\lambda (S - \E S)}{\VP} \Bigr\}
	& \leq &
	\lambda^{2}	,
	\qquad
	\lambda \leq \log(2) \VP/\weights.
\label{llEelY1tsel}
\end{EQA}
Furthermore, suppose that given \( \xx \geq 0 \), 
\begin{EQA}
	\VP 
	& \geq &
	\frac{3}{2} \, \xx^{1/2} \weights \, .
\label{Het2Fetetpzz}
\end{EQA}
Then
\begin{EQA}
	\P\bigl( \VP^{-1} |S - \E S| \geq 2 \sqrt{\xx} \bigr)
	& \leq &
	2 \ex^{-\xx} .
\label{PHetm1SESs2x}
\end{EQA}
Without \eqref{Het2Fetetpzz}, the bound \eqref{PHetm1SESs2x} applies
with \( \VP \) replaced by \( \VP_{\xx} = \VP \vee (3 \, \xx^{1/2} \weights/2) \).
\end{lemma}

\begin{proof}
Without loss of generality assume \( \weights = 1 \), otherwise just rescale all the weights by the factor \( 1/\weights \).
We use that 
\begin{EQA}
	f(u)
	\eqdef
	\log \E \exp\Bigl\{ u (S - \E S) \Bigr\}
	&=&
	\sum_{i=1}^{\nbin} \Bigl[ \log\bigl( \thetas_{i} \ex^{u \weight_{i}} + 1 - \thetas_{i} \bigr) 
	- u \weight_{i} \thetas_{i} \Bigr] \, .
\label{llEelY1tsela}
\end{EQA}
This is an analytic function of \( u \) for \( |u| \leq \log 2 \) satisfying 
\( f(0) = 0 \), \( f'(0) = 0 \), and, with \( \upss_{i} = \log \thetas_{i} - \log(1-\thetas_{i}) \),
\begin{EQA}
	f''(u)
	&=&
	\sum_{i=1}^{\nbin} 
	\frac{\weight_{i}^{2} \, \thetas_{i}(1-\thetas_{i}) \, \ex^{u \weight_{i}}}
		 {( \thetas_{i} \ex^{u \weight_{i}} + 1 - \thetas_{i})^{2}}
	=
	\sum_{i=1}^{\nbin} \frac{\weight_{i}^{2} \, \ex^{\upss_{i}+u \weight_{i}}}{(\ex^{\upss_{i}+u \weight_{i}} + 1)^{2}} 
	=
	\sum_{i=1}^{\nbin} \theta_{i}(u) \bigl\{ 1-\theta_{i}(u) \bigr\} \weight_{i}^{2}
\label{nu2He2fpput}
\end{EQA}
for \( \theta_{i}(u) = \ex^{\upss_{i}+u \weight_{i}}/(\ex^{\upss_{i}+u \weight_{i}} + 1) \).
Clearly \( \theta_{i}(u) \) and thus, \( \theta_{i}(u) \bigl\{ 1-\theta_{i}(u) \bigr\} \) monotonously increases with \( u \) and 
it holds for \( \thetas_{i} = \theta_{i}(0) \) 
\begin{EQA}
	\theta_{i}(u) \bigl\{ 1-\theta_{i}(u) \bigr\}
	& \leq &
	\ex^{|u|} \, \thetas_{i} \, (1 - \thetas_{i})
	\leq 
	2 \, \thetas_{i} \, (1 - \thetas_{i}) ,
	\qquad
	|u| \leq \log 2.
\label{2ts1tseu1tb}
\end{EQA}
This yields 
\begin{EQA}
	f(u)
	& \leq & 
	\VP^{2} \, u^{2}
	\qquad
	|u| \leq \log 2.
\label{flaHetnu2la22}
\end{EQA}
As \( \xx \leq 4\VP^{2}/9 \), the value \( \lambda = \sqrt{\xx} \) fulfills 
\( \lambda/\VP = \sqrt{\xx}/\VP \leq  \log 2 \leq 2^{-1/2} \).
Now by the exponential Chebyshev inequality 
\begin{EQA}
	\P\Bigl( \VP^{-1}(S - \E S) \geq 2 \sqrt{\xx} \Bigr)
	& \leq &
	\exp\bigl\{ - 2 \lambda \sqrt{\xx} + f(\lambda/\VP) \bigr\}
	\\
	& \leq &
	\exp\bigl( - 2 \lambda \sqrt{\xx} + \lambda^{2} \bigr)
	=
	\ex^{-\xx} .
\label{exmxxsq2xnu2la22}
\end{EQA}
Similarly one can bound \( \E S - S \).
\end{proof}

\Section{Deviation bounds for weighted sums of Bernoulli r.v.'s in a \( \ell_{2} \)-norm}

Now we present a bound for the norm of the normalized vector \( \xiv = \Psiv (\Yv - \E \Yv) \),
where \( \Psiv \) is a linear mapping \( \Psiv \colon \R^{\dimp} \to \R^{\dimq} \).
Let 
\begin{EQA}
	\VP^{2} 
	=  
	\Var(\xiv)
	&=&
	\Var(\Psiv \Yv)
	=
	\Psiv \Var(\Yv) \Psiv^{\T} .
\label{HP2VxiWWYT}
\end{EQA}
We aim at bounding the norm \( \| \HP^{-1} \xiv \| \), where an operator \( \HP \) satisfies
\( \HP^{2} \geq \VP^{2} \).

\begin{lemma}
\label{LHDxibound}
Let \( Y_{i} \sim \Bernoulli(\thetas_{i}) \), \( i = 1,\ldots, n \), and let 
\( (\weight_{i,\lsh}) \), \( \lsh \in \WS \),
be a collection of weighting schemes. 
Consider \( \xiv = \Psiv (\Yv - \E \Yv) \), and let \( \HV^{2} \geq \Var(\xiv) \).
Define 
\begin{EQA}
	\weights
	&=&
	\max_{i} \| \HP^{-1} \Psiv_{i} \| \, .
\label{WHm1uinfmiwi}
\end{EQA} 
Then with 
\( \dimA = \tr\bigl( \HV^{-2} \VP^{2} \bigr) \), it holds
\begin{EQA}
    \P\bigl( \| \HV^{-1} \xiv \| \ge \sqrt{\dimH} + \sqrt{2\xx} \bigr)
    & \le &
    3 \ex^{-\xx} ,
\label{PxivbzzBBroB3m}
\end{EQA}    
provided that
\begin{EQA}
	1 / \weights
	& \geq &
	\sqrt{\xx}/2 + (\dimA \xx/4)^{1/4} .
\label{g2x1214px14b}
\end{EQA}
\end{lemma}

\begin{proof}
We apply the general result of Theorem~\ref{Tdevboundgm} under conditions \eqref{expgamgm} and \eqref{g2x1214px14}.
For any vector \( \uv \), consider the scalar product 
\( \langle \HP^{-1} \xiv,\uv \rangle = \langle \HP^{-1} \Psiv (\Yv - \E \Yv), \uv \rangle \).
It is obviously a weighted centered sum of the Bernoulli r.v.'s \( Y_{i} - \thetas_{i} \) with
\begin{EQA}
	\Var \langle \HP^{-1} \xiv,\uv \rangle
	& \leq &
	\| \uv \|^{2} .
\label{VlHm1xiuu22}
\end{EQA}
One can write with \( \eps_{i} = Y_{i} - \thetas_{i} \) and \( \epsv = (\eps_{i}) \)
\begin{EQA}
	\langle \HP^{-1} \xiv,\uv \rangle
	&=&
	\bigl\langle \epsv, \Psiv^{\T} \HP^{-1} \uv \bigr\rangle .
\label{lHm1xiuepWH1u}
\end{EQA}
By the Cauchy-Schwarz inequality, it holds 
\begin{EQA}
	\| \Psiv^{\T} \HP^{-1} \uv \|_{\infty}
	=
	\max_{i} \bigl| (\HP^{-1} \Psiv_{i})^{\T} \uv \bigr|
	& \leq &
	\weights \| \uv \| .
\label{WHm1uinfmiwi}
\end{EQA} 
The bound \eqref{llEelY1tsel} of Lemma~\ref{LBeBvM} on the exponential moments of \( \langle \HV^{-1} \xiv,\uv \rangle \) implies 
\begin{EQA}
	\log \E \exp\bigl\{ \langle \HV^{-1} \xiv,\uv \rangle \bigr\}
	& \leq &
	\| \uv \|^{2}	,
	\qquad
	\| \uv \| \leq \log(2) / \weights \, .
\label{llEelY1tsel2}
\end{EQA}
Therefore, \eqref{expgamgm} is fulfilled with \( \HVB^{2} = 2 \HV^{2} \) and 
\( \gm = \sqrt{2} \, \log(2) / \weights \), and
\eqref{g2x1214px14} follows from \eqref{g2x1214px14b}. 
The deviation bound \eqref{PxivbzzBBroB3} of Theorem~\ref{Tdevboundgm} yields the assertion.
\end{proof}

\bibliography{exp_ts,listpubm-with-url}

\end{document}